\documentclass[12pt]{amsart}

\usepackage[utf8]{inputenc}
\usepackage[margin=1in]{geometry}
\usepackage{todonotes}
\usepackage{xcolor}
\usepackage{bbm}
\usepackage{tikz}
\usetikzlibrary{arrows}
\usetikzlibrary{arrows.meta}
\usepackage{tikz-cd}
\usepackage{pgffor}
\usepackage{amssymb,amsmath,amsfonts,mathrsfs,amsthm}
\usepackage{enumitem}
\usepackage{adjustbox}
\usepackage{hyperref}
\usepackage{relsize}
\usepackage{graphicx}
\usepackage[font={small}]{caption}
\usepackage{float}
\usepackage{subcaption}
\usepackage{changepage}
\usepackage[toc]{appendix}
\usepackage{pgfplots}
\pgfplotsset{compat=1.7}
\usepackage{adjustbox}
\usepackage{varwidth}
\usepackage{chngcntr}
\usepackage{standalone}
\usepackage{tcolorbox}
\usepackage{mathtools}
\usepackage{relsize}
\usepackage{wasysym}
\usepackage{tensor} 
\usepackage{longtable}
\usepackage{makecell}
\usepackage{tabu}
\usepackage{tabulary}
\usetikzlibrary{calc}
\usetikzlibrary{patterns}
\usetikzlibrary{patterns.meta}
\usetikzlibrary{knots}
\usetikzlibrary{arrows.meta,calc,decorations.markings,math,arrows.meta}
\usepackage{esvect}
\usepackage{array}   
\newcolumntype{L}{>{$}l<{$}}
\newcolumntype{C}{>{$}c<{$}}

\newlist{alphalist}{enumerate}{1}
\setlist[alphalist,1]{label=\textbf{\Alph*.}}

\newcommand{\R}{\mathbb{R}}

\newcommand{\Z}{\mathbb{Z}}

\renewcommand{\S}{\mathbb{S}}
\newcommand{\F}{\mathcal{F}}
\newcommand{\Sq}{\S[G]}
\newcommand{\Zq}{\Z[G]}
\newcommand{\ZG}{\Z[G]}
\newcommand{\ZGmod}{\ggRmod[{{\ZG}}]}
\newcommand{\Hom}{\operatorname{Hom}}

\newcommand{\hocolim}{\operatorname{hocolim}}
\newcommand{\id}{\operatorname{id}}

\renewcommand{\d}{\partial}
\newcommand{\q}{\mathfrak{q}}

\newcommand{\LR}[1]{\lVert{#1}\rVert}

\newcommand{\D}{\mathbb{D}}  
\newcommand{\A}{\mathbb{A}}  
\newcommand{\Aq}{\A_q} 

\newcommand{\BN}{\mathcal{BN}}
\newcommand{\BNq}[1][]{\BN^{#1}(\Aq)}
\newcommand{\divcob}{\operatorname{Cob}_d} 
\newcommand{\divcobq}[1][]{\divcob^{#1}(\Aq)}

\newcommand{\FKh}[1][\D]{\F_{#1}} 
\newcommand{\qFKh}[1][]{\prescript{#1}{q}{\FKh}} 
\newcommand{\FA}{\F_{\A}} 
\newcommand{\FCKbar}{\til{\F}^k_{\mathsmaller{CK}}}
\newcommand{\FCK}{\F_{\mathsmaller{CK}}} 
\newcommand{\qFCK}[1][]{\prescript{#1}{q}{\FCK}} 
\newcommand{\FCKk}{\FCK^k}
\newcommand{\qFCKk}[1][]{\prescript{#1}{q}{\FCKk}}
 
\newcommand{\FAq}[1][]{\prescript{#1}{}{\F}_{\Aq}} 
\newcommand{\choiceFAq}[1][] {\FAq[#1]^{\choice}} 

\newcommand{\VHKK}[1][\D]{V_{#1}}
\newcommand{\VHKKlax}[1][\D]{\VHKK[#1]^{\text{lax}}}
\newcommand{\VHKKann}{V_{\A}}
\newcommand{\VHKKannlax}{\VHKKann^{\text{lax}}}
\newcommand{\qVHKK}[1][]{\prescript{#1}{}{V}_{\Aq}}
\newcommand{\qVHKKchoice}[1][]{\prescript{#1}{}{V}_{\Aq,\Psi}}
\newcommand{\qVHKKlax}[1][]{\qVHKKchoice[#1]^{\text{lax}}}
\newcommand{\modifiedqHMcombo}{\func{G}_T^{H^\omega}}
\newcommand{\modifiedqHMcomborestricted}{\func{G}_T^{\omega}}
\newcommand{\rectifiedqHMcombo}{\rect_T^{H^\omega}}

\newcommand{\rect}{\mathfrak{R}}

\newcommand{\KC}[1][\D]{C_{#1}}
\newcommand{\aKC}{C_{\A}}
\newcommand{\qaKC}[1][]{\prescript{#1}{}{C}_{\Aq}}
\newcommand{\qaKCchoice}[1][]{\qaKC[#1]^\choice}

\newcommand{\aKh}{Kh_{\A}}
\newcommand{\qaKh}[1][]{\prescript{#1}{}{Kh}_{\Aq}}
\newcommand{\qaKhchoice}[1][]{\qaKh[#1]^\choice}
\newcommand{\CCK}{C_{\mathsmaller{CK}}} 
\newcommand{\CCKnk}{\CCK^{n-k,k}}
\newcommand{\qCCK}[1][]{\prescript{#1}{q}{C}_{\mathsmaller{CK}}} 
\newcommand{\qCCKnk}{\qCCK^{n-k,k}}

\newcommand{\IMC}{\mathrm{IMC}} 
\newcommand{\X}{\mathcal{X}_{\mathsmaller{CK}}} 
\newcommand{\Xk}{\X^{n-k,k}} 
\newcommand{\Xfunc}[1][\mathsmaller{CK}]{\mathscr{X}_{#1}}  
\newcommand{\Xfunck}[1][\mathsmaller{CK}]{\Xfunc[#1]^{n-k,k}}  

\newcommand{\qXk}[1][]{\prescript{#1}{q}{\Xk}}
\newcommand{\XKh}[1][\D]{\mathcal{X}_{#1}} 
\newcommand{\aX}{\mathcal{X}_{\A}} 
\newcommand{\qaX}[1][]{\prescript{#1}{}{\mathcal{X}}_{\Aq}} 
\newcommand{\qaXchoice}[1][]{\prescript{#1}{}{\mathcal{X}}^\choice_{\Aq}}
\newcommand{\qakwX}[1][]{\prescript{#1}{}{\mathcal{X}}_{\mathsmaller{AKW}}} 
\newcommand{\hatX}[1][]{\prescript{#1}{}{\h{\mathcal{X}}}} 

\newcommand{\Top}{\operatorname{Top}}
\newcommand{\til}[1]{\widetilde{#1}}
\renewcommand{\b}[1]{\overline{#1}}
\newcommand{\h}[1]{\widehat{#1}}

\renewcommand{\k}{\Bbbk}

\newcommand{\adeg}{\operatorname{adeg}}
\newcommand{\gr}{\mathrm{gr}} 
\newcommand{\lr}[1]{\vert {#1} \vert}
\newcommand{\Kom}{\operatorname{Kom}}
\newcommand{\Func}{\operatorname{Func}}
\newcommand{\Map}{\operatorname{Map}}
\newcommand{\EnrFunc}{\operatorname{EnrFunc}}

\renewcommand{\u}[1]{\underline{#1}}
\renewcommand{\o}{\otimes}

\newcommand{\sh}{\operatorname{sh}}
\newcommand{\THH}{\mathrm{THH}}
\newcommand{\qTHH}{\mathrm{qTHH}}
\newcommand{\qhh}{\operatorname{qHH}}
\newcommand{\qch}{\operatorname{qCH}}
\newcommand{\hh}{\mathrm{HH}}
\newcommand{\ch}{\mathrm{CH}}
\newcommand{\Id}{\id}
\newcommand{\eps}{\varepsilon}

\newcommand{\Ab}{\operatorname{Ab}}
\newcommand{\Aut}{\operatorname{Aut}}

\newcommand{\Spectra}{\mathscr{S}}
\newcommand{\gSpectra}{\Spectra^{gr}}
\newcommand{\ggSpectra}{\Spectra^{gg}}

\newcommand{\sset}{\operatorname{SSet}}
\newcommand{\gAb}{\Ab^{gr}}
\newcommand{\ggAb}{\Ab^{gg}}

\newcommand{\Rmod}[1][R]{#1\text{-}\mathrm{Mod}}
\newcommand{\gRmod}[1][R]{\Rmod[#1]^{gr}}
\newcommand{\ggRmod}[1][R]{\Rmod[#1]^{gg}}

\newcommand{\bimod}[2]{\tensor[_{#1}]{\mathrm{Mod}}{_{#2}}}

\newcommand{\cat}[1]{\mathcal{#1}} 

\newcommand{\opcat}[1]{\cat{#1}^{\mathrm{op}}}
\newcommand{\spcat}[1]{\cat{#1}} 

\newcommand{\func}[1]{\mathscr{#1}} 
\newcommand{\shape}{\cat{S}} 

\newcommand{\HMshape}{\shape_{HM}}

\newcommand{\qHMfunc}[1]{\func{F}^{qH}_{#1}}
\newcommand{\chainsfunc}{C_h}
\newcommand{\ring}{R}
\newcommand{\ringspec}{\mathcal{R}}

\newcommand{\platalg}[1][n]{\spcat{A}^{#1}}
\newcommand{\platalgsummand}[2][n]{\spcat{A}^{#1-#2,#2}}
\newcommand{\platalgk}{\platalgsummand[n]{k}}

\newcommand{\qplatalgk}[1][] {\prescript{#1}{q}{\spcat{A}}^{n-k,k}}

\newcommand{\matchingssummand}[2][n]{B^{#1-#2,#2}}
\newcommand{\matchingsk}{\matchingssummand[n]{k}}

\newcommand{\Burn}{\mathscr{B}}
\newcommand{\gBurn}{\Burn^{gr}}
\newcommand{\ggBurn}{\Burn^{gg}}

\newcommand{\Permu}{\mathrm{Permu}}
\newcommand{\Sets}{\mathrm{Sets}}

\newcommand{\Arr}{\operatorname{Arr}}

\newcommand{\BtoAb}{\func{U}} 
\newcommand{\KEM}{\func{K}} 

\newcommand{\sq}{s_q}
\newcommand{\tq}{t_q}

\newcommand{\cobfunc}{\func{C}} 

\newcommand{\ob}{\operatorname{ob}} 
\newcommand{\lar}[1]{\xleftarrow{#1}}
\newcommand{\rar}[1]{\xrightarrow{#1}}
\newcommand{\squig}[1]{\overset{\rightsquigarrow}{#1}}
\newcommand{\nerve}{\mathscr{N}}

\newcommand{\match}[1]{B^{#1}} 
\newcommand{\cupdiag}[2]{B^{{#1}-{#2},{#2}}} 
\newcommand{\closure}[1]{\left[{#1}\right]} 
\newcommand{\flip}[1]{\widehat{#1}} 
\newcommand{\strands}[1]{\underline{\overline{#1}}} 

\newcommand{\qXi}[1][]{\prescript{#1}{q}{\Xi}} 
\newcommand{\nmshift}[2]{\phi_{#1,#2}}

\newcommand{\Plat}[1]{{A}^{#1}} 
\newcommand{\qPlat}[2][]{\prescript{#1}{q}{{A}^{#2}}} 
\newcommand{\qPlatzero}[2][]{\prescript{#1}{q}{{A}^{#2}_0}} 
\newcommand{\Platk}[2]{{A}^{{#1}-{#2},{#2}}}
\newcommand{\Platnk}{\Platk{n}{k}}
\newcommand{\qPlatk}[2]{\prescript{}{q}{A}^{{#1}-{#2},{#2}}}
\newcommand{\qPlatnk}[1][]{\prescript{#1}{q}{A}^{n-k,k}}

\newcommand{\Platbark}[2]{{\til{A}}^{{#1}-{#2},{#2}}}
\newcommand{\Ideal}[2]{\mathcal{I}^{{#1}-{#2}, {#2}}} 
\newcommand{\mq}{\mu_q} 

\newcommand{\scext}[2][]{\prescript{#1}{q}{#2}}  

\newcommand{\choice}{\Psi}

\newtheorem{theorem}{Theorem}[section]
\newtheorem{lemma}[theorem]{Lemma}
\newtheorem{proposition}[theorem]{Proposition}
\newtheorem{corollary}[theorem]{Corollary}

\theoremstyle{definition}

\theoremstyle{remark}
\newtheorem{example}[theorem]{Example}

\theoremstyle{remark}
\newtheorem{remark}[theorem]{Remark}

\theoremstyle{remark}
\newtheorem{question}[theorem]{Question}

\theoremstyle{definition}
\newtheorem{definition}[theorem]{Definition}

\newcommand{\tiltimes}{\widetilde{\times}}
\newcommand{\cube}{\underline{2}} 
\newcommand{\B}{\mathscr{B}}

\newcommand{\Mq}{\mathcal{M}_q}

\newcommand{\Diff}{\operatorname{Diff}}
\newcommand{\st}{\circ}
\newcommand{\cut}{\operatorname{cut}} 
\newcommand{\Gen}{\operatorname{Gen}}

\newcommand{\sphere}{\mathbb{S}}

\title{Quantum topological Hochschild homology and annular Khovanov spectra}
\author[R. Akhmechet]{Rostislav Akhmechet}
\address{Department of Mathematics, Columbia University, New York, NY, 10027}
\email{\href{mailto:akhmechet@math.columbia.edu}{akhmechet@math.columbia.edu}}

\author[T. Gerhardt]{Teena Gerhardt}
\address{Department of Mathematics, Michigan State University, East Lansing, MI, 48824}
\email{\href{mailto:teena@math.msu.edu}{teena@math.msu.edu}}

\author[M. Willis]{Michael Willis}
\address{Department of Mathematics, Texas A\&M University, College Station, TX, 77840}
\email{\href{mailto:msw188@tamu.edu}{msw188@tamu.edu}}

\date{}

\begin{document}

\begin{abstract}
    Topological Hochschild homology is a topological analogue of classical Hochschild homology of algebras and bimodules. Beliakova, Putyra, and Wehrli introduced quantum Hochschild homology ($\qhh$) and used it to define a quantization of annular Khovanov homology as $\qhh$ of the tangle bimodules of Chen-Khovanov and Stroppel.  After introducing quantum topological Hochschild homology ($\qTHH$), we construct a new stable homotopy refinement of quantum annular Khovanov homology and show that it agrees with $\qTHH$ of the spectral Chen-Khovanov tangle bimodules of Lawson, Lipshitz, and Sarkar.  We also show that this new stable homotopy refinement recovers the construction introduced in earlier work of Krushkal together with the first and third authors.
\end{abstract}
\maketitle
\tableofcontents

\section{Introduction}
Given a fixed commutative ground ring $R$ with chosen unit $q\in R$, Beliakova, Putyra, and Wehrli in \cite{BPW} define a notion of \emph{quantum Hochschild homology} $\qhh(A;M)$ for a graded $R$-algebra $A$ and graded $(A,A)$-bimodule $M$, or more generally for a chain complex of graded bimodules.  They use this to define the \emph{quantum annular homology} $\qaKh(L)$ of a link $L$ in the thickened annulus as follows.  We let $T$ denote any $(n,n)$-tangle whose annular closure recovers $L$.  Chen and Khovanov  \cite{Chen-Khovanov} and Stroppel \cite{Stroppel} define the \emph{platform algebra} $\Plat{n}$ and a complex $\CCK(T)$ of $(\Plat{n},\Plat{n})$-bimodules, all over the ground ring $\Z$.  The quantum annular homology can be defined as $\qaKh(L):=\qhh(\Plat{n}\otimes R; \CCK(T)\otimes R)$.

Moreover, \cite{BPW} provides an effective method of computing $\qaKh(L)$ by first considering the \emph{quantum annular complex} $\qaKC(L)$, defined by applying the functor 
\[
\qhh(\Plat{n}\otimes R ; (-) \otimes R)
\]
term-wise to $\CCK(T)$, and showing that the higher quantum Hochschild homology groups of each chain group $\CCK^i(T)\otimes R$ vanish \cite[Proposition 6.6]{BPW}. Moreover, the chain groups and differentials in $\qaKC(L)$ can be computed without reference to $\qhh$ \cite[Theorem 6.3]{BPW}, and it is shown in \cite[Theorem C]{BPW} that 
\begin{equation}
\label{eq:total cx vs termise iso}
\qaKh(L) \cong H_*(\qaKC(L)).
\end{equation}
In the case that $R=\Z$ and $q=1$ (so that $\qhh$ is identified with Hochschild homology $\hh$), $\qaKC(L)$ recovers the \emph{annular Khovanov complex} $\aKC(L)$ of Asaeda, Przytycki, and Sikora \cite{APS}. The above isomorphism then specializes to
\begin{equation}
    \label{eq:akh cong hh}
    \hh(\Plat{n};\CCK(T))\cong \aKh(L),
\end{equation}
where $\aKh(L)$ is the \emph{annular Khovanov homology}; this establishes a conjecture of Auroux, Grigsby, and Wehrli \cite{AGW}. 

Recent interesting developments concerning quantum Hochschild homology include work of Beliakova, Hogancamp, Putyra, and Wehrli \cite{BHPW} and Beliakova, Putyra, Robert, and Wagner \cite{BPRW}, whose methods rely on $q$ being sufficiently generic.

Many of the Khovanov link and tangle invariants have spectral lifts beginning with the work of Lipshitz and Sarkar \cite{LS-Kh-spectra} and Hu, Kriz, and Kriz \cite{HKK}, which are known to contain more information than their homological counterparts \cite{LS-Steenrod,Seed,SZ-localization,SW_projectors}.  A spectral refinement $\aX(L)$ of annular Khovanov homology is discussed in \cite{SZ-localization, AKW, LLS_CK_Spectra, Transverseinvtspectra}.  Lawson, Lipshitz, and Sarkar  \cite{LLS_CK_Spectra} define the \emph{spectral platform algebras} $\platalg$ and associate to an $(n,n)$-tangle $T$ the \emph{spectral Chen-Khovanov bimodule} $\X(T)$, building on earlier work of the same authors \cite{LLS_Khovanov_Spectra}.  These spectral invariants lift the corresponding algebraic invariants in the sense that $H_*(\aX(L))\cong \aKh(L)$,\footnote{Note, elsewhere in the literature the \emph{cohomology} of the spectrum recovers the Khovanov invariant; see Remark \ref{rmk:homology vs cohomology} for a discussion.} and similarly for $\platalg$ and $\X(T)$. The isomorphism \eqref{eq:akh cong hh} was lifted to the level of spectra in \cite[Theorem 8]{LLS_CK_Spectra}, which establishes a weak equivalence 
\begin{equation}\label{eq:aX cong THH}
\THH(\platalg;\X(T)) \simeq \aX(L),
\end{equation}
where again $T$ is any tangle whose annular closure recovers $L$, and $\THH$ denotes \emph{topological Hochschild homology}.  In addition, the first and third authors together with Krushkal in \cite{AKW} defined a quantum annular spectrum $\qakwX(L)$, without reference to any Hochschild constructions, in cases where $R=\Z[G]$ for some finite cyclic group $G$ with fixed generator $q\in G$. The generator $q$ is viewed as part of a $G$-action on $\qakwX(L)$, and it is shown that $H_*(\qakwX(L))\cong \qaKh(L)$ as $\Z[G]$-modules.  

In this paper, we begin by presenting a definition for quantum topological Hochschild homology $\qTHH(A;M)$ of a graded $\Sq$-algebra $A$ and graded $(A,A)$-bimodule $M$, where $\Sq$ is the spherical group ring of a (potentially infinite) cyclic group $G$.  After presenting some of the general properties of the construction, we compute the quantum topological Hochschild homology of the spectral Khovanov invariants discussed above in the spirit of the earlier discussion on computing $\qaKh(L)$.  To that end we define the \emph{quantum annular spectrum} $\qaX(L)$, a spectrum with $G$-action, without reference to Hochschild constructions, which  satisfies $H_*(\qaX(L)) \cong \qaKh(L)$ as $\Z[G]$-modules.  Our first main result is the following. 

\begin{theorem}
\label{thm: qthh recovers bpw spectrum}
Let $T$ be an $(n,n)$-tangle whose annular closure is denoted $L$.  Then we have a weak equivalence of spectra with $G$-action
\[ \qTHH(\platalg\wedge \Sq;\X(T) \wedge \Sq ) \simeq \qaX(L).\]
\end{theorem}

Note that taking the smash product with $\Sq$ corresponds to extending scalars from $\Z$ to $R$ as in the definition of $\qaKh(L)$.  Indeed there is a more precise version of Theorem \ref{thm: qthh recovers bpw spectrum} which takes grading information into account; see Theorem \ref{thm:qTHH recovers qAX}.  Theorem \ref{thm: qthh recovers bpw spectrum} gives an effective method of computing $\qTHH$ for spectral Chen-Khovanov bimodules and generalizes Equation \eqref{eq:aX cong THH}, which can be viewed as the case where $G=\{1\}$.  Our second main result shows that we can recover the version in \cite{AKW}.

\begin{theorem}
\label{thm:main thm recover AKW}
Let $L$ be an annular link, and let $G$ be finite cyclic.  Then we have a weak equivalence of spectra with $G$-action
\[\qaX(L) \simeq \qakwX(L).\]
\end{theorem}

An analogue of $\qakwX(L)$ when $G \cong \Z$ is infinite was introduced in \cite{bigger_Burnside}; we expect their construction to agree with $\qaX(L)$ but do not pursue it here.

This work was motivated in part by the fact that, while the various algebraic Khovanov invariants have strong ties with representation theory, there are almost no known analogues for the spectral invariants (see \cite{DGML} for some recent work on the subject; earlier constructions such as those in \cite{BCSteenrod} should be related).  In particular, in \cite{BPW} it is shown that the quantum annular homology $\qaKh(L)$ carries an action of the quantum group $U_q(\mathfrak{sl}_2)$, generalizing the $\mathfrak{sl}_2$ action of Grigsby, Licata, and Wehrli \cite{GLW} on $\aKh(L)$.  In \cite{AKWsl2} the first and third authors, together with  Krushkal, constructed corresponding endomorphisms of $\aX(L)$ for each of the standard generators of $\mathfrak{sl}_2$ which lift the action on homology, although the $\mathfrak{sl}_2$ relations between these generators have not yet been verified at the spectral level.  In fact, the $\mathfrak{sl}_2$ action in \cite{GLW} was motivated by the isomorphism \eqref{eq:akh cong hh} (which was conjectured but unproven at the time), which would give a conceptual explanation for an $\mathfrak{sl}_2$ action on $\aKh(L)$ from the point of view of categorified quantum groups; see \cite[Section 1.3]{GLW} for a more detailed discussion. With these facts in mind we ask:

\begin{question}
Is there any sense in which $\qTHH(\platalg\wedge \Sq;\X(T) \wedge \Sq )$ naturally carries an action of $U_q(\mathfrak{sl}_2)$?
\end{question}

\subsection{An overview of the proof of Theorem \ref{thm: qthh recovers bpw spectrum}}
The proof of Theorem \ref{thm: qthh recovers bpw spectrum} is a generalization of the proof of \cite[Theorem 8]{LLS_CK_Spectra}.  As before, $T$ will denote an $(n,n)$-tangle with annular closure $L$ in the thickened annulus.

To begin, we present the situation more carefully.  The platform algebras $\Plat{n}$ and tangle bimodules $\CCK(T)$ discussed above carry a \emph{weight grading}; we use the notation $\Platnk$ and $\CCKnk(T)$ to denote the summands in weight $k$.  The spectral lifts are similar, with notations $\platalgk$ and $\Xk(T)$.  To move into the quantum setting, we fix a cyclic group $G$ and introduce notations for the scalar extensions 
\[
\qPlatnk=\Platnk\otimes \Z[G] \ \text{ and } \ \qplatalgk=\platalgk \wedge \Sq,
\]
and similarly $\qCCKnk(T)$ and $\qXk(T)$.  The (quantum) annular Khovanov invariants meanwhile carry an \emph{annular grading}.  We use the notation $\qaKC(L;w)$ for the quantum annular complex of $L$ in annular grading $w$, and similarly for $\qaKh(L;w)$, $\qaX(L;w)$, and when $G$ is finite, $\qakwX(L;w)$.

\begin{figure}
\[\begin{tikzpicture}[xscale=2mm, yscale=-.3mm, baseline=(current bounding box.center)]
    \node (tang) at (0,0) {\{$(n,n)$-tangle $T$\}};
    \node (closure) at (2,0) {\{annular link $L$\}};
    \node (spec bimod) at (0,2) {$\qXk(T)$};
    \node[red] (THH) at (.8,2) {$\qTHH(\qplatalgk;\qXk(T))$};
    \node[red] (ann spec) at (1.65,1.5) {$\qaX(L;n-2k)$};
    \node (AKW) at (2.35,2.5) {$\qakwX(L;n-2k)$};
    \node (bimod) at (0,4) {$\qCCKnk(T)$};
    \node (HC) at (.8,4) {$\qch_*(\qPlatnk;\qCCKnk(T))$};
    \node (ann cx) at (2,4) {$\qaKC(L;n-2k)$};
    \node (hom) at (1,6) {$\qaKh(L;n-2k) \cong H_*(\qaKC(L;n-2k))$};
    
    \draw[|->] (tang)--(closure) node[midway,above]{annular closure};
    \draw[|->] (tang)--(spec bimod) node[midway,left]{\cite{LLS_CK_Spectra} $\wedge\, \Sq$};
    \draw[|->] (spec bimod)--(bimod) node[midway,left]{$\chainsfunc$};
    \draw[|->,thick,red] (closure)--(ann spec);
    \draw[|->] (closure)--(AKW) node[midway,right]{$\,$\cite{AKW}};
    \draw[|->] (ann spec)--(ann cx) node[midway,left]{$\chainsfunc$};
    \draw[|->] (AKW)--(ann cx) node[midway,right]{$\,\chainsfunc$};
    \draw[|->] (bimod)--(HC) node[midway,below]{$\qch_*$};
    \draw[->,red,thick] (HC)--(ann cx) node[midway,above]{quasi-iso} node[midway,below]{$\qXi$};
    \draw[|->,red] (spec bimod)--(THH) node[midway,above]{$\qTHH$};
    \draw[<->,red,thick] (THH)--(ann spec) node[midway,above]{$\simeq$};
    \draw[|->] (THH)--(HC) node[midway,left]{$\chainsfunc$};
    \draw[|->] (bimod) to[out=90,in=-135] node[midway,left]{$\qhh_*$} (hom);
    \draw[|->] (HC)--(hom.165) node[midway,left]{$H_*$};
    \draw[|->] (ann cx) to[out=90,in=-45] node[midway,right]{$H_*$} (hom);
    \draw[<->,red,thick] (ann spec)--(AKW) node[midway,red, above]{$\simeq$};
\end{tikzpicture}\]
\caption{A schematic of various Khovanov invariants.  The first column consists of scalar extended bimodule invariants associated to tangles.  The second column consists of quantum Hochschild constructions.  The third column consists of annular link invariants.  We use $\longmapsto$ to denote an assignment of one object to another (such as taking homology of a chain complex), $\longrightarrow$ to denote an actual morphism between objects, and $\overset{\simeq}{\longleftrightarrow}$ to denote a zig-zag of weak equivalences. The portions in red indicate new constructions in this paper.}
\label{fig:Khovanov qTHH schematic}
\end{figure}
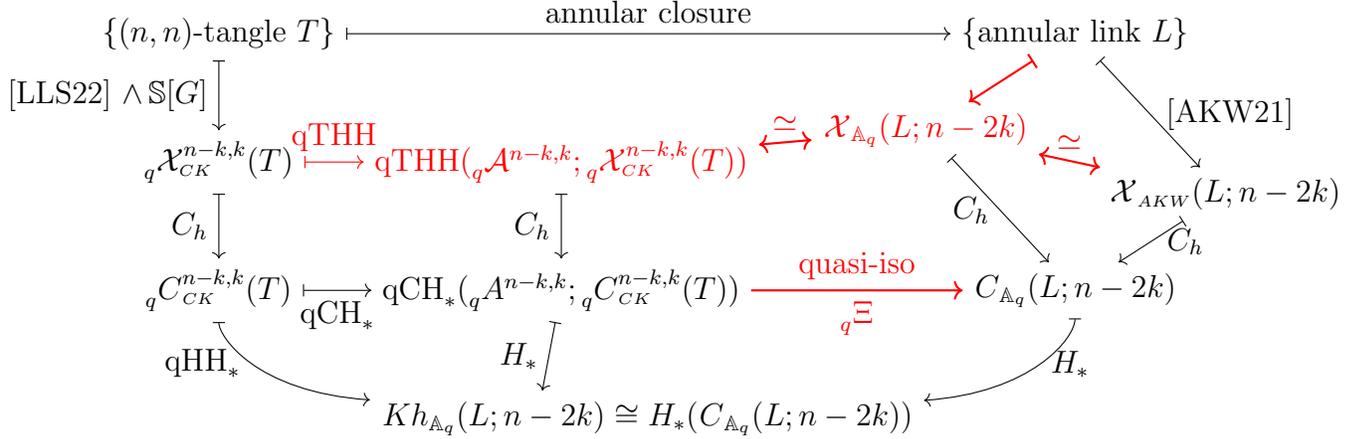

Our overall goal then is to complete Figure \ref{fig:Khovanov qTHH schematic}.  In this figure we use the notation $\qch$ for the quantum Hochschild complex, and $\chainsfunc$ for the \emph{chains functor} described in \cite[Section 2.7]{LLS_Khovanov_Spectra}.  Our first task will be to construct the chain map $\qXi$ which induces the isomorphism \eqref{eq:total cx vs termise iso}; this will be done in Section \ref{sec:explicitly constructing the quasi-iso}.  Our main task from there will be to lift $\qXi$ to a weak equivalence of spectra.

The process of lifting $\qXi$ will involve a number of categories which we summarize below:

\begin{itemize}
    \item The \emph{cube category} $\cube^N =\{0\to 1\}^N$, where $N$ is the number of crossings of $T$.
    \item The \emph{augmented Hochschild-Mitchell category} $\HMshape^\omega(\matchingsk)$, which tracks the semi-simplicial shape of $\qch(\qplatalgk;\qCCKnk(T))$, together with one extra terminal object $\omega$.  See Section \ref{sec:qTHH via shape cats} and \cite[Section 5]{LLS_CK_Spectra} where $\HMshape^\omega(\matchingsk)$ is denoted $\mathscr{C}^T$, with terminal object $T$.
    \item The \emph{quantum divided cobordism category of the annulus} $\divcobq$, a 2-category of cobordisms in the thickened annulus with some extra information (see Section \ref{sec:defining qdivcob}).
    \item The \emph{free $G$-Burnside category} $\Burn_G$ of sets with free $G$-action, acting as a 2-categorical lift of free $\Z[G]$-modules (see Section \ref{sec:Burnside}).
    \item The category $\Spectra_G$ of \emph{symmetric spectra with $G$-action} (see Section \ref{sec:defining qTHH}).
\end{itemize}

The construction can then be summarized by the schematic of Equation \eqref{eq:qHM composition via divcobq}:
\begin{equation}\label{eq:qHM composition via divcobq}
\begin{tikzcd}[row sep=tiny]
\big[(\cube^N \times \HMshape^\omega(\matchingsk))' \ar[r,"\cobfunc_T^{H^\omega}"] &
\divcobq' \ar[r,"\qVHKK'"] &
\Burn_G \big] & \\
& \Big\Downarrow & & \\
(\cube^N \times \HMshape^\omega(\matchingsk))'
\ar[rr, "\scext{\func{G}'}"] & &
\Burn_G \ar[r,"K_G"] &
\Spectra_G
\end{tikzcd}
\end{equation}
    
Roughly speaking, the tangle $T$ prescribes the topological data for the functor $\cobfunc_T^{H^\omega}$ (see Section \ref{sec:func HMshape to qdivcob}), while the functor $\qVHKK$ (see Section \ref{sec:qVHKKlax}) translates topological data in $\divcobq$ into combinatorial data in $\Burn_G$.  The prime notation $(-)'$ indicates \emph{canonical thickening} (see Section \ref{sec:Groupoid enrichments}).  The composition $\qVHKK'\circ \cobfunc_T^{H^\omega}$ is then modified to reflect that platform algebras and bimodules are subquotients (see Section \ref{sec:modifying qVHKK}) to give rise to a new functor $\scext{\func{G}'}$, which can then be composed with a further functor $K_G$ (based upon Elmendorf-Mandell $K$-theory \cite{Elmendorf-Mandell}) to arrive at spectra with $G$-action.  From here, a rectification result (see Section \ref{sec:Groupoid enrichments} and Appendix \ref{sec:simplicial stuff}) and some formal arguments are enough to complete the proof.

For the reader's convenience, we provide a summary of notation in Appendix \ref{sec:summary of notation}.

\subsection*{Acknowledgments} The authors would like to thank Anna Marie Bohmann, David Chan, Matt Hogancamp, Mikhail Khovanov, Slava Krushkal,  Robert Lipshitz, Cary Malkiewich, and Lucy Yang  for helpful discussions. The first author was supported by NSF grants DMS-2105467 and RTG DMS-1839968. The second author was supported by NSF grants DMS-2104233 and DMS-2404932. This material is also in part based on work supported by the National
Science Foundation under DMS-1928930, while the second and third authors were in residence at the Simons Laufer
Mathematical Sciences Institute (previously known as MSRI) in Berkeley, California, during the
Fall 2022 semester.

\section{Algebraic background}
In this section we review the definition of Hochschild homology, as well as its recently defined variant, \emph{quantum} Hochschild homology from \cite{BPW}. We then discuss the extensions of these theories to linear categories. 

\subsection{Quantum Hochschild homology} \label{sec:qHH}
In \cite{BPW}, Beliakova, Putyra, and Wehrli define a notion of quantum Hochschild homology (qHH). Before we recall their definition of qHH, we first recall the definition of classical Hochschild homology. 

\begin{definition}
\label{def:HH cx}
Let $R$ be a commutative ring, $A$ an $R$-algebra, and $M$ an $(A,A)$-bimodule. The \emph{Hochschild homology} of $A$ with coefficients in $M$, denoted $\hh_*(A;M)$, is the homology of the Hochschild complex $\ch_n(A;M ) = M \o_{R} A^{\o_{R} n}$ with differential $\partial = \sum_{i=0}^n (-1)^id_i$, where
\begin{equation}
\label{eq:HH differential}
d_i(m \otimes a_1 \otimes \ldots \otimes a_n) = \left\{\begin{array}{ll} ma_1 \otimes a_2 \ldots \otimes a_n &: i=0 \\
m\otimes a_1 \ldots \otimes a_ia_{i+1} \otimes \ldots \otimes a_n &: 0< i<n \\ a_nm\otimes a_1 \otimes \ldots \otimes a_{n-1}  &: i=n.
\end{array}\right.
\end{equation}
\end{definition}

 \begin{definition}
\label{def:qHH cx}
 Let $R$ be a commutative ring, $A$ a graded $R$-algebra, and $q\in R$ a unit. For a graded $A$-bimodule $M$,  the \emph{quantum Hochschild complex} of $A$ with coefficients in $M$ and parameter $q$ is 
 \[
 \qch_n(A;M) = M \otimes_R A^{\otimes_R n}.
 \]
The differential $\partial: \qch_n(A;M) \to \qch_{n-1}(A;M)$ is given by the alternating sum $\partial = \sum_{i=0}^n (-1)^id_i$
where 
\begin{equation}
\label{eq:qHH differential}
d_i(m \otimes a_1 \otimes \ldots \otimes a_n) = \left\{\begin{array}{ll} ma_1 \otimes a_2 \ldots \otimes a_n &: i=0 \\
m\otimes a_1 \ldots \otimes a_ia_{i+1} \otimes \ldots \otimes a_n &: 0< i<n \\ q^{-|a_n|}a_nm\otimes a_1 \otimes \ldots \otimes a_{n-1}  &: i=n.
\end{array}\right.
\end{equation}
We use $\vert a \vert$ to denote the degree of a homogeneous element $a\in A$. The quantum Hochschild homology of $A$ with coefficients in $M$ and parameter $q$ is defined to be the homology of this chain complex,
\[
\qhh_k(A;M) := H_k(\qch_{*}(A;M)).
\]
\end{definition}

\begin{remark}
   It is sometimes convenient to repackage the above notions as simplicial objects. Specifically, the simplicial $R$-module $\ch(A;M)_\bullet$ has $n$-simplices $\ch(A;M)_n = M \o_R A^{\o_R n}$, face maps given by $d_i$ as in \eqref{eq:HH differential}, and degeneracy maps given by $s_i = \id^{\o (i+1)} \o \eta \o \id^{\o (n-i)}$, where $\eta : R\to A$ is the unit map. In the graded setting, the simplicial $R$-module $\qch(A;M)_\bullet$ is defined analogously, except using \eqref{eq:qHH differential} for the face maps. 
\end{remark}
\begin{remark}\label{rmk:qHH as relative HH of twist}
As noted by Lipshitz in \cite{Lipshitz}, quantum Hochschild homology can also be viewed as ordinary Hochschild homology with coefficients in a twisted bimodule. Let $f_q: A \to A$ denote the $R$-algebra homomorphism defined on homogeneous elements by
\[
f_q(a) = q^{-|a|}a.
\]
From the $A$-bimodule $M$, we define an $A$-bimodule ${}_{q}M$ with the left module action twisted by $f_q$. In other words, ${}_{q}M$ has the same right $A$-action as $M$, and if $\lambda$ is the left action map for $M$, the left action map ${}_{q}\lambda$ for ${}_{q}M$ is given by
\[
\begin{tikzcd}
    A \otimes_R M \ar[d,"f_q \otimes 1"'] \ar[dr,"{}_{q}\lambda"] & \\ A \otimes_R M \ar[r,"\lambda"] & M. 
\end{tikzcd}
\]
Then $\qch_*(A;M) = \ch_*(A;{}_{q}M)$ and $\qhh_*(A;M) = \hh_*(A; {}_{q}M)$. 
\end{remark}

Suppose $A$ and $B$ are graded $R$-algebras, $M$ is a graded $(A,B)$-bimodule, and $N$ is a graded $(B,A)$-bimodule. Suppose further that $N$ and $M$ are finitely generated and projective as left modules. As discussed in \cite[Section 3.8.5]{BPW} there is a chain homotopy equivalence 
\[
\qch_*(A; M\o_B N) \simeq \qch_*(B; N\o_A M).
\]

\begin{remark}
\label{rem:trace property of qHH on chain level}
There is an explicit isomorphism $\qhh_0(A; M\o_B N) \rar{\sim} \qhh_0(B; N\o_A M)$ on the level of homology given by $[x\o y] = q^{-\lr{y}}[y\o x]$. However, the map on the chain level is not immediately described; note that $M\o_B N \to  N\o_A M$ given by $x\o y\mapsto y\o x$ is not necessarily well-defined. 
\end{remark}

If $M$ is a complex of $(A,A)$-bimodules with differential $\partial_M$, then $\hh_n(A;M)$ is defined to be the homology of \[
\ch_n(A;M) = \bigoplus_{i+j=n} M_i \o_R A^{\o_R j}
\]
with differential given on the summand $M_i \o_R A^{\o_R j}$ by $\partial_M + (-1)^i\partial$, where $\partial :  M_i \o_R A^{\o_R j} \to  M_i \o_R A^{\o_R (j-1)}$ is as in Definition \ref{def:HH cx}. Similarly, if $q\in R$ is a choice of unit, $A$ is a graded $R$-algebra, and $M$ is a complex of graded $(A,A)$-bimodules with degree-preserving differential, then $\qhh(A;M)$ is the homology of the complex with the same terms as above but now with $\partial$ as in Definition \ref{def:qHH cx}. 

\begin{remark}
\label{rmk:ground ring is trivially graded}
    The ground ring $R$ will always be trivially graded (that is, concentrated in degree zero), both in this algebraic setting and also later in the spectral setting.  While working over a nontrivially graded ground ring (such as the ring of polynomials $\Z[x_1, \ldots, x_n]$ where each $x_i$ is in degree $2$, or its subring of symmetric polynomials) is common in link homology, it is not clear how to extend $\qhh$ to the case where the ground ring is nontrivially graded.
\end{remark}

\subsection{Algebras with idempotents and linear categories}
\label{sec:linear cats}

In many cases of interest arising from geometric topology and knot theory, we have algebras (and bimodules over them) which contain extra structure.  Algebras such as the arc algebras of Khovanov \cite{Kh-ring_H} and platform algebras of Chen-Khovanov and Stroppel \cite{Chen-Khovanov,Stroppel} (which we will describe in Section \ref{sec:Platform algebras and bimodules}), as well as the strands algebras and pong algebras appearing in bordered Heegaard Floer theories (see for instance \cite{LOT_bordered,OS_pong,zarev_suturedFloer,MMW_strands}), and many others, come equipped with a decomposition via idempotents (see below for a precise statement), which makes them amenable to study using \emph{linear categories}.  This point of view in turn aids in the study of their Hochschild constructions  as well as their `topological lifts' via \emph{spectral categories} (see Section \ref{sec:qTHH for spectral cats}).

Consider a commutative ring $R$ and an $R$-algebra $A$ with a finite set $S_A:=\{a_1,\dots,a_m\}\subset A$ of mutually orthogonal idempotents that sum to the identity:
\begin{equation*}
    a_i\cdot a_j = \begin{cases}
        a_i & : i=j\\
        0 & : i\neq j
        \end{cases},
    \quad\quad
    \sum_{i=1}^m a_i = 1_A.
\end{equation*}
We call $S_A$ a \emph{complete set of idempotents} for $A$. The set $S_A$ endows $A$ with a direct sum decomposition into $R$-submodules
\[A=\bigoplus_{a_i,a_j\in S_A} a_iAa_j.\]
Suppose $B$ is another $R$-algebra with a complete set of idempotents $S_B$, and let $\bimod AB$ denote the category of $(A,B)$-bimodules. Any $M\in \bimod{A}{B}$ has a similar direct sum decomposition
\[M=\bigoplus_{a_i\in S_A,\, b_j\in S_B} a_iMb_j.\]
Such algebras and bimodules can be recast in a more categorical language which we review here as a definition.
Let $\Rmod$ denote the category of $R$-modules. 

\begin{definition}\label{def:Alg as lin Cat}
An algebra $A$ over $R$ with a complete set of idempotents $S_A$ is equivalent to the data of an \emph{$R$-linear category} $\cat{A}$ with
\begin{itemize}
    \item a finite object set $S_A$ recording the idempotents, 
    \item morphism sets $\cat{A}(a_i,a_j)=\Hom_{\cat{A}}(a_i,a_j):=a_i A a_j\in \Rmod$ encoding the idempotent decomposition of $A$, 
    \item special morphisms $R\hookrightarrow\cat{A}(a,a)$ encoding the unital structure, and
    \item composition encoding multiplication in $A$.
\end{itemize}

In this language, a bimodule $M\in\bimod AB$ is equivalent to the data of an enriched\footnote{Here \emph{enriched} means that the assignment on Hom-sets is $R$-linear.} functor $\opcat{A}\times \cat{B} \xrightarrow{\cat{M}} \Rmod$ with 
\begin{itemize}
    \item $\cat{M}(a_i,b_j):= a_i M b_j\in\Rmod$ encoding the idempotent decomposition of $M$, and
    \item values on Hom-sets encoding the bimodule actions as in
    \[\Hom_{\opcat{A}\times \cat{B}}((a_i,b_k),(a_j,b_\ell)) \xrightarrow{\cat{M}} \Hom_{\Rmod}(\cat{M}(a_i,b_k), \cat{M}(a_j,b_\ell))\]
    being reinterpreted as a map
    \[\cat{A}(a_j,a_i) \otimes_R \cat{M}(a_i,b_k) \otimes_R \cat{B}(b_k,b_\ell) \rightarrow \cat{M}(a_j,b_\ell).\]
\end{itemize}
We shall refer to such a functor $\cat{M}$ as an $(\cat{A},\cat{B})$-bimodule, written as $\cat{M}\in\bimod{\cat{A}}{\cat{B}}$.

Similarly, a $\Z$-graded algebra is equivalent to a \emph{graded $R$-linear category}, meaning that Hom-sets are in the category of graded $R$-modules and homogeneous $R$-linear maps (of any degree), which we denote by $\gRmod$; a $\Z$-graded bimodule is equivalent to a functor $\opcat{A}\times \cat{B} \rightarrow \gRmod$ whose assignment on Hom-sets is $R$-linear and degree-preserving. Analogous equivalences hold for bigraded ($\Z^2$-graded) bimodules and algebras, where we denote by $\ggRmod$ the category of $\Z^2$-graded $R$-modules. 

Finally, a chain complex of $(A,B)$-bimodules is equivalent to a functor $\opcat{A}\times\cat{B} \rightarrow \Kom(\Rmod)$ where $\Kom(\Rmod)$ denotes the category of chain complexes of $R$-modules, and similarly for a complex of graded or bigraded bimodules.
\end{definition}

Linear categories have their own notions of classical and quantum Hochschild homology, which we review here.

\begin{definition}\label{def:HM complex for lin cat and bimod}
Let $\cat{A}$ be an $R$-linear category with finite object set $S_A$, and fix an $(\cat{A},\cat{A})$-bimodule $\opcat{A}\times\cat{A}\rar{\cat{M}}\Rmod$ as in Definition \ref{def:Alg as lin Cat}.  The \emph{Hochschild-Mitchell complex of $\cat{A}$ with coefficients in $\cat{M}$} consists of chain groups
\[\ch_n(\cat{A};\cat{M}):= \bigoplus_{a_0,\dots,a_n\in S_A} \cat{M}(a_n,a_0)\otimes_R \cat{A}(a_0,a_1) \otimes_R \cdots \otimes_R \cat{A}(a_{n-1},a_n)\]
with differential $\partial=\sum (-1)^i d_i$, where
\begin{equation}
\label{eq:HM differential}
d_i(m \otimes \alpha_1 \otimes \ldots \otimes \alpha_n) = \left\{\begin{array}{ll} m\alpha_1 \otimes \alpha_2 \ldots \otimes \alpha_n &: i=0 \\
m\otimes \alpha_1 \ldots \otimes \alpha_i\alpha_{i+1} \otimes \ldots \otimes \alpha_n &: 0< i<n \\ \alpha_n m\otimes \alpha_1 \otimes \ldots \otimes \alpha_{n-1}  &: i=n,
\end{array}\right.
\end{equation}
and each of the products here is defined via the interpretations of Definition \ref{def:Alg as lin Cat}.

Then the \emph{Hochschild homology of $\cat{A}$ with coefficients in $\cat{M}$}, denoted $\hh_*(\cat{A};\cat{M})$, is defined to be the homology of the Hochschild-Mitchell complex
\[\hh_n(\cat{A};\cat{M}):=H_n(\ch_*(\cat{A};\cat{M})).\]

\end{definition}

\begin{definition}
\label{def:qHM complex for lin cat and bimod}
Let $R$ be a commutative ring with a chosen unit $q\in R$. Let $\cat{A}$ be a graded $R$-linear category with finite object set $S_A$, and fix a graded bimodule $(\cat{A},\cat{A}) \rar{\cat{M}} \gRmod $. The \emph{quantum Hochschild-Mitchell complex of $\cat{A}$ with coefficients in $\cat{M}$} consists of the same chain groups as in Definition \ref{def:HM complex for lin cat and bimod}, except where the last face map $d_n$ in the differential is given by 
\[
d_n(m\o \alpha_1 \o \cdots \o \alpha_n) = q^{-\lr{\alpha_n}}\alpha_n m\otimes \alpha_1 \otimes \ldots \otimes \alpha_{n-1} 
\]
Its homology is denoted $\qhh_*(\cat{A};\cat{M})$.  
\end{definition}

More generally if $\opcat{A}\times\cat{A}\rar{\cat{M}}\Kom(\Rmod)$ is a complex of bimodules with differential $\partial_M$, then we define $HH_*(\cat{A};\cat{M})$ as the homology of the complex
\[
\ch_n(\cat{A};\cat{M}) = \bigoplus_{i+j=n} \ch_j(\cat{A};\cat{M}_i),
\]
where now $\cat{M}_i$ is a bimodule and by $\ch_j(\cat{A};\cat{M}_i)$ we mean the module in Definition \ref{def:HM complex for lin cat and bimod}, with respect to the differential $\partial_M + (-1)^i\partial$, where $\partial :\ch_j(\cat{A};\cat{M}_i) \to \ch_{j-1}(\cat{A};\cat{M}_i)$ is as in Definition \ref{def:HM complex for lin cat and bimod}. If $\cat{A}$ is a graded $R$-linear category and $\cat{M}$ is a complex of graded $(A,A)$-bimodules, we define $\qhh_*(\cat{A};\cat{M})$ completely analogously, except using the differential in Definition \ref{def:qHM complex for lin cat and bimod}.

Suppose we are viewing the linear category $\cat{A}$ as equivalent to the algebra with idempotent decomposition $A$, with bimodule functor $\cat{M}\in\bimod{\cat{A}}{\cat{A}}$ equivalent to bimodule $M\in\bimod{A}{A}$ (or the graded versions). Then the following shows that the Hochschild-Mitchell constructions are equivalent to the Hochschild constructions of Section \ref{sec:qHH}, justifying the use of the same notation in this setting.

\begin{proposition}
Given the equivalences set out above (utilizing Definition \ref{def:Alg as lin Cat}), we have
\[\hh_*(A;M) \cong \hh_*(\cat{A};\cat{M}), \quad \qhh_*(A;M)\cong \qhh_*(\cat{A};\cat{M}),\]
where the symbol $\cong$ indicates isomorphism in whichever category our coefficients are in ($\Rmod$ or $\gRmod$). 
\end{proposition}
\begin{proof}
This is a classical result.  See for instance \cite{CR_HH_reference}.
\end{proof}

\section{Quantum topological Hochschild homology}

Classical Hochschild homology of associative algebras has a topological analogue: topological Hochschild homology (THH) for ring spectra. In this section we develop a topological analogue of quantum Hochschild homology, called quantum topological Hochschild homology (qTHH). As noted in Section \ref{sec:qHH}, the ground ring for the algebraic theory of qHH was a commutative ring $R$ with an invertible element $q\in R$. For our lift to spectra, we choose to focus on the simplest such rings, namely the group rings $\Z[G]$ where $G$ is a cyclic group with designated generator $q\in G$.
The ground ring for our topological theory of qTHH  will thus be the spherical group ring $\mathbb{S}[G]$. We begin this section by recalling the necessary spectral algebra definitions, and then define quantum topological Hochschild homology.   

\subsection{Spectral algebra}\label{sec:qTHH input definitions}

Recall that for a ring spectrum $R$, one can define $R$-module spectra and $R$-algebra spectra, which are topological analogues of the classical algebraic notions. We will begin by recalling these definitions. Later in this section we will consider the case where $R$ is the spherical group ring $\Sq$. For a more extensive treatment of spectral algebra, see, for example, \cite{ekmm}.

\begin{definition}\label{def:ring spectrum}
A \emph{ring spectrum} is a monoid in the symmetric monoidal category of spectra $\Spectra$.\footnote{We take $\Spectra$ to be Hovey-Shipley-Smith's category of symmetric spectra \cite{symmetric_spectra}.} In other words, a ring spectrum $R$ is a spectrum together with multiplication and unit maps \[\mu: R\wedge R \to R \hspace{.4cm} \textup{and} \hspace{.4cm} \eta: \S \to R
\]
such that the multiplication is unital and associative. The ring spectrum $R$ is \emph{commutative} if $\mu = \mu \circ \tau$, where $\tau : R\wedge R \to R\wedge R$ swaps the two factors.
\end{definition}

\begin{definition}\label{def:spectral Sq module} 
Let $R$ be a ring spectrum. A \emph{left} $R$-\emph{module} is a spectrum $M$ with a map
\[
\lambda: R \wedge M \to M
\]
which is associative and unital.  A \emph{map of left $R$-modules} is a map of spectra $M\xrightarrow{f}N$ such that the following diagram commutes:
\[\begin{tikzcd}
R\wedge M \ar[r,"\id\wedge f"] \ar[d,"\lambda_M"] & R\wedge N \ar[d,"\lambda_N"]\\
M \ar[r,"f"] & N
\end{tikzcd}.\]
A \emph{graded ($\Z$-graded) left $R$-module} is of the form $M= \bigvee_{m\in \Z} M_m$ with each $M_m$ a left $R$-module. A \emph{map of $\Z$-graded left $R$-modules of degree $d$} is a sequence of maps of left $R$-modules  $M_m\rar{f_m}N_{m+d}$.  
We will also consider \emph{bigraded} ($\Z^2$-graded) $R$-modules, allowing maps of bigraded degree $(c,d)$.
\end{definition}

The definition of a \emph{(graded or bigraded) right $R$-module} is analogous, where we use the notation $\rho : M\wedge R \to M$ for the right action map. We will also need to consider bimodules in this spectral setting.

\begin{definition}\label{def:bimodules}
For ring spectra $R$ and $R'$, an \emph{$(R, R')$-bimodule} is a left $R$- and right $R'$-module $M$, with action maps 
\[
\lambda: R \wedge M \to M \textup{   and  } \rho: M \wedge R' \to M
\]
such that the following diagram commutes:
\[
\begin{tikzcd}
    R \wedge M \wedge R' \ar[r,"\lambda \wedge \id"] \ar[d,"\id \wedge \rho"] & M\wedge R' \ar[d,"\rho"] \\
R \wedge M \ar[r,"\lambda"] & M.
\end{tikzcd}
\]
Graded (and bigraded) $(R,R')$-bimodules are defined analogously.
\end{definition}

Observe that for a \emph{commutative} ring spectrum $R$, a left $R$-module $M$ also has a right $R$-module structure via
\[
\lambda \circ \tau: M\wedge R  \to M,
\]
such that $M$ is an $(R, R)$-bimodule.
\begin{definition}\label{def:relative smash product}
Let $R$ be a ring spectrum, let $M$ be a right $R$-module, and let $N$ a left $R$-module. The \emph{relative smash product}, $M\wedge_R N$, is the coequalizer in spectra:
\[\begin{tikzcd}
M \wedge R \wedge N \ar[r,shift left=.75ex,"\rho \wedge \id"]
  \ar[r,shift right=.75ex,"\id \wedge \lambda"']
&
M \wedge N \ar[r] & M\wedge_{R}N.
\end{tikzcd}\]
Here $\rho$ and $\lambda$ denote the actions of $R$ on $M$ and $N$ respectively. If $R$ is a commutative ring spectrum, $M\wedge_{R}N$ has a canonical $R$-module structure induced from the $R$-module structure of $M$, or equivalently $N$.  If $M$ and $N$ are graded, then $M\wedge_R N$ is also graded via
\[
(M\wedge_{R}N)_m = \bigvee_{i+j=m} M_i \wedge_{R} N_j, 
\]
and similarly in the bigraded case.
\end{definition}

For a commutative ring spectrum $R$, and $\bullet\in\{\varnothing,gr,gg\}$, we let $\Rmod[R]^\bullet$ denote the symmetric monoidal category of $R$-modules and maps between them with relevant grading (see Definition \ref{def:spectral Sq module}); the monoidal structure is given by the relative smash product (see \cite[Theorem 7.1]{ekmm} and \cite[IV.1]{schwedesymmetric}).  The case where $R$ is the sphere spectrum $\S$ recovers the category $\Spectra^\bullet$ of (possibly graded or bigraded) spectra.  There is a forgetful functor $(-)^{un} : \Rmod[R]^\bullet \to \Spectra^\bullet$ sending an object $M$ to its \emph{underlying spectrum}, where the $R$-module structure is forgotten.

\begin{definition}
    A map $f: M\to N$ in $\Rmod[R]^\bullet$ is a \emph{weak equivalence} if $f^{un}$ is a weak equivalence in $\Spectra^\bullet$.
\end{definition}

\begin{theorem}[Whitehead's Theorem]
\label{them:Whitehead thm}
A map of $R$-modules $f:M\to N$ which induces an isomorphism on homology is a weak equivalence of $R$-modules if both $M$ and $N$ are $n$-connected for some $n$ (and similarly for graded or bigraded $R$-modules).
\end{theorem}

In addition to considering modules in the spectral setting, we will also need the notions of spectral algebras and modules over them. We recall those definitions now. 

\begin{definition}\label{def:spectral algebra over Sq}
Let $R$ be a commutative ring spectrum. An \emph{$R$-algebra} is a monoid in the symmetric monoidal category of $R$-modules. In other words, an $R$-algebra $A$ is an $R$-module with multiplication and unit maps \[\mu: A\wedge_{R} A \to A \hspace{.4cm} \textup {and} \hspace{.4cm} \eta: R \to A
\]
of $R$-modules, satisfying unitality and associativity.  A \emph{graded} (respectively \emph{bigraded}) $R$-algebra is one whose underlying $R$-module is graded (respectively bigraded), with multiplication and unit maps being degree zero (respectively bidegree $(0,0)$) maps. Explicitly, $A = \bigvee_j A_m$ where each $A_m$ is an $R$-module, and $\mu$ and $\eta$ are assembled from maps 
\[
A_i \wedge_R A_j \to A_{i+j} \hspace{.4cm} \textup {and} \hspace{.4cm} R \to A_0.
\]
\end{definition}

Observe that when $R=\S$, a (graded) $R$-algebra $A$ is simply a (graded) ring spectrum. We can also consider modules over $R$-algebras. 
\begin{definition}\label{def:modules over an algebra}
For $R$ a commutative ring spectrum, and $A$ an $R$-algebra, a \emph{left $A$-module} is an $R$-module $M$ with action map
\[A\wedge_{R} M \xrightarrow{\lambda} M \]
which is associative and unital. Similarly, a \emph{right $A$-module} is an $R$-module $N$ with action map
\[
N\wedge_{R} A \xrightarrow{\rho} N
\]
which is associative and unital. For $R$-algebras $A$ and $B$, the definition of an $(A,B)$-bimodule is analogous to Definition \ref{def:bimodules} above.  Graded (or bigraded) $A$-modules require the action maps to be degree zero (or bidegree $(0,0)$) maps, and similarly for graded (or bigraded) $(A,B)$-bimodules.
\end{definition}

\begin{remark}
    One could allow the ground ring $R$ to itself be graded and consider graded algebras over $R$ and graded bimodules over these algebras. In our case, we view $R$ as trivially graded in light of Remark \ref{rmk:ground ring is trivially graded}. 
\end{remark}

\subsection{Defining quantum topological Hochschild homology}\label{sec:defining qTHH}

In this section we will give the definition of \emph{quantum topological Hochschild homology}, qTHH. We begin first by recalling the definition of $R$-relative topological Hochschild homology with coefficients.
\begin{definition}
Let $R$ be a commutative ring spectrum and $A$ an $R$-algebra, with product and unit maps 
\[
\mu: A \wedge_R A \to A \hspace{.4cm}\textup{   and   } \hspace{.4cm} \eta: R\to A.
\]
Let $M$ be an $(A,A)$-bimodule with left and right action maps
\[
\lambda: A \wedge_R M \to M \hspace{.4cm}\textup{  and  }\hspace{.4cm} \rho: M \wedge_R A \to M, 
\]
respectively. The simplicial $R$-module $\THH_R(A;M)_{\bullet}$ has $n$-simplices
\[
\THH_R(A;M)_n = M \wedge_R A^{\wedge_R n},
\]
where $A^{\wedge_R n}$ denotes the $n$-fold $R$-relative smash product of $A$. The face maps are given by:
\[
d_i = \left\{\begin{array}{ll}
\rho \wedge \id^{\wedge n} &: i =0 \\
\id^{\wedge i} \wedge \mu \wedge \id^{\wedge (n-i-1)} &: 0< i<n \\
(\lambda \wedge \id^{\wedge (n-1)}) \circ \tau & i=n,
\end{array}\right.
\]
where $\tau$ rotates the last factor to the front, and the degeneracy maps are
\[
s_i =  \id^{\wedge (i+1)} \wedge \eta \wedge \id^{\wedge (n-i)}.
\]
This simplicial object is a cyclic bar construction with coefficients, in the sense of \cite{BHM}. The \emph{$R$-relative THH of $A$ with coefficients in $M$} is the geometric realization of this simplicial $R$-module:
\[
\THH_R(A;M) = |\THH_R(A;M)_{\bullet}|.
\]
\end{definition}
When $M$ is $A$, we denote $\THH_R(A;A)$ by simply $\THH_R(A)$. When $R= \S$, the sphere spectrum, the $R$-algebra $A$ is a ring spectrum and we omit the $R$ from the notation, denoting $\THH_{\S}(A;M)$ by $\THH(A;M).$

The input for quantum topological Hochschild homology will be a graded algebra over the spherical group ring $\Sq$ for a cyclic group $G$. The cyclic group $G$ can be finite or infinite, and we moreover fix a distinguished generator $q\in G$.  Before defining $\qTHH$, we first recall the definition of spherical group rings.

\begin{definition}\label{def:spherical group ring} Let $H$ be a topological group with multiplication map $H\times H \xrightarrow{m}H$ and unit $1\in H$. The \emph{spherical group ring} $\S[H]$ is a spectrum defined by 
\[
(\S[H])_n = H_+ \wedge S^n,
\]
where $H_+$ is the space $H$ with disjoint basepoint added. In other words, $\S[H]$ is the suspension spectrum $\Sigma^{\infty}_+ H$. The spherical group ring becomes a ring spectrum with multiplication maps defined by the compositions
\[
 H_+ \wedge S^n \wedge H_+ \wedge S^m \xrightarrow{\id \wedge \tau \wedge \id} (H \times H)_+ \wedge S^n \wedge S^m  \xrightarrow{m_+ \wedge \mu_{n,m}} H_+ \wedge S^{n+m},
\]
and unit map given by 
\[
1 \wedge (-): S^n \to H_+ \wedge S^n.
\]
\end{definition}

In general, $\pi_0$ of a spherical group ring recovers the corresponding classical group ring: 
\[\pi_0(\S[H]) \cong \Z[H].\]
In the case where $G$ is cyclic, the spherical group ring $\Sq$ can be identified as a wedge sum
\[\Sq:=\bigvee_{q^i\in G} q^i\S,\]
where the notation $q^i$ is used to distinguish different copies of $\S$. From this perspective, the multiplication $\Sq\wedge\Sq\xrightarrow{\mu} \Sq$ is given by the usual ring spectrum structure on the sphere spectrum via
\[q^i\S \wedge q^j\S \cong \S \rar{\id} q^{i+j} \S,\]
 and is commutative.  The unit map $\S\xrightarrow{\eta}\Sq$ is then $\S\xrightarrow{\id}q^0\S\hookrightarrow \Sq$.

Recall that in the classical algebra setting, for an abelian group $H$ the category of $\Z[H]$-modules and the functor category $\Func(BH, \Ab)$ are isomorphic, where $\Ab$ is the category of abelian groups and $BH$ is the category with one object $*$ and $\Hom_{BH}(*,*)=H$.  We have the following analogous results in the topological setting.

\begin{proposition}\label{prop:Sq module automorphism}
Let $G$ be a cyclic group with generator $q\in G$. An $\Sq$-module $X$ comes equipped with designated  automorphisms $X \xrightarrow{q_X^i} X$ for each $q^i\in G$ which satisfy:
\begin{itemize}
    \item the map $X\xrightarrow{q_X^0}X$ is the identity map,
    \item the maps satisfy $q_X^i\circ q_X^j = q_X^{i+j }$, and
    \item for any $\Sq$-module map $X\rar{f}Y$, we have $f\circ q_X^i = q_Y^i \circ f$. 
\end{itemize}
Further, the functor $\Rmod[\Sq] \rar{(-)^{un}} \Func(BG, \Spectra)$ sending an $\Sq$-module to its underlying spectrum equipped with the above $G$-action via $q^i_X$ and sending a map of $\Sq$-modules to the map on underlying spectra is an isomorphism of categories. 

\end{proposition}
\begin{proof}
Define $q^i_X$ as the composition $
X \cong q^i \sphere \wedge X \hookrightarrow \Sq\wedge X \rar{\lambda} X$. The first and second bullet points follow from unitality and associativity, respectively. 
Suppose $X \rar{f} Y$ is a map of $\Sq$-modules, and consider the diagram
\[\begin{tikzcd}
X \ar[r,"\cong"] \ar[d,"f"] &
q^i\S \wedge X \ar[r,hook] \ar[d,"\id\wedge f"] &
\Sq \wedge X \ar[r,"\lambda"] \ar[d,"\id\wedge f"] &
X \ar[d,"f"] \\
Y \ar[r,"\cong"] &
q^i\S \wedge Y \ar[r,hook] &
\Sq \wedge Y \ar[r,"\lambda"] &
Y
\end{tikzcd}\]
where the left and middle rectangles commute trivially,  and the right rectangle commutes since $f$ is an $\Sq$-module map. This verifies the third bullet point. 

The inverse of the functor $(-)^{un}$ is given as follows. Given $F\in \Func(BG, \Spectra)$, define $X'\in \Rmod[\Sq]$ by $X' = F(*)$ with left $\Sq$-action map 
\[
\Sq \wedge X' \cong \bigvee_{q^i \in G} q^i \sphere \wedge X' \cong \bigvee_{q^i \in G} q^i X' \to X'
\]
where $\bigvee_{q^i \in G} q^i X' \to X'$ is given  by the automorphism $F(q^i) : X' \to X'$ on the wedge summand corresponding to $q^i\in G$. That this defines an $\Sq$-module and provides an inverse to $(-)^{un}$ is straightforward to verify.
\end{proof}
  We note that Proposition \ref{prop:Sq module automorphism} works just as well to give an isomorphism $\Rmod[\Sq]^\bullet \cong \Func(BG, \Spectra^\bullet)$ for $\bullet \in \{gr, gg\}$.

\begin{definition}
\label{def:S_G}
  For a cyclic group $G$ and $\bullet \in \{\varnothing, gr, gg\}$, let $\Spectra^\bullet_G$ denote the category $\Func(BG, \Spectra^\bullet)$ of \emph{spectra with $G$-action}.  In light of Proposition \ref{prop:Sq module automorphism}, we will identify $\Spectra^\bullet_G$ with 
  $\Rmod[\Sq]^\bullet$ and will sometimes specify an object of  $\Rmod[\Sq]^\bullet$ by giving a spectrum $X$ together with a group homomorphism $G \to \Aut_{\Spectra}(X)$.   
\end{definition}

These automorphisms will play the role of the multiplication by $q^{-\vert a_n\vert}$ in the last face map \eqref{eq:qTHH cyclic face map} described in the definition of $\qTHH$ given below.

\begin{definition}\label{def:chains functor}
Let $\chainsfunc:\Spectra\rightarrow\Kom(\Ab)$ denote the chains functor on spectra described in \cite[Definition 2.30]{LLS_Khovanov_Spectra}.
\end{definition}

We note that $\chainsfunc$ can be upgraded to a functor (with the same notation) 
\begin{equation}
\label{eq:chains func ZG modules}
\chainsfunc : \Spectra_G \to \Kom(\Rmod[{\Z[G]}]).
\end{equation}
One way to see this is to use Proposition \ref{prop:Sq module automorphism} to get a functor 
\[
\Spectra_G \cong \Func(BG, \Spectra) \rar{ C_h \circ (-)} \Func(BG, \Kom(\Ab)) \cong \Kom(\Rmod[{\Z[G]}]). 
\]
In particular, we do not modify the outputs of $\chainsfunc$ but rather note that if the inputs are in $\Spectra_G$ then the outputs can naturally be viewed as landing in $\Kom(\Rmod[{\Z[G]}])$.

For $G$ a cyclic group, the input for qTHH will be a graded $\Sq$-algebra $A=\bigvee_j A_j$, together with a graded $(A,A)$-bimodule $M=\bigvee_j M_j$. Although the group ring $\Sq$ can be given the structure of a $G$-graded spectrum, here we are considering it as a $\Z$-graded spectrum with trivial $\Z$-grading. Observe that in this graded case, the automorphism $A\xrightarrow{q}A$ described in Proposition \ref{prop:Sq module automorphism} decomposes over the wedge sum, and the map on each term will also be denoted $A_j\xrightarrow{q}A_j$.

\begin{definition}
\label{def:qTHH as a simplicial object}
Let $G$ be a cyclic group and $A$ be a graded $\Sq$-algebra with product and unit maps
\[
\mu: A \wedge_{\Sq} A \to A \textup{ and  } \eta: \Sq \to A.
\]
Let $M$ be a graded $(A,A)$-bimodule with left and right action maps
\[
\lambda: A \wedge_{\Sq} M \to M \hspace{.4cm}\textup{  and  }\hspace{.4cm} \rho: M \wedge_{\Sq} A \to M,
\]
respectively. The simplicial $\Sq$-module $\qTHH(A;M)_{\bullet}$ has $n$-simplices 
\[
\qTHH(A;M)_n= M \wedge_{\Sq} A^{\wedge_{\Sq} n}.
\]
As in the standard cyclic bar construction with coefficients, the degeneracy maps $s_i$ are given by 
\[
s_i = \Id^{(i+1)} \wedge \eta \wedge \Id^{(n-i)}, \textup{    for  } 0\leq i \leq n. 
\]
The first $n$ face maps on $\qTHH(X;M)_{n}$ are given by
\begin{equation}\label{eq:qTHH noncyclic face maps}
d_i = \left\{\begin{array}{ll}
\rho \wedge \Id^{\wedge n} &: i =0 \\
\Id^{\wedge i} \wedge \mu \wedge \Id^{\wedge (n-i-1)} &: 0< i<n
\end{array}\right.
\end{equation}
as usual for the cyclic bar construction with coefficients. Let $\beta_q$ denote the map that cyclically permutes the last factor of $A$ to the front, then decomposes as a wedge sum, 
\[
A \wedge_{\Sq} M\wedge_{\Sq} A^{\wedge_{\Sq} (n-1)} = \bigvee_i \left(A_i\wedge_{\Sq} M \wedge_{\Sq} A^{\wedge_{\Sq} (n-1)}\right),
\] and applies the automorphism $q^{-i}$ of Proposition \ref{prop:Sq module automorphism} to the first factor $A_i$:
\begin{align*}
\beta_q: M\wedge_{\Sq} A^{\wedge_{\Sq} (n-1)} \wedge_{\Sq} A & \xrightarrow{} \bigvee_i \left(A_i\wedge_{\Sq} M \wedge_{\Sq} A^{\wedge_{\Sq} (n-1)}\right) \\
& \xrightarrow{\bigvee_i  q^{-i} \wedge \id^{\wedge n} } \bigvee_i\left( A_i\wedge_{\Sq} M \wedge_{\Sq} A^{\wedge_{\Sq} (n-1)}\right). 
\end{align*}
The last face map, then, is given by:
\begin{equation}\label{eq:qTHH cyclic face map}
d_n = (\lambda \wedge \id^{\wedge(n-1)}) \circ \beta_q.
\end{equation}
We refer to the simplicial $\Sq$-module $\qTHH(A;M)_{\bullet}$ as the \emph{quantum cyclic bar construction} of $A$ with coefficients in $M$. 
\end{definition}

The verification that this is indeed a simplicial $\Sq$-module proceeds as usual for the cyclic bar construction.

\begin{remark} 
    There have recently been other generalizations of topological Hochschild homology using variations on the cyclic bar complex, such as the dihedral bar complex, or the twisted cyclic bar complex (see, e.g. \cite{DMPR21}, \cite{ABGHLM18}, \cite{AGHKK23}). While the quantum cyclic bar construction defined above looks somewhat similar to the twisted cyclic bar construction of \cite{ABGHLM18} used to define twisted topological Hochschild homology for $C_n$-ring spectra, the two constructions are distinct.
\end{remark}

\begin{definition}\label{def:qTHH}
Let $G$ be a cyclic group, $A$ a graded $\Sq$-algebra, and $M$ a graded $(A,A)$-bimodule. The \emph{quantum topological Hochschild homology of $A$ with coefficients in $M$} is given by the geometric realization in $\Sq$-module spectra
\[
\qTHH(A;M) = |\qTHH(A;M)_{\bullet}|.
\]
\end{definition}

\begin{remark}\label{rmk:qTHH is relative THH of twist}
We note that, analogous to the algebraic case described in Remark \ref{rmk:qHH as relative HH of twist}, quantum topological Hochschild homology can also be viewed as a special case of relative topological Hochschild homology with coefficients in a twisted bimodule. 
\end{remark}

Classically, Hochschild homology and topological Hochschild homology are related via a map called \emph{linearization}. For $R$ a $(-1)$-connected ring spectrum, the linearization map takes the form 
\[
\pi_n(\THH(R)) \to \hh_n(\pi_0(R)),
\]
and is an isomorphism when $n=0$. In particular, when $A$ is a ring and $HA$ its Eilenberg-MacLane spectrum, there is such a linearization map
\[
\pi_n(\THH(HA)) \to \hh_n(A).
\]
The linearization map on topological Hochschild homology is induced by the linearization map $R \to H\pi_0(R)$ of a $(-1)$-connected ring spectrum.\footnote{A spectrum $X$ is $n$-connected if $\pi_k(X) = 0$ for $k\leq n$.}

For our construction of quantum topological Hochschild homology, we prove that there is similarly a linearization map to the algebraic theory.
\begin{theorem}
For $G$ a cyclic group, $A$ a $(-1)$-connected graded $\Sq$-algebra, and $M$ a $(-1)$-connected graded $(A,A)$-bimodule, there is a natural homomorphism
\[
\pi_n(\qTHH(A;M)) \to \qhh_n(\pi_0(A); \pi_0(M)),
\]
which is an isomorphism when $n=0.$
\end{theorem}
\begin{proof}
For a simplicial spectrum $X_{\bullet}$, by \cite[X.2.9]{ekmm}, there is a spectral sequence 
\[
E^2_{s,t} = H_s(\pi_t(X_{\bullet})) \Rightarrow \pi_{s+t}(|X_{\bullet}|).
\]
from filtering by the skeleton. The spectral sequence in the case at hand will be of the form \[
E^2_{s,t} = H_s(\pi_t(\qTHH(A;M)_{\bullet})) \Rightarrow \pi_{s+t}(\qTHH(A;M)).
\]
The edge homomorphism of this spectral sequence is a map
\[
\pi_s(\qTHH(A;M)) \to H_s(\pi_0(\qTHH(A;M)_{\bullet})). 
\]
We rewrite this as 
\[
\pi_s(\qTHH(A;M)) \to H_s(\pi_0(M\wedge_{\Sq} A^{\wedge_{\Sq}\bullet})). 
\]
Since $A$ and $M$ are $(-1)$-connected, the K\"unneth spectral sequence collapses in degree 0, and we can identify the right hand side as 
\[
H_s(\pi_0(M) \otimes_{\mathbb{Z}[G]} \pi_0(A)^{\otimes_{\mathbb{Z}[G]}\bullet}).
\]
 Since $A$ is an $\Sq$-module, it follows that $\pi_0(A)$ is a $\mathbb{Z}[G]$-module.  The right hand side of the edge homomorphism is the homology of the quantum Hochschild complex $\qch_{\bullet}(\pi_0A; \pi_0M)$, where the unit $q$ in the group ring  $\mathbb{Z}[G]$ is a generator of the cyclic group. Hence, we get a linearization map 
\[
\pi_s(\qTHH(A;M)) \to \qhh_s(\pi_0(A);\pi_0(M)).
\]
We now show that this is an isomorphism in degree 0. Observe that in the spectral sequence from the skeletal filtration the only contribution to $s+t=0$ on the $E^2$-page is 
\[
E^2_{0,0} = H_0(\pi_0(M) \otimes_{\mathbb{Z}[G]} \pi_0(A)^{\otimes_{\mathbb{Z}[G]}\bullet}),
\]
which is $\hh_0(\pi_0(A);\pi_0(M)).$ This is a first quadrant spectral sequence so we conclude that $E^{\infty}_{0,0} = E^2_{0,0}$ and hence 
\[
\pi_0(\qTHH(A;M)) \cong \qhh_0(\pi_0(A); \pi_0(M)).
\]

\end{proof}

\begin{remark}
One of the foundational tools for computing classical topological Hochschild homology is the B\"okstedt spectral sequence. 
We have a similar result for the quantum setting. Let $A$ be a graded $\Sq$-algebra, $M$ a graded $(A,A)$-bimodule, and $k$ a commutative ring. If $H_*(A;k)$ is flat over $H_*(\Sq;k)$, then there is a quantum B\"okstedt spectral sequence
\[
E^2_{s,*} = \qhh_s(H_*(A;k);H_*(M;k)) \Rightarrow H_{s+*}(\qTHH(A;M);k)
\]
that converges strongly.
\end{remark}

\subsubsection{Scalar extensions}\label{sec:sc ext}
In our applications of $\qTHH$, the $\Sq$-algebras and modules over these algebras are obtained as scalar extensions (see \cite[Chapter IV.1]{schwedesymmetric}), which we briefly review here.  Let $R$ be a commutative ring spectrum with multiplication $\mu_R$.  If $A$ is any other ring spectrum (viewed as an $\S$-algebra) with multiplication and unit $\mu_A$ and $\eta_A$, then the \emph{$R$-scalar extension} of $A$ is the $R$-algebra $\prescript{}{R}{A}:=A\wedge R$ with $R$-module structure\footnote{Note that, since $R$ is commutative, an $R$-module is naturally an $R$-bimodule as usual.} induced by $\mu_R$ and multiplication and unit induced by $\mu_A$ and $\eta_A$ as follows:
\begin{align*}
(A\wedge R) \wedge_R (A\wedge R) & \cong (A\wedge A) \wedge R \rar{\mu_A \wedge \id_R} A\wedge R , \\
R & \cong \S \wedge R \rar{\eta_A \wedge \id_R} A\wedge R.
\end{align*}

Similarly, if $M$ is an $A$-bimodule, then the \emph{$R$-scalar extension} $\prescript{}{R}{M}:=M\wedge R$ is a bimodule over $\prescript{}{R}{A}$ with actions induced by the actions of $A$ on $M$.  For another $A$-bimodule $N$, there is an isomorphism of $R$-modules
\begin{equation}\label{eq:scext commutes with smash}
\left(\prescript{}{R}{M}\right) \wedge_{R} \left( \prescript{}{R}{N} \right) \cong \prescript{}{R}{(M\wedge N)}.
\end{equation}

We focus on the case of interest where $R=\Sq$.  In this setting we use a special notation for the $\Sq$-scalar extensions, which take the following special forms
\begin{equation}\label{eq:scalar ext for Sq}
\scext{A}:=A\wedge \Sq\cong \bigvee_\ell q^\ell A,\quad \scext{N}:=N\wedge \Sq\cong \bigvee_m q^m N.
\end{equation}

As with $\Sq$, we write $q^\ell A$ here to denote different copies of $A$ as $q^\ell$ ranges over $G$, and similarly for $q^m N$.  The designated automorphisms of $\scext{A}$ and $\scext{N}$, as in Proposition \ref{prop:Sq module automorphism}, are then permutations of wedge summands, induced by identity maps
\begin{equation} 
\label{eq:automorphism on extension}
q^\ell A \xrightarrow{} q^{\ell+1} A, \quad q^m N \xrightarrow{} q^{m+1} N.
\end{equation}
The following lemma is then immediate from the definitions.

\begin{lemma}\label{lem:scalar ext properties}
Let $A$ be a ring spectrum ($\sphere$-algebra) with multiplication map $A\wedge A\xrightarrow{\mu_0} A$.  Then the scalar extended multiplication $\scext{A}\wedge_{\Sq}\scext{A} \xrightarrow{\scext{\mu}} \scext{A}$ respects the equivalence of Equation \eqref{eq:scext commutes with smash} as in the following commuting diagram.
\[\begin{tikzcd}
\scext{A} \wedge_{\Sq} \scext{A} \ar[rr,<->,"\cong"] \ar[dr,"\scext{\mu}"'] & &
\scext{(A\wedge A)}= \bigvee_\ell q^\ell (A\wedge A) \ar[dl,"\bigvee_\ell \mu_0"]\\
 & \scext{A} = \bigvee_\ell q^\ell A & 
\end{tikzcd}.\]
Similarly, if $N$ is an $(A,A)$-bimodule with left action $A\wedge N \xrightarrow{\lambda_0} N$, then the scalar extended left action $\scext{A}\wedge_{\Sq}\scext{N} \xrightarrow{\scext{\lambda}} \scext{N}$ satisfies:
\[\begin{tikzcd}
\scext{A} \wedge_{\Sq} \scext{N} \ar[rr,<->,"\cong"] \ar[rd,"\scext{\lambda}"'] & &
\scext{(A\wedge N)}= \bigvee_\ell q^\ell (A\wedge N) \ar[dl,"\bigvee_\ell \lambda_0"]\\
 & \scext{N} = \bigvee_\ell q^\ell N & 
\end{tikzcd}.\]
The right action behaves similarly.
\end{lemma}

Lemma \ref{lem:scalar ext properties} can be used to give a simplified version of Definition \ref{def:qTHH} in the setting where $A$ and $M$ are $\Sq$-scalar extensions of a graded ring spectrum and bimodule respectively; rather than spell this out precisely, we will instead apply the lemma in a more specific context in Section \ref{sec:proof of main qTHH thm}.

\subsection{$\qTHH$ for spectral categories}
\label{sec:qTHH for spectral cats}
Recall from Definition \ref{def:Alg as lin Cat} that we may view an algebra as a linear category.  We begin by reviewing the definition of an $R$-spectral category for $R$ a commutative ring spectrum, which is the topological analogue of a linear category. For further details on spectral categories and their topological Hochschild homology see \cite{Blumberg-Mandell}.  
\begin{definition}\label{def:spectral cat}
Let $R$ be a commutative ring spectrum. An \emph{$R$-spectral category} is a category enriched in $R$-module spectra. More explicitly, an $R$-spectral category $\mathcal{C}$ consists of:
\begin{enumerate}
    \item a collection of objects $\ob(\mathcal{C})$,
    \item an $R$-module spectrum $\mathcal{C}(a,b)$ for each pair of objects $a, b \in $ob$(\mathcal{C})$,
    \item a unit map $R \to \mathcal{C}(a,a)$ for each $a\in $ob$(\mathcal{C})$, and
    \item a composition map
    \[
    \mathcal{C}(b,c) \wedge_R \mathcal{C}(a,b) \to \mathcal{C}(a,c)
    \]
    for each triple of objects, $a,b,c \in $ob$(\mathcal{C}$), which is associative and unital.
\end{enumerate}
\end{definition}

\begin{definition}
Let $\mathcal{C}$ and $\mathcal{D}$ be $R$-spectral categories. An \emph{$R$-spectral functor} $F: \mathcal{C} \to \mathcal{D}$ is an enriched functor. More explicitly an $R$-spectral functor consists of:
\begin{enumerate}
    \item a function on objects $F:$ ob($\mathcal{C}) \to $ob($\mathcal{D})$, and 
    \item a map of $R$-module spectra 
    \[
    F_{a,b}: \mathcal{C}(a,b) \to \mathcal{D}(Fa, Fb)
    \]
    for each pair of objects $a,b \in $ob($\mathcal{C}),$
\end{enumerate}
compatible with composition and units.
\end{definition}

\begin{definition}\label{def:bimodule over spectral cats}
Let $\mathcal{C}$ be an $R$-spectral category. We denote by $\cat{C}^{op}$ the \emph{opposite} spectral category, defined by $\cat{C}^{\textup{op}}(c, c') = \cat{C}(c',c)$. A \emph{$\mathcal{C}$-module} is an $R$-spectral functor from $\cat{C}^{\textup{op}}$ to the spectral category of $R$-module spectra.
\end{definition}
\begin{definition}
Let $\mathcal{C}$ and $\cat{D}$ be $R$-spectral categories. Their \emph{smash product} $\cat{C}\wedge \cat{D}$ is defined to have objects $\ob(\cat{C})\times \ob(\cat{D})$ and $(\cat{C}\wedge \cat{D})((c, d), (c', d')) = \cat{C}(c, d) \wedge_R \cat{D}(c', d')$, with units and composition defined in the natural way. A \emph{$(\cat{C}, \cat{D})$-bimodule} is an $R$-spectral functor from $\cat{C}^{\textup{op}}\wedge \cat{D}$ to the spectral category of $R$-module spectra.
\end{definition}
Observe that, explicitly, a $(\cat{C}, \cat{D})$-bimodule consists of a choice of $R$-module $\cat{M}(c,d)$ for each pair of objects $c \in \textup{ob}(\cat{C})$ and $d\in \textup{ob}(\cat{D})$, along with maps
\begin{align*}
\cat{C}(c, c') \wedge_R \cat{M}(c',d) &\to \cat{M}(c,d) \\
\cat{M}(c,d) \wedge_R  \cat{D}(d,d') &\to \cat{M}(c,d')
\end{align*}
for each $c'\in \textup{ob}(\cat{C})$ and $d'\in \textup{ob}(\cat{D})$, which are compatible with composition in $\cat{C}$ and $\cat{D}$.

\begin{definition}
Let $R$ be a cofibrant commutative ring spectum. An $R$-spectral category $\mathcal{C}$ is said to be \emph{pointwise cofibrant} if for each pair $a,b \in $ ob($\mathcal{C})$, $\mathcal{C}(a,b)$ is a cofibrant $R$-module spectrum in the standard model structure on the category of $R$-modules \cite{Shipley}. If $\cat{D}$ is an $R$-spectral category, a $(\cat{C}, \cat{D})$-bimodule $\cat{M}$ is \emph{pointwise} cofibrant if $\cat{M}(c,d)$ is a cofibrant $R$-module for every $c\in \ob(\cat{C})$ and $d\in \ob(\cat{D})$.
\end{definition}

Graded $R$-spectral categories, functors, and bimodules are defined analogously. For instance, the graded analogue of Definition \ref{def:spectral cat} each $\cat{C}(a,b)$ is graded, the unit map lands in the degree-zero piece of $\cat{C}(a,a)$, and composition is degree-preserving. In the setting of quantum topological Hochschild homology, we are interested in graded spectral categories over the spherical group ring $\Sq$, where $G$ is a cyclic group. We now define quantum topological Hochschild homology for graded $\Sq$-spectral categories. 

\begin{definition}\label{def:qTHH for spectral cat}
Given a cyclic group $G$, a pointwise cofibrant graded $\Sq$-spectral category $\mathcal{C}$, and a graded $(\mathcal{C},\mathcal{C})$-bimodule $\mathcal{M},$ we define a simplicial $\Sq$-module $\qTHH(\mathcal{C};\mathcal{M})_{\bullet}$ with $n$ simplices
\[
\qTHH(\mathcal{C};\mathcal{M})_n = \bigvee_{(c_0,\dots,c_n)\in \textup{ob}\mathcal{C}} \mathcal{M}(c_n,c_0)\wedge_{\Sq} \cat{C}(c_0,c_1) \wedge_{\Sq} \cdots \wedge_{\Sq} \cat{C}(c_{n-1},c_n).
\]
Degeneracy maps are induced by the unit maps in $\cat{C}$. The face maps $d_i$ for $0<i<n$ are induced by the composition
\[
\cat{C}(c_i, c_{i+1}) \wedge _{\Sq} \cat{C}(c_{i+1}, c_{i+2}) \to \cat{C}(c_{i}, c_{i+2}),
\]
and the map $d_0$ is induced by the module action
\[
\cat{M}(c_n, c_0) \wedge_{\Sq} \cat{C}(c_0, c_1) \to \cat{M}(c_n, c_1).
\]
The final face map $d_n$ is given by the composite 
\begin{align*}
\mathcal{M}(c_n,c_0)\wedge_{\Sq}  \cdots \wedge_{\Sq} \cat{C}(c_{n-1},c_n) \cong  
& \bigvee_{j\in\Z} \cat{C}(c_{n-1},c_n)_j \wedge_{\Sq} \mathcal{M}(c_n,c_0)\wedge_{\Sq} \cdots   \mathcal{C}(c_{n-2}, c_{n-1})
\\
\rar{\bigvee_j q^{-j} \wedge \id } & \bigvee_{j\in\Z}\cat{C}(c_{n-1},c_n)_j \wedge_{\Sq} \mathcal{M}(c_n,c_0)\wedge_{\Sq} \cdots \mathcal{C}(c_{n-2}, c_{n-1})
\\
 \to & \mathcal{M}(c_n,c_0)\wedge_{\Sq} \cdots   \mathcal{C}(c_{n-2}, c_{n-1})
,
\end{align*}
where the second map applies the automorphisms $q^{-j}$ from Proposition \ref{prop:Sq module automorphism} only on the first factor, and the third map is induced by the left module action.
\end{definition}

\section{Multifunctors and shape multicategories}
\label{sec:spectral invariants via multifunctors}

We will focus on applying the definitions above to the Chen-Khovanov spectral invariants of Lawson-Lipshitz-Sarkar \cite{LLS_CK_Spectra}.  These invariants take the form of spectral categories associated to points on a line, and spectral bimodules associated to tangles between these points, which lift the platform algebras and bimodules of Chen-Khovanov \cite{Chen-Khovanov} (suitably reinterpreted via Definition \ref{def:Alg as lin Cat}) to the stable homotopy category.  However, in order to build and study these invariants, Lawson-Lipshitz-Sarkar take advantage of a further reformulation of linear/spectral categories and bimodules via multifunctors out of so-called \emph{shape} multicategories which we describe in Section \ref{sec:linear and spectral cats via multifunctors}.

In \cite[Section 5]{LLS_CK_Spectra}, Lawson-Lipshitz-Sarkar study the topological Hochschild homologies of their invariants using functors out of different sorts of shape categories (which are ordinary categories, not multicategories).  In Section \ref{sec:qTHH via shape cats} we generalize this approach to handle our quantum case.
Throughout this section we will  use $R$ to denote a commutative ring, and $\ringspec$ to denote a commutative ring spectrum.

\subsection{Linear and spectral categories via shape multicategories}
\label{sec:linear and spectral cats via multifunctors}

This section presents a review of the material for describing linear and spectral categories via multifunctors out of shape multicategories, based upon the presentation in \cite{LLS_Khovanov_Spectra}.  As we will not be making direct use of the majority of these concepts, we will only present brief summaries of the main ideas, together with references to the relevant discussions in \cite{LLS_Khovanov_Spectra}.  We begin with an outline of the notion of a multicategory as described there.

\begin{definition}\label{def:multicategory}
A \emph{multicategory} $\cat{C}$ (also referred to as a \emph{colored operad}) consists of a set of objects $\ob(\cat{C})$ together with, for any sequence of objects $x_1,\dots,x_n,y\in\ob(\cat{C})$, ($n\geq 0$), a multimorphism set
\[\cat{C}(x_1,\dots,x_n;y):= \Hom_{\cat{C}}(x_1,\dots,x_n;y).\]
The multimorphism sets come equipped with distinguished identity elements $\id_x\in \cat{C}(x;x)$ and  compositions of the form
\[\cat{C}(y_1,\dots, y_m;z) \times \cat{C}(x_{1,1},\dots,x_{1,n_1};y_1) \times\cdots\times \cat{C}(x_{m,1},\dots,x_{m,n_m};y_m) \xrightarrow{\circ} \cat{C}(x_{1,1},\dots,x_{m,n_m};z).\]
Together these data are required to satisfy associativity and identity axioms similar to those of a usual category. 

A \emph{multifunctor} $\cat{C}\xrightarrow{\func{F}}\cat{D}$ is then a function on the object sets and various multimorphism sets that respects identities and compositions.  See \cite[Definition 2.1]{LLS_Khovanov_Spectra} for more details.
\end{definition}

Any monoidal category $(\cat{C},\otimes)$ can be viewed as a multicategory with the help of the monoidal product $\otimes$ by declaring multimorphism sets to be
\[\cat{C}(x_1,\dots,x_n;y):=\Hom_{\cat{C}}(x_1\otimes\cdots\otimes x_n;y),\]
for any choice of parenthesis placement in the monoidal product. In this way, we may view $\Rmod$ ($R$-modules with $\otimes_R$), $\Kom(\Rmod)$ (complexes of $R$-modules with $\o_R)$, and $\Rmod[\ringspec]$ (spectral $\ringspec$-modules with $\wedge_{\ringspec}$) as multicategories.

There are also some specific types of multicategories which are particularly well-suited to studying linear/spectral categories.

\begin{definition}[{\cite[Definition 2.2]{LLS_Khovanov_Spectra}}]
\label{def:shape(X) multicat}
For a finite set $X$, let $\shape(X)$ denote the \emph{shape multicategory of $X$}. Objects of $\shape(X)$ are pairs $(x_1, x_2) \in X\times X$, with unique multimorphisms
\[
(x_1, x_2), (x_2, x_3), \cdots ,(x_{m-1}, x_m) \to (x_1, x_m)
\]
and all other multimorphism sets empty. In the above $m=0$ is allowed, producing a unique morphism $\varnothing \to (x,x)$ for each $x\in X$.
\end{definition}

More generally, we have a notion of shape category for any number of sets as follows.

\begin{definition}[{\cite[Definition 2.3]{LLS_Khovanov_Spectra}}]
\label{def:shape gluing multicat}
Given finite sets $X^1, \ldots, X^k$, $\shape(X^1, \ldots, X^k)$ denotes the \emph{shape multicategory of $X^1, \ldots, X^k$}. Its objects are $\coprod_{i\leq j} X^i \times X^j$, there is a unique multimorphism  $(x_1, x_2), (x_2, x_3), \cdots ,(x_{m-1}, x_m) \to (x_1, x_m)$, and all other multimorphism sets are empty. The case $m=0$ is again allowed, producing a unique multimorphism $\varnothing \to (x,x)$ for every $x\in \coprod_i X^i$.
\end{definition}

In particular, the shape category $\shape(X,Y)$ of two sets $X,Y$ has objects of the form $(x_1,x_2),(x,y),$ and $(y_1,y_2)$ where $x_1, x_2, x\in X$, $y_1, y_2, y\in Y$. In this setting, we use the notation
\begin{equation}
\label{eq:x vert y notation}
    (x\mid y) 
\end{equation}
to denote objects of the second type, since we will often consider the case $X=Y$ where we need to distinguish between objects of the form $(x,y)$ and $(x\mid y)$ when $x,y\in X$. There are unique morphisms in $S(X,Y)$ of the form
\begin{align*}
    (x_1,x_2),(x_2,x_3),\dots,(x_{m-1},x_m) &\rar{} (x_1,x_m)\\
    (x_1,x_2),(x_2,x_3),\ldots,(x_{m-1},x_m),(x_m \mid y_1),(y_1,y_2),\ldots, (y_{\ell-1},y_\ell) &\rar{} (x_1 \mid y_\ell)\\
    (y_1,y_2),(y_2,y_3),\ldots(y_{\ell-1},y_\ell), &\rar{} (y_1,y_\ell)
\end{align*}
where $x_1, \ldots, x_m \in X$, $y_1, \ldots, y_\ell \in Y$, together with a unique morphism $\varnothing \to (z,z)$ for each $z\in X\sqcup Y$.  Note that $\shape(X,Y)$ contains both $\shape(X)$ and $\shape(Y)$ as full submulticategories.

With these definitions in place, we follow \cite[Section 2.3]{LLS_Khovanov_Spectra} to recast linear and spectral categories, as well as their bimodules, into a multicategorical language.

\begin{proposition}\label{prop:Alg as multifunctor}
The data of an $\ring$-linear category $\cat{A}$ with finite object set $S_A$ is equivalent to the data of a multifunctor $\shape(S_A)\rar{\func{A}}\Rmod$. Under this equivalence, an $(\cat{A},\cat{B})$-bimodule $\opcat{B}\times\cat{A}\rar{\cat{M}} \Rmod$ is equivalent to the data of a multifunctor $\shape(S_A,S_B)\rar{\func{M}}\Rmod$ which restricts to $\func{A}$ and $\func{B}$ on the submulticategories $\shape(S_A)$ and $\shape(S_B)$ respectively.

Similarly, a complex of $(\cat{A},\cat{B})$-bimodules $\opcat{B}\times\cat{A}\rar{\cat{M}} \Kom(\Rmod)$ is equivalent to a multifunctor $\shape(S_A,S_B)\rar{\func{M}}\Kom(\Rmod)$ 
which restricts to $\func{A}$ and $\func{B}$ on  $\shape(S_A)$ and $\shape(S_B)$, where we view $\Rmod$ as a full subcategory of $\Kom(\Rmod)$ as complexes concentrated in homological degree zero.

All of the above holds for $\ringspec$-spectral categories in place of $R$-linear categories, with all relevant multifunctors having codomain $\Rmod[\ringspec]$ in place of $\Rmod$.  Similarly, graded versions in the categorical language are equivalent to multifunctors from relevant shape categories to $\Rmod^\bullet$ (or $\Rmod[\ringspec]^\bullet$) where $\bullet \in \{gr, gg\}$. 

\end{proposition}

The reader should consult \cite[Section 2.3]{LLS_Khovanov_Spectra} for more precise details.

\subsection{Quantum topological Hochschild homology via shape categories}
\label{sec:qTHH via shape cats}

Similar to the case of their spectral algebras and bimodules, Lawson-Lipshitz-Sarkar in \cite{LLS_CK_Spectra} study the topological Hochschild homology of their Khovanov invariants by presenting it via a functor from a special sort of shape category.  We present a generalization of their ideas here geared towards studying the quantum topological Hochschild homology of their invariants (see Remark \ref{rmk:recovering LLS THH constructions} for the relationship between our arrangement and theirs).  The reader should think of the set $S$ in the following definition as the object set of some linear or spectral category.

\begin{definition}\label{def:HM shape category}
The \emph{Hochschild-Mitchell shape category} of a set $S$, denoted $\HMshape(S)$, is the category consisting of the following.
\begin{itemize}
    \item Objects are finite sequences $(a_0,\dots,a_\alpha)$ of elements of  $S$ (here $\alpha\geq 0$).
    \item For each (non-empty) subsequence $(a_\beta,\dots,a_\gamma)$ of $(a_0,\dots,a_\alpha)$ (including $\beta=0,\gamma=\alpha$), there is a unique morphism $(a_0,\dots,a_\alpha)\rightarrow (a_\beta,\dots,a_\gamma)$ corresponding to a deletion of some number of elements from the sequence.
    \item Compositions are defined in the natural way for sub-subsequences, corresponding to multiple deletions.
\end{itemize}
The \emph{augmented} Hochschild-Mitchell shape category of $S$, denoted $\HMshape^\omega(S)$, is the category obtained by adding one \emph{terminal object} $\omega$ to $\HMshape(S)$, together with one \emph{terminal morphism} $x\rightarrow\omega$ for all objects $x\in\HMshape^\omega(S)$ (including $x=\omega$).  Thus $\omega$ plays the role of the empty sequence.  Note that the augmented $\HMshape^\omega(S)$ contains $\HMshape(S)$ as a full subcategory.
\end{definition}

\begin{remark}\label{rmk:HMshape is cat not multicat}
Note that $\HMshape(S)$ and $\HMshape^\omega(S)$ are only ordinary categories, having no non-empty multimorphism sets with more than one input.
\end{remark}

\begin{definition}\label{def:qHM functor}
Let $\cat{A}$ be a graded $\Sq$-spectral category with object set $S_A$, and let $\func{M}$ be a graded $(\cat{A}, \cat{A})$-bimodule.  The \emph{quantum Hochschild-Mitchell functor} associated $\func{M}$ is the functor

\[\HMshape(S_A) \rar{\qHMfunc{\func{M}}} \gSpectra_G\]
defined as follows. 
\begin{itemize}
    \item (Simplices) On objects $(a_0,\dots,a_\alpha)$ we define
    \[
    \qHMfunc{\func{M}}(a_0,\dots,a_\alpha) = \func{M}(a_\alpha,a_0)\wedge_{\Sq} \cat{A}(a_0,a_1)\wedge_{\Sq} \cdots \wedge_{\Sq} \cat{A}(a_{\alpha-1},a_{\alpha}).
    \]
    \item (Ordinary face maps) Each morphism of the form 
    \[
    (a_0,\dots,a_i,\dots,a_\alpha)\rar{} (a_0,\dots,a_{i-1},a_{i+1},\dots,a_\alpha)
    \]
    for $i\neq\alpha$ is assigned a multiplication or action map of the form 
    \begin{align*}
    \cat{A}(a_{i-1},a_i)\wedge_{\Sq} \cat{A}(a_i,a_{i+1}) & \rightarrow \cat{A}(a_{i-1},a_{i+1}) & & (0<i <\alpha), \\
    \func{M}(a_\alpha,a_0)\wedge_{\Sq} \cat{A}(a_0,a_1) & \rightarrow \func{M}(a_\alpha,a_1) & & (i=0),
    \end{align*}
    smashed with identity maps on the other factors, mimicking the usual Hochschild face maps as in Equation \eqref{eq:qTHH noncyclic face maps} of Definition \ref{def:qTHH as a simplicial object}.
    \item (Last face map) Each morphism of the form $(a_0,\dots,a_\alpha) \rar{} (a_0,\dots,a_{\alpha-1})$ is assigned the composite
    \begin{align*}
    \func{M}(a_\alpha,a_0)\wedge_{\Sq}\cdots \wedge_{\Sq}\cat{A}(a_{\alpha-1},a_\alpha)_j  \cong & \cat{A}(a_{\alpha-1},a_\alpha)_j \wedge_{\Sq} \func{M}(a_\alpha,a_0)\wedge_{\Sq}\cdots  \cat{A}(a_{\alpha - 2}, a_{\alpha-1})\\
    \xrightarrow{q^{-j}_{\cat{A}}\wedge \id} &
    \cat{A}(a_{\alpha-1},a_\alpha)_j \wedge_{\Sq} \func{M}(a_\alpha,a_0)\wedge_{\Sq}\cdots \cat{A}(a_{\alpha - 2}, a_{\alpha-1})  \\
    \rightarrow & \func{M}(a_{\alpha-1},a_0)\wedge_{\Sq}\cdots \cat{A}(a_{\alpha - 2}, a_{\alpha-1})
    \end{align*}
    on each graded summand of $\cat{A}(a_{\alpha-1},a_\alpha)= \bigvee_j\cat{A}(a_{\alpha-1},a_{\alpha})_j$ as in Equation \eqref{eq:qTHH cyclic face map} of Definition \ref{def:qTHH as a simplicial object}.
\end{itemize}

Then an \emph{augmented} quantum Hochschild-Mitchell functor associated to $\func{M}$ is a functor
\[\func{F}:\HMshape^\omega(S_A) \rightarrow \gSpectra_G\]
which restricts to a functor naturally isomorphic to the quantum Hochschild-Mitchell functor,
\[\func{F}|_{\HMshape(S_A)} \cong \qHMfunc{\func{M}}.\]
\end{definition}

\begin{remark}\label{rmk:recovering LLS THH constructions}
The categories $\HMshape(S_A)$ and $\HMshape^\omega(S_A)$ correspond to the categories denoted $\mathscr{C}$ and $\mathscr{C}^T$, respectively, in \cite[Section 5]{LLS_CK_Spectra}.   Then the functors $F$ and $\overline{K\circ L^T}$ in \cite[Proof of Theorem 8]{LLS_CK_Spectra} are (trivially quantum, for $G=\{1\}$) Hochschild-Mitchell and augmented Hochschild-Mitchell functors for their Khovanov tangle invariants (in the case that the cube $\cube^N$ described there is $(N=0)$-dimensional).
\end{remark}

\begin{proposition}\label{prop:qTHH via qHM functor}
The data of the semi-simplicial $\Sq$-module $\qTHH(\cat{A};\func{M})_\bullet$ (Definition \ref{def:qTHH for spectral cat}) is equivalent to the data of the quantum Hochschild-Mitchell functor
\[\HMshape(S_A) \rar{\qHMfunc{\func{M}}} \gSpectra_G\]
associated to $\func{M}$.  In particular, $\qTHH(\cat{A};\func{M})\simeq \hocolim\left( \qHMfunc{\func{M}}\right).$
\end{proposition}
\begin{proof}
The equivalence of data is clear from the definitions; compare to the argument in \cite[proof of Theorem 8]{LLS_CK_Spectra}.
\end{proof}

\begin{corollary}\label{cor:aug qHM gives map on qTHH}
An augmented quantum Hochschild-Mitchell functor $\func{F}$ associated to $\func{M}$ induces a map of $\Sq$-module spectra
\[\qTHH(\cat{A};\func{M}) \rar{F} \func{F}(\omega).\]
\end{corollary}
\begin{proof}
This is an immediate abstraction of the argument in \cite[proof of Theorem 8]{LLS_CK_Spectra}.  In short, the category $\HMshape^\omega(S_A)$ contains $\HMshape(S_A)$ as a full subcategory with no incoming morphisms, with complementary subcategory $\{\omega\}$ having no outgoing morphisms, giving rise to a cofibration sequence among the homotopy colimits
\[\func{F}(\omega) \rar{} \hocolim \func{F} \rar{} \hocolim \func{F}|_{\HMshape(S_A)}.\]
The Puppe construction then gives rise to a map
\[\Sigma \hocolim \func{F}|_{\HMshape(S_A)} \rar{} \Sigma\func{F}(\omega),\]
and $F$ is the desuspension of this map via Proposition \ref{prop:qTHH via qHM functor}. 
\end{proof}

\section{An overview of the relevant link and tangle homology theories}
\label{sec:quantum annular homological background}

\subsection{The Khovanov and annular Khovanov functors}
\label{sec:The Khovanov and annular Khovanov functors}
We briefly review the Bar-Natan categories \cite{Bar-Natan} $\BN(\D)$ and $\BN(\A)$ of the open disk $\D:=(0,1)\times(-1,1)$ and annulus $\A:=S^1\times (-1,1)$. We set $I := [0,1]$.

\begin{figure}
\centering
\subcaptionbox{Sphere.
\label{fig:sphere}}[.45\linewidth]
{\begin{tikzpicture}
\draw (0,0) circle [radius=.75];
\draw (.75,0) arc(0:-180:.75 and .3);
\draw[dashed] (-.75,0) arc(180:-0:.75 and .3);
\node at (1.25,0) {$=\ 0$};
\end{tikzpicture}
}
\subcaptionbox{Dotted sphere.
\label{fig:dotted sphere}}[.45\linewidth]
{\begin{tikzpicture}
\draw (0,0) circle [radius=.75];
\draw (.75,0) arc(0:-180:.75 and .3);
\draw[dashed] (-.75,0) arc(180:-0:.75 and .3);
\node at (1.25,0) {$=\ 1$};

\node at (0,0) {$\bullet$};
\end{tikzpicture}
}
\\
\vskip2ex
\subcaptionbox{Neck-cutting.
\label{fig:neck-cutting}}[.5\linewidth]
{\begin{tikzpicture}
\draw (0,0) arc (180:360:.5 and .2);
\draw[densely dashed] (1,0) arc (0:180:.5 and .2);

\draw (0,2) arc (180:360:.5 and .2);
\draw(1,2) arc (0:180:.5 and .2);

\draw (0,0) -- (0, 2);
\draw (1,0) -- (1, 2);




\begin{scope}[shift={(2,0)}]
\draw (0,0) arc (180:360:.5 and .2);
\draw[densely dashed] (1,0) arc (0:180:.5 and .2);

\draw (0,2) arc (180:360:.5 and .2);

\draw (1,2) arc (0:180:.5 and .2);

\draw (1,0) arc (0:180:.5 and .8);
\draw (0,2) arc (180:360:.5 and .8);

\node at (0.5,1.5) {$\bullet$};

\end{scope}


\begin{scope}[shift={(4,0)}]
\draw (0,0) arc (180:360:.5 and .2);
\draw[densely dashed] (1,0) arc (0:180:.5 and .2);

\draw (0,2) arc (180:360:.5 and .2);

\draw (1,2) arc (0:180:.5 and .2);

\draw (1,0) arc (0:180:.5 and .8);
\draw (0,2) arc (180:360:.5 and .8);

\node at (.5,.5) {$\bullet$};
\end{scope}


\node at (1.5,1) {$=$};
\node at (3.5,1) {$+$};

\end{tikzpicture}
}
\subcaptionbox{Two dots.
\label{fig:two dots}}[.4\linewidth]
{\begin{tikzpicture}[scale=.7]
\draw (0,0) -- (2,0);
\draw (0,0) -- (1,2);
\draw (1,2) -- (3,2);
\draw (2,0) -- (3,2);
\node at (1.25,1) {$\bullet$};
\node at (1.75,1) {$\bullet$};

\node at (3.75,1) {$=\ 0$};

\end{tikzpicture}
}
\caption{The local Bar-Natan relations.}\label{fig:BN relations}
\end{figure}
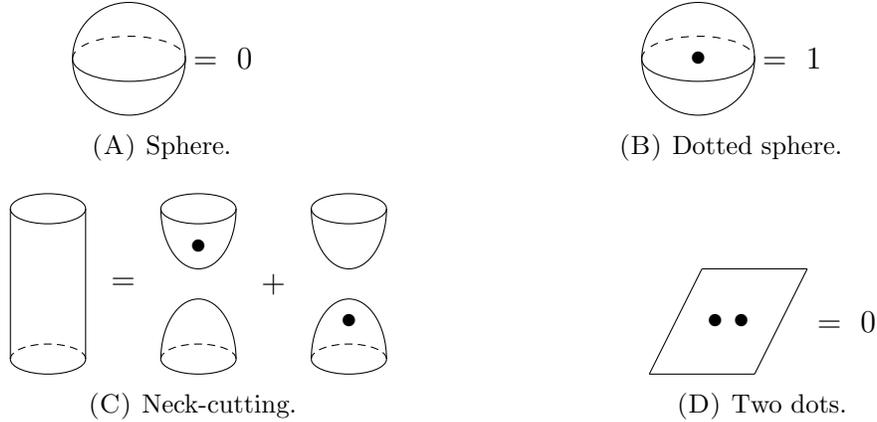

\begin{definition}
\label{def:BN category}
For $X=\D$ or $X=\A$, the \emph{Bar-Natan category} of $X$, denoted $\BN(X)$, is the category whose objects are formal direct sums of formally graded collections of disjoint simple closed curves in $X$, and whose morphisms are matrices whose entries are formal $\Z$-linear combinations of dotted cobordisms in $I\times X$, considered up to ambient isotopy fixing the boundary pointwise and modulo the local Bar-Natan relations shown in Figure \ref{fig:BN relations}.   The \emph{topological degree} of a dotted cobordism $\Sigma$ is defined to be
\begin{equation}\label{eq:top deg of dotted cob}
\chi(\Sigma):= \chi(u(\Sigma)) - 2d
\end{equation}
where $\chi(u(\Sigma))$ is the Euler characteristic of the underlying surface of $\Sigma$ (ignoring the dots), and $d$ is the number of dots.
\end{definition}

\begin{definition}\label{def:FKh on BN(D)}
The \emph{Khovanov functor} $\BN(\D)\rar{\FKh}\gAb$ is the representable functor $\FKh(-):= \Hom_{\BN(\D)}(\varnothing,-)$.  For a collection of disjoint circles $Z \subset \D$, let $\Gen(Z)$ denote the set of labellings of each circle by $1$ or $X$. Then $\FKh(Z)$ is identified with the free abelian group generated by $\Gen(Z)$, where a label of $1$ (respectively of $X$) on a circle corresponds to an undotted (respectively dotted) cup cobordism in $\Hom(\varnothing,Z)$.  We refer to such basis elements $z\in\Gen(Z)$ as \emph{Khovanov generators}.  The \emph{quantum grading} $|\cdot|$ on $\FKh(Z)$ is defined on Khovanov generators $z\in\Gen(Z)$ by
\begin{equation}
\label{eq:qdeg}
|z|:=\chi(z) = \#\{\text{circles labelled $1$}\} - \#\{\text{circles labelled $X$}\}.
\end{equation}
Throughout the paper, we will use the brackets $\{-\}$ to denote an upwards shift in quantum grading, so that $z\in\FKh(Z)\{m\}$ is in quantum grading $|z|=\chi(z)+m$.  Moreover, morphisms $\Sigma$ in $\BN(\D)$ give rise to graded maps $\FKh(\Sigma)$ of degree $\chi(\Sigma)$ in $\gAb$.  
\end{definition}

\begin{remark}
    Two quantum grading conventions differing by an overall sign exist in the literature. The one used in this paper is consistent with \cite{Khovanov, BPW} but opposite that of \cite{Kh-ring_H, Chen-Khovanov, LLS_Khovanov_Spectra, LLS_CK_Spectra}. In particular, the Chen-Khovanov platform algebras described in Section \ref{sec:Platform algebras and bimodules} are non-positively graded. 
\end{remark}

The annular version of the theory is similar, but we differentiate between circles in $Z$ which are trivial or essential (zero or nonzero in $H_1(\A;\Z)$, respectively), allowing us to alter the grading information.

\begin{definition}
\label{def:FA on BN(A)}
Given an object $Z\in\BN(\A)$, the set of \emph{Khovanov generators} $\Gen(Z)$ is the set of labellings of circles in $Z$ such that any trivial circle is labelled by $1$ or $X$, while any essential circle is labelled by $v_-$ or $v_+$.  The \emph{annular Khovanov functor} $\BN(\A)\rar{\FA}\ggAb$ is then defined as follows.  On an object $Z$, $\FA(Z)$ is the free abelian group generated by $\Gen(Z)$, with \emph{quantum grading} $|\cdot|$ and \emph{annular grading} $|\cdot|_\A$ defined on Khovanov generators $z\in\Gen(Z)$ as
\begin{align}
    \begin{aligned}
     \label{eq:degrees in QAH}
     |z| & = \#\{\text{trivial circles labelled $1$}\} - \#\{\text{trivial circles labelled $X$}\}, \\
     |z|_\A & = \#\{\text{essential circles labelled $v_+$}\} - \#\{\text{essential circles labelled $v_-$}\}.
    \end{aligned}
\end{align}

\begin{remark}
    The quantum grading above is natural from several perspectives but differs slightly from the quantum grading elsewhere in the literature. It agrees with the grading denoted $j'$ in \cite[Remark 2]{GLW} and with the quantum grading in \cite{AKW, BPW}. See also \cite[Remark 2.3]{AKW}.
\end{remark}

 We may embed $\A\hookrightarrow\D$ and view any morphism $Z_0\rar{\Sigma}Z_1$ in $\BN(\A)$ as a morphism in $\BN(\D)$. We define $\FA(\Sigma)$ to be the annular degree-preserving part of $\FKh(\Sigma): \FA(Z_0) \to \FA(Z_1)$, where we use the natural identification $\FA(Z_i) \cong \FKh(Z_i)$ as free abelian groups given by setting $v_+$ and $v_-$ to correspond to $1$ and $X$, respectively.  It was shown in \cite[Section 2]{Roberts} that $\FKh(\Sigma)$ is non-increasing with respect to the annular grading.
\end{definition}

We elaborate slightly on the definition of $\FA$ on morphisms here.  If all boundary circles of $\Sigma$ are trivial, then $\FA(\Sigma) = \FKh(\Sigma)$.  There are four types of saddle cobordisms in $I\times \A$ where at least one boundary circle is essential; the corresponding maps under $\FA$ are recorded below in tensor product notation. We depict the annulus $\A$ as the punctured plane, where the puncture is depicted by the symbol $\times$.

\noindent\begin{minipage}{.5\linewidth}
\begin{equation*}\label{eq:formula1}
  \begin{split}
  \begin{tikzpicture}[baseline=-2pt]
      \node at (0,0) {$\mathsmaller{\times}$};
      \draw (0,0) circle (8pt);
      \draw (.8,0) circle (8pt);
  \end{tikzpicture} & \to   
  \begin{tikzpicture}[baseline=-2pt]
  \node at (0,0) {$\mathsmaller{\times}$};
      \draw (0,0) circle (8pt);
  \end{tikzpicture}
  \\
   v_{\pm} \o 1 & \mapsto v_{\pm} \\ 
    v_\pm \o X & \mapsto 0
  \end{split}
\end{equation*}
\end{minipage}%
\begin{minipage}{.35\linewidth}
\begin{equation*}\label{eq:formula2}
  \begin{split}
   \begin{tikzpicture}[baseline=-2pt]
  \node at (0,0) {$\mathsmaller{\times}$};
      \draw (0,0) circle (5pt);
      \draw (0,0) circle (9pt);
  \end{tikzpicture}
  &\to
     \begin{tikzpicture}[x=1pt, y=1pt,baseline=-2pt]
  \node at (0,0) {$\mathsmaller{\times}$};
      \draw (0,9) arc (90:-90:2 and 2); 
      \draw (0,-5) arc (90:-90:2 and 2); 
      \draw (0,9) arc (90:270:9 and 9); 
      \draw (0,5) arc (90:270:5 and 5); 
  \end{tikzpicture}
  \\
  v_\pm \o v_\mp &\mapsto X \\
v_\pm \o v_\pm & \mapsto 0
  \end{split}
\end{equation*}
\end{minipage}

\vskip1em
 
\noindent\begin{minipage}{.49\linewidth}
\begin{equation*}\label{eq:formula3}
  \begin{split}
   \begin{tikzpicture}[baseline=-2pt]
  \node at (0,0) {$\mathsmaller{\times}$};
      \draw (0,0) circle (8pt);
  \end{tikzpicture}
  & \to 
   \begin{tikzpicture}[baseline=-2pt]
      \node at (0,0) {$\mathsmaller{\times}$};
      \draw (0,0) circle (8pt);
      \draw (.8,0) circle (8pt);
  \end{tikzpicture}
  \\
v_- & \mapsto v_- \o X \\
v_+ & \mapsto v_+ \o X
  \end{split}
\end{equation*}
\end{minipage} 
\begin{minipage}{.4\linewidth}
\begin{equation*}\label{eq:formula4}
  \begin{split}
      \begin{tikzpicture}[x=1pt, y=1pt,baseline=-2pt]
  \node at (0,0) {$\mathsmaller{\times}$};
      \draw (0,9) arc (90:-90:2 and 2); 
      \draw (0,-5) arc (90:-90:2 and 2); 
      \draw (0,9) arc (90:270:9 and 9); 
      \draw (0,5) arc (90:270:5 and 5); 
  \end{tikzpicture}
  &\to
  \begin{tikzpicture}[baseline=-2pt]
  \node at (0,0) {$\mathsmaller{\times}$};
      \draw (0,0) circle (5pt);
      \draw (0,0) circle (9pt);
  \end{tikzpicture}
  \\
1 &\mapsto v_- \o v_+ + v_+ \o v_-\\
X & \mapsto 0
  \end{split}
\end{equation*}
\end{minipage}
\vskip1em

We now discuss the cube of resolutions and the (annular) Khovanov chain complex. 

\begin{definition}
\label{def:cube cat}
The $N$-dimensional \emph{cube category} $\cube^N$ has objects $\{0,1\}^N$ and a unique morphism $v=(v_1,\dots,v_N)\rightarrow w=(w_1,\dots,w_N)$ whenever $v_i\leq w_i$ for all $1\leq i \leq N$. We write $v\leq w$ if there is a morphism from $v$ to $w$. Moreover, setting $\lr{v} = \sum_i v_i$, we write $v \leq_k w$ if $v\leq w$ and $\lr{w} = \lr{v}+k$. An \emph{edge} is a morphism $v\to w$ such that $v\leq_1 w$. 
\end{definition}

Let $X$ denote either $\D$ or $\A$, and let $L \subset X$ be a diagram of an oriented link in $I\times X$.  Order its crossings $1,\ldots, N$, and let  $N_+$ and $N_-$ denote the number of positive and negative crossings, respectively. Following \cite{Bar-Natan}, we recall the \emph{cube of resolutions} of $L$, which is a functor $\cube^N \to \BN(X)$, and the Bar-Natan complex $[[L]]$, which is a chain complex over $\BN(X)$. 

\begin{figure}
    \centering
    \begin{tikzpicture}[scale=.5]
\draw[ ] (0,2) -- (2,0);

\draw[ , preaction={draw, line width=12pt, white}] (0,0) -- (2,2); 

\draw[->] (-.2,1) -- (-1,1);
\draw[->] (2.2,1) -- (3,1);

\node[draw=none, anchor = south] at (-.6,1) {$0$};
\node[draw=none, anchor = south] at (2.6,1) {$1$};

\draw[ ] (-1.2,0) .. controls (-2.2,1) .. (-1.2,2);
\draw[ ] (-3.2,0) .. controls (-2.2,1) .. (-3.2,2); 

\draw[ ] (3.2,0) .. controls (4.2,1) .. (5.2,0);
\draw[ ] (3.2,2) .. controls (4.2,1) .. (5.2,2); 

\end{tikzpicture}
    \caption{The two resolutions of a crossing.}
    \label{fig:resolutions}
\end{figure}
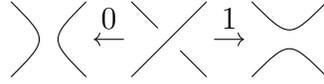

For $v\in \cube^N$, resolve each crossing according to Figure \ref{fig:resolutions} to obtain an object $L_v$ in $\BN(X)$. To each edge $v\leq_1 w$ we assign the cobordism $S_{v,w}: L_v\to L_w$ which is the evident saddle cobordism near the crossing where $v$ and $w$ differ and the identity cobordism elsewhere. To a general morphism $v\to w$,  decompose it in some way as a composition of edges and assign to it the composition of each edge saddle cobordism; the result is independent of this decomposition, and we have a commutative cube in $\BN(X)$. There is a way to assign $\eps_{v,w} \in \{\pm 1\}$ to each edge  $v\leq_1 w$ so that multiplying the edge map $S_{v,w}$ by $\eps_{v,w}$ results in a cube with anticommuting square faces; that is, for every $v \leq_1, u,u' \leq_1 w$ in $\cube^N$, we have $\eps_{u,w}\eps_{v,u} S_{u,w} S_{v,u} = -\eps_{u',w}\eps_{v,u'} S_{u',w} S_{v,u'}$
The \emph{Bar-Natan complex of $L$} is 
\[
[[L]]^i = \bigoplus_{\lr{v} = -i + N_-} L_v\{-i + N_+ - N_-\}.
\]

The differential on the summand corresponding to $v$ is the sum of edge maps $\eps_{v,w} S_{v,w}$ over all $w$ with $v\leq_1 w$. Anti-commutativity ensures that $[[L]]$ is a chain complex. By \cite[Theorem 1]{Bar-Natan}, if $L$ and $L'$ represent isotopic links in $I\times X$, then $[[L]]$ and $[L']]$ are chain homotopy equivalent. 
If $X=\D$, the \emph{Khovanov complex} $\KC(L)$ is obtained by applying the TQFT $\FKh$ to $[[L]]$. Likewise, if $X=\A$, then the \emph{annular Khovanov complex} $\aKC(L)$ is obtained by applying $\FA$ to $[[L]]$; in this case, the grading shifts $\{- \}$ are applied only to the quantum grading.

\begin{remark}
\label{rmk:homology vs cohomology}
    Let us explain the homological grading convention; see \cite[Section 2.10.1 and Remark 2.58]{LLS_Khovanov_Spectra} for an analogous discussion. The homological grading above follows the conventions in \cite{LLS_Khovanov_Spectra, LLS_CK_Spectra} and is obtained by negating the usual homological grading. In particular, the differential  decreases the homological grading. 

    This change is due to the fact that the \emph{homologies} of the spectra considered in this paper will recover the Khovanov invariants of links and Chen-Khovanov invariants of tangles. This is in contrast to other constructions \cite{LS-Kh-spectra, LLS-Burnside, AKW} where the \emph{cohomologies} of spectra recover the Khovanov invariants. Likewise, the lifts to the Burnside category discussed in Section \ref{sec:VHKKlax for D,A} are opposite to those in \cite{LLS-Burnside, HKK}. This will be relevant in Section \ref{sec:comparison} when comparing with \cite{AKW}. 
\end{remark}

\subsection{Platform algebras and tangle bimodules}
\label{sec:Platform algebras and bimodules}
In this section we use the Khovanov functor $\FKh$ to describe the platform algebras and tangle bimodules of Chen-Khovanov \cite{Chen-Khovanov} and Stroppel \cite{Stroppel} at the heart of the constructions to follow. These algebras were further extensively studied in a series of papers of Brundan-Stroppel, beginning with \cite{Brundan-Stroppel}.

Our convention is that tangles are read from left to right. In particular, an $(n,m)$-tangle has $n$ left endpoints and $m$ right endpoints. The composition $T_1 T_2$ of an $(n,m)$-tangle $T_1$ and an $(m,\ell)$-tangle $T_2$ is given by gluing $T_2$ to the right of $T_1$.

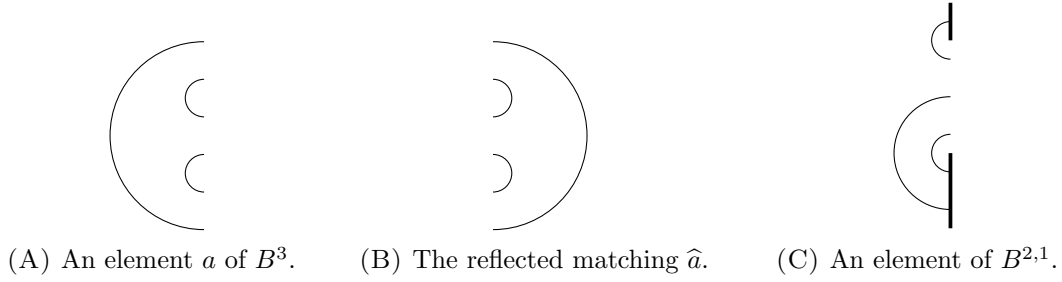
\begin{figure}
\centering
\subcaptionbox{An element $a$ of $\match{3}$.
\label{fig:matching}}[.3\linewidth]
{\begin{tikzpicture}[scale=.5]
    \coordinate (1) at (0,0);
    \coordinate (2) at (0,1);
    \coordinate (3) at (0,2);
    \coordinate (4) at (0,3);
    \coordinate (5) at (0,4);
    \coordinate (6) at (0,5);

  \draw[] (1) arc[start angle=-90, end angle=-270, radius=2.5];
  \draw[] (2) arc[start angle=-90, end angle=-270, radius=.5];
  \draw[] (4) arc[start angle=-90, end angle=-270, radius=.5];
\end{tikzpicture}
}
\subcaptionbox{The reflected matching $\flip{a}$.
\label{fig:flip}}[.3\linewidth]
{\begin{tikzpicture}[xscale=-1,scale=.5]
    \coordinate (1) at (0,0);
    \coordinate (2) at (0,1);
    \coordinate (3) at (0,2);
    \coordinate (4) at (0,3);
    \coordinate (5) at (0,4);
    \coordinate (6) at (0,5);

  \draw[] (1) arc[start angle=-90, end angle=-270, radius=2.5];
  \draw[] (2) arc[start angle=-90, end angle=-270, radius=.5];
  \draw[] (4) arc[start angle=-90, end angle=-270, radius=.5];
\end{tikzpicture}
}
\subcaptionbox{An element of $B^{2,1}$.
\label{fig:platform matching}}[.3\linewidth]
{\begin{tikzpicture}[yscale=-1,scale = .5]
\begin{scope}[rotate=-90]
    \coordinate (1) at (0,0);
    \coordinate (2) at (1,0);
    \coordinate (3) at (2,0);
    \coordinate (4) at (3,0);
    \coordinate (5) at (4,0);
    \coordinate (6) at (5,0);

  \draw[]  (1) arc[start angle=-180, end angle=0, radius=1.5];
  \draw[]  (2) arc[start angle=-180, end angle=0, radius=.5];
  \draw [] (5) arc[start angle=-180, end angle=0, radius=.5];
  
  \coordinate (l1) at (-.5,0); 
  \coordinate (r1) at (1.5,0);
  \coordinate (l2) at (4.5,0);
  \coordinate (r2) at (5.5,0);
    
    \draw[line width = .5mm] (l1) -- (r1);
    \draw[line width = .5mm] (l2) -- (r2);
  \end{scope}
\end{tikzpicture}
}
\caption{A depiction of Definition \ref{def:matchings and platforms}}\label{fig:matchings and platforms}
\end{figure}

\begin{definition}
\label{def:matchings and platforms}
Let $n\geq 0$. 
\begin{enumerate}
    \item  Let $\match{n}$ denote the set of crossingless matchings between $2n$ points on a vertical line. Matchings will be drawn in the left half plane. We can (and will) view elements of $\match{n}$ as planar $(0,2n)$-tangles.  See Figure \ref{fig:matching}.
    \item For $a\in \match{n}$, we let $\flip{a}$ denote the reflection of $a$ across the vertical line. See Figure \ref{fig:flip}. 
    \item For $0\leq k \leq n$,  $\cupdiag{n}{k}$ denotes the subset of $\match{n}$ in which there are no matchings between the top $k$ points and between the bottom $n-k$ points. Bold vertical lines called \emph{platforms} are drawn to indicate that two points on the same platform cannot be matched. Elements of $\cupdiag{n}{k}$ are called \emph{platform matchings}. See Figure \ref{fig:platform matching}.
\end{enumerate}
\end{definition}

We now define the \emph{platform algebra} $\Plat{n}$. First, consider the algebra
\begin{equation}
\label{eq:subalgebra of H^n}
\Platbark{n}{k} = \bigoplus_{a,b\in \cupdiag{n}{k}}\FKh\left(a \flip{b}\right)\{-n\}.
\end{equation}
Multiplication on $\FKh(a\flip{b})\{-n\} \otimes \FKh(c\flip{d})\{-n\}$ is zero unless $b=c$, in which case it is the map 
\[
\FKh(a\flip{b})\{-n\} \otimes \FKh(b\flip{d})\{-n\} \to \FKh(a\flip{d})\{-n\}
\]
induced by the minimal cobordism from $\flip{b} b$ to $2n$ parallel horizontal strands, consisting of $n$ saddle cobordisms. This is a subalgebra of Khovanov's arc algebra $H^n$\cite{Kh-ring_H}. 

Next, given $a, b\in \cupdiag{n}{k}$, distinguish $3$ types of circles in $a\flip{b}$:
\begin{enumerate}
    \item Type I circles do not intersect any platform. 
    \item Type II circles intersect at least one platform and intersect each platform at most once. 
    \item Type III circles intersect one of the platforms at least twice.
\end{enumerate}

Let $\Ideal{n}{k}_{a,b} \subset  \FKh(a\flip{b})$ be the following subgroup. If $a\flip{b}$ contains a type III circle, then $\Ideal{n}{k}_{a,b} = \FKh(a\flip{b})$. Otherwise, if $a\flip{b}$ does not contain any type III circles, then $\Ideal{n}{k}_{a,b} $ is generated by standard Khovanov generators in which every type II circle is labeled $X$. Set
\[
\Ideal{n}{k} = \bigoplus\limits_{a,b\in \cupdiag{n}{k}} \Ideal{n}{k}_{a,b}  \subset \Platbark{n}{k}.
\]
By \cite[Lemma 1]{Chen-Khovanov}, $\Ideal{n}{k}$ is a two-sided ideal of $\Platbark{n}{k}$. 
\begin{definition}\label{def:platform algebra}
Set $\Platk{n}{k} = \Platbark{n}{k}/\Ideal{n}{k}$,  and define the \emph{platform algebra} to be
\[
\Plat{n}:=\prod_{0\leq k \leq n} \Platk{n}{k}.
\]
\end{definition}

Note that $\Platk{n}{k}$ is a free abelian group with a distinguished basis consisting of the Khovanov generators of the diagrams $a\flip{b}$ (with $a, b\in \cupdiag{n}{k}$) which do not contain type III circles and such that each type II circle is labeled by $1$. Grading shifts ensure that $\Plat{n}$ is a graded ring supported in non-positive gradings.  Platform algebras also carry a \emph{weight}: an element in $\Platk{n}{k}$ is of weight $n-2k$.

The algebra $\Plat{n}$ is equipped with a complete set of idempotents, as defined in Section \ref{sec:linear cats}, which is in bijection with the platform matchings. Namely, for $0\leq k \leq n$ and $a\in \cupdiag{n}{k}$, the idempotent 
\[
e_a \in \Platk{n}{k}
\]
is (the image in the quotient of) the element of $\FKh(a\flip{a})$ in which each circle is labeled $1$. 

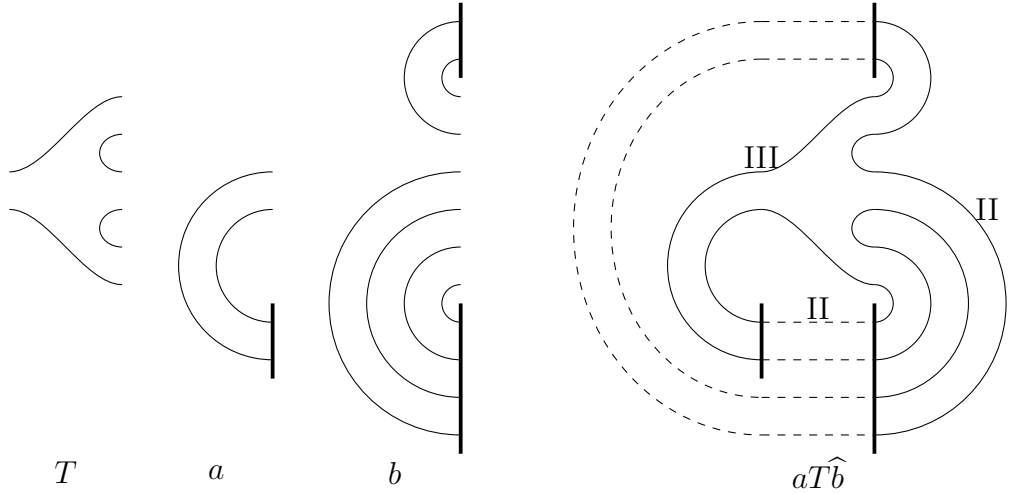
\begin{figure}
    \centering
    \begin{tikzpicture}

\draw (0,.5) .. controls (.5,.5)  and (1,-.5) .. (1.5,-.5); 

\draw (0,1) .. controls (.5,1)  and (1,2) .. (1.5,2); 

\draw (1.5,1.5) arc (90:270:.3 and .25);

\draw (1.5,.5) arc (90:270:.3 and .25);

\node at (.75,-3) {$T$};


\begin{scope}[shift={(3.5,0)}]
\draw[line width = .5mm] (0,-.75) --++ (0,-1); 


\draw (0,.5) arc (90:270:{1.5/2}  and {1.5/2});

\draw (0,1) arc (90:270: {2.5/2} and {2.5/2});

\node at (-.75,-3) {$a$};

\end{scope}

\begin{scope}[shift={(6,0)}]
\draw(0,-.5) arc  (90:270: {.5/2} and {.5/2});
\draw (0,-0) arc( 90:270: {1.5/2} and {1.5/2});
\draw (0,.5) arc( 90:270: {2.5/2} and {2.5/2});
\draw (0,1) arc( 90:270: {3.5/2} and {3.5/2});

\draw[line width = .5mm] (0,-.75) --++ (0,-2); 
\draw[line width = .5mm] (0,2.25) --++ (0,1); 

\draw (0,2.5) arc (90:270: .25 and {.5/2}); 
\draw (0,3) arc( 90:270: .75 and {1.5/2});

\node at (-{1.75/2},-3) {$b$};
\end{scope}


\begin{scope}[shift= {(10,0)}]
\draw (0,.5) .. controls (.5,.5)  and (1,-.5) .. (1.5,-.5); 

\draw (0,1) .. controls (.5,1)  and (1,2) .. (1.5,2); 

\draw (1.5,1.5) arc (90:270:.3 and .25);

\draw (1.5,.5) arc (90:270:.3 and .25);

\node at (.75,-3) {$aT\flip{b}$};

\draw[line width = .5mm] (0,-.75) --++ (0,-1); 
\draw[line width = .5mm] (1.5,-.75) --++ (0,-2); 
\draw[line width = .5mm] (1.5,2.25) --++ (0,1); 

\draw (1.5,2) arc (-90:90: .25 and {.5/2}); 
\draw (1.5,1.5) arc( -90:90: .75 and {1.5/2});

\draw(1.5,-1) arc  (-90:90: {.5/2} and {.5/2});
\draw (1.5,-1.5) arc( -90:90: {1.5/2} and {1.5/2});
\draw (1.5,-2) arc( -90:90: {2.5/2} and {2.5/2});
\draw (1.5,-2.5) arc( -90:90: {3.5/2} and {3.5/2});


\draw (0,.5) arc (90:270:{1.5/2}  and {1.5/2});

\draw (0,1) arc (90:270: {2.5/2} and {2.5/2});


\draw[dashed] (0,-1) --++ (1.5,0);
\draw[dashed] (0,-1.5) --++ (1.5,0);
\draw[dashed] (0,-2) --++ (1.5,0);
\draw[dashed] (0,-2.5) --++ (1.5,0);

\draw[dashed] (0,2.5) --++ (1.5,0);
\draw[dashed] (0,3) --++ (1.5,0);

\draw[dashed] (0,2.5) arc (90:270: 2 and {4.5/2});
\draw[dashed] (0,3) arc (90:270: 2.5 and {5.5/2});


\node[above] at (.75,-1.1) {II};

\node at (0,1.2) {III};

\node at (3,.5) {II};
\end{scope}
\end{tikzpicture}
    \caption{A planar $(2,6)$-tangle $T$, elements $a \in B^{2,0}$ and $b\in B^{4,2}$, and the closure $a T \flip{b}$. The dashed arcs in $a T \flip{b}$ are added to close the picture to a union of circles in the plane.  Here $k=0$ and $h(k)=2$ in the notations of Equation \eqref{eq:FCKBar(T,k)}. The type of each circle is indicated.}
    \label{fig:capping off platform tangle}
\end{figure}

Now let $T$ be a planar $(n,m)$-tangle. For  $k$ satisfying
\begin{equation}
    \label{eq:k range}
    \max\left(0, \frac{n-m}{2}\right) \leq k \leq \min\left(n,\frac{n+m}{2}\right),
\end{equation}
set $h(k)=k+\frac{m-n}{2}$ and define an $(\Platk{n}{k},\Platk{m}{h(k)})$-bimodule $\FCKk(T)$ as follows.\footnote{The notation $h(k)$ here follows \cite{LLS_CK_Spectra}; in \cite{Chen-Khovanov} the notation $l$ is used for the quantity $h(k)-k=\frac{m-n}{2}$.} First, consider the $(\Platbark{n}{k},\Platbark{m}{h(k)})$-bimodule $\FCKbar(T)$, which as an abelian group is given by
\begin{equation}\label{eq:FCKBar(T,k)}
\FCKbar(T) = \bigoplus_{\substack{a\in \cupdiag{n}{k} \\ b\in  \cupdiag{m}{h(k)} }} \FKh(a T \flip{b})\{\nmshift{n}{m}\}
\end{equation}
where
\begin{equation}
    \nmshift{n}{m}= -n - \max\left(0, \frac{n-m}{2} \right).
\end{equation}
Note that $\nmshift{n}{n}=-n$, which coincides with the shift for the corresponding platform algebra if $T$ is the identity tangle.  Here $a T \flip{b}$ is a collection of circles and arcs which by condition \eqref{eq:k range} can be closed to a collection of circles in a natural way; we continue to denote this closure by $a T \flip{b}$. See Figure \ref{fig:capping off platform tangle} for an illustration.  Then the left action by $\Platbark{n}{k}$ on the summand $\FKh(a T \flip{b})$ of $\FCKbar(T)$ is trivial unless we are acting by an element of $\Platbark{n}{k}e_a$, in which case the action is induced by the natural minimal saddle cobordism; the right action is analogous. 

For $a\in\cupdiag{n}{k} $ and  $b\in \cupdiag{m}{h(k)}$, consider the subgroup $\mathcal{I}(T;a,b) \subset \FKh(a T \flip{b})$ which is equal to $\FKh(a T \flip{b})$ if $a T \flip{b}$ contains a type III circle, and otherwise is generated by all elements of $\FKh(a T \flip{b})$ which label any type II circle by $X$ (see Figure \ref{fig:capping off platform tangle}). By \cite[Lemma 2, Lemma 3]{Chen-Khovanov},  
\[
\mathcal{I}^k(T) := \bigoplus_{\substack{a\in \cupdiag{n}{k} \\ b\in \cupdiag{m}{h(k)}}} \mathcal{I}(T;a,b)
\]
is a sub-bimodule of $\FCKbar(T)$, and the left (resp. right) action of $\Ideal{n}{k}$  (resp. $\Ideal{m}{h(k)}$) is trivial on $\FCKbar(T) /\mathcal{I}^k(T) $. 
\begin{definition}
\label{def:platform bimodule}
Set 
$\FCKk(T) = \FCKbar(T) /\mathcal{I}^k(T)$ and define  the $(\Plat{n}, \Plat{m}$)-bimodule
\[
\FCK(T): = \bigoplus_k \FCKk(T)
\]
where $k$ is in the range given in \eqref{eq:k range}.
\end{definition}

In \cite{Chen-Khovanov}, the gradings on platform algebras and bimodules are inherited from the gradings on arc algebras and bimodules. However, as noted in \cite[Section 5.5]{BPW}, this grading is not coherent with respect to tensor products of tangle bimodules.  The extra grading shift $\nmshift{n}{m}$ introduced into Equation \eqref{eq:FCKBar(T,k)} can be regarded as a fix for this issue as we show below (see also \cite[Equation (5.35)]{BPW}).

\begin{lemma}\label{lem:gluing planar tangles}
Let $T_1$ denote a planar $(p,n)$-tangle and $T_2$ denote a planar $(n,m)$-tangle.  Then there is an isomorphism of graded $(\Plat{p},\Plat{m})$-bimodules
\[\FCK(T_1)\otimes_{\Plat{n}}\FCK(T_2) \cong \FCK(T_1T_2).\]
\end{lemma}
\begin{proof}
As in \cite[Proposition 4]{Chen-Khovanov}, this isomorphism is given by the map
\[\FCK(T_1)\otimes_{\Plat{n}}\FCK(T_2) \to \FCK(T_1T_2)\]
induced by the minimal cobordisms $a T_1 \flip{b} \sqcup b T_2 \flip{c} \to a T_1T_2 \flip{c}$ consisting of one saddle for each component of $b$ as well as for any dashed arc added to both $b$ and $\flip{b}$.  The number of dashed arcs added to $b$ (respectively $\flip{b}$) is $\max(p-n,0)$ (respectively $\max(m-n,0)$), and thus the total number of saddles is
\[n+\min(\max(p-n,0),\max(m-n,0))= \min(\max(p,n),\max(n,m)).\]
After performing these saddles, any circles which are formed entirely from dashed arcs are removed from the diagram.  We may assume that each such circle is labelled by $1$, and so its removal is a degree $-1$ map; the total number of such circle removals is
\[\min\left(\max\left(\frac{n-p}{2},0\right),\max\left(\frac{n-m}{2},0\right)\right).\]
Example \ref{ex:tangle concatenation degree shifts} illustrates the above counts. It is therefore enough to compute
\[\nmshift{p}{n}+\nmshift{n}{m}+\min(\max(p,n),\max(n,m)) + \min\left(\max\left(\frac{n-p}{2},0\right),\max\left(\frac{n-m}{2},0\right)\right)\]
and ensure that it matches $\nmshift{p}{m}$ for all $p,n,m\geq 0$ of equal parity.  There are six cases to check depending on the relative ordering of $p,n,m$, which we leave to the reader.
\end{proof}

\begin{example}
\label{ex:tangle concatenation degree shifts}
We illustrate the proof of Lemma \ref{lem:gluing planar tangles} by two small examples. First, consider the following $(1,3)$-tangle $T_1$ and $(3,1)$-tangle $T_2$, with the closures $a T_1\flip{b}$ and $b T_2 \flip{a}$ as depicted below. 
\begin{center}
    \begin{tikzpicture}[scale=.7]

\draw (0,.5) .. controls (.5,.5)  and (1,-.5) .. (1.5,-.5); 

\draw (1.5,.5) arc (90:270:.25 and .25);

\node at (.75,-1) {$T_1$};

\begin{scope}[shift={(5,0)}, xscale=-1]
\draw (0,.5) .. controls (.5,.5)  and (1,-.5) .. (1.5,-.5); 

\draw (1.5,.5) arc (90:270:.25 and .25);

\node at (.75,-1) {$T_2$};
\end{scope}

\begin{scope}[shift={(8,0)}]
\draw (0,.5) .. controls (.5,.5)  and (1,-.5) .. (1.5,-.5); 

\draw (1.5,.5) arc (90:270:.25 and .25);

\draw[line width = .5mm] (1.5,-.75) --++ (0,-.5); 

\draw[line width = .5mm] (1.5,.75) --++ (0,1); 

\draw[line width = .5mm] (0,.75) --++ (0,.5); 


\draw (1.5,0) arc (-90:90:.25 and .25);

\draw (1.5,-.5) arc (-90:90:.75 and .75);

\draw (1.5,-1) arc (-90:90:1.25 and 1.25);

\draw (0,1) arc (90:270:.25 and .25);


\draw[dashed] (0,1) --++ (1.5,0);

\draw[dashed] (0,1.5) --++ (1.5,0);

\draw[dashed] (0,-1) --++ (1.5,0);

\draw[dashed] (0,1.5) arc (90:270:1.25 and 1.25);

\end{scope}


\begin{scope}[shift={(15,0)}, xscale=-1]

\draw (0,.5) .. controls (.5,.5)  and (1,-.5) .. (1.5,-.5); 

\draw (1.5,.5) arc (90:270:.25 and .25);

\draw[line width = .5mm] (1.5,-.75) --++ (0,-.5); 

\draw[line width = .5mm] (1.5,.75) --++ (0,1); 

\draw[line width = .5mm] (0,.75) --++ (0,.5);


\draw (1.5,0) arc (-90:90:.25 and .25);

\draw (1.5,-.5) arc (-90:90:.75 and .75);

\draw (1.5,-1) arc (-90:90:1.25 and 1.25);

\draw (0,1) arc (90:270:.25 and .25);


\draw[dashed] (0,1) --++ (1.5,0);

\draw[dashed] (0,1.5) --++ (1.5,0);

\draw[dashed] (0,-1) --++ (1.5,0);

\draw[dashed] (0,1.5) arc (90:270:1.25 and 1.25);

\end{scope}

\end{tikzpicture}
\end{center}
\noindent
The cobordism from $\flip{b} \sqcup b$ to the six-strand identity tangle $\id_6$ consists of three saddles and results in the diagram   
\begin{center}    \begin{tikzpicture}[scale=.7]

\begin{scope}
\draw (0,.5) .. controls (.5,.5)  and (1,-.5) .. (1.5,-.5); 

\draw (1.5,.5) arc (90:270:.25 and .25);



\draw[line width = .5mm] (0,.75) --++ (0,.5);

\draw (0,1) arc (90:270:.25 and .25);


\draw[dashed] (0,1) --++ (1.5,0);

\draw[dashed] (0,1.5) --++ (1.5,0);

\draw[dashed] (0,-1) --++ (1.5,0);

\draw[dashed] (0,1.5) arc (90:270:1.25 and 1.25);

\end{scope} 

\begin{scope}[shift={(3,0)},xscale=-1]

\draw (0,.5) .. controls (.5,.5)  and (1,-.5) .. (1.5,-.5); 

\draw (1.5,.5) arc (90:270:.25 and .25);



\draw[line width = .5mm] (0,.75) --++ (0,.5);


\draw (0,1) arc (90:270:.25 and .25);


\draw[dashed] (0,1) --++ (1.5,0);

\draw[dashed] (0,1.5) --++ (1.5,0);

\draw[dashed] (0,-1) --++ (1.5,0);

\draw[dashed] (0,1.5) arc (90:270:1.25 and 1.25);

\end{scope}

\end{tikzpicture}
\end{center}
\noindent
whose outermost dashed circle (necessarily labeled by $1$) must be removed. The total degree of this map before the grading shifts $\nmshift{*}{*}$ are applied is $2$. We have $\nmshift{1}{3} + \nmshift{3}{1} = -5$, $\nmshift{3}{3} = -3$, so that this map is degree preserving after introducing grading shifts.

As another example, consider the following tangles $T_1, T_2$ and the depicted choice of their closures. 
\begin{center}\begin{tikzpicture}[scale=.7]


\draw (0,.5) .. controls (.5,.5)  and (1,-.5) .. (1.5,-.5); 

\draw (1.5,.5) arc (90:270:.25 and .25);

\node at (.75,-1.5) {$T_1$};

\begin{scope}[shift={(4,0)}]

\draw (0,-.5) .. controls (.5,-.5)  and (1,0) ..  (1.5,0); 

\draw (0,.5) arc (90:-90:.25 and .25);

\draw (1.5,-.5) arc (90:270:.25 and .25);

\draw (1.5,1) arc (90:270:.25 and .25);

\node at (.75,-1.5) {$T_2$};
\end{scope}


\begin{scope}[shift={(9,0)}]

\draw (0,.5) .. controls (.5,.5)  and (1,-.5) .. (1.5,-.5); 

\draw (1.5,.5) arc (90:270:.25 and .25);

\draw[line width = .5mm] (1.5,-1.25) --++ (0,-.5); 

\draw[line width = .5mm] (1.5,1.25) --++ (0,1); 

\draw[line width = .5mm] (0,1.25) --++ (0,.5);


\draw (1.5,-.5) arc (-90:90:.25 and .25);

\draw (1.5,.5) arc (-90:90:.5 and .5);

\draw (1.5,-1.5) arc (-90:90:1.5 and 1.75);

\draw (0,1.5) arc (90:270:.5 and .5);


\draw[dashed] (0,1.5) --++ (1.5,0);

\draw[dashed] (0,2) --++ (1.5,0);

\draw[dashed] (0,-1.5) --++ (1.5,0);

\draw[dashed] (0,2) arc (90:270:1.5 and 1.75);

\end{scope}

\begin{scope}[shift={(16,0)}]

\draw (0,-.5) .. controls (.5,-.5)  and (1,0) ..  (1.5,0); 

\draw (0,.5) arc (90:-90:.25 and .25);

\draw (1.5,-.5) arc (90:270:.25 and .25);

\draw (1.5,1) arc (90:270:.25 and .25);

\draw[line width = .5mm] (1.5,-1.25) --++ (0,-1); 

\draw[line width = .5mm] (1.5,1.25) --++ (0,1.5); 

\draw[line width = .5mm] (0,1.25) --++ (0,1); 

\draw[line width = .5mm] (0,-1.25) --++ (0,-.5); 


\draw (1.5,.5) arc (-90:90:.25 and .25);

\draw (1.5,-1) arc (-90:90:.25 and .25);

\draw (1.5,0) arc (-90:90:.75 and .75);

\draw (1.5,-1.5) arc (-90:90:1.75 and 1.75);

\draw (1.5,-2) arc (-90:90:2.25 and 2.25);

\draw (0,1.5) arc (90:270:.5 and .5);

\draw (0,0) arc (90:270:.25 and .25);

\draw (0,2) arc (90:270:1.5 and 1.75);

\draw[dashed] (0,2.5) --++ (1.5,0);

\draw[dashed] (0,2) --++ (1.5,0);

\draw[dashed] (0,1.5) --++ (1.5,0);

\draw[dashed] (0,-1.5) --++ (1.5,0);

\draw[dashed] (0,-2) --++ (1.5,0);

\draw[dashed] (0,2.5) arc (90:270:2 and 2.25);

\end{scope} 

\end{tikzpicture}
\end{center}
The cobordism realizing the isomorphism of Lemma \ref{lem:gluing planar tangles} in this case again consists of three saddles, resulting in 
\begin{center}
\begin{tikzpicture}[scale=.7]

\draw (0,-.5) .. controls (.5,-.5)  and (1,0) ..  (1.5,0); 

\draw (0,.5) arc (90:-90:.25 and .25);

\draw (1.5,-.5) arc (90:270:.25 and .25);

\draw (1.5,1) arc (90:270:.25 and .25);

\draw[line width = .5mm] (1.5,-1.25) --++ (0,-1); 

\draw[line width = .5mm] (1.5,1.25) --++ (0,1.5); 



\begin{scope}[shift={(-1.5,0)}]
    \draw (0,.5) .. controls (.5,.5)  and (1,-.5) .. (1.5,-.5); 

\draw (1.5,.5) arc (90:270:.25 and .25);

\draw (0,1.5) arc (90:270:.5 and .5);
\draw[line width = .5mm] (0,1.25) --++ (0,.5); 

\draw[dashed] (0,2.5) --++ (1.5,0);

\draw[dashed] (0,2) --++ (1.5,0);

\draw[dashed] (0,1.5) --++ (1.5,0);

\draw[dashed] (0,-1.5) --++ (1.5,0);

\draw[dashed] (0,-2) --++ (1.5,0);

\draw[dashed] (0,2.5) arc (90:270:2.25 and 2.25);

\draw[dashed] (0,2) arc (90:270:1.75 and 1.75);
\end{scope}

\draw (1.5,.5) arc (-90:90:.25 and .25);

\draw (1.5,-1) arc (-90:90:.25 and .25);

\draw (1.5,0) arc (-90:90:.75 and .75);

\draw (1.5,-1.5) arc (-90:90:1.75 and 1.75);

\draw (1.5,-2) arc (-90:90:2.25 and 2.25);

\draw[dashed] (0,2.5) --++ (1.5,0);

\draw[dashed] (0,2) --++ (1.5,0);

\draw[dashed] (0,1.5) --++ (1.5,0);

\draw[dashed] (0,-1.5) --++ (1.5,0);

\draw[dashed] (0,-2) --++ (1.5,0);

\end{tikzpicture},
\end{center}
and no circles are removed in the process. The total degree of this map before degree shifts is $3$, and the shifts $\nmshift{1}{3} + \nmshift{3}{5} = -4$, $\nmshift{1}{5} = -1$ ensure that this map is degree preserving.  
\end{example}

By \cite[Proposition 1]{Chen-Khovanov}, if $T_0$ and $T_1$ are planar $(n,m)$-tangles and $S:T_0 \to T_1$ is a (possibly dotted) tangle cobordism, then there is an induced bimodule map $\FCK(S) : \FCK(T_0) \to \FCK(T_1)$ of degree
\begin{equation}
    \label{eq:degree of tangle cob}
     \chi(S) -  \frac{n+m}{2},
\end{equation}
where $\chi(S)$ is as defined in \eqref{eq:top deg of dotted cob}. We call \eqref{eq:degree of tangle cob} the \emph{degree} of $S$. A dot on $S$ corresponds to multiplication by $X$.   Degree is additive under composition of cobordisms. 

Now let $T$ be a (not necessarily planar) $(n,m)$-tangle diagram. Analogously to the construction in Section \ref{sec:The Khovanov and annular Khovanov functors}, we form the cube of resolutions of $T$, collapse to a chain complex, and apply $\FCK$ to obtain a chain complex $\CCK(T)$ of  graded $(\Plat{n}, \Plat{m})$-bimodules. By \cite{Chen-Khovanov}, if $T$ and $T'$ represent isotopic $(n,m)$-tangles, then  $\CCK(T)$ and $\CCK(T')$ are chain homotopy equivalent. 

For $0\leq k \leq n$ and $a,b\in \cupdiag{n}{k}$, we set $\Plat{n}(a,b) = e_a \Plat{n} e_b$. More generally, consider an $(n,m)$-tangle diagram $T$, $k$ in the range \eqref{eq:k range}, and  $a\in \cupdiag{n}{k}$, $b\in  \cupdiag{m}{h(k)}$. Let $\CCKnk(T;a,b)$ denote the subcomplex $e_a \CCK(T) e_b$ of $\CCK(T)$.
When $a=b$ (so necessarily $n=m$) we will write $\CCK(T;a)$ rather than $\CCKnk(T;a,a)$. We also define 
\[
\CCKnk(T) = \bigoplus_{\substack{a\in \cupdiag{n}{k} \\ b\in  \cupdiag{m}{h(k)} }} \CCKnk(T;a,b).
\]
We will also need versions of the above constructions with scalars extended to $\Zq$.  We set the convention that, if $X$ is any of the algebraic constructions in this section, then 
\begin{equation}
\label{eq:quantum extensions of CK}
    _qX=X\otimes \Zq
\end{equation}
denotes the extension of scalars to $\Zq$, analogous to the notation in \eqref{eq:scalar ext for Sq}. For instance,  $\qPlat{n} = \Plat{n} \o \Zq$.

\subsection{Quantum annular homology}
\label{sec:quantum annular homology defs}
We now outline the construction of Beliakova-Putyra-Wehrli's quantum annular link homology \cite{BPW}. The theory is built over a commutative ring $\k$ and a choice of unit $q\in \k$. We set $\k=\Zq$ for some fixed cyclic group $G=\langle q\rangle$, with the generator $q$ providing the distinguished unit $q\in \Zq$.  We begin with a definition and notational convention.

\begin{definition}\label{def:quantum annulus}
The notation $\Aq$ will indicate the data of the annulus $\A$ together with a fixed ray $\mq$ emanating from the puncture, called the \emph{seam}.  We fix the co-orientation on $\mq$, indicated with a small arrow, to be induced by the clockwise orientation on the core circle of $\A$.  The seam then determines the \emph{membrane} $\Mq:= I \times \mq \subset I \times \Aq$ with similar co-orientation.  By convention, closed 1-manifolds $Z\subset\Aq$ and link diagrams $L\subset\Aq$ are assumed to be transverse to $\mq$, with crossings disjoint from $\mq$.  Similarly, cobordisms $\Sigma\subset I\times \Aq$ are assumed to be transverse to $\Mq$, and dots on $\Sigma$ are disjoint from $\Mq$.  We will continue to use the notation $\Gen(Z)$ for the set of \emph{Khovanov generators} of a collection of circles $Z\subset\Aq$ as in Definition \ref{def:FA on BN(A)}.
\end{definition}

\begin{figure}
    \centering
    \includegraphics[scale=1]{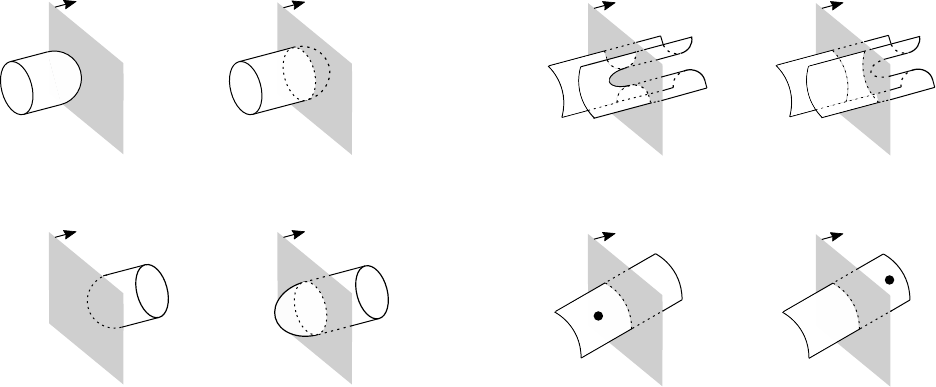}
        \put(-380,140){$=\,\q^{-1}$}
    \put(-360,40){$=\,\q$}
    \put(-102,140){$=\,\q$}
    \put(-102,40){$=\,\q^2$}
    \caption{The trace relations in $\BNq$.}\label{fig:trace moves}
\end{figure}

We will draw $\mq$ as a dashed line and $\Mq$ as a gray plane; Figures \ref{fig:trace moves} and \ref{fig:saddles through seam} illustrate some of these notions.  We now recall a deformation of the Bar-Natan category of the annulus $\BN(\A)$ from \cite{BPW}.

\begin{definition}
\label{def:qBN cat}
The \emph{quantum Bar-Natan category of the annulus}, denoted $\BNq$, is defined as follows.
\begin{itemize}
    \item Objects of $\BNq$ are formal direct sums of formally graded collections of circles $Z\subset\Aq$  transverse to $\mq$.
    \item Morphisms in $\BNq$ are matrices whose entries are $\Zq$-linear combinations of cobordisms in $I\times \Aq$ (each of which is transverse to $\Mq$).  Isotopic cobordisms are identified if the isotopy fixes the membrane $\Mq$ set-wise. Otherwise, the cobordisms are scaled by a power of $q$ according to the degree (Equation \eqref{eq:degree of tangle cob}) of the part of the cobordism that passes through the membrane during the isotopy, accounting also for the coorientation of the membrane. The generating relations are depicted  in Figure $\ref{fig:trace moves}$; for details see \cite[Section 6.2]{BPW}. Moreover, the Bar-Natan relations of Figure \ref{fig:BN relations} are imposed, where the local pictures are understood to be disjoint from $\Mq$. 
\end{itemize}
\end{definition}

By general position, if two cobordisms $\Sigma,\Sigma'\subset I\times\Aq$ are isotopic, then they are related by a sequence of isotopies fixing the membrane and the generating moves in Figure \ref{fig:trace moves}, and thus $\Sigma = q^k \Sigma'$ as morphisms in $\BNq$, for some $k\in \Z$; see also \cite[Proposition 6.2]{BPW}.

We are now ready to define the quantum annular TQFT 
\[
\FAq : \BNq \to \ggRmod[\Zq].
\]
Let $Z\subset \Aq$ be an object in $\BNq$, and let $n$ denote the number of intersection points of $Z$ with $\mq$.  Cut $Z$ open along $\mq$ to obtain a planar $(n,n)$-tangle $Z^{\rm cut}$, and form the $(\qPlat{n},\qPlat{n})$-bimodule $\qFCK(Z^{cut})$ (recall the scalar extension notation of Equation \eqref{eq:quantum extensions of CK}).  Then $\FAq$ is defined on objects by
\[
\FAq(Z) : = \qhh_*(\qPlat{n}; \qFCK(Z^{cut})).
\]
By \cite[Proposition 6.6]{BPW},  $\qhh_i(\qPlat{n}; \qFCK(Z^{cut})) = 0$
for $i>0$. The $\Zq$-module $\FAq(Z)$ carries a bigrading. The quantum grading $\lr{\cdot}$ is inherited from the quantum grading on Chen-Khovanov bimodules, while the annular grading $\lr{\cdot}_\A$ is inherited from the weight grading of Chen-Khovanov bimodules; that is, the image in $\FAq(Z)$ of an element in $\qFCKk(Z^{cut})$ is in annular grading $n-2k$.

The full construction of $\FAq$ in \cite{BPW} follows from a general theory of (twisted) horizontal traces of bicategories, which we will not describe. The definition of $\FAq$ on morphisms follows from this general theory as well. The main computational tool is \cite[Theorem 6.3]{BPW}, which allows one to bypass this categorical machinery. Instead of describing \cite[Theorem 6.3]{BPW} in detail, we record key properties of $\FAq$, beginning with a more explicit description of the value of $\FAq$ on certain objects.

\begin{definition}
\label{def:standard circles}
A collection of circles $Z\subset \Aq$ is \emph{standard} if every component of $Z$ intersects $\mq$ at most once. 
\end{definition}

If $Z\subset \Aq$ is standard, with $n$ essential and $t$ trivial circles, then $Z^{cut}$ consists of the identity $(n,n)$-tangle together with $t$ trivial circles which can be delooped, resulting in a distinguished isomorphism\footnote{We often do not keep track of bigradings to simplify the notation.}  of $\Zq$-modules 
\[
\FAq(Z) \cong \qhh_0(\qPlat{n}; \qPlat{n}) \otimes  (\Zq\oplus \Zq)^{\o t}
\]
where we write  $\o$ to denote $\o_{\Zq}$ for simplicity. Let $\qPlatzero{n} \subset \qPlat{n}$ denote the subalgebra consisting of elements of quantum degree $0$. By \cite[Proposition 6.6]{BPW}, the inclusion $\qPlatzero{n} \hookrightarrow \qPlat{n}$ induces an isomorphism 
\[
\qhh_0(\qPlatzero{n}; \qPlatzero{n}) \cong \qhh_0(\qPlat{n}; \qPlat{n})
\]
Moreover, $\qPlatzero{n}$ is freely generated over $\Zq$ by the $2^n$ distinguished mutually orthogonal idempotents $\{e_a\}$ parametrized by platform matchings. It follows from the definition of $\qhh_0$ that $\qhh_0(\qPlatzero{n};\qPlatzero{n})$ is a free $\Zq$-module of rank $2^n$, with basis given by the images of the $e_a$. In fact, there is a distinguished isomorphism
\begin{equation}
\label{eq:distinguished iso for all standard}
\FAq(Z) \cong \qhh_0(\qPlatzero{n}; \qPlatzero{n}) \otimes  (\Zq\oplus \Zq)^{\o t} \cong \Zq\langle \Gen(Z) \rangle,
\end{equation}
where $\Zq\langle\Gen(Z)\rangle$ denotes the free $\Zq$-module on the set $\Gen(Z)$ of Khovanov generators of Definition \ref{def:FA on BN(A)}.  We will freely make use of Equation \eqref{eq:distinguished iso for all standard} when $Z$ is standard; in particular the elements of $\FAq(Z)$ which correspond to elements of $\Gen(Z)$ will also be referred to as Khovanov generators of $\FAq(Z)$.\footnote{Note that a Khovanov generator of $Z$ corresponds to a linear combination of idempotents $e_a\in\qPlatzero{n}$, see \cite[Section 6]{Chen-Khovanov}.}

As explained in the proof of \cite[Corollary 6.7]{BPW}, if $Z\subset \Aq$ is any collection of $n$ circles (not necessarily standard), then $\FAq(Z)$ is a free $\Zq$-module of rank $2^n$; however, there is not a canonical choice of basis in general.  We now describe how we fix a choice of bases, as will be necessary for the lift to the spectral setting; see also \cite[Section 2.3]{AKW}.  Note that any $Z\subset\Aq$ is isotopic to a standard collection of circles, which we denote by $Z^\st$, that is unique up to planar isotopy of the cut-open planar tangle $(Z^\st)^{cut}$.

\begin{definition}
\label{def:choice of generators}
A \emph{regular} cobordism in $I\times\Aq$ is one which contains no critical points with respect to the projection to $I$. Note, every component of a regular cobordism is an annulus with one boundary component in $\{0\}\times \Aq$ and the other boundary component in $\{1\}\times \Aq$.

A \emph{choice of generators}, denoted $\choice$, consists of a regular cobordism $\choice_Z$ from  $Z^\st$ to $Z$ for each collection of circles  $Z\subset \Aq$ and satisfies the following.  
\begin{itemize}
    \item If $Z$ is standard then $\choice_Z$ is the identity cobordism.
    \item If a component of $Z$ is a standard essential circle $C$, then $\choice_Z$ is the disjoint union of $\choice_X$, where $X = Z \setminus C$,  and the identity cobordism on $C$. 
\end{itemize}
\end{definition}

\begin{remark}
    In \cite[Definition 2.3]{AKW} the above two bullet points were not required. We include them in the present article  for convenience, but what follows could be restated without them. We also note that a choice of generators can be obtained by taking the cobordism formed by a choice of isotopy from $Z^\st$ to $Z$; if $Z$ contains a standard essential circle $C$ then this isotopy can be chosen to fix $C$ pointwise.
\end{remark}

By \cite[Lemma 2.6]{AKW}, the cobordism $\choice_Z$ is an isomorphism in $\BNq$. In fact, if we let $\b{\choice_Z}$ denote its reflection in the $I$ coordinate, then in $\BNq$ we have $(\choice_Z)^{-1} = q^k \b{\choice_Z}$ for some $k\in \Z$. In particular, $\choice_Z$ induces an isomorphism 
\[
\FAq(\choice_Z) : \FAq(Z^\st) \rar{\sim} \FAq(Z).
\]

We now use a choice of generators $\choice$ to construct a functor $\choiceFAq : \BNq \to \ggRmod[\Zq]$ that is naturally isomorphic to $\FAq$ but which will assign a $\Z[G]$-module equipped with a choice of basis to each collection of circles $Z\subset \Aq$. This will be used later in Section \ref{sec:qVHKKlax} when lifting to the Burnside category; see also Section \ref{sec:Burnside}.

\begin{definition}
\label{def:FAq choice}
Let $\choice$ be a choice of generators. For $Z\subset \Aq$, let $\choiceFAq(Z)$ be the free $\Z[G]$-module generated by $\Gen(Z)$.  The  quantum and annular degree of $z\in \Gen(Z)$ is as in Equation \eqref{eq:degrees in QAH}.  Before defining $\choiceFAq$ on cobordisms, note that $\choice$ provides an isomorphism 
\[
\beta^\choice_Z: \choiceFAq(Z) \to \FAq(Z)
\]
as follows. Let $y\in \choiceFAq(Z)$ be a basis element (i.e. a Khovanov generator). The cobordism $\choice_Z$ provides a bijection between circles of  $Z^\st$ and of $Z$, so that $y$ corresponds to a Khovanov generator $y^\st\in \FAq(Z^\st)$ using Equation \eqref{eq:distinguished iso for all standard}. Then $\beta^\choice_Z(y)$ is defined to be the image of $y^\st$ under the isomorphism $\FAq(\choice_Z)$.  Note that, if $Z$ is already standard, then $\choice_Z=\id_Z$ by definition, so that $\beta^\choice_Z$ is the identification of Equation \eqref{eq:distinguished iso for all standard}.

For a cobordism $\Sigma : Z_0 \to Z_1$ in $I\times\Aq$, we define $\choiceFAq(\Sigma)$ to be the unique map fitting into the diagram 
\begin{equation}
\label{eq:FAq choices square}
\begin{tikzcd}[column sep = large]
    \choiceFAq(Z_0) \ar[d,"\beta^\choice_{Z_0}"]  \ar[r, "\choiceFAq(\Sigma)"] & \choiceFAq(Z_1) \ar[d, "\beta^\choice_{Z_1}"] \\
    \FAq(Z_0) \ar[r,"\FAq(\Sigma)"] & \FAq(Z_1).
\end{tikzcd}
\end{equation}
The maps $\FAq(\Sigma)$ and $\choiceFAq(\Sigma)$ have quantum and annular degrees $\chi(\Sigma)$ and $0$, respectively, where $\chi(\Sigma)$ is as defined in \eqref{eq:top deg of dotted cob}. 
\end{definition}

In other words, $\choice$ fixes a basis for each $\FAq(Z)$, and $\choiceFAq(\Sigma)$ computes the map assigned to $\Sigma$ by  $\FAq$ in terms of these fixed bases. In practice, $\choiceFAq(\Sigma)$ can be computed using the Bar-Natan relations in Figure \ref{fig:BN relations} (understood to occur away from the membrane), the relations in Figure \ref{fig:trace moves}, and \cite[Theorem 6.3]{BPW}. For some explicit examples, we refer the reader to \cite[Section 2.4]{AKW}. For an annular link diagram $L \subset \Aq$, let $\qaKC(L)$ and $\qaKCchoice(L)$ denote the chain complex of bigraded free $\Zq$-modules obtained by applying $\FAq$ and $\choiceFAq$, respectively, to the Bar-Natan complex of $L$ (see Section \ref{sec:The Khovanov and annular Khovanov functors}). For $j\in \Z$, we denote by  $\qaKC(L;j)$ and $\qaKCchoice(L;j)$ the subcomplexes generated by annular degree $j$ Khovanov generators. We let $\qaKh(L)$, $\qaKhchoice(L)$, $\qaKh(L;j)$, and $\qaKhchoice(L;j)$ denote the homology of these complexes. Note that $\FAq$ and $\choiceFAq$ are naturally isomorphic by construction, so $\qaKC(L)$ and $\qaKCchoice(L)$ are isomorphic chain complexes.

The remainder of this subsection will focus on key properties of $\FAq$ and of $\choiceFAq$ on morphisms.  To begin, if $\Sigma \subset I \times \Aq$ is a closed surface, define its \emph{winding number} $w(\Sigma) \in \Z_{\geq 0}$ to be the (non-negative) generator of the image of the map $H_1(\Sigma;\Z ) \to H_1(I\times \A;\Z) \cong \Z.$ The following may be viewed as a generalization of \cite[Corollary 2.18]{AKW}.

\begin{lemma}
 \label{lem:torus}
 Let $\Sigma \subset I\times \Aq$ be an embedded torus with no dots. Then in $\BNq$, we have $\Sigma = q^\ell (1 + q^{2 w(\Sigma)})$ for some $\ell\in \Z$. 
 \end{lemma}
\begin{proof}
  Consider the number of circles in $\Sigma \cap \Mq$. If there are none, then $w(\Sigma)=0$ and $\Sigma = 2$.  Otherwise, let $C$ be an innermost circle in $\Sigma \cap \Mq$, in the sense that $C$ bounds a disk $D^2 \subset \Mq$ such that $D^2\cap \Sigma = C$. If $C$ also bounds a disk in $\Sigma$, then we can perform an isotopy to remove $C$ from $\Sigma \cap \Mq$ at the cost of a factor of $q^{k_1}$ for some $k_1\in \Z$. We continue this procedure until all circles in the intersection are removed or until we reach an innermost circle $C$ in $\Sigma \cap \Mq$ which is nontrivial in $\Sigma$. In the former case,  $\Sigma$ lies within a ball in $I\times \A$ so that $w(\Sigma) =0$ and $\Sigma$ evaluates to $2q^k$ for some $k\in \Z$.   In the latter case, perform neck-cutting along a small pushoff of $C$ to write $\Sigma = \Sigma_1 + \Sigma_2$ as in \eqref{eq:neck cutting winding number}.
  \begin{equation}
      \begin{aligned}
      \label{eq:neck cutting winding number}
\begin{tikzpicture}[xscale = .9]

\begin{scope}[shift={(0,.25)}]
\draw[>=Stealth,->]  (.5,1.7) -- (.8,1.75); 

\draw[white,fill=gray!40] (-.5, 2) -- (2,1.25) -- (2,-1) -- (-.5, -.3) --  (-.5, -.15); 

\end{scope}

\draw[white, fill=white] (-.5,.5) -- (.5,.5) -- (.5,1.5) -- (-.5,1.5) -- cycle;

\draw[densely dashed] (.5,.5) --++ (1,0);
\draw[densely dashed] (.5,1.5) --++ (1,0);

\begin{scope}[shift={(-1,0)}]
\end{scope}

\draw[] (.5,.5) --++ (-2,0); 
\draw[] (.5,1.5) --++ (-2,0);  
\draw[fill=white] (.5,1.5) arc (90:-90:.5 and .5); 
\draw[densely dashed]  (.5,.5) arc (-90:--270:.5 and .5);

\draw[]  (-1.5,.5) arc (-90:--270:.5 and .5);

\node at (.5,-1) {$\Sigma$};


\begin{scope}[shift={(5,0)}]
    \begin{scope}[shift={(0,.25)}]
\draw[>=Stealth,->]  (.5,1.7) -- (.8,1.75); 

\draw[white,fill=gray!40] (-.5, 2) -- (2,1.25) -- (2,-1) -- (-.5, -.3) --  (-.5, -.15); 

\end{scope}

\draw[white, fill=white]  (.5,1.5) arc (90:270:.7 and .5);

\draw (.5,1.5) arc (90:270:.7 and .5);

\draw[densely dashed] (.5,.5) --++ (1,0);
\draw[densely dashed] (.5,1.5) --++ (1,0);

\begin{scope}[shift={(-1,0)}]
\end{scope}

\draw[fill=white] (.5,1.5) arc (90:-90:.5 and .5); 
\draw[densely dashed]  (.5,.5) arc (-90:--270:.5 and .5);

\draw[]  (-1.5,.5) arc (270:450:.8 and .5);
\draw[]  (-1.5,.5) arc (-90:--270:.5 and .5);

\node at (-.85,1) {$\bullet$};

\node at (.5,-1) {$\Sigma_1$};

\end{scope}


\begin{scope}[shift={(10,0)}]
    \begin{scope}[shift={(0,.25)}]
\draw[>=Stealth,->]  (.5,1.7) -- (.8,1.75); 

\draw[white,fill=gray!40] (-.5, 2) -- (2,1.25) -- (2,-1) -- (-.5, -.3) --  (-.5, -.15); 

\end{scope}

\draw[white, fill=white]  (.5,1.5) arc (90:270:.7 and .5);

\draw (.5,1.5) arc (90:270:.7 and .5);

\draw[densely dashed] (.5,.5) --++ (1,0);
\draw[densely dashed] (.5,1.5) --++ (1,0);

\begin{scope}[shift={(-1,0)}]
\end{scope}

\draw[fill=white] (.5,1.5) arc (90:-90:.5 and .5); 
\draw[densely dashed]  (.5,.5) arc (-90:--270:.5 and .5);

\draw[]  (-1.5,.5) arc (270:450:.8 and .5);
\draw[]  (-1.5,.5) arc (-90:--270:.5 and .5);

\node at (.5,1) {$\bullet$};

\node at (.5,-1) {$\Sigma_2$};

\end{scope}

\node at (2.5,1) {$=$};

\node at (7.5,1) {$+$};
\end{tikzpicture}
      \end{aligned}
  \end{equation}
  Since $\Sigma_1$ is a once-dotted $2$-sphere, we have $\Sigma_1= q^{k_2}$ in $\BN(\Aq)$ for some $k_2\in \Z$. Pick a curve $\gamma \subset \Sigma_2$ which intersects $C$ in one point and which connects the dot on $\Sigma_2$ to the point on $\Sigma_2$ corresponding to the dot on $\Sigma_1$. We may also view $\gamma$ as a curve on $\Sigma$, which we complete to a closed loop $\b{\gamma}\subset \Sigma$ that intersects $C$ in one point. Then $\b{\gamma}$ is sent to $ w(\Sigma)$ under $H_1(\Sigma;\Z ) \to H_1(I\times \A;\Z)$. From the dot sliding relation in Figure \ref{fig:trace moves}, by moving the dot on $\Sigma_2$ along $\gamma$, we see that $\Sigma_2 = q^{2 s} \Sigma_1$ where $s$ is the oriented intersection number between $\gamma$ and $\Mq$. On the other hand, $s = w(\Sigma)$ since oriented intersection is preserved under homotopy, which completes the proof. 
\end{proof}

\begin{proposition}
\label{prop:FAq properties}
Fix a choice of generators $\choice$, let $\Sigma : Z_0 \to Z_1$ be a connected cobordism in $I\times \Aq$, and consider a generator $x\in \Gen(Z_0)$.  We may view $\Sigma$ as a cobordism in $I\times\A$ and $x$ as a Khovanov generator of $\FA(Z_0)$ to write $\displaystyle \FA(\Sigma)(x) = \sum_{y\in\Gen(Z_1)} n_y y$, with each $n_y\in \Z_{\geq 0}$. The following are satisfied.
    \begin{enumerate} 
\item \label{item:g neq 1} If $g(\Sigma) \neq 1$, then
$\choiceFAq(\Sigma)(x) = \sum\limits_{y\in \Gen(Z_1)} n_y q^{\ell_y} y$ 
where each $n_y\in\{0,1\}$ and  $\ell_y\in \Z$. In particular, if $g(\Sigma)>1$ then $\choiceFAq(\Sigma) = 0$. 
\item \label{item:g = 1} If $g(\Sigma) = 1$, then we have the following two cases. 
\begin{itemize}
    \item If there is an essential circle in $Z_0$ or $Z_1$ then $\choiceFAq(\Sigma) = 0$. 
    \item Otherwise, when $Z_0$ and $Z_1$ consist of contractible circles, let $\Sigma'$ denote the torus obtained by capping off the boundary circles of $\Sigma$. Let $x_1 \in \Gen(Z_0)$ and $y_X \in\Gen(Z_1) $ denote the generators which label each circle in $Z_0$ by $1$ and each circle in $Z_1$ by $X$, respectively. Then for some $\ell\in \Z$,
    \begin{equation}\label{eq:FAq on genus 1 with wrapping}
    \choiceFAq(\Sigma)(x)=
    \begin{cases}
    q^{\ell}\left(1+q^{2 w(\Sigma') }\right) y_X& \text{if } x = x_1, \\
    0 & \text{otherwise.}
    \end{cases}
    \end{equation}
\end{itemize}
    \end{enumerate}
\end{proposition}

\begin{proof}
We may assume that $Z_0$ and $Z_1$ are standard. By\footnote{This can also be deduced from \cite[Theorem 6.3]{BPW}.}  \cite[Lemma 2.17]{AKW}, 
\begin{equation}
\label{eq:FAq is q deformed of FA}
\choiceFAq(\Sigma)(x) = \sum_{y\in \Gen(Z_1)} f_y(q) y
\end{equation}
where each $f_y(q) \in \Z_{\geq 0}[G]$ such that $f_y(1) = n_y$. If $g(\Sigma) > 1$ then $\FA(\Sigma) = 0$ since $\FKh(\Sigma)=0$, so $\choiceFAq(\Sigma)= 0$ as well. If $g(\Sigma) = 0$, then it is straightforward to see that each $n_y\in \{0,1\}$, which by the above implies that each $f_y(q) = n_y q^{\ell_y}$ for some $\ell_y\in \Z$. This completes item (1).

We now address item (2). Suppose $Z_0$ or $Z_1$ contains an essential circle. By neck-cutting near all the trivial circles,\footnote{We may need to push them off the membrane at the cost of introducing a power of $q$.} we write $\Sigma$ as a linear combination where each summand consists of a disjoint union of: a (dotted) cap on each trivial circle in $Z_0$, a (dotted) cup on each trivial circle in $Z_1$, and a (dotted) genus $g(\Sigma)$ cobordism from the essential circles in $Z_0$ to the essential circles in $Z_1$. Each cobordism of this last type induces the zero map, since it must decrease the quantum grading while preserving the annular grading, which establishes the first bullet point of item (2). 

For the second bullet point, by degree reasons we see that $\choiceFAq(\Sigma)(x) = f(q) y_X$ for some $f(q) \in \Z_{\geq 0}[G]$ if $x=x_1$ and $\choiceFAq(\Sigma)(x) = 0$ if $x\neq x_1$. That $f(q) = q^\ell (1+ q^{2 w(\Sigma')}) y_X$ follows from Lemma \ref{lem:torus}.

\end{proof}

\begin{lemma}
\label{lem:bundt pan}
    Let $Z\subset \Aq$ consist of two standard essential circles, and let $\Sigma : \varnothing \to Z$ be the annular cobordism which is topologically an annulus and whose intersection with $\Mq$ is a single arc.
    Then $\choiceFAq(\Sigma) : \Zq \to \choiceFAq(Z)$ is given by 
    \begin{equation}
    \label{eq:bundt cobordism}
        1 \mapsto v_+ \o v_-  + q^{-1} v_- \o v_+. 
    \end{equation}
    Here, the left-to-right ordering of tensor factors corresponds to the innermost-to-outermost ordering on the circles in $Z$. 
\end{lemma}

\begin{proof}
      This follows from \cite[Theorem 6.3]{BPW} and the coevaluation map in \cite[Equation (A.6)]{BPW}. 
\end{proof}

\subsection{Explicitly constructing the quasi-isomorphism $\qXi^\choice$}
\label{sec:explicitly constructing the quasi-iso}

Let $T$ be a planar $(n,n)$-tangle. 
\begin{enumerate}
    \item For a fixed $0\leq k \leq n$, let $\strands{T}$ denote the planar $(2n, 2n)$-tangle obtained by adding $k$ horizontal strands above $T$ and $n-k$ horizontal strands below $T$. See Figure \ref{fig:tangle strands}.
    \item Let $\closure{T}$ denote the closure of $T$ in $\Aq$. See Figure \ref{fig:tangle closure}. Note, the $n-k$ lower and $k$ upper horizontal strands in $\strands{T}$ become inner and outer circles in $\closure{\strands{T}}$, respectively. 
\end{enumerate}
\begin{figure}
    \centering
\subcaptionbox{A planar $(3,3)$-tangle $T$.
\label{fig:planar tangle ex}}[.27\linewidth]
{\begin{tikzpicture}[yscale=.7, xscale=.8]
\draw[densely dashed] (0,0) -- (0,2);
\draw[densely dashed] (2,0) -- (2,2);

\draw (.9,1.2) circle [radius=.2];

\draw(0,.75) arc (-90:90:.5 and .5);

\draw (0,.25) .. controls (1,.25) and (1.5,1.75) .. (2,1.75);

\draw (2,.25) arc (-90:-270:.5 and .5); 

\node at (1,2.5) {};
\node at (1,-1) {};

\end{tikzpicture}
}
\subcaptionbox{The planar $(6,6)$ tangle $\strands{T}$, when $k=1$. 
\label{fig:tangle strands}}[.27\linewidth]
{\begin{tikzpicture}[yscale = .7, xscale=.8]

\draw (.9,1.2) circle [radius=.2];

\draw(0,.75) arc (-90:90:.5 and .5);

\draw (0,.25) .. controls (1,.25) and (1.5,1.75) .. (2,1.75);

\draw (2,.25) arc (-90:-270:.5 and .5);

\draw[line width = .5mm] (0,2.5) -- (0,2); 
\draw[line width = .5mm] (2,2.5) -- (2,2); 
\draw (0,2.25) -- (2,2.25);

\draw[line width = .5mm] (0,-1) -- (0,0); 
\draw[line width = .5mm] (2,-1) -- (2,0); 
\draw (0,-.25) -- (2,-.25); 
\draw (0,-.75) -- (2,-.75); 

\end{tikzpicture}
}
\subcaptionbox{The annular closure $\closure{T} \subset \Aq$ of $T$. 
\label{fig:tangle closure}}[.42\linewidth]
{\begin{tikzpicture}[yscale = .7, xscale=.8]
\draw (.9,1.2) circle [radius=.2];

\draw(0,.75) arc (-90:90:.5 and .5);

\draw (0,.25) .. controls (1,.25) and (1.5,1.75) .. (2,1.75);

\draw (2,.25) arc (-90:-270:.5 and .5); 

\node at (1,2.5) {};
\node at (1,-1) {};

\node at (1,-.3) {$\times$}; 
\draw[densely dashed] (1,-.3) --++ (0,-1.75); 

\draw (0,.25) arc (90:270:.5 and .5); 
\draw (0,.75) arc (90:270:1 and 1); 
\draw (0,1.75) arc (90:270:1.5 and 1.75); 


\draw (2,.25) arc (90:-90:.5 and .5); 
\draw (2,1.25) arc (90:-90:1 and 1.25); 
\draw (2,1.75) arc (90:-90:1.5 and 1.75);


\draw (0,-.75) --++ (2,0);
\draw (0,-1.25) --++ (2,0);
\draw (0,-1.75) --++ (2,0);

\end{tikzpicture}
}
\caption{}
\end{figure}

Let $L\subset \A$ be an annular link diagram intersecting $\mq$ transversely in $n$ points, and let $L^{\cut}$ denote the $(n,n)$-tangle diagram obtained by cutting $\A$ along $\mq$. By \cite[Theorem C]{BPW}, for each $0\leq k \leq n$, there is an isomorphism 
\begin{equation}
\label{eq:qHH of CK is AKH}
\qhh(\qPlatnk;\qCCKnk(L^{cut})) \cong \qaKh(L;n-2k).
\end{equation}
The case where $G=\{1\}$ proves a conjecture of \cite{AGW}. In \cite[Section 5]{LLS_CK_Spectra}, the authors give an explicit chain map that realizes\footnote{As pointed out in \cite{LLS_CK_Spectra}, the authors use the statement of \cite[Theorem C]{BPW} and thus do not give an independent proof of \eqref{eq:qHH of CK is AKH}.} the isomorphism \eqref{eq:qHH of CK is AKH} when $G=\{1\}$. In this subsection we provide an analogous chain map for general cyclic $G$ and also record some additional results which will be used later.

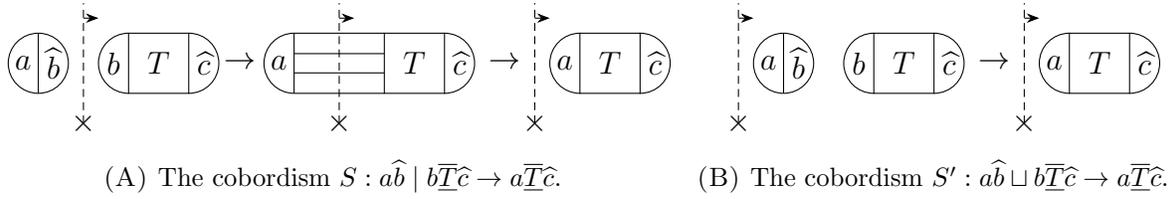
\begin{figure}
\centering
\subcaptionbox{The cobordism $S : a\flip{b} \mid b \strands{T} \flip{c} \to a \strands{T} \flip{c}$. 
\label{eq:qHH differential through seam}}[.55\linewidth]
{\begin{tikzpicture}[scale=.8]

\node at (0,0) {$T$};

\node at (-.75,0) {$b$}; 

\node at (.75,0) {$\widehat{c}$};

\draw (.5,-.5) -- (.5,.5) -- (-.5,.5) -- (-.5,-.5) -- cycle;

\draw (-.5,.5) arc (90:270:.5 and .5); 

\draw (.5,.5) arc (90:-90:.5 and .5); 


\draw (-2,.5) arc (90:270:.5 and .5); 

\draw (-2,.5) arc (90:-90:.5 and .5); 
\draw (-2,.5) -- (-2,-.5); 
\node at (-1.75,0) {$\widehat{b}$}; 
\node at (-2.25,0) {$a$};

\node at (-1.25,-1) {$\times$}; 
\draw[densely dashed] (-1.25,-1) -- (-1.25,1);

\draw[>=Stealth,->] (-1.25,.75) -- (-1,.75);


\begin{scope}[shift= {(4.25,0)}]

\node at (0,0) {$T$};

\node at (.75,0) {$\widehat{c}$};

\draw (.5,-.5) -- (.5,.5) -- (-.5,.5) -- (-.5,-.5) -- cycle;

\draw (.5,.5) arc (90:-90:.5 and .5); 


\draw (-2,.5) arc (90:270:.5 and .5); 
 
\draw (-2,.5) -- (-2,-.5); 
 
\node at (-2.25,0) {$a$};

\draw (-2,.5) -- (-.5,.5); 
\draw (-2,.16) -- (-.5,.16); 
\draw (-2,-.16) -- (-.5,-.16); 
\draw (-2,-.5) -- (-.5,-.5);

\node at (-1.25,-1) {$\times$}; 
\draw[densely dashed] (-1.25,-1) -- (-1.25,1);

\draw[>=Stealth,->] (-1.25,.75) -- (-1,.75); 
\end{scope}


\begin{scope}[shift= {(7.5,0)}]

\node at (0,0) {$T$};

\node at (-.75,0) {$a$}; 

\node at (.75,0) {$\widehat{c}$};

\draw (.5,-.5) -- (.5,.5) -- (-.5,.5) -- (-.5,-.5) -- cycle;

\draw (-.5,.5) arc (90:270:.5 and .5); 

\draw (.5,.5) arc (90:-90:.5 and .5); 


\node at (-1.25,-1) {$\times$}; 
\draw[densely dashed] (-1.25,-1) -- (-1.25,1);

\draw[>=Stealth,->] (-1.25,.75) -- (-1,.75); 
\end{scope}

\node at (1.35,0) {$\to$}; 
\node at (5.75,0) {$\to$}; 
\end{tikzpicture}
}
\subcaptionbox{The cobordism $S ': a\flip{b} \sqcup b \strands{T} \flip{c} \to a \strands{T} \flip{c}$. 
\label{eq:qHH differential no seam}}[.4\linewidth]
{\begin{tikzpicture}[scale=.8]

\node at (0,0) {$T$};

\node at (-.75,0) {$b$}; 

\node at (.75,0) {$\widehat{c}$};

\draw (.5,-.5) -- (.5,.5) -- (-.5,.5) -- (-.5,-.5) -- cycle;

\draw (-.5,.5) arc (90:270:.5 and .5); 

\draw (.5,.5) arc (90:-90:.5 and .5); 


\draw (-2,.5) arc (90:270:.5 and .5); 

\draw (-2,.5) arc (90:-90:.5 and .5); 
\draw (-2,.5) -- (-2,-.5); 
\node at (-1.75,0) {$\widehat{b}$}; 
\node at (-2.25,0) {$a$};

\node at (-2.75,-1) {$\times$}; 
\draw[densely dashed] (-2.75,-1) -- (-2.75,1);

\draw[>=Stealth,->] (-2.75,.75) -- (-2.5,.75);


\begin{scope}[shift= {(3.25,0)}]

\node at (0,0) {$T$};

\node at (-.75,0) {$a$}; 

\node at (.75,0) {$\widehat{c}$};

\draw (.5,-.5) -- (.5,.5) -- (-.5,.5) -- (-.5,-.5) -- cycle;

\draw (-.5,.5) arc (90:270:.5 and .5); 

\draw (.5,.5) arc (90:-90:.5 and .5); 


\node at (-1.25,-1) {$\times$}; 
\draw[densely dashed] (-1.25,-1) -- (-1.25,1);

\draw[>=Stealth,->] (-1.25,.75) -- (-1,.75); 
\end{scope}

\node at (1.5,0) {$\to$}; 
\end{tikzpicture}
}
\caption{The cobordisms in Lemma \ref{lem:FAq recovers deformed face map}. }\label{fig:saddles through seam}
\end{figure}

Fix a planar $(n,n)$-tangle $T$, and let $a,b,c\in \cupdiag{n}{k}$. Place $a\flip{b}$ and $b \strands{T} \flip{c}$ in $\Aq$ as shown in the left-most picture of Figure \ref{eq:qHH differential through seam}, and denote this configuration by $ a\flip{b} \mid b \strands{T} \flip{c}$. Consider the cobordism $S : a\flip{b} \mid b \strands{T} \flip{c} \to a \strands{T} \flip{c}$ given as the composition of the minimal saddle from $\flip{b} b$ to $2n$ horizontal lines, followed by moving $a$ to the other side of $\mq$, as shown in Figure \ref{eq:qHH differential through seam}. 

On the other hand, let $a\flip{b} \sqcup b \strands{T} \flip{c}$ denote the configuration shown in Figure \ref{eq:qHH differential no seam}, and let $S' : a\flip{b} \sqcup b \strands{T} \flip{c} \to a \strands{T} \flip{c}$ be the minimal saddle cobordism. Given a homogeneous $x\in \FAq(a \flip{b})$ and any $m\in \FAq(b \strands{T} \flip{c})$, each represented by (dotted) cup cobordisms, we write $x\mid m \in \FAq(a\flip{b} \mid b \strands{T} \flip{c})$ and $x\o m \in \FAq(a\flip{b} \sqcup b \strands{T} \flip{c})$ to be the natural elements, represented by the disjoint union of the aforementioned cup cobordisms. 

The following lemma will be used in Section \ref{sec:aug quant HM functor}. 
\begin{lemma}
\label{lem:FAq recovers deformed face map}
With the notation established above, we have \[
\FAq(S)(x\mid m) = q^{-\lr{x}} \FAq(S')(x\otimes m).
\]
\end{lemma}

\begin{proof}
Represent $x$ as a dotted cup cobordism $\varnothing \to a \flip{b}$. We can relate the two maps by moving the cobordism schematically depicted in Equation \eqref{eq:qHH differential through seam cob} to the other side of the membrane. 
\begin{equation}
    \label{eq:qHH differential through seam cob}
    \begin{aligned}
    \begin{tikzpicture}

\begin{scope}[shift={(0,.25)}]
\draw[>=Stealth,->]  (.5,1.7) -- (.8,1.75); 

\draw[white,fill=gray!40] (-.5, 2) -- (2,1.25) -- (2,-1.75) -- (-.5, -1) --  (-.5, -.85); 

\end{scope}

\draw[white, fill=white] (-.5,.5) -- (.5,.5) -- (.5,1.5) -- (-.5,1.5) -- cycle;

\draw[densely dashed] (.5,.5) --++ (1,0);
\draw[densely dashed] (.5,1.5) --++ (1,0);
\begin{scope}[shift={(-1,0)}]

\draw[densely dashed] (0,0) arc (0:180:.8 and .4); 
\draw (0,0) arc (0:-180:.8 and 1);

\draw[] (0,0) arc (0:-180:.8 and .4); 
\node at (-.5,0) {$\widehat{b}$};
\node at (-1.2,0) {$a$};
\draw (-.8,.4) -- (-.8,-.4);

\node at (-.8,-.75) {$x$}; 

\draw[] (.5,.5) arc(90:180:.5 and .5); 
 
\draw[]( -.6,1.5) arc (90:180:1 and 1);
\draw[] (-1.6,.5) -- (-1.6,0); 
\end{scope}

\draw (.5,.5) --++ (-1,0); 
\draw (.5,1.5) --++ (-2.1,0);  
\draw[fill=white] (.5,1.5) arc (90:-90:.5 and .5); 
\draw[densely dashed]  (.5,.5) arc (-90:--270:.5 and .5);

\end{tikzpicture}
    \end{aligned}
\end{equation}
This cobordism has degree equal to $\lr{x}$, which completes the proof.
\end{proof}

The remainder of this subsection is devoted to providing an explicit chain map yielding the isomorphism in \eqref{eq:qHH of CK is AKH}. Much of what follows is an analogue of \cite[Section 5]{LLS_CK_Spectra} for quantum annular homology; in particular, we will now define two key maps $A_a$ and $C$.

For a planar $(n,n)$-tangle $T$, $0\leq k \leq n$, and $a\in \cupdiag{n}{k}$, define a map 
\begin{equation}
\label{eq:A map} 
A_a: \FAq(a \strands{T} \flip{a}) \to \FAq( \closure{\strands{T}})
\end{equation}
by placing $a\strands{T} \flip{a}$ into $\Aq$ and applying $\FAq$ to the minimal saddle cobordism from $a\strands{T} \flip{a}$ to $\closure{\strands{T}}$
around the puncture, as shown in Figure \ref{fig:A map}. For a choice of generators $\choice$, we will denote by $A_a^\choice:\choiceFAq(a \strands{T} \flip{a}) \to \choiceFAq( \closure{\strands{T}})$ the natural analogue of \eqref{eq:A map}, defined using $\choiceFAq$ in place of $\FAq$.

\begin{figure}
    \centering
       \begin{tikzpicture}

\node at (0,0) {$T$};

\node at (-.75,0) {$a$}; 

\node at (.75,0) {$\widehat{a}$};

\draw (.5,-.5) -- (.5,.5) -- (-.5,.5) -- (-.5,-.5) -- cycle;

\draw (-.5,.5) arc (90:270:.5 and .5); 

\draw (.5,.5) arc (90:-90:.5 and .5); 

\node at (0,-1) {$\times$}; 
\draw[densely dashed] (0,-1)-- (-2.25,-1);

\draw[>=Stealth,->] (-1.5,-1) -- (-1.5,-.75); 

\draw[very thick, red] (-1,0) arc (90:270:.75 and 1); 
\draw[very thick, red] (1,0) arc (90:-90:.75 and 1); 
\draw[very thick, red] (-1,-2) -- (1,-2); 
\end{tikzpicture}
    \caption{The cobordism $a \strands{T} \flip{a} \to \closure{\strands{T}}$ inducing the map  $A_a: \FAq(a \strands{T} \flip{a}) \to \FAq( \closure{\strands{T}})$.}
    \label{fig:A map}
\end{figure}

  Now we define the $C$ map. First, for a Khovanov generator $y\in\Gen(\closure{\strands{T}})$, let $y|_{\closure{T}}\in\Gen(\closure{T})$ denote the corresponding generator whose labels are the restriction of the labels assigned by $y$, and define $C^\choice$ on Khovanov generators as
\begin{equation}\label{eq:C map on generators}
C^\choice(y) := \begin{cases}
    y|_{\closure{T}} & \text{if $y$ labels the inner $n-k$ (resp. the outer $k$) essential circles}\\
    & \text{ of $\closure{\strands{T}}$ by $v_-$ (resp. by $v_+$), and}\\
    0 & \text{otherwise.}
\end{cases}
\end{equation}

\begin{lemma}
\label{lem:C map naturality}
    Let $T$ and $T'$ be planar $(n,n)$-tangles and let $S:T\to T'$ be a cobordism. We denote by $\strands{S} : \strands{T}\to \strands{T'}$ the cobordism obtained by adding to $S$ the identity cobordism on $n-k$ (resp. $k$) horizontal lower (resp. upper) strands. We also denote by $\closure{\strands{S}} : \closure{\strands{T}} \to \closure{\strands{T'}}$ and by $\closure{S}: \closure{T} \to \closure{T'}$ the annular cobordisms obtained by closing $\strands{S}$ and $S$, respectively. Then the diagram \eqref{eq:C map commutative diagram} commutes. 
    \begin{equation}
    \label{eq:C map commutative diagram}
\begin{tikzcd}
    \choiceFAq(\closure{\strands{T}}) \ar[r,"C^\choice"] \ar[d,"\choiceFAq\left(\closure{\strands{S}}\right)"'] & \choiceFAq(\closure{T}) \ar[d, "\choiceFAq(\closure{S})"] \\
    \choiceFAq(\closure{\strands{T'}}) \ar[r,"C^\Psi"] & \choiceFAq(\closure{T'})
\end{tikzcd}.
    \end{equation}
\end{lemma}

\begin{proof} This follows from the second bullet point in Definition \ref{def:choice of generators}. 
\end{proof}

\begin{remark}
In \cite[Section 5]{LLS_CK_Spectra}, the authors consider maps  which we will denote by $A_{LLS}, B_{LLS}$, and $C_{LLS}$ to avoid confusion with our $A_a$ and $C$. The Khovanov TQFT $\FKh$ rather than the annular Khovanov TQFT $\FA$ is used to define $A_{LLS}$. The map $A_{LLS}$ lands in non-positive annular degree since $\FKh$ does not increase $\adeg$, and  $B_{LLS}$ projects to the annular degree zero summand. So $B_{LLS} \circ A_{LLS}$ recovers the annular TQFT $\FA$. Since we use the quantum annular TQFT to define $A_a$, we do not have an analogue of the map $B_{LLS}$. 
\end{remark}

\begin{figure}
    \centering
    \begin{tikzpicture}

\draw (0,.5) .. controls (.5,.5)  and (1,-.5) .. (1.5,-.5); 

\draw (1.5,.5) arc (90:270:.25 and .25);

\node at (.75,-1.6) {$T$};

\begin{scope}[shift={(3,0)}]

\draw (0,.5) .. controls (.5,.5)  and (1,-.5) .. (1.5,-.5); 

\draw (1.5,.5) arc (90:270:.25 and .25);

\draw[dashed] (0,1) --++ (1.5,0);

\draw[dashed] (0,-1) --++ (1.5,0);

\draw[dashed] (0,1.5) --++ (1.5,0);

\node at (.75,-1.6) {$\underline{\overline{T}}$};
\end{scope}

\begin{scope}[shift={(6,0)}, xscale=-1]
\draw (0,.5) arc (-90:90:.25 and .25);
\draw[line width = .5mm] (0,.75) --++ (0,.5); 

\node at (0,-1.6) {$a$};
\end{scope}

\begin{scope}[shift={(8.5,0)}, xscale=-1]
\draw (0,.5) arc (-90:90:.25 and .25);

\draw[line width = .5mm] (0,.75) --++ (0,.5); 

\draw[dashed] (0,-1) arc (-90:90:1.25 and 1.25);

\node at (.5,-1.6) {$a^*$};
\end{scope}

\begin{scope}[shift={(10,0)}, xscale=1]
\draw (0,0) arc (-90:90:.25 and .25);

\draw (0,-.5) arc (-90:90:.75 and .75);

\draw (0,-1) arc (-90:90:1.25 and 1.25);

\draw[line width = .5mm] (0,-.75) --++ (0,-.5); 

\draw[line width = .5mm] (0,.75) --++ (0,1); 

\node at (.5,-1.6) {$\widehat{b} = \widehat{b}^*$};
\end{scope}

\begin{scope}[shift={(13.5,0)}]
\draw (0,.5) .. controls (.5,.5)  and (1,-.5) .. (1.5,-.5); 

\draw (1.5,.5) arc (90:270:.25 and .25);

\draw[line width = .5mm] (1.5,-.75) --++ (0,-.5); 

\draw[line width = .5mm] (1.5,.75) --++ (0,1); 

\draw[line width = .5mm] (0,.75) --++ (0,.5); 


\draw (1.5,0) arc (-90:90:.25 and .25);

\draw (1.5,-.5) arc (-90:90:.75 and .75);

\draw (1.5,-1) arc (-90:90:1.25 and 1.25);

\draw (0,1) arc (90:270:.25 and .25);


\draw[dashed] (0,1) --++ (1.5,0);

\draw[dashed] (0,1.5) --++ (1.5,0);

\draw[dashed] (0,-1) --++ (1.5,0);

\draw[dashed] (0,1.5) arc (90:270:1.25 and 1.25);

\node at (.75,-1.6) {$aT\widehat{b} = a^* \underline{\overline{T}} \widehat{b}^*$};
\end{scope} 

\end{tikzpicture}
    \caption{An example of some of the notation in Lemma \ref{lem:trace property for qXi} with $k=1$. Arcs and horizontal strands that are added are dashed. }
    \label{fig:adding arcs and strands}
\end{figure}
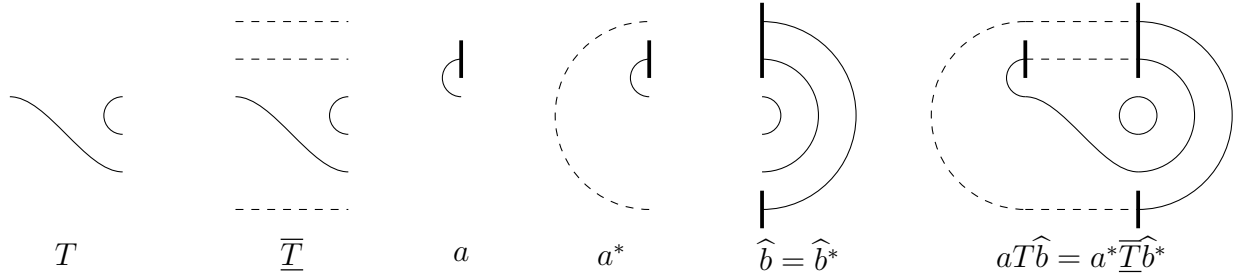

Next, consider a planar $(n,m)$-tangle  $T_1$, a planar $(m,n)$-tangle $T_2$, $a\in \cupdiag{n}{k}$, and $b\in \cupdiag{m}{h(k)}$ where  $k$ is in the range \eqref{eq:k range}. Recall that arcs and horizontal strands are added when forming $a T_1 \flip{b}$ and $b T_2 \flip{a}$  (see Figure \ref{fig:capping off platform tangle}). For the following lemma, for $T\in \{T_1, T_2\}$, define $\strands{T}$ to mean $T$ with $\max(k,h(k))$ horizontal strands added above and $\max(n,m)-\max(k,h(k))$ strands added below.  In this way, $a T \flip{b}$ is precisely the diagram $a^*\strands{T} \flip{b}^*$ with 
\[
a^*\in B^{n + \max ( (m-n)/2,0
)} \text{ and } b^* \in B^{m +  \max ( (n-m)/2,0) }
\]
obtained from $a$ and $b$ by potentially adding outermost arcs as in Figure \ref{fig:adding arcs and strands}. Similarly, $b T_2 \flip{a} = b^* \strands{T_2} \flip{a}^*$.

Placing $a \strands{T_1} \flip{b}$ to the left and right, respectively, of $b \strands{T_2} \flip{a}$ and applying saddle cobordisms, we obtain two maps
\begin{align}
\label{eq:mu_b}
\mu_{b^*} &: \FAq(a^* \strands{T_1} \flip{b}^*) \o_{\Zq} \FAq( b^* \strands{T_2} \flip{a}^*) \cong  \FAq(a^* \strands{T_1} \flip{b}^* \sqcup b^* \strands{T_2} \flip{a}^*) \to \FAq( a^* \strands {T_1}\, \strands{T_2} \flip{a}^*) \\
\label{eq:mu_a}
\mu_{a^*} &:   \FAq( b^* \strands{T_2} \flip{a}^*)  \o_{\Zq} \FAq(a^* \strands{T_1} \flip{b}^*)  \cong  \FAq(b^* \strands{T_2} \flip{a}^* \sqcup a^* \strands{T_1} \flip{b}^*) \to \FAq( b^* \strands {T_2}\,  \strands{T_1} \flip{b}^*).
\end{align}
Note that $\mu_{b^*}$ is not the map realizing the isomorphism in Lemma \ref{lem:gluing planar tangles}, in that in $\mu_b$, dashed circles are not removed after applying saddles (similarly for $\mu_{a^*}$). Note also that $ \strands {T_1}\, \strands{T_2} $  is different from $\strands{T_1 T_2}$ if $n< m$, for the same reason.  We also have maps
\begin{align}
\label{eq:mu_b gluing}
\mu_b &: \FAq(a^* \strands{T_1} \flip{b}^*) \o_{\Zq} \FAq( b^* \strands{T_2} \flip{a}^*) \cong  \FAq(a^* \strands{T_1} \flip{b}^* \sqcup b^* \strands{T_2} \flip{a}^*) \to \FAq( a \strands {T_1 T_2} \flip{a}) \\
\label{eq:mu_a gluing}
\mu_a &:   \FAq( b^* \strands{T_2} \flip{a}^*)  \o_{\Zq} \FAq(a^* \strands{T_1} \flip{b}^*)  \cong  \FAq(b^* \strands{T_2} \flip{a}^* \sqcup a^* \strands{T_1} \flip{b}^*) \to \FAq( b \strands {T_2 T_1} \flip{b}).
\end{align}
in which outermost dashed circles are removed after performing saddles (for instance, the maps $\mu_b$ as $b$ varies assemble to form the isomorphism of Lemma \ref{lem:gluing planar tangles}).

\begin{lemma}
\label{lem:trace property for qXi}
With the notation established above, the following diagram commutes
    \[
\begin{tikzcd}
    \choiceFAq( a^* \strands{T_1} \flip{b}^*)\{\nmshift{n}{m}\}\o_{\Zq}  \choiceFAq( b^* \strands{T_2} \flip{a}^*) \{\nmshift{m}{n}\}  \ar[r, "\tau"] \ar[d, "C^\choice \circ A_{a}^\choice \circ \mu_b"]&    \choiceFAq( b^* \strands{T_2} \flip{a}^*) \{\nmshift{m}{n}\}  \o_{\Zq} \choiceFAq( a^* \strands{T_1} \flip{b}^*) \{\nmshift{n}{m}\}  \ar[d, "C^\choice \circ A_{b}^\choice \circ \mu_a"] \\
    \choiceFAq( \closure{T_1 T_2} ) \ar[r, "\choiceFAq(\Sigma)"] & \choiceFAq( \closure{T_2 T_1}) 
\end{tikzcd}
\]
where $\tau (x\o y) = q^{-\lr{y}}y\o x$ and $\Sigma : \closure{T_1 T_2}  \to \closure{T_2 T_1}$ is the annular cobordism formed by the minimal isotopy that moves $T_2$ across the seam, as shown in \eqref{eq:moving T2 across the seam}.

\begin{equation}
\label{eq:moving T2 across the seam}
\begin{aligned}
    \begin{tikzpicture}

\node at (0,0) {$T_1$};

\draw (.5,-.5) -- (.5,.5) -- (-.5,.5) -- (-.5,-.5) -- cycle;

\node at (1.5,0) {$T_2$};

\draw (2,-.5) -- (2,.5) -- (1,.5) -- (1,-.5) -- cycle;

\draw (.5,.15) --++ (.5,0);

\draw (.5,-.15) --++ (.5,0);

\draw[densely dashed] (2.5,-1)--++ (0,2);

\draw[>=Stealth,->] (2.5,.75) --++ (.3,0); 

\draw (2,.3) --++ (1,0);
\draw (2,.1) --++ (1,0);
\draw (2,-.1) --++ (1,0);
\draw (2,-.3) --++ (1,0);

\draw (-.5,.3) --++ (-.5,0);
\draw (-.5,.1) --++ (-.5,0);
\draw (-.5,-.1) --++ (-.5,0);
\draw (-.5,-.3) --++ (-.5,0);

\draw[->] (3.75,0) -- (4.75,0) node[midway,above] {$\Sigma$};


\begin{scope}[shift={(6.5,0)}]
\node at (0,0) {$T_1$};

\draw (.5,-.5) -- (.5,.5) -- (-.5,.5) -- (-.5,-.5) -- cycle;

\node at (1.5,0) {$T_2$};

\draw (2,-.5) -- (2,.5) -- (1,.5) -- (1,-.5) -- cycle;

\draw (.5,.15) --++ (.5,0);

\draw (.5,-.15) --++ (.5,0);

\draw[densely dashed] (.75,-1)--++ (0,2);

\draw[>=Stealth,->] (.75,.75) --++ (.3,0); 

\draw (2,.3) --++ (.5,0);
\draw (2,.1) --++ (.5,0);
\draw (2,-.1) --++ (.5,0);
\draw (2,-.3) --++ (.5,0);

\draw (-.5,.3) --++ (-.5,0);
\draw (-.5,.1) --++ (-.5,0);
\draw (-.5,-.1) --++ (-.5,0);
\draw (-.5,-.3) --++ (-.5,0);

\end{scope}

 
\end{tikzpicture}
\end{aligned}
\end{equation}
\end{lemma}

\begin{proof}

Throughout the proof we fix a choice of generators $\choice$ but omit it from the notation to avoid clutter. 
We will actually show that the following diagram commutes
\[
\begin{tikzcd}
    \FAq( a^* \strands{T_1} \flip{b}^*)\{\nmshift{n}{m}\}\o_{\Zq}  \FAq( b^* \strands{T_2} \flip{a}^*) \{\nmshift{m}{n}\}  \ar[r, "\tau"] \ar[d, " A_{a^*} \circ \mu_{b^*}"] \ar[bend right = 70, dd, "C \circ A_{a} \circ \mu_b"'] &    \FAq( b^* \strands{T_2} \flip{a}^*) \{\nmshift{m}{n}\}  \o_{\Zq} \FAq( a^* \strands{T_1} \flip{b}^*) \{\nmshift{n}{m}\}  \ar[d, "A_{b^*}\circ \mu_{a^*}"]  \ar[bend left = 70, dd, "C\circ A_{b} \circ \mu_a"] \\
    \choiceFAq( \closure{\, \strands{T_1}\, \strands{T_2} \, } ) \ar[r, "q^{\frac{m-n}{2}}  \FAq(\strands{\Sigma})"] \ar[d,"q^{s_{m,n}}C"] & \FAq( \closure{\,\strands{T_2} \, \strands{T_1}\, }) \ar[d,"q^{s_{n,m}}C"] \\
    \FAq(\closure{T_1 T_2}) \ar[r, "\FAq(\Sigma)"] & \FAq(\closure{T_2 T_1}) 
\end{tikzcd}
\]
where $s_{m,n} = \max((m-n)/2,0)$ and $\strands{\Sigma} : \closure{\, \strands{T_1}\, \strands{T_2} \, }  \to \closure{\, \strands{T_2} \, \strands{T_1}\, }$ is obtained from $\Sigma$ in the natural way by adding the identity cobordism on the $\max(k,h(k))$ outermost and $\max(n,m) - \max(k,h(k))$ innermost standard essential circles.

Let us begin with the top middle square.  Assume that $x$ and $y$ are standard generators, represented as dotted cup cobordisms $S_x: \varnothing \to a^* \strands{T_1} \flip{b}^*$ and $S_y: \varnothing \to b^* \strands{T_2} \flip{a}^*$ which are both disjoint from $\Mq$. Consider the cobordism $S': \varnothing \to \closure{\,\strands{T_2}\, \strands{T_1}\, }$ formed by starting at the top left of the square, performing the saddle cobordisms defining $A_{a^*} \circ \mu_{b^*}$, and finally performing $\strands{\Sigma}$. We show a portion of $S'$ around $S_y$ in \eqref{eq:moving y to other side}, which we view as a tangle cobordism $S :b^* \to \strands{T_2} a^*$ read from left to right. 
 \begin{equation}
 \label{eq:moving y to other side}
     \begin{aligned}
\begin{tikzpicture}[yscale=1.5]

\draw[>=Stealth,->]  (.5,1.7) -- (.8,1.75);

\begin{scope}
\clip (-.5, 2) -- (-.5,1.7)  -- (0,1.7) -- (0,1.4)  -- (1,1.4) -- (1,1.3) -- (.99,1.2) -- (.97,1.1) -- (.955,1) -- (.925,.9) -- (.9,.85)-- ( .8,.69) -- (.7,.57) -- (.6,.52) -- (.5,.5) -- (.4,.5) -- (-.5,.5) -- (-.5,-1) -- (2,-1.75) -- (2,1.25) -- cycle;
\filldraw[color = green] (-.5, 2) -- (-.5,1.3)  -- (0,1.3) -- (.05,1) -- (1,1) -- (1,.9) -- (.9,.7) -- (.8,.6) -- (.7,.55)-- (.5,.5) -- (.4,.5) -- (-.5,.5) -- (-.5,-1) -- (2,-1.75) -- (2,1.25) -- cycle;
\fill[blue] (-.5, 2) -- (-.5,1.3)  -- (0,1.3) -- (.1,1) -- (1,1) -- (1,.9) -- (.9,.65) -- (.8,.57) -- (.7,.52)-- (.5,.48) -- (.4,.5) -- (-.5,.5) -- (-.5,-1) -- (2,-1.75) -- (2,1.25) -- cycle;
\draw[white,fill=gray!40] (-.5, 2) -- (2,1.25) -- (2,-1.75) -- (-.5, -1) --  (-.5, -.85); 
\end{scope}

\begin{scope}[shift={(-1,0)}]
\draw[] (0,0) arc (0:-180:1.5 and .4); 
\draw[densely dashed] (0,0) arc (0:180:1.5 and .4); 
\draw[] (0,0) arc (0:-180:1.5 and 1.2); 

\node at (-.5,0) {$\widehat{a}^*$};
\node at (-2.5,0) {$b^*$};
\node at (-1.5,0) {$\underline{\overline{T_2}}$};
\draw (-1,.38) -- (-1,-.38); 
\draw (-2,.38) -- (-2,-.38);

\node at (-1.5,-.75) {$S_y$};

\draw[] (.5,.5) arc(90:180:.5 and .5); 

 \draw [ double distance=3pt,
             arrows = {-Stealth[length=0pt 2.2 0,open]}] (-1.5,.5) -- (1,1) node[midway, above] {$\Sigma$} ;

\begin{scope}[shift={(-1,0)}]
\draw[] (-2.5,.5) arc (-90:-180:.5 and 1.2); 
\draw[] (-2,1.4) arc (0:-90:.5 and .9); 

\draw[] (-2.5,.5) arc (90:0:.5 and .5); 

\node[rotate=90] at (-2.5,1) {$b^*$};
\end{scope}
\end{scope}

\begin{scope}
\clip (-.5,1.7) -- (-.5,1.4) -- (1,1.4) -- (1,1.7) -- cycle;
\draw[]  (.5,.5) arc (-90:-180:.5 and 1.2);
\end{scope}

\draw[] (.5,.5) --++ (-1,0); 
\draw[] (1,1.4) --++ (-5,0); 
\draw[] (0,1.7) --++ (-5,0); 
\draw[] (1,1.4) arc (0:-90:.5 and .9); 
\draw[densely dashed]  (.5,.5) arc (-90:-180:.5 and 1.2);

\node[rotate=-90] at (.5,.9) {$ \underline{\overline{T_2}} \widehat{a}^*$};

\end{tikzpicture}
     \end{aligned}
 \end{equation}
The degree $s$ of $S$, in the sense of \eqref{eq:degree of tangle cob}, is
\[
s := \underbrace{\left( \lr{y} - \nmshift{m}{n} \right)}_{\chi(S)} - \underbrace{\left(m + s_{n,m} \right)}_{\text{contribution from } \partial S }.
\]
We can move $S$ to the other side of the membrane at the cost of introducing a factor of $q^{-s}$. Note that 
 \[
\frac{m-n}{2} = s_{m,n} - s_{n,m},
 \]
 so $
-s + \frac{m-n}{2} = - \lr{y}$, precisely the power of $q$ appearing in $\tau$. This shows that the top middle square commutes.

Next, let us address the bottom middle square. First, consider this square with all powers of $q$ omitted, which commutes by the condition in the second bullet point of Definition \ref{def:choice of generators}. Moreover, $s_{m,n} = \frac{m-n}{2} + s_{n,m}$, so the original square with powers of $q$ included also commutes. 

Finally, we show that the left triangle commutes; the argument for the right triangle is similar. For the left triangle, observe that $\mu_b$ is equal to $\mu_{b^*}$ followed by removing $s_{m,n}$ dashed outermost circles.  By \cite[Lemma 5.4]{LLS_CK_Spectra}, we may assume that each such circle is labeled $1$. Then by Lemma \ref{lem:bundt pan}, under the $ A_{a^*}$ map each such circle contributes a sum of two generators, precisely one of which survives the $C$ map and contributes a factor of $q^{-1}$. The factor of $q^{s_{m,n}}$ then guarantees that the left triangle commutes. 
\end{proof}

We establish some notation to record  a special case of Lemma \ref{lem:trace property for qXi} which will be used later in the proof of Proposition \ref{prop:qXi is a chain map}. Given a planar $(n,n)$-tangle $T$ and $a,b\in \cupdiag{n}{k}$, we have special cases of Equations \eqref{eq:mu_a} and \eqref{eq:mu_b} (where $T_2 = \id_n$)
\[
 \FAq(a \strands{T} \flip{b}) \o_{\Zq} \FAq(b\flip{a})   \to \FAq(a \strands{T} \flip{a}) \text{ and } \FAq(b\flip{a}) \o_{\Zq} \FAq(a \strands{T} \flip{b})    \to \FAq(b \strands{T} \flip{b})
\]
which we denote by $x\o y \mapsto x\cdot y$ and $y\o x \mapsto y\cdot x$ for $x\in  \FAq(a \strands{T} \flip{b}), y\in  \FAq(b\flip{a})$.

\begin{corollary}
\label{cor:A is a chain map}
    For a planar $(n,n)$-tangle $T$, platform matchings $a,b\in \cupdiag{n}{k}$, and homogeneous elements $x\in \FAq(a \strands{T} \flip{b})\{-n\}$ and $y\in \FAq(b \flip{a})\{-n\}$, we have
\begin{equation}
\label{eq:A is a chain map}
A_a( x \cdot y) =  q^{-\lr{y}} A_b(y \cdot x). 
\end{equation}
This is illustrated in Equation \eqref{eq:A is a chain map figure}.
\begin{equation}
\label{eq:A is a chain map figure}
    \begin{aligned}
    \begin{tikzpicture}

\begin{scope}[shift={(-.5,0)}]
\node at (0,0) {$T$};

\node at (-.75,0) {$a$}; 

\node at (.75,0) {$\widehat{b}$};

\draw (.5,-.5) -- (.5,.5) -- (-.5,.5) -- (-.5,-.5) -- cycle;

\draw (-.5,.5) arc (90:270:.5 and .5); 

\draw (.5,.5) arc (90:-90:.5 and .5);

\draw (2,.5) arc (90:270:.5 and .5); 
\draw (2,.5) arc (90:-90:.5 and .5); 
\draw (2,.5) -- (2,-.5); 
\node at (2.25,0) {$\widehat{a}$}; 
\node at (1.75,0) {$b$};

\draw[very thick, red] (1.5,0) -- (1,0);
\end{scope}

\node at (0,-1) {$\times$}; 
\draw[densely dashed] (0,-1)-- (-2.25,-1);

\draw[>=Stealth,->] (-1.5,-1) -- (-1.5,-.75); 
 
\draw[very thick, red] (-1.5,0) arc (90:270:.6 and 1); 

\draw[very thick, red] (2,0) arc (90:-90:.6 and 1); 

\draw[very thick, red] (-1.5,-2) -- (2,-2); 


\begin{scope}[shift={(7.75,0)}]

\begin{scope}[shift={(1,0)}]
\node at (0,0) {$T$};

\node at (-.75,0) {$a$}; 

\node at (.75,0) {$\widehat{b}$};

\draw (.5,-.5) -- (.5,.5) -- (-.5,.5) -- (-.5,-.5) -- cycle;

\draw (-.5,.5) arc (90:270:.5 and .5); 

\draw (.5,.5) arc (90:-90:.5 and .5);

\draw (-2,.5) arc (90:270:.5 and .5); 

\draw (-2,.5) arc (90:-90:.5 and .5); 
\draw (-2,.5) -- (-2,-.5); 
\node at (-1.75,0) {$\widehat{a}$}; 
\node at (-2.25,0) {$b$};

 
\draw[very thick, red] (-1.5,0) -- (-1,0);
\end{scope}

\node at (0,-1) {$\times$}; 
\draw[densely dashed] (0,-1)-- (-2.25,-1); 

\draw[>=Stealth,->] (-1.5,-1) -- (-1.5,-.75); 
 
\draw[very thick, red] (-1.5,0) arc (90:270:.6 and 1); 

\draw[very thick, red] (2,0) arc (90:-90:.6 and 1); 

\draw[very thick, red] (-1.5,-2) -- (2,-2);

\end{scope}


\node at (4,-1) {$\ =\ q^{-(\vert y \vert)}$}; 
\end{tikzpicture}
    \end{aligned}
\end{equation}
\end{corollary}

\begin{remark}
    \label{rem:A to A choice}
    Lemma  \ref{lem:FAq recovers deformed face map} and Corollary \ref{cor:A is a chain map} hold with $\FAq$ and $A_a$ replaced with $\choiceFAq$ and $A_a^\choice$, for some choice of generators $\Psi$.  Note that many of the diagrams are already standard. 
\end{remark}

Consider now an $(n,n)$-tangle diagram $T$ with $N$ crossings. The chain complex $\qCCK(T)$ splits as a direct sum 
\[
\qCCK(T) = \bigoplus_{0\leq k \leq n} \qCCKnk(T),
\]
with each $\qCCKnk(T)$ a complex of $\qPlatnk$-bimodules; we refer to its differential as the \emph{internal} differential. Moreover, $\qhh_*(\qPlatnk;\qCCKnk(T))$ can be computed as the homology of $\qch_*(\qPlatnk;\qCCKnk(T))$, the latter of which is the total complex of the bicomplex 
\begin{align}
\begin{aligned}
\label{eq:qHH bicomplex}
& \cdots \longrightarrow \bigoplus_{a,b,c\in \cupdiag{n}{k}} \qCCKnk(T;a, b)  \o_{\Zq} \qPlat{n}(b,c) \o_{\Zq} \qPlat{n}(c,a) \\
& \longrightarrow \bigoplus_{a, b\in \cupdiag{n}{k}} \qCCKnk(T;a, b)  \o_{\Zq} \qPlat{n}(b,a)  \longrightarrow \bigoplus_{a\in \cupdiag{n}{k}} \qCCKnk(T;a).
\end{aligned}
\end{align}
We refer to the maps written in \eqref{eq:qHH bicomplex} as the \emph{Hochschild differential}. 

For any $a,b\in \cupdiag{n}{k}$ and $v\in \cube^N$, there are canonical identifications $\qFKh(a\flip{b}) \cong \choiceFAq(a \flip{b})$ and $\qFKh(a\strands{T_v}\flip{b}) \cong \choiceFAq(a \strands{T}_v \flip{b})$, where on the right-hand side of these isomorphisms we view $a\flip{b}$ and $a\strands{T}_v \flip{b}$ as diagrams in $\Aq$, as explained above Lemma \ref{lem:FAq recovers deformed face map}. 
Chain groups in the right-most term of \eqref{eq:qHH bicomplex} consist of direct sums of $\qPlatnk$-bimodules of the form
\[
\qFCKk(T_v) = \bigoplus_{a\in \cupdiag{n}{k}}  \left(\choiceFAq(a\strands{T_v}\flip{a}) / _q\mathcal{I}(T_v;a)\right) \{-n\}.
\]
By \eqref{eq:FAq is q deformed of FA} and \cite[Lemma 5.4]{LLS_CK_Spectra}, $C^\choice\circ A_a^\choice$ vanishes on each sub-bimodule $_q\mathcal{I}(T_v;a)$. Letting $A^\choice$ denote the direct sum of the maps $A^\choice_a$ over $a\in \matchingsk$, the composition $C^\choice \circ A^\choice$ induces a map
\[
\qXi_0^\choice: \bigoplus_{a\in \cupdiag{n}{k}} \qCCKnk(T;a) \to \qaKC^\choice(\closure{T};n-2k).
\]
There is then an induced map
\begin{equation}
\label{eq:qXi}
\qXi^\choice : \qch_*(\qPlatnk;\qCCKnk(T)) \to \qaKC^\choice(\closure{T};n-2k)
\end{equation}
given by precomposing $\qXi_0^\choice$ with projection onto $\bigoplus\limits_{a\in \cupdiag{n}{k}} \qCCKnk(T;a)$.
\begin{proposition}\label{prop:qXi is a chain map}
    The map $\qXi^\choice$ is a chain map. 
\end{proposition}

\begin{proof}
    By Corollary \ref{cor:A is a chain map}, $\qXi_0^\choice$ vanishes on the image of the Hochschild differential. On the other hand, $\qXi_0^\choice$ commutes with the internal differential by Lemma \ref{lem:C map naturality} and far-commutativity of the saddle maps involved in $A^\choice$ and in the internal differential. 
\end{proof}

\begin{lemma}
Let $T$ be a planar $(n,n)$-tangle diagram. For $0\leq k \leq n$, the map 
\[
\qXi^\choice : \qch_*(\qPlatnk;\qCCKnk(T)) \to \qaKC^\choice(\closure{T};n-2k)
\]
is a quasi-isomorphism. 
\end{lemma}

\begin{proof}
We first prove the lemma when $T=\id_n$, the identity $(n,n)$-tangle diagram. 
For $a\in\cupdiag{n}{k}$, the image of its corresponding idempotent $e_a$ under $\qXi^\choice$ can be computed using Lemma \ref{lem:bundt pan} to be precisely the element\footnote{There are minor differences in conventions: our $v_+$ and $v_-$ correspond to their $v_1$ and $v_{-1}$, respectively; our $q$ is their $q^{-1}$; our annular degree corresponds to their weight where the roles of $k$ and $n-k$ are interchanged.} $p_a$ from \cite[Section 6]{Chen-Khovanov}.  As explained in \cite[Section 6]{Chen-Khovanov} (see also  the proof of \cite[Lemma 5.9]{LLS_CK_Spectra})  $\{ p_a \mid a\in \cupdiag{n}{k} \}$ form a basis for the annular degree $n-2k$ part of $\choiceFAq(\closure{\id_n})$. Therefore the map induced by $\qXi^\choice$ on homology is surjective. Since $\qhh(\qPlatnk;\qCCKnk(\id_n))$ and $\qaKC^\choice(\closure{\id_n};n-2k)$  are isomorphic by \cite[Proposition 6.6]{BPW} and $\Z[G]$ is a Noetherian ring, we conclude that  $\qXi^\choice$ is a quasi-isomorphism when $T= \id_n$.

For a general $T$, one can reduce to the above case of the identity tangle by arguing as in the proof of \cite[Theorem 7]{LLS_CK_Spectra}; we present it here for completeness.  If $n=0$ then the statement is immediate. For $n>0$, if $T$ is not the identity tangle then it can be decomposed as $T= T_1 T_2$ where $T_1$ is an $(n,m)$-tangle, $T_2$ is an $(m,n)$-tangle, and $m<n$. We will establish the existence of a commutative square 
\begin{equation}
\label{eq:Xi trace relation}
\begin{tikzcd}
    \qhh_0(\qPlatk{n}{k}, \qCCKnk(T_1 T_2)) \ar[d, "\qXi^\choice "] \ar[r, "\til{\tau}"] &   \qhh_0(\qPlatk{m}{h(k)}, \qCCK^{m-h(k),h(k)}(T_2 T_1)) \ar[d,"\qXi^\choice"] \\
    \qaKC^\choice(\closure{T_1 T_2};n-2k) \ar[r, "\sigma"] & \qaKC^\choice(\closure{T_2T_1};n-2k = m - 2h(k))
\end{tikzcd} 
\end{equation}
where $\til{\tau}$ and $\sigma$ are isomorphisms; this would complete the proof by induction on $n$. For simplicity, we set 
$A_{12} = \qPlatnk$, $A_{21} = \qPlatk{m}{h(k)}$, $M_1 = \qCCKnk(T_1)$, $M_2 = \qCCK^{m-h(k),h(k)}(T_2)$, $M_{12} = \qCCKnk(T_1 T_2)$, and $M_{21} = \qCCK^{m-h(k),h(k)}(T_2 T_1)$.  Consider the following diagram

\[
\begin{tikzcd}[row sep = small]
    M_1\o_{\Zq} M_2 \ar[r,"\tau"]  \ar[d, two heads, "\pi_{12}"] 
    & M_2 \o_{\Zq} M_1 \ar[d, two heads, "\pi_{21}"] 
    \\
    \qhh_0(A_{12}; M_1\o_{A_{21}} M_2) \ar[r, "\tau'"] \ar[d, "\cong"] & \qhh_0(A_{21};  M_2 \o_{A_{12}} M_1) \ar[d, "\cong"] \\
    \qhh_0(A_{12}; M_{12}) \ar[r, dashed, "\til{\tau}"] \ar[d, "\qXi^\choice"]  & \qhh_0(A_{21}; M_{21}) \ar[d, "\qXi^\choice "] \\
    \qaKC^\choice(\closure{T_1 T_2};n-2k) \ar[r, "\sigma"] & \qaKC^\choice(\closure{T_2T_1};n-2k)
\end{tikzcd}
\]
where 
\begin{itemize}
    \item $\pi_{ij}$ is the natural projection.
    \item $\tau(x\o y) = q^{-\lr{y}} y\o x$ (note that tensor products in the top row are over $\Zq$ rather than over $A_{ij}$). 
    \item $\tau'([x\o y]) = q^{-\lr{y}} [y\o x]$ (see Remark \ref{rem:trace property of qHH on chain level}).  
    \item The vertical maps marked $\cong$ are induced by the isomorphisms in Lemma \ref{lem:gluing planar tangles} (see \eqref{eq:mu_b gluing} and \eqref{eq:mu_a gluing}). 
    \item The dashed map $\til{\tau}$ is the unique map so that the middle square commutes.
    \item $\sigma =\choiceFAq(\Sigma)$, where $\Sigma : \closure{T_1 T_2} \to \closure{T_2 T_1}$ is depicted in \eqref{eq:moving T2 across the seam}. 
\end{itemize}
It is straightforward to see that
the top square commutes. Moreover, the outermost square diagram commutes 
by Lemma \ref{lem:trace property for qXi}. It follows that the bottom square commutes, which establishes \eqref{eq:Xi trace relation} and completes the proof.
\end{proof}

\begin{proposition}\label{prop:qXichoice is a quasi-iso}
    Let $T$ be an $(n,n)$-tangle diagram. For $0\leq k \leq n$, the map 
\[
\qXi^\choice : \qch_*(\qPlatnk;\qCCKnk(T))\to \qaKC^\choice(\closure{T};n-2k)
\]
is a quasi-isomorphism. 
\end{proposition}

\begin{proof}
The argument in the proof of \cite[Theorem 7]{LLS_CK_Spectra} applies with no changes. 
\end{proof}

\section{Lax lifts of the Khovanov TQFT's}
\label{sec:lax lifts of Khov TQFTs}

In Section \ref{sec:quantum annular homological background} we described various homological link invariants via TQFTs of the form $\BN(X)\rar{\FKh[X]}\Ab^\bullet$ for $X=\D,\A$ (the bullet $\bullet$ indicates grading $gr$ for $X=\D$ and bigrading $gg$ for $X=\A$), as well as the quantum annular link homology via the functor $\BNq\rar{\choiceFAq}\ZGmod$.  Here we summarize the functors $V_X^{\text{lax}}$ as described in \cite{LLS_Khovanov_Spectra,LLS_CK_Spectra}
which lift the Khovanov TQFT's $\FKh[X]$ in a suitable sense, and define our corresponding $\qVHKKlax$ to lift $\choiceFAq$.  All of these lifts require extra data in both the domain and the codomain. These will be used  to construct spectral link and tangle invariants in Section \ref{sec:spectral link invts} and Section  \ref{sec:spectral platalg and tangle bimod}, respectively. 

\subsection{The Burnside category and lax 2-functors}
\label{sec:Burnside}
In this section we will consider the codomain and the types of functors to be used in our lifts.  We begin by describing the (free) Burnside category, giving an informal description first (which yields a lax $2$-category) before giving the formal definition (giving a strict $2$-category), as in \cite[Section 2.11]{LLS_Khovanov_Spectra}.  Note that, in order to accommodate the infinite cyclic group $G=\Z$, we will allow infinite sets in our definition (unlike in \cite{LLS_Khovanov_Spectra,LLS_CK_Spectra}).  In order to have our formulas still make sense, we do impose a certain finiteness condition on the morphisms involved; see also \cite{bigger_Burnside} for a related discussion.

Let $G$ be a cyclic group. The \emph{free Burnside category} associated to $G$, denoted $\B_G$, is the $2$-category defined informally as follows. 
\begin{itemize}
    \item Objects are at most countable, free $G$-sets. 
    \item A morphism from $X$ to $Y$ is a correspondence, consisting of a free $G$-set $A$ which is at most countable, as well as $G$-equivariant \emph{source} and \emph{target} maps $X\lar{s} A \rar{t} Y$. Composition is given by fiber product.
    \item A $2$-morphism from $X\lar{s} A \rar{t} Y$ to $X\lar{s'} B\rar{t'} Y$ is a bijection $\alpha: A\to B$ which is $G$-equivariant and commutes with source and target maps:
    \[
 \begin{tikzcd}[row sep = tiny]
     & A \ar[dl, "s"'] \ar[dr, "t"] \ar[dd, "\alpha"] & \\
     X & & Y \\
        & B \ar[ul, "s'"] \ar[ur, "t'"'] & 
 \end{tikzcd}
    \]
\end{itemize}

\begin{definition}
    \label{def:G-Burn}
    Formally, $\B_G$ is the $2$-category whose objects are as above.  A $1$-morphism from $X$ to $Y$ consists of an integer $n\geq 0$ and a collection of subsets $A_{y,x} \subset \R^n$ for each $x\in X, y\in Y$ satisfying the following disjointness, equivariance, and finiteness conditions 
{
\begin{enumerate}[label= (\Alph*)]
\setcounter{enumi}{3}
    \item \label{item:D}  $A_{y,x} \cap A_{y',x} = \varnothing$ if $y\neq y'$ and $A_{y,x} \cap A_{y,x'} =\varnothing$ if $x\neq x'$.  
    \item \label{item:E} $A_{gy,gx} = A_{y,x}$ for all $x\in X, y\in Y$, and $g\in G$. 
    \item \label{item:F} For each $x\in X$, $\bigcup_{y\in Y} A_{y,x}$ is finite.
\end{enumerate}
}
 The integer $n$ will almost always be omitted from the notation. Note, the union in \ref{item:F} is disjoint by \ref{item:D}, and moreover \ref{item:F} implies that each $A_{y,x}$ is finite and that for fixed $x\in X$, $A_{y,x}$ is empty except for finitely many $y\in Y$. We also note that properties \ref{item:D} and \ref{item:E} are compatible, since freeness of the actions implies $gy=y$ if and only if $g = 1$, and similarly $gx=x$ if and only if $g=1$. 
 
  Given $1$-morphisms $A = (A_{y,x})_{x\in X, y\in Y}$ from $X$ to $Y$ and $B = (B_{z,y})_{y\in Y, z\in Z}$ from $Y$ to $Z$, their composition $(B\circ A)_{x\in X, z\in Z}$ is defined to be
 \[
 (B \circ A)_{z,x} := \bigcup_{y\in Y} B_{z,y} \times A_{y,x}.
\]

The composition satisfies \ref{item:D} as explained in \cite[Definition 2.56]{LLS_Khovanov_Spectra}. It also satisfies \ref{item:F} since for a fixed $x\in X$, there are finitely many $y\in Y$ such that $A_{y,x}$ is nonempty, and for each such $y$ there are finitely many $z$ such that $B_{z,y}$ is nonempty. Then 
\[
\bigcup_{z\in Z} \bigcup_{y\in Y} B_{z,y} \times A_{y,x}
\]
is finite. Finally, it satisfies \ref{item:E} since
\[
(B\circ A)_{gz,gx} = \bigcup_{y'\in Y} B_{gz,y'} \times A_{y',gx} = \bigcup_{y\in Y} B_{z,g^{-1} y} \times A_{g^{-1}y,x} = (B\circ A)_{z,x}.
\]
Given an object $X$, the \emph{identity} $1$-morphism $1_X$ from $X$ to $X$ is defined by $(1_X)_{x,x'} = \R^0$ if $x=x'$ and $(1_X)_{x,x'} = \varnothing$ if $x\neq x'$. With this, composition of $1$-morphisms in $\B_G$ is strictly unital and associative.

If $A$ and $B$ are $1$-morphisms from $X$ to $Y$ in $\B_G$, a $2$-morphism $\alpha: A \Rightarrow B$ is a collection of bijections $\alpha_{y,x} : A_{y,x} \to B_{y,x}$ such that $\alpha_{gy,gx} = \alpha_{y,x}$ for all $x\in X, y\in Y, g\in G$.  
\end{definition}

\begin{remark}
   The sets $A_{y,x}$ in the strict version correspond to $s^{-1}(x) \cap t^{-1}(y)$ in the lax version. To upgrade from the weak version to the above strict version the embeddings can be chosen arbitrarily except in some special cases (such as the identity morphism above) which we will carefully point out when needed. Moreover, we will often use the same notation as for the lax version; namely, we write a $1$-morphism as $X \lar{s} A \rar{t} Y$, sometimes with $s$ or $t$ (or both) omitted, and write composition of $1$-morphisms as 
    \[
\begin{tikzcd}[row sep = small, column sep =small]
    & A \ar[dl, "s"'] \ar[dr, "t"] & & B \ar[dl, "s'"'] \ar[dr, "t'"] & \\
    X & & Y & & Z 
\end{tikzcd}
\]
 with the understanding that we mean the strict versions. 
\end{remark}

\begin{definition}
\label{def: finite Burnside}
    When $G$ is finite, let $\til{\B}_G$ denote the full sub-$2$-category of $\B_G$ generated by objects which are finite sets. 
\end{definition}

Note that in the case of the trivial group, $\til{\B}_{\{1\}}$ is exactly the $2$-category in \cite[Definition 2.56]{LLS_Khovanov_Spectra}; we denote ${\B}_{\{1\}}$ and $\til{\B}_{\{1\}}$ by $\B$ and $\til{\B}$, respectively.  When $G$ is finite, our functors into $\B_G$ will land in $\til{\B}_G$, and for the most part we can work in this smaller $2$-category. However, it will occasionally be helpful to have the additional flexibility to consider infinite sets. For instance, we would like to have a $2$-functor $\B_G \to \B$ that forgets the $G$-action even when $G$ is infinite. 

We may view $\B$ as a multicategory enriched in groupoids via 
\[
\Hom_{\B}(X_1, \ldots, X_m; Y) = \Hom_{\B}(X_1 \times \cdots \times X_m, Y).
\]
Some care is needed to make Cartesian product strictly associative and unital, as explained in \cite[Section 3.2.1]{LLS_Khovanov_Spectra}.

\begin{remark}Note, we do not view $\B_G$ for $G\neq \{1\}$ as a multicategory, since Cartesian product does not reflect the correct structure. One could potentially modify the monoidal product, but we will not need this. 
\end{remark}

The \emph{graded} free Burnside category $\gBurn_G$ consists of graded free $G$-sets\footnote{A \emph{graded} (resp. bigraded) set is a set $X$ equipped with a function $\gr: X\to \Z$ (resp. $\gr:X\to \Z^2$).} $X$ and correspondences $A$ which \emph{respect} the grading, meaning that if $A_{y,x}\neq \varnothing$ then $\gr(x) = \gr(y)$.  Analogously, we also have the \emph{bigraded} free Burnside category $\ggBurn_G$.  We will use $\Burn_G^\bullet$ to denote either $\gBurn_G$ or $\ggBurn_G$  when it is clear from context which one is meant, and similarly we write $\til{\B}_G^\bullet$ for the finite versions. 

There is a forgetful functor
\[
\Burn^\bullet_G \rar{\BtoAb}
\Rmod[{\Z[G]}]^\bullet
\]
for $\bullet\in\{\varnothing,gr,gg\}$. It sends an object  $X$ to the free abelian group with basis $X$,
\[\BtoAb(X):= \bigoplus_{x\in X} \Z.\] Letting $1_x$ denote the generator of the summand $\Z$ corresponding to $x\in X$, the $\Z[G]$-module structure is given by $g\cdot 1_x = 1_{gx}$. 
The functor sends a $1$-morphism $A$ from $X$ to $Y$ to the map $\BtoAb(A):\BtoAb(X)\rightarrow \BtoAb(Y)$ defined by
\[
1_x \mapsto \sum_{y\in Y} \vert A_{y,x} \vert 1_{y}.
\] 
Condition \ref{item:F} ensures that this sum is finite, and condition \ref{item:E} ensures that this map is a map of $\Z[G]$-modules. Non-identity 2-morphisms are forgotten.

In the case of the trivial group, the above assembles into a multifunctor $\B^\bullet \rar{\BtoAb} \Ab^\bullet$ landing in abelian groups, where the multicategory structure is given via tensor product. 

\begin{remark}
\label{rem:BurnG as a lift of free modules}
If we view $\ZG$-modules as abelian groups with $G$-action, then $\Burn_G$ informally behaves as a sort of categorical lift of free $\ZG$-modules with designated ($G$-equivariant) bases. Sets in $\Burn_G$ are the bases, and a correspondence presents a map between two such modules as a (possibly infinite) $Y\times X$-matrix of finite sets $A_{y,x}$.  The cardinalities of these sets determine the corresponding coefficients in the matrix of numbers which would represent the ($G$-equivariant) map between the underlying abelian groups (thus this version is meant to lift matrices of \emph{non-negative} integers only; see \cite[Section 3]{odd_khovanov_homotopy} for some generalizations).  All of these ideas for $G=\{1\}$ are discussed in more detail in \cite[Section 2.11]{LLS_Khovanov_Spectra}.
\end{remark}

Finally, a lax 2-functor can be summarized as a 2-functor which is strictly functorial on 2-morphisms, but is only functorial on 1-morphisms up to a further 2-morphism.  In \cite[Section 4.1]{LLS-Burnside}, Lawson-Lipshitz-Sarkar present a general definition of a lax 2-functor before giving a simplified version for the cases of interest.  Since we do not need the general definition of a lax $2$-functor, we present an immediate generalization of their simplified version in the lemma below.

\begin{lemma}[{\cite[Lemma 4.4]{LLS-Burnside}}]
\label{def:lax 2-functor}
A \emph{strictly unitary lax 2-functor} $\func{F}$ from some strict 2-category $\cat{C}$ with only identity $2$-morphisms to $\Burn^\bullet_G$ consists of the following data.
\begin{itemize}
    \item For each object $x\in\cat{C}$, an object $\func{F}(x)\in\Burn^\bullet_G$.
    \item For each 1-morphism $x\rar{\phi} y$ in $\cat{C}$, a $1$-morphism $\func{F}(x) \lar{s}\func{F}(\phi)\rar{t}\func{F}(y)$ in $\Burn^\bullet_G$, such that $\func{F}(\id_x)$ is the identity $1$-morphism in $\B^\bullet_G$.
    \item For each pair of composable 1-morphisms $x\rar{\phi}y\rar{\psi}z$ in $\cat{C}$, a \emph{composition 2-morphism}
    \[\func{F}(\psi)\times_{\func{F}(y)}\func{F}(\phi)
    \rar{\Phi_{\phi,\psi}} \func{F}(\psi\circ\phi)\]
   in $\Burn^\bullet_G$.
\end{itemize}
The composition 2-morphisms are required to satisfy the following compatibility condition: for any triplet of composable 1-morphisms
\[w \rar{\phi} x \rar{\psi} y \rar{\xi} z\]
in $\cat{C}$, the following diagram of bijections commutes.
\[\begin{tikzcd}
\func{F}(\xi) \times_{\func{F}(y)} \func{F}(\psi) \times_{\func{F}(x)} \func{F}(\phi) \ar[r," \Phi_{\psi,\xi} \times \Id"] \ar[d,"\id \times \Phi_{\phi,\psi}"] &
\func{F}(\xi\circ\psi)\times_{\func{F}(x)}  \func{F}(\phi)  \ar[d,"\Phi_{\phi,\xi\circ\psi}"] \\
 \func{F}(\xi)  \times_{\func{F}(y)} \func{F}(\psi\circ\phi)  \ar[r,"\Phi_{\psi\circ\phi,\xi}"] &
\func{F}(\xi\circ\psi\circ\phi)
\end{tikzcd}\]

\end{lemma}
\begin{proof}
    The proof is entirely analogous to the proof of \cite[Lemma 4.4]{LLS-Burnside}.
\end{proof}

\subsection{Divided cobordism categories and translating Khovanov TQFT's}
\label{sec:old divcob}
Here we recall the divided cobordism categories $\divcob(X)$ of \cite{LLS_Khovanov_Spectra,LLS_CK_Spectra} for $X\in\{\D,\A\}$.  The first definition below is used to define equivalences among objects (1-manifolds) in $\divcob(X)$.  Note that we use the interval $[0,2\pi]$ to define $S^1$ to avoid confusion with the interval $I=[0,1]$ to appear shortly.

\begin{definition}[{\cite[Definition 2.46]{LLS_Khovanov_Spectra} and \cite[Definition 5.11]{LLS_CK_Spectra}}]
\label{def:Diff' groups for objects}
Consider the open interval $(0,1)$ (that will contribute to $X=\D=(0,1)\times(-1,1)$).  We let $\Diff'((0,1))$ denote the group of orientation-preserving diffeomorphisms of $(0,1)$ which are the identity in some $\eps$-neighborhood $B_{\eps}(0\sqcup 1)$ of the boundary.

Similarly, let $S^1:=[0,2\pi]/\partial$ (that will contribute to $X=\A=S^1\times(-1,1)$).  We let $\Diff'(S^1)$ denote the group of orientation-preserving diffeomorphisms of $S^1$ which are the identity in some $\eps$-neighborhood $B_\eps(0)$ of the point $0\in S^1$.
\end{definition}

The next definition is used to define equivalences among morphisms (2-dimensional cobordisms) in $\divcob(X)$.  The closed interval $I=[0,1]$ will correspond to the `cobordism-direction'.

\begin{definition}[{\cite[Definition 2.47]{LLS_Khovanov_Spectra} and \cite[Definition 5.11]{LLS_CK_Spectra}}]
\label{def:Diff' groups for morphisms}
Let $W$ be either $(0,1)$ or $S^1$ as in Definition \ref{def:Diff' groups for objects}, and let $I=[0,1]$.  We define $\Diff'(I \times W)$ to be the group of orientation-preserving diffeomorphisms $ \psi: I \times W \to I \times W$ such that there exists some $\eps >0$ and $\phi_0, \phi_1\in \Diff'(W)$ satisfying the following:
\begin{itemize}
    \item $\psi(p,q)  = (p,\phi_0(q))$ for all $p\in [0,\eps)$; 
    \item $\psi(p,q) = (p,\phi_1(q))$ for all $p\in (1-\eps,1]$;
    \item for the case $W=(0,1)$, $\psi\vert_{I\times B_\eps(0\sqcup 1)} = \id$; and
    \item for the case $W=S^1$, $\psi\vert_{I\times B_\eps(0)}=\xi\times\id_{B_\eps(0)}$ for some orientation preserving diffeomorphism $\xi$ of $[0,1]$ (note that, by the previous bullet points, $\xi$ will be the identity near $0$ and $1$).
\end{itemize}
\end{definition}

We now present the definition of the divided cobordism categories $\divcob(X)$.

\begin{definition}[{\cite[Sections 3.1 and 3.23]{LLS_Khovanov_Spectra} and \cite[Definition 5.12]{LLS_CK_Spectra}}]
\label{def:divcob}
Let $W$ denote either $(0,1)$ or $S^1$ as in Definition \ref{def:Diff' groups for objects}.  Let $X=W\times (-1,1)$ denote either $\D$ or $\A$ accordingly.  The \emph{divided cobordism category of $X$}, denoted $\divcob(X)$, is a multicategory enriched in groupoids defined as follows.
\begin{itemize}
    \item A \emph{divided 1-manifold} (object) is an equivalence class of smooth closed embedded 1-manifolds $Z\subset X$ divided into \emph{active} arcs $A\subset Z$ and \emph{inactive} arcs $Z\setminus A$ (when $W=S^1$, we require $Z$ be transverse to $\{0\} \times (-1,1)$).
    Two such 1-manifolds are equivalent if there is $\phi\in\Diff'(W)$ such that $\left( \phi\times\id_{(-1,1)}\right)$ takes one to the other.
    \item A \emph{divided cobordism} (morphism) from $(Z_1,A_1)$ to $(Z_2,A_2)$ is an equivalence class of pairs $(\Sigma,\Gamma)$ consisting of
    \begin{itemize}
        \item a smoothly embedded cobordism $\Sigma\subset I\times X$ from $Z_1$ to $Z_2$ (when $W=S^1$, we further require $\Sigma$ to be transverse to $I\times \{0\} \times (-1,1)$), and
        \item a collection of smooth, properly embedded arcs $\Gamma\subset \Sigma$ with $\partial \Gamma = (\partial A_1 \sqcup \partial A_2)$,
     so that every component of $\Sigma\setminus \Gamma$ has one of the following forms.
        \begin{enumerate}
            \item A rectangle, with two sides components of $\Gamma$ and two sides components of $A\sqcup A'$.
            \item A $(2\ell+2)$-gon, with $(\ell+1)$ sides components of $\Gamma$, one side a component of $Z'\setminus A'$, and the other $\ell$ sides components of $Z\setminus A$. (The integer $\ell$ is allowed to be zero.)
        \end{enumerate}
    \end{itemize}
    Two such pairs are equivalent if there is $\psi\in\Diff'(I\times W)$ such that $\left( \psi\times\id_{(-1,1)}\right)$ takes one to the other.
    \item In the case of $W=(0,1)$, more general multimorphisms in $\divcob(\D)$ from 
    \[
    ((Z_1,A_1),\dots,(Z_{b-1},A_{b-1}))
    \]
    to $(Z_b,A_b)$ are defined as divided cobordisms from $(Z_1, A_1)\sqcup \cdots \sqcup (Z_{b-1}, A_{b-1})$ to $(Z_b, A_b)$; in the case of $W=S^1$, there are no multimorphisms in $\divcob(\A)$ with more than one input, so $\divcob(\A)$ is an ordinary category (which we may view as a multicategory).
    \item Composition of (multi-)morphisms is induced by concatenating divided cobordisms.
    \item There is a unique 2-morphism between divided cobordisms $(\Sigma,\Gamma)\rar{}(\Sigma',\Gamma')$ if and only if they are isotopic rel boundary.\footnote{For us this means an ambient isotopy rel boundary which is the identity outside of some compact subspace of $I\times X$.}  
\end{itemize}
\end{definition}

Note that in $\divcob(\D)$ we quotient our embeddings by diffeomorphisms which can affect the $(0,1)$-direction to allow disjoint union to be sensible, and the $I$-direction so that composition is strictly unital and associative.  We do not quotient by more general diffeomorphisms which could potentially interchange connected components of our divided 1-manifolds (conditions on the divides ensure that we will not encounter closed divided cobordisms which potentially could be interchanged). Meanwhile, when two divided cobordisms are isotopic (rel boundary), we record this fact with a unique 2-morphism rather than keeping track of the (isotopy class) of the isotopy.  This loss of data is justified by the presence of the divides \cite[Lemma 3.1]{LLS_Khovanov_Spectra}. Similarly, in $\divcob(\A)$ we quotient in the $I$-direction for the same reason as above and in the $S^1$-direction via $\Diff'(S^1)$ which fixes $0\in S^1$, allowing Hochschild-Mitchell constructions to be well-defined, while preventing diffeomorphisms from interchanging components. Moreover, $\divcob(\A)$ is a category rather than a multicategory (see Remark \ref{rmk:HMshape is cat not multicat}). 

\begin{remark}\label{rmk:LLS Divcob vs our Divcob}
Our definition for $\Diff'(I\times(0,1))$ matches that of \cite[Definition 2.47]{LLS_Khovanov_Spectra} (denoted $\Diff^2$ there), and thus our $\divcob(\D)$ matches the definition of $\rm{Cob}_d$ given in \cite[Section 3.1]{LLS_Khovanov_Spectra}.  

However, in \cite[Definition 5.12]{LLS_CK_Spectra} for $\mathrm{Cob}_d(\A)$, cobordisms are required to be vertical (invariant in the $I$-direction) near  $I\times \{0\}\times (-1,1)$; with this requirement, their definition of $\Diff'(I\times S^1)$ \cite[Definition 5.11]{LLS_CK_Spectra} demands $\psi$ to be the identity on $I\times B_\eps(0)$, and this is enough to guarantee strict unitality and associativity in their $\mathrm{Cob}_d(\A)$.  In our definition for $\divcob(\A)$, we only require our cobordisms to be transverse to $I\times\{0\}\times (-1,1)$, but as a result our allowable equivalences on $1$-morphisms via $\Diff'(I\times S^1)$ need to allow stretching in the vertical direction near $I\times\{0\}$ to maintain strict unitality and associativity.

Moreover, in \cite[Definition 5.12]{LLS_CK_Spectra}, an isotopy realizing a $2$-morphism is required to fix $I\times\{0\}\times (-1,1)$ pointwise. Here we do not impose this condition.  Later in Definition \ref{def:qdivcob of annulus}, $I\times\{0\}\times (-1,1)$ will play the role of the membrane $\Mq$. 

\end{remark}

Finally, for $X\in\{\D,\A\}$, we translate the data of our Khovanov functors $\FKh[X]$ out of Bar-Natan categories $\BN(X)$ from Section \ref{sec:Platform algebras and bimodules} into functors with same notation $\FKh[X]$ out of the divided cobordism categories $\divcob(X)$ instead, where the objects are equivalence classes of 1-manifolds.

\begin{lemma}\label{lem:FKh and FA work on divcob}
For $X\in\{\D,\A\}$ there exists a multifunctor 
\[\divcob(X)\rar{\FKh[X]}\Ab^\bullet\]
which forgets 2-morphisms, and is otherwise determined by the functor $\BN(X)\rar{\FKh[X]}\Ab^\bullet$ as follows (here the bullet $\bullet$ indicates  $gr$ for $X=\D$ and $gg$ for $X=\A$).
\begin{itemize}
    \item On objects, $\FKh[X](Z,A) = \FKh[X](\tilde{Z})$, where $\tilde{Z}\in\BN(X)$ is any 1-manifold in the equivalence class of $Z$ (and the arc data $A$ is ignored).
    \item On morphisms, $\FKh[X](\Sigma,\Gamma) = \FKh[X](\tilde{\Sigma})$, where $\tilde{\Sigma}$ is any morphism in $\BN(X)$ (cobordism) in the equivalence class of $\Sigma$ (and the divide data $\Gamma$ is ignored).
\end{itemize}
\end{lemma}
\begin{proof}
    Equivalent objects and morphisms in $\divcob(X)$ are related by isotopies which cannot exchange components, ensuring that $\FKh[X]$ is well-defined on $\divcob(X)$.
\end{proof}

\subsection{Defining the lax lifts $V_X^{\text{lax}}$ for $X\in\{\D,\A\}$}
\label{sec:VHKKlax for D,A}
In \cite[Proposition 3.2]{LLS_Khovanov_Spectra}, building on \cite{HKK}, the authors define a strictly unitary lax 2-functor (see Lemma \ref{def:lax 2-functor}) $\divcob(\D)\rar{\VHKKlax}\gBurn$ which fits in to the following diagram \begin{equation}\label{eq:VHKKlax lifts Kh}
\begin{tikzcd}
 & \gBurn \ar[d,"\BtoAb"]\\
\divcob(\D) \ar[ur,"\VHKKlax"] \ar[r,"\FKh"] & \gAb.
\end{tikzcd}
\end{equation} 
We now describe $\VHKKlax$,  where we view $\divcob(\D)$ as an ordinary 2-category, ignoring the multicategorical/monoidal structure.  Large portions of this description are taken nearly verbatim from \cite[Section 2.11]{LLS_Khovanov_Spectra}; however, we include some extra details about the composition 2-isomorphisms which will be relevant when describing the quantum case.

\begin{remark}\label{rem:VHKK largely determined by FKh}
As in Remark \ref{rem:BurnG as a lift of free modules}, for objects $(Z,A)\in\divcob(\D)$, $\VHKKlax(Z,A)$ must coincide with a basis for $\FKh(Z)$ if we are to satisfy the diagram \eqref{eq:VHKKlax lifts Kh}.  Similarly, if $(\Sigma,\Gamma)$ is a 1-morphism from $(Z,A)$ to $(Z',A')$, then for each $x\in\VHKKlax(Z,A)$ and $y\in\VHKKlax(Z',A')$, the cardinality of $s^{-1}(x)\cap t^{-1}(y) \subset \VHKKlax(\Sigma,\Gamma)$ is equal to the coefficient of $y$ in $\FKh(\Sigma)(x)$; in particular, if all such coefficients are 0 or 1, then the value of $\VHKKlax$ on any 2-morphism $(\Sigma, \Gamma) \Rightarrow (\Sigma', \Gamma')$ is uniquely determined. In light of this, many of the descriptions below are determined by $\FKh$; the main new ingredient, provided by \cite{HKK, LLS_Khovanov_Spectra} is in prescribing the two element correspondences for connected genus 1 morphisms $\Sigma$, and compatible 2-morphisms between them.
\end{remark}

For objects $(Z,A)\in\divcob(\D)$ (up to equivalence, see Definition \ref{def:divcob}), we ignore the active arc data $A$ and define $\VHKKlax(Z,A)$ to be the set $\Gen(Z)$ of Definition \ref{def:FKh on BN(D)}, i.e. the set of functions $\pi_0(Z)\rar{}\{1,X\}$.  Note that the allowable equivalences cannot interchange components of $Z$, and so this functor is well-defined on objects. Since the value of $\VHKKlax$ on an object $(Z,A)$ does not depend on $A$, we will often simply write $\VHKKlax(Z)$. 

To define $\VHKKlax$ on 1-morphisms $(\Sigma,\Gamma)$ in $\divcob(\D)$, we ignore the divides $\Gamma$ and fix a checkerboard coloring of the complement of $\Sigma$ in $I\times \D$ with the outer region colored white. Since $\VHKKlax(\Sigma,\Gamma)$ will not depend on $\Gamma$, we will typically simply write $\VHKKlax(\Sigma)$. The value of $\VHKKlax(\Sigma)$ as a 1-morphism in $\gBurn$ (correspondence with source and target maps $s$ and $t$) is the product over the components $\Sigma'$ of $\Sigma$ of $\VHKKlax(\Sigma')$ with respect to the checkerboard coloring of the complement of $\Sigma'$ induced from the $\Sigma$-complement coloring by declaring that the two colorings agree in a neighborhood of $\Sigma'$.  

Thus we may assume that $\Sigma$ is a connected cobordism between 1-manifolds $Z_0,Z_1$, but now the checkerboard coloring is arbitrary.  If $\Sigma$ has genus $> 1$ then $\VHKKlax(\Sigma):=\varnothing$.  Otherwise, we fix $x\in\VHKKlax(Z_0)$ and $y\in\VHKKlax(Z_1)$.  If $\Sigma$ has genus zero, then $s^{-1}(x)\cap t^{-1}(y) \subset \VHKKlax(\Sigma)$ has at most one element and $\VHKKlax(\Sigma)$ is determined by the diagram \eqref{eq:VHKKlax lifts Kh}.  If $\Sigma$ has genus 1, then $s^{-1}(x)\cap t^{-1}(y)$ is empty unless $x$ labels each circle in $Z_0$ by $1$, and $y$ labels each circle in $Z_1$ by $X$.  In this case, $\VHKKlax(\Sigma)$ has two elements, which we describe as follows.  Let $S^2$ denote the one point compactification of $\D$.  Let 
\[B(\Sigma)\subset ([0,1]\times S^2) \setminus \Sigma\]
denote the black region in the checkerboard coloring (possibly extended to the new points at infinity), and let
\[B^\partial(\Sigma):= B(\Sigma) \cap (\{0,1\}\times S^2)\]
denote the black region at the top and bottom boundary.    Then  define $\VHKKlax(\Sigma)$ to be the set of generators of
\begin{equation}
\label{eq:two homology generators} 
H_1(B(\Sigma)) / H_1(B^\partial(\Sigma)) \cong \Z.
\end{equation}

If there is a 2-morphism $(\Sigma,\Gamma)\rar{f}(\Sigma',\Gamma')$ in $\divcob(\D)$, then \cite[Lemma 3.1]{LLS_Khovanov_Spectra} ensures that the isotopy $\Sigma\rar{f}\Sigma'$ is unique (up to isotopy of isotopies).  Thus after extending to the thickened 1-point compactification $[0,1]\times S^2$,\footnote{Such an extension exists since our isotopies were compactly supported.} $f$ induces a well-defined bijection on the homology generators $\VHKKlax(\Sigma)\rar{\VHKKlax(f)}\VHKKlax(\Sigma')$.  (This is the content of \cite[Proposition 3.2]{LLS_Khovanov_Spectra}.)

Finally, the composition 2-isomorphism is obvious except when composing two genus 0 components
\[ Z_0 \rar{\Sigma_0} Z_1 \rar{\Sigma_1} Z_2 \]
to obtain a genus 1 component
\[ Z_0 \rar{\Sigma} Z_2\]
(again we are ignoring the divides).  In this case, it again suffices to assume that $\Sigma$ is connected.  Let
\[B^{1/2}(\Sigma):= B(\Sigma)\cap (\{1/2\}\times S^2), \quad W^{1/2}(\Sigma):=
(\{1/2\}\times S^2) \setminus B^{1/2}(\Sigma)\]
denote the black and white regions at the `middle time slice' of the cobordism (where we have glued $\Sigma_0$ to $\Sigma_1$ along $Z_1$).  For any curve $C$ on $\Sigma$, let $C_B$ and $C_W$ be its pushoffs into $B^{1/2}(\Sigma)$ and $W^{1/2}(\Sigma)$ respectively.

Now consider an element
\[ (b,a) \in \VHKKlax(\Sigma_1) \times_{\VHKKlax(Z_1)} \VHKKlax(\Sigma_0).\]

Such an element determines a unique component (curve) $C$ of $Z_1=\partial\Sigma_0\cap\partial\Sigma_1$ that is non-separating in $\Sigma$ and is labelled $1$ by $y:=t(a)=s(b)\in\VHKKlax(Z_1)$.  We orient $C$ as the boundary of $B^{1/2}(\Sigma)$ and consider its two pushoffs $C_B$ and $C_W$.  One of these is a generator of $H_1(([0,1]\times S^2)\setminus \Sigma) / H_1((\{0,1\}\times S^2)\setminus \partial\Sigma) \cong \Z^2$, and the other is zero.  If $C_B$ is the generator, then the composition 2-isomorphism is defined by sending
\[(b,a) \mapsto [C_B] \in \VHKKlax(\Sigma).\]
If $C_W$ is the generator, then we let $D$ be another curve on $\Sigma$ oriented so that the algebraic intersection number $D\cdot C=-1$ (with $\Sigma$ being oriented as the boundary of $B(\Sigma)$) and define the composition 2-isomorphism by sending
\[(b,a) \mapsto [D_B] \in \VHKKlax(\Sigma).\]

\begin{remark}
\label{rmk:intersection is -1}
    In \cite[Section 2.11]{LLS_Khovanov_Spectra} the curve $D$ is chosen so that $D\cdot C = 1$. As explained in \cite[Remark 8.4]{LLS-Burnside}, the two choices are essentially the two global choices of ladybug matchings in \cite{LS-Kh-spectra}, called the left pair and the right pair. For the comparison with \cite{AKW} in Section \ref{sec:comparison}, we need to choose the left pair in light of \cite[Corollary 2.21]{AKW}, which corresponds to $D\cdot C = -1$; see also  Lemma \ref{lem:AKW functor agrees with HKK}. 
\end{remark}

The description above gives rise to the lax 2-functor $\VHKKlax$ fitting into diagram \eqref{eq:VHKKlax lifts Kh}.  There is an analogous construction for $\VHKKannlax$ fitting into the diagram

\begin{equation}\label{eq:VA lifts AKh}
\begin{tikzcd}
 & \ggBurn \ar[d,"\BtoAb"]\\
\divcob(\A) \ar[ur,"\VHKKannlax"] \ar[r,"\FA"] & \ggAb
\end{tikzcd},
\end{equation}
using the fact that $\FA$ is the annular-degree preserving part of $\FKh$.  To avoid excessive repetition, we will not spell out all of the details here.  Instead we point out some key properties of $\VHKKannlax$.

\begin{itemize}
    \item For an object $(Z,A)$ of $\divcob(\A)$, $\VHKKannlax(Z,A)=\Gen(Z)$, i.e. the set of ways to label each essential circle of $Z$ with either $v_+$ or $v_-$ and each trivial circle with either $1$ or $X$. The bigradings are as in \eqref{eq:degrees in QAH}. 
    \item Like $\VHKKlax$, $\VHKKannlax$ also splits as a product over disjoint unions of $1$-morphisms.
    \item For a $1$-morphism $(\Sigma, \Gamma)$ of $\divcob(\A)$, we continue to write $\VHKKannlax(\Sigma)$ since the value does not depend on $\Gamma$.  If $\Sigma$ has contractible boundary then $\VHKKannlax(\Sigma) = \VHKKlax(\Sigma)$, where on the right-hand side we view $\Sigma$ as a cobordism in $I\times \D$ via some fixed embedding of $\A$ into $\D$. On the other hand, the value of $\VHKKannlax$ on any connected annular cobordism with at least one non-contractible boundary component is uniquely determined by the diagram \eqref{eq:VA lifts AKh}, as in Remark \ref{rem:VHKK largely determined by FKh}.
    \item The value of $\VHKKannlax$ on $2$-morphisms as well as the compatibility requirements to form a lax $2$-functor are entirely analogous to the definition of $\VHKKlax$ (via embedding $\A\hookrightarrow S^2$ as the complement of two points).  
\end{itemize}

\subsection{The quantum divided cobordism category $\divcobq$}
\label{sec:defining qdivcob}
We now turn our attention to lifting the quantum annular functor $\BNq\rar{\choiceFAq}\ZGmod$, beginning with adjusting the domain accordingly.  For each cyclic group $G=\langle q \rangle$, we will define a corresponding quantum divided cobordism category $\divcobq$; we will write  $\divcobq[G]$  when the dependence on $G$ needs to be specified. 

Recall Definition \ref{def:quantum annulus} for the quantum annulus $\Aq$ equipped with a seam $\mq$ determining a membrane $\Mq=I\times\mq \subset I\times \Aq$, to which all embeddings are assumed transverse.  From this point onwards we will fix coordinates so that $\mq=\{0\}\times(-1,1)\subset S^1\times (-1,1)$. Recall also $\Diff'(S^1)$ and $\Diff'(I\times S^1)$ of Definitions \ref{def:Diff' groups for objects} and \ref{def:Diff' groups for morphisms}.

\begin{definition}
\label{def:qdivcob of annulus}
The \emph{quantum divided cobordism category of the annulus}, denoted $\divcobq$, is a category enriched in groupoids defined as follows.

\begin{itemize}
    \item A \emph{divided 1-manifold} (object) is an equivalence class of smooth closed embedded 1-manifolds $Z\subset \Aq$ divided into \emph{active} arcs $A\subset Z$ and \emph{inactive} arcs $Z\setminus A$.  We declare $(Z,A)$ and $(Z',A')$ to be equivalent if there is $\phi \in \Diff'(S^1)$ such that $\left( \phi \times \id_{(-1,1)}\right) ( Z,A) = (Z',A')$. 
    \item A morphism (\emph{divided cobordism}) from $(Z_1,A_1)$ to $(Z_2,A_2)$ is an equivalence class of triples $(\Sigma, \Gamma,q^r)$ (sometimes written $q^r(\Sigma, \Gamma)$, or just  $q^r\Sigma$, with the divides $\Gamma$ understood) consisting of
    \begin{itemize}
        \item a cyclic group element $q^r\in G$,
        \item a smoothly embedded cobordism $\Sigma \subset I \times \A$ from $Z_1$ to $Z_2$, which is transverse to $\Mq$ and invariant in the $I$-direction near $\{0,1\}\times \A$, and
        \item a collection of smooth, properly embedded arcs $\Gamma\subset\Sigma$ with $\partial \Gamma = (\partial A \sqcup \partial A')$
        so that every component of $\Sigma\setminus \Gamma$ is of the form specified in Definition \ref{def:divcob}. 
        We declare $(\Sigma,\Gamma,q^r)$ and $(\Sigma',\Gamma',q^{r'})$ to be equivalent if $q^r=q^{r'}$ in $G$ and there is $\psi \in \Diff'(I \times S^1)$ such that $\left(\psi \times \id_{(-1,1)} \right)(\Sigma,\Gamma) = (\Sigma',\Gamma')$.
       
        \item Composition of morphisms is induced by concatenating divided cobordisms and group multiplication, so that $q^r\Sigma \circ q^{r'}\Sigma' = q^{r+r'}(\Sigma\circ\Sigma')$.  Note that this involves choosing composable representatives of each equivalence class; the nature of $\Diff'(I\times S^1)$ ensures that this operation is well-defined.
        \end{itemize}
        \item There is a unique 2-morphism $(\Sigma,\Gamma,q^r)\rar{}(\Sigma',\Gamma',q^{r'})$ between divided cobordisms  if and only if, for some representatives $(\Sigma_0, \Gamma_0)$ and $(\Sigma_0', \Gamma_0')$ of $(\Sigma,\Gamma)$ and $(\Sigma',\Gamma')$, respectively, there exists a sequence of isotopies of $I\times \Aq$ fixing $\Mq$ set-wise and trace moves (Figure \ref{fig:trace moves}) from $q^r \Sigma_0$ to $q^{r'} \Sigma_0'$.  (In particular, $\FAq(q^r\Sigma) = \FAq(q^{r'}\Sigma')$.)
    \end{itemize}
\end{definition}

We also have quotient functors relating various quantum divided cobordism categories.
\begin{proposition}\label{prop:qdivcob G to G/H}
Let $H=\langle q^h \rangle$ be a (cyclic) subgroup of $G=\langle q \rangle$.  Then there exists a \emph{quotient functor} (surjective on objects, 1-morphisms, and 2-morphisms)
\[ \divcobq[G] \rar{/q^h} \divcobq[G/H]\]
which:
\begin{itemize}
    \item is the identity on objects;
    \item sends any 1-morphism $q^r\Sigma$ to $[q^r]\Sigma$ (i.e. sends the coefficient $q^r\in G$ to its equivalence class in $G/H$); and
    \item sends any 2-morphism $(\Sigma,\Gamma,q^r)\rar{}(\Sigma',\Gamma',q^{r'})$ to the corresponding 2-morphism $(\Sigma,\Gamma,[q^r])\rar{}(\Sigma',\Gamma',[q^{r'}])$.
\end{itemize}
\end{proposition}
\begin{proof}
This is a simple consequence of the definitions.
\end{proof}

\subsection{Viewing $\choiceFAq$ as a functor $\divcobq[\Z]\rar{\choiceFAq}  \ggRmod[{{\Z[G]}}]$} \label{sec:FAq with choices summary}
As discussed in Remark \ref{rem:BurnG as a lift of free modules}, the free $G$-Burnside category is designed to lift free $\ZG$-modules with designated bases.  It is for this reason that we focus on the functor $\BNq \rar{\choiceFAq} \ZGmod$ of Definition \ref{def:FAq choice}, rather than the `basis-free' version $\FAq$.  As in Section \ref{sec:old divcob}, we seek to reinterpret $\choiceFAq$ as a functor (with same notation) $\divcobq\rar{\choiceFAq}\ZGmod$.  In order to do so, we will require further restrictions on our choice of bases.

Recall the groups $\Diff'(S^1)$ and $\Diff'(I\times S^1)$ from Definitions \ref{def:Diff' groups for objects} and \ref{def:Diff' groups for morphisms}.  Note that, as is implicit in the proof of Lemma \ref{lem:FKh and FA work on divcob}, if $Z,Z'\subset \Aq$ are two collections of circles such that $Z'=(\phi\times \id_{(-1,1)})(Z)$ for some $\phi\in\Diff'(S^1)$, then the sets of Khovanov generators $\Gen(Z)$ and $\Gen(Z')$ (see Definition \ref{def:FA on BN(A)}) are canonically identified since components cannot be interchanged.

The following definition and lemma will make heavy use of the notions of standard circles, choices of generators, and regular cobordisms introduced in Definitions \ref{def:standard circles} and \ref{def:choice of generators}.

\begin{definition}\label{def:coherent choice of generators}
    A \emph{coherent choice of generators} is a choice of generators $\choice$ such that, for any two cobordisms $\Sigma,\Sigma'\subset I\times\Aq$ with $\Sigma'=(\psi\times \id_{(-1,1)})(\Sigma)$ for some $\psi\in\Diff'(I\times S^1)$, we have $\choiceFAq(\Sigma)=\choiceFAq(\Sigma')$, where we canonically identify the domains and codomains of $\choiceFAq(\Sigma)$ and $\choiceFAq(\Sigma')$ as above.
\end{definition}

\begin{lemma}\label{lem:coherent generators exist}
There exist coherent choices of generators.
\end{lemma}
\begin{proof}
Define an equivalence relation on the set of all embedded circles in $\Aq$ by declaring that $Z\sim Z'$ if and only if there exists $\phi\in\Diff'(S^1)$ such that $Z'=(\phi\times\id_{(-1,1)})(Z)$.  For every such equivalence class, choose a representative $Z_r$ and a regular cobordism $\choice_{Z_r}$ from $Z_r^\st$ to $Z_r$  satisfying the following.
\begin{itemize}
    \item If $Z_r^\st = Z_r$, we choose the identity cobordism $\choice_{Z_r}:= I\times Z_r$.
    \item If $Z_r$ contains any standard essential circle $C$, we choose the corresponding component of $\choice_{Z_r}$ to be the identity cobordism $I\times C$. 
\end{itemize}

Now let $Z$ be any collection of circles.  If $Z$ is standard, we continue to use the cylinder $\choice_Z=I\times Z$ as required.  Otherwise, let $Z_r$ be the chosen representative in the equivalence class of $Z$, so that there exists $\phi\in\Diff'(S^1)$ with $Z=(\phi\times\id_{(-1,1)})(Z_r)$.  We can pick a homotopy from  $\id_{S^1}$  to $\phi$ through elements of $\Diff'(S^1)$.  By extending via the identity across $(-1,1)$, this homotopy forms a regular cobordism $\Sigma_{Z_r,Z} \subset I\times \Aq$ from $Z_r$ to $Z$ which is invariant in the $I$-direction near $\Mq$ and $\{0,1\}\times \Aq$. We note that $\Sigma_{Z_r,Z}$ does not depend on $\phi$ or the chosen homotopy up to ambient isotopy of $I\times \Aq$ fixing a neighborhood of $\Mq \sqcup (\{0,1\}\times \Aq)$.  Then we define $\choice_{Z}$ to be the concatenation $\Sigma_{Z_r,Z}\choice_{Z_r}$ (note that we implicitly have chosen $Z^\st$ to be $Z_r^\st$ in this way).  This gives a choice of generators $\choice$ and we seek to show that it is coherent.

Suppose $\Sigma\subset I\times\Aq$ is a cobordism from $Z_0$ to $Z_1$, with $\Sigma'=(\psi \times \id_{(-1,1)})(\Sigma)$ for some $\psi\in\Diff'(I\times S^1)$ a cobordism from $Z'_0$ to $Z'_1$.  In particular there are some $\phi_0,\phi_1\in\Diff'(S^1)$ such that $Z'_i=(\phi_i\times \id_{(-1,1)})(Z_i)$ so that $Z_i\sim Z'_i$ for $i=0,1$.  Given a pair of basis elements (Khovanov generators) $x\in \Gen(Z_0)\subset \choiceFAq(Z_0)$ and $y\in \Gen(Z_1)\subset \choiceFAq(Z_1)$, let $c_{x,y}$ denote the coefficient of $y$ in the value of
\[\choiceFAq(\Sigma)(x)=\cdots + c_{x,y} y + \cdots.\]
Then as described above, $\phi_0$ identifies $x$ with some $x'\in\Gen(Z_0')\subset \choiceFAq(Z_0')$, and similarly for $\phi_1$, so that we can let $c_{x',y'}$ denote the corresponding coefficient in
\[\choiceFAq(\Sigma')(x') = \cdots + c_{x',y'} y' + \cdots.\]
Our goal is to show that $c_{x,y}=c_{x',y'}$.

Let $Z_{i,r}$ denote the chosen representative for the equivalence class of $Z_i$ and $Z'_i$, with corresponding standard $Z_i^\st$.  Altogether we have the following diagram in $\BNq$ for each $i=0,1$.
\[
\begin{tikzcd}
 & & & Z_i \ar[dd,"\Sigma_{\phi_i}"] \\
Z_i^\st \ar[r,"\choice_{Z_{i,r}}"] \ar[urrr,bend left, "\choice_{Z_i}"] \ar[drrr,bend right,"\choice_{Z_i'}",swap] & Z_{i,r} \ar[urr,"\Sigma_{Z_{i,r},Z_i}"] \ar[drr,"\Sigma_{Z_{i,r},Z'_i}",swap] & & \\
 & & & Z'_i
\end{tikzcd}
\]
The two curved triangles commute by our definition for $\choice$, while the rightmost triangle commutes in $\BNq$ since all of the cobordisms involved are regular cobordisms tracing out homotopies of $\phi_i\in\Diff'(S^1)$.  Moreover, since these cobordisms in the rightmost triangle cannot permute components, we have that $x\in\Gen(Z_0)$ and $x'\in\Gen(Z_0')$ both correspond to the same Khovanov generator $x^\st\in\FAq(Z_0^\st)$ and thus to the same element $x_r=\FAq(\choice_{Z_{0,r}})(x^\st)$, and similarly for $y,y'$ corresponding to one Khovanov generator $y^\st\in\FAq(Z_1^\st)$ and one element $y_r=\FAq(\choice_{Z_{1,r}})(y^\st)$.  In the following computations, we use the notation $\b{\Sigma}$ to indicate the cobordism $\Sigma$ viewed in the reverse direction.

By definition $c_{x,y}$ is determined by the value $c_{x^\st,y^\st}$ in
\[\FAq(\choice_{Z_0}\Sigma\b{\choice_{Z_1}})(x^\st) = \cdots + c_{x^\st,y^\st}y^\st + \cdots\]
which is then identical to the value $c_{x_r,y_r}$ in
\[\FAq( \Sigma_{Z_{0,r},Z_0}\Sigma\b{\Sigma_{Z_{1,r},Z_1}} ) (x_r) = \cdots + c_{x_r,y_r}y_r + \cdots,\]
and likewise $c_{x',y'}$ is determined by the value $c'_{x_r,y_r}$ in
\[\FAq( \Sigma_{Z_{0,r},Z'_0}\Sigma'\b{\Sigma_{Z_{1,r},Z'_1}} ) (x_r) = \cdots + c'_{x_r,y_r}y_r + \cdots.\]
But the cobordisms $\Sigma_{Z_{0,r},Z_0}\Sigma\b{\Sigma_{Z_{1,r},Z_1}}$ and $\Sigma_{Z_{0,r},Z'_0}\Sigma'\b{\Sigma_{Z_{1,r},Z'_1}}$ are equal in $\BNq$ since they are related by diffeomorphisms in $\Diff'(I\times S^1)$, and thus
\[ \FAq(\Sigma_{Z_{0,r},Z_0}\Sigma\b{\Sigma_{Z_{1,r},Z_1}}) = \FAq(\Sigma_{Z_{0,r},Z'_0}\Sigma'\b{\Sigma_{Z_{1,r},Z'_1}})\]
implying $c_{x_r,y_r}=c'_{x_r,y_r}$ as required.
\end{proof}

 \begin{lemma}
     \label{lem:FAq choice well-defined on divcob}
For any coherent choice of generators $\choice$ there exists a functor
\[\divcobq\rar{\choiceFAq}\ggRmod[{{\Z[G]}}]\]
which forgets 2-morphisms, and is otherwise determined by the functor $\BNq\rar{\choiceFAq}\ggRmod[{{\Z[G]}}]$ as follows.
\begin{itemize}
    \item On objects, $\choiceFAq(Z,A) = \choiceFAq(\tilde{Z})$, where $\tilde{Z}\in\BNq$ is any 1-manifold in the equivalence class of $Z$ (and the arc data of $A$ is ignored).
    \item On morphisms, $\choiceFAq(q^r\Sigma,\Gamma) = q^r\choiceFAq(\tilde{\Sigma})$, where $\tilde{\Sigma}$ is any morphism in $\BNq$ (cobordism) in the equivalence class of $\Sigma$ (and the divide data $\Gamma$ is ignored).
\end{itemize}
 \end{lemma}

 \begin{proof}
The fact that $\choiceFAq$ is well-defined on objects follows just as in Lemma \ref{lem:FKh and FA work on divcob}.  Then $\choice$ being coherent by definition implies that $\choiceFAq$ is well-defined on morphisms.
 \end{proof}

\subsection{Defining $\qVHKKlax$}\label{sec:qVHKKlax}
\subsubsection{The case where $G=\Z$}
In this section we set $G=\Z$ and provide the construction for our lax 2-functor $\qVHKKlax[\Z]$ which lifts $\choiceFAq[\Z]$, in the sense that the diagram \eqref{eq:qVHKK commutes with FAq} below commutes. Unlike in Section \ref{sec:FAq with choices summary}, here we include $\Z$ in the notation for emphasis. 

\begin{equation}
    \label{eq:qVHKK commutes with FAq} 
    \begin{tikzcd}
     & \ggBurn_\Z \ar[d, "\BtoAb"] \\
     \divcobq[\Z]  \ar[ur, "{\qVHKKlax[\Z]}"] \ar[r, "{\choiceFAq[\Z]}"']& \ggRmod[{{\Z[G]}}]
    \end{tikzcd}
\end{equation}

Just as in Remark \ref{rem:VHKK largely determined by FKh}, large parts of the following definition are determined by the functor $\choiceFAq[\Z]$, but we include all of the details for completeness.  As before, the main new ingredients center around genus 1 connected cobordisms.  The construction will utilize the lax 2-functors $\VHKKlax,\VHKKannlax$ described in Section \ref{sec:VHKKlax for D,A}.

\begin{definition}\label{def:qVHKK Z lax}
Given a coherent choice of generators $\Psi$ as in Definition \ref{def:coherent choice of generators}, the lax 2-functor $\divcobq[\Z] \rar{\qVHKKlax[\Z]} \ggBurn_\Z$ is defined as follows.  Note that all of the sets below will be of the form $Y \times \Z$ for some set $Y$; in all cases, $G=\Z$ is taken to act on the second factor.
\begin{itemize}
    \item Given an object $(Z,A)\in\divcobq[\Z]$ (a divided 1-manifold up to equivalence), we define
    \begin{align*}
    \qVHKKlax[\Z](Z,A) &:= \VHKKannlax(Z,A) \times \Z\\
    &= \Gen(Z)\times \Z.
    \end{align*}
    Diagram \eqref{eq:qVHKK commutes with FAq} is then satisfied on objects: for $(x,r) \in \qVHKKlax[\Z](Z,A)$, the generator 
    \[
    1_{(x,r)} \in \BtoAb\left(\qVHKKlax[\Z](Z,A)\right)
    \]
    corresponds to $q^r x \in \choiceFAq(Z,A)$. 
    \item Given a 1-morphism $(Z_0,A_0)\rar{q^r\Sigma}(Z_1,A_1)$ in $\divcobq[\Z]$ (a divided cobordism up to equivalence together with an integer), we define a correspondence 
    \[\begin{tikzcd}
    & \qVHKKlax[\Z](q^r\Sigma) \ar[ld, "\sq"'] \ar[rd,"\tq"] & \\
    \VHKKannlax(Z_0) \times \Z &  & \VHKKannlax(Z_1) \times \Z
    \end{tikzcd}\]
    as follows.
    \begin{itemize}
        \item If $\Sigma$ is connected with genus $g(\Sigma)\neq 1$, then we define
            \[\qVHKKlax[\Z](q^r\Sigma) := \VHKKannlax(\Sigma) \times \Z.\]
        (In particular, if $g(\Sigma)>1$, $\qVHKKlax[\Z](q^r\Sigma)=\varnothing$.)  Then for each $(a,k)\in \VHKKannlax(\Sigma) \times \Z$, we use the source map $s$ of $\VHKKannlax(\Sigma)$ to define
            \[\sq(a,k):=(s(a),k)\]
        and the target map $t$ of $\VHKKannlax(\Sigma)$, together with the integer $\ell=\ell_{t(a)}$ of Property \eqref{item:g neq 1} of Proposition \ref{prop:FAq properties}, to define
            \[\tq(a,k):=(t(a), k+r+\ell).\]
        In other words, when writing $\choiceFAq(q^r\Sigma)(q^k s(a))$ as a linear combination of the designated basis of $\choiceFAq(Z_1)$, we have  $\choiceFAq(q^r\Sigma)(q^k s(a)) = q^{k+ r + \ell} t(a) + \text{other terms}$.
        \item If $\Sigma$ is connected with genus $g(\Sigma)=1$, and $\partial\Sigma=Z_0\sqcup Z_1$ contains a homologically essential circle, then
            \[\qVHKKlax[\Z](q^r\Sigma):=\varnothing.\]
        
        \item If $\Sigma$ is connected with $g(\Sigma)=1$, $\partial\Sigma$ containing no essential circles, and the geometric winding number $w(\Sigma)\neq 0$, we let $z_1 \in \Gen(Z_0)$ and $z_X\in \Gen(Z_1)$ be the generators that label each circle by $1$ and $X$, respectively, and define
            \[\qVHKKlax[\Z](q^r\Sigma):=\{o,w\}\times \Z\]
        with
            \[\sq(o,k)=\sq(w,k):=(z_1,k)\in\Gen(Z_0)\times \Z.\]
        Then property \eqref{item:g = 1} of Proposition \ref{prop:FAq properties} provides an integer $\ell$ allowing us to define
            \[\tq(o,k):= (z_X,k+r+\ell)\in \Gen(Z_1)\times \Z\]
        and
            \[\tq(w,k):= (z_X,k+r+\ell+2w(\Sigma))\in\Gen(Z_1)\times \Z.\]

        \item If $\Sigma$ is connected with $g(\Sigma)=1$, $\partial\Sigma$ containing no essential circles, and $w(\Sigma)=0$, then we define
            \[\qVHKKlax[\Z](q^r\Sigma):=\VHKKannlax(\Sigma)\times \Z\]
        where $\VHKKannlax(\Sigma)$ is a two element set of homology generators as in \eqref{eq:two homology generators}.  Then property \eqref{item:g = 1} of Proposition \ref{prop:FAq properties} again provides an integer $\ell$ allowing us to define the map
            \[\sq(a,k):=(z_1,k)\in\Gen(Z_0)\times \Z,\]
            \[\tq(a,k):=(z_X,k+r+\ell)\in\Gen(Z_1)\times \Z\]
            for all $a\in\VHKKannlax(\Sigma)$.
        \item If $\Sigma$ consists of several connected components $\Sigma_i$, then as a set $\qVHKKlax[\Z](q^r\Sigma)$ is a product over the various components
        \[\qVHKKlax[\Z](q^r\Sigma) := \left(\prod_i \VHKKannlax(\Sigma_i)\right) \times \Z.\]
        If the various source and target maps are denoted $s_i,t_i$, then we define
        \[\sq((\cdots,a_i,\cdots),k):= ((\cdots,s(a_i),\cdots),k),\]
        and
        \[\tq((\cdots,a_i,\cdots),k):= ((\cdots, t(a_i),\cdots),k+r+\sum_i\ell_i)\]
        where the various $\ell_i$ are defined for each individual connected component as above (note that these $\ell_i$ are taken to include any winding degree shifts $2w(\Sigma_i)$ if present).
    \end{itemize}
    In all of these cases, the correspondences can be defined strictly to satisfy property \ref{item:E} of Definition \ref{def:G-Burn} using the fact that $\FAq$ is a functor to $\Z[\Z]$-modules, while properties \ref{item:D} and \ref{item:F} are immediate. Proposition \ref{prop:FAq properties} ensures that diagram \eqref{eq:qVHKK commutes with FAq} commutes on 1-morphisms.

    \item Given a 2-morphism $q^{r_1}\Sigma_1 \rar{f} q^{r_2}\Sigma_2$ in $\divcobq[\Z]$ realized by some isotopy $f$ from $\Sigma_1$ to $\Sigma_2$ with the same notation, we define the bijection $\qVHKKlax[\Z](f)$ fitting into the diagram 
    \[
    \begin{tikzcd}
    \qVHKKlax[\Z](q^{r_1}\Sigma_1) \ar[rr,"{\qVHKKlax[\Z](f)}", scale=1.2] \ar[d,"s_1"] \ar[drr,"t_1"' near end] &
    & \qVHKKlax[\Z](q^{r_2}\Sigma_2) \ar[dll,"s_2" near end] \ar[d,"t_2"] \\
    \VHKKannlax(Z_0)\times \Z & & \VHKKannlax(Z_1)\times \Z
    \end{tikzcd}
    \]
    as follows.  First we note that the isotopy $f$ must preserve the various topological cases outlined for assigning the 1-morphisms $\qVHKKlax[\Z](q^{r_i}\Sigma_i)$.  We then define $\qVHKKlax[\Z](f)$ in these cases, noting first that if $\qVHKKlax[\Z](q^{r_i}\Sigma_i)=\varnothing$ then the bijection is null.
    \begin{itemize}
        \item In connected cases where $\qVHKKlax[\Z](q^{r_i}\Sigma_i)= \VHKKannlax(\Sigma_i)\times \Z$ (i.e. all cases with genus not equal to one, as well as the genus one case with contractible boundary and winding number zero), we define
        \[\qVHKKlax[\Z](f):=\VHKKannlax(f)\times \Id_{\Z}.\]
        \item In the connected genus 1 case where $\qVHKKlax[\Z](q^{r_i}\Sigma_i)=\{o_i,w_i\}\times \Z$, we define
        \[\qVHKKlax[\Z](f)(o_1,k):=(o_2,k),\quad \qVHKKlax[\Z](f)(w_1,k):=(w_2,k)\]
        for all $k\in\Z$.
        \item For disconnected $\Sigma_i$, $f$ must maintain each connected component and we define $\qVHKKlax[\Z](f)$ to be the product over each component $\Sigma_{i,j}$ of the corresponding bijection $f_{i,j}$.
    \end{itemize}
    In all of these cases, the commutation with the source maps is clear.  Verifying commutation with the target maps utilizes the uniqueness of the various integers $\ell$ used to define them together with the fact that $q^{r_1}\FAq[\Z](\Sigma_1) = q^{r_2}\FAq[\Z](\Sigma_2)$ as maps of free $\Z[\Z]$-modules.
    \item Non-trivial composition 2-isomorphisms are built in precisely the same way as for $\VHKKannlax$ (which themselves were built in the same way as for $\VHKKlax$) in all cases other than the case where two genus zero cobordisms $\Sigma_1,\Sigma_2$ compose to give a genus one cobordism $\Sigma=\Sigma_1\circ\Sigma_2$ with non-zero winding number, so that $\qVHKKlax[\Z](\Sigma)=\{o,w\}\times \Z$.  In this case however,  the composition 2-isomorphism is uniquely determined by the source and target maps, since if $x,y\in \qVHKKlax[\Z](\Sigma)$ satisfy $s_q(x)=s_q(y)$, then the integer component of $t_q(x)$ differs\footnote{This 2-isomorphism requiring $2w(\Sigma)\neq 0$ is the reason why we begin with $G=\Z$ rather than an arbitrary cyclic group.} from that of $t_q(y)$ by $2w(\Sigma)\neq 0$.
\end{itemize}
\end{definition}

In summary, the lax 2-functor $\divcobq[\Z]\rar{\qVHKKlax[\Z]} \ggBurn_\Z$ can be thought of as $\VHKKannlax \times \Z$ in almost all situations.  The only new situations to consider are those involving 1-morphisms (cobordisms) $\Sigma$ which are genus 1 and \emph{not} contractible in $[0,1]\times \A$.  If such a $\Sigma$ has an essential circle on its boundary, $\qVHKKlax[\Z](\Sigma)=\varnothing$ and no further constructions are necessary.  If such a $\Sigma$ has only trivial circles on its boundary, then $\qVHKKlax[\Z](\Sigma)$ involves a two element set $\{o,w\}$ where the elements are distinguished by how they correspond to the summands in Equation \eqref{eq:FAq on genus 1 with wrapping}.  This distinguishing of elements then uniquely prescribes how 2-morphisms (bijections) involving this case must behave.  Then the various relationships between the $\Z$ factors are entirely determined by the functor $\choiceFAq[\Z]$.

\subsubsection{The case where $G$ is any cyclic group}
\label{sec:qVHKK lax}

Let $G=\langle q \rangle$ be an arbitrary cyclic group, so $G=\Z/H$ for some subgroup $H\subset \Z$. Proposition \ref{prop:qdivcob G to G/H} and its proof indicate how all of the data in $\divcobq$ is determined by the data in $\divcobq[\Z]$.  Thus we may use the data of $\qVHKKlax[\Z]$ to define the corresponding lax 2-functor $\qVHKKlax$, essentially by replacing every occurrence of $\Z$ in Definition \ref{def:qVHKK Z lax} by $G$.

\begin{proposition}\label{prop:qVHKK lax}
Given a coherent choice of generators $\Psi$ as in Definition \ref{def:coherent choice of generators}, the data of Definition \ref{def:qVHKK Z lax} determines a lax 2-functor $\divcobq \rar{\qVHKKlax} \ggBurn_G$ as follows.

\begin{itemize}
    \item If $(Z,A)\in\ob(\divcobq)=\ob(\divcobq[\Z])$, we define
    \[\qVHKKlax(Z,A):=\Gen(Z) \times G.\]
    \item  If $[q^r]\Sigma$ is a 1-morphism in $\divcobq$, we consider any lift $q^a\Sigma$ in $\divcobq[\Z]$ so that
    \[\qVHKKlax[\Z](q^a\Sigma) = A \times \Z\]
    for some set $A$, allowing us to define
    \[\qVHKKlax([q^r]\Sigma):= A \times G,\]
    with source and target maps defined on the $G$ factor accordingly.  It is clear that this definition is independent of lift.
    \item If $[q^r]\Sigma \rar{f} [q^{r'}]\Sigma'$ is a 2-morphism in $\divcobq$, then there must be some lifts with corresponding 2-morphism $q^a\Sigma \rar{\widetilde{f}} q^{a'}\Sigma'$ in $\divcobq[\Z]$, for which the bijection of the form
    \[A\times \Z \rar{\qVHKKlax[\Z](\widetilde{f})} A'\times \Z\]
    determines a bijection defining
    \[A \times G \rar{\qVHKKlax(f)} A'\times G.\]
\end{itemize}

The resulting lax 2-functor lifts $\choiceFAq$, as in the following diagram:

\begin{equation}
    \label{eq:qVHKK G commutes with FAq G} 
    \begin{tikzcd}
     & \ggBurn_G \ar[d, "\BtoAb"] \\
    \divcobq  \ar[ru, "{\qVHKKlax}"] \ar[r, "{\choiceFAq}"']& \ggRmod[{{\Z[G]}}]
    \end{tikzcd}.
\end{equation}

\end{proposition}

\begin{proof}
All of these statements are straightforward consequences of the definitions and Proposition \ref{prop:qdivcob G to G/H}. In particular, the composition 2-isomorphism involving genus one $\Sigma$ with $w(\Sigma)\neq 0$ distinguishes the two elements via their target values differing by $2w(\Sigma)\neq 0$ in $\Z$, despite the fact that $q^{2w(\Sigma)}$ may be the identity in $G$.
\end{proof}

\section{Spectral link invariants}
\label{sec:spectral link invts}
In this section we use the constructions of Section \ref{sec:lax lifts of Khov TQFTs} to define various Khovanov spectra for links which lift the corresponding Khovanov complexes of Sections \ref{sec:The Khovanov and annular Khovanov functors} and \ref{sec:quantum annular homology defs}.  Sections \ref{sec:Groupoid enrichments} and \ref{sec:XKh for links in D or A} mainly review material from \cite{LLS_Khovanov_Spectra}.  In Section \ref{sec:qaXKh} which is of most immediate interest and contains new material, we will be lifting the quantum annular complex $\qaKCchoice(L)$ rather than $\qaKC(L)$.

\subsection{Groupoid enrichments, Elmendorf-Mandell $K$-theory, and rectification}
\label{sec:Groupoid enrichments}
In this subsection we establish several key categorical constructions that will be used heavily in what follows.

Given any multicategory $\cat{C}$, the \emph{canonical thickening} $\cat{C}'$ is a groupoid-enriched multicategory with the same objects.  Informally, the multimorphism groupoids of $\cat{C}'$ consist of trees of composable multimorphisms in $\cat{C}$.  Then a multimorphism in $\cat{C}$ can be viewed as a multimorphism in $\cat{C}'$ with one internal vertex.  Such multimorphisms in $\cat{C}'$ will be called \emph{basic multimorphisms}.  There is a unique $2$-morphism in $\cat{C}'$ between any two trees whose multimorphisms compose to be equal in $\cat{C}$. For the precise definition, see \cite[Section 2.4.1]{LLS_Khovanov_Spectra}.  
More generally, if $\cat{C}$ was already enriched in groupoids (for instance, $\divcob(X)$), then there is a unique $2$-morphism in $\cat{C}'$ between two trees if there is a $2$-morphism in $\cat{C}$ between their respective compositions.\footnote{While the canonical thickening of an already groupoid-enriched multicategory is not explicitly defined in \cite{LLS_Khovanov_Spectra}, this is implicitly the thickening procedure used to define $\widetilde{\underline{\mathrm{Cob}_d}}$ in \cite[Section 3.2.3]{LLS_Khovanov_Spectra}.} There is a \emph{canonical projection}  $\pi: \cat{C}'\rightarrow\cat{C}$ which is the identity on objects and composes the multimorphisms along any tree.

If $\cat{D}$ is a multicategory enriched in groupoids, then there is a \emph{strictification} $\cat{D}^0$, which is a non-enriched multicategory where Ob($\cat{D}^0) = $ Ob$(\cat{D})$, and the multimorphism sets for $\cat{D}^0$ are given by the set of path components in the groupoid $\Hom_{\cat{D}}(x_1, \ldots, x_m;y)$. Viewing $\cat{D}^0$ as trivially enriched in groupoids, there is a projection $2$-functor $\pi : \cat{D} \to \cat{D}^0$.

More generally, let $\cat{M}$ be a simplicial multicategory; that is, a multicategory enriched in  simplicial sets. In this case we write $\Map_{\cat{M}}(x_1, \ldots, x_m;y)$ (rather than $\Hom$) for the corresponding simplicial set of multimorphisms.
For instance, if $\cat{D}$ is a multicategory enriched in groupoids (or more generally, in categories), then taking the nerve of each Hom-category gives a simplicial multicategory $\nerve(\cat{D})$. 
There is a \emph{strictificiation} $\cat{M}^0$, which is a (non-enriched) multicategory, given by replacing each Hom simplicial set by its set of path components. There is also a projection functor $\pi : \cat{M} \to \cat{M}^0$, viewing $\cat{M}^0$ as a constant simplicial multicategory. 

Any category can be viewed as a multicategory which has no multimorphisms with more than one input. The discussion and constructions in this subsection so far apply equally well if all instances of ``multicategory'' are replaced with ``category''. In particular, we may consider the canonical thickening of a category (potentially enriched in groupoids) which is a category enriched in groupoids. In this case, a tree of composable multimorphisms is just a sequence of composable morphisms in the original category. We will elaborate on the relevant canonical thickenings such as $(\cube^N)'$ and $\divcob(X)'$ as they are required.

We now discuss \emph{rectification}, which is a procedure that turns a simplicially enriched functor into a strict (non-enriched) functor out of the strictification.  This is largely due to Elmendorf-Mandell \cite{Elmendorf-Mandell}, but we follow the exposition in \cite[Section 2.9]{LLS_Khovanov_Spectra}. Rectification uses the language of (multi)categories enriched in simplicial sets and (multi)functors between them; we refer the reader to \cite{Kelly} for necessary background on enriched categories. We also discuss some of the details in Appendix \ref{sec:simplicial stuff}. 

The category $\Spectra$ is enriched in (pointed) simplicial sets. We may then view $\Spectra$ as a simplicial multicategory by using the smash product. The following is a combination of \cite[Definition 2.42 and Lemma 2.43]{LLS_Khovanov_Spectra}. See \cite[Definition 12.1]{Elmendorf-Mandell} for the definition of a weak equivalence. 

\begin{lemma}\label{lem:rectification}
    Let $\cat{M}$ be simplicial multicategory. Given an enriched multifunctor $J: \cat{M} \to \Spectra$, there exists a multifunctor $\rect_J : \cat{M}^0 \to \Spectra$ such that $\rect_J \circ \pi$ is naturally weakly equivalent\footnote{This means there is a zig-zag of natural transformations where each one is an object-wise weak equivalence.} to $J$ if $\pi:\cat{M} \to \cat{M}^0$ is a weak equivalence. 
\end{lemma}

Similarly, $\Spectra_G$ is simplicially enriched as we describe in Appendix \ref{sec:simplicial stuff}. We have the following analogue of Lemma \ref{lem:rectification}. 

\begin{lemma}
\label{lem:G rectification}
    Let $\cat{C}$ be a simplicial category. Given an enriched functor $J: \cat{C} \to \Spectra_G$, there exists a functor $\rect_J : \cat{C}^0 \to \Spectra_G$ such that $\rect_J \circ \pi$ is naturally weakly equivalent to $J$ if $\pi:\cat{C} \to \cat{C}^0$ is a weak equivalence.
\end{lemma}

 We prove Lemma \ref{lem:G rectification} in Appendix \ref{sec:simplicial stuff}. Note that $\cat{C}$ in the statement is a category rather than a multicategory. We do not have (nor do we need) an analogous statement if $\cat{C}$ is a multicategory. The functor $\rect_J$ (in both cases) is called the \emph{rectification} of $J$. 

We now discuss Elmendorf-Mandell $K$-theory. Let $\Permu$ denote the simplicial\footnote{Strictly speaking, the multicategory $\Permu$ in \cite{Elmendorf-Mandell} is enriched in categories. We apply the nerve construction but do not include it in the notation, as explained in the introduction in \cite{Elmendorf-Mandell}.} multicategory of all small permutative multicategories \cite[Definition 3.1]{Elmendorf-Mandell}.  We do not need the precise definition of this category.  Instead we focus on the existence of two (simplicially enriched) multifunctors into and out of $\Permu$,
\[\nerve(\Burn^\bullet) \rar{\func{O}}\Permu\rar{\KEM}\Spectra^\bullet,\]
where the bullet $\bullet$ stands in for either graded or bigraded as throughout Section \ref{sec:lax lifts of Khov TQFTs}.  The multifunctor $\func{O}$ sends an object $X$ of $\Burn^\bullet$ to the overcategory $\Sets/X\in\Permu$ of finite sets over $X$, and sends a correspondence $X\lar{s}A\rar{t}Y$ to the functor $(-)\times_X A$ between such overcategories (again there are some details to track in order to ensure strict identities and associativity; see \cite[Section 2.8]{LLS_Khovanov_Spectra}).  The functor $\KEM$ is the $K$-theory functor of Elmendorf-Mandell \cite{Elmendorf-Mandell}.  The gradings are not difficult to track - see \cite[Section 6]{LLS_Khovanov_Spectra} for a formal procedure.

The functor $\KEM$ and the composition $K:=\KEM\circ\func{O}$ satisfy a number of important properties, some of which we list below. 

\begin{enumerate}[label= (K\arabic*)]
\item \label{item:K(CxD) = K(C) x K(D)} For two permutative categories $\cat{C}$ and $\cat{D}$, there is an isomorphism $\KEM(\cat{C} \times \cat{D}) \to \KEM(\cat{C}) \times \KEM(\cat{D})$ and a weak equivalence $\KEM(\cat{C}) \vee \KEM(\cat{D}) \to \KEM(\cat{C})\times \KEM(\cat{D})$ which are natural in both factors. 
    \item\label{item:K(X) is wedge of sphere} For a set $X\in\Burn^\bullet$, there is a weak equivalence $\bigvee_{x\in X} \S \rar{} K(X)$ (here $X$ is an at most countable set). 
    \item\label{item:K(A) is the right map} Under the weak equivalence above, given a $1$-morphism $X\lar{s}A\rar{t}Y$, the induced map on spectra $K(X)\rar{K(A)} K(Y)$ sends a sphere $\S_x$ associated to the element $x\in X$ to the sphere $\S_y$ by a map of degree $\vert A_{y,x} \vert$.
\end{enumerate}

We will lift $K$ to an enriched functor\footnote{Note, not multifunctor.}  $K_G : \nerve(\B_G) \to \Spectra_G$.  Precisely,  view $\B$ as a category and $K : \nerve(\B) \to \Spectra$ as a simplicially enriched functor (rather than as a multicategory and as a multifunctor). Consider the diagram \eqref{eq:lifting K}. 
\begin{equation}
\label{eq:lifting K} 
\begin{tikzcd}
    \nerve(\B) \ar[r, "K"] & \Spectra \\
    \nerve(\B_G) \ar[u] \ar[ur] \ar[r,dashed] & \Spectra_G \ar[u]
\end{tikzcd}
\end{equation}
The vertical maps are forgetful functors (forgetting the $G$-action), and the diagonal map is defined to be the composition of the left vertical map and the top horizontal map. The following says that the dashed arrow exists.

\begin{proposition}
\label{prop:K_G}
There is an enriched functor $K_G: \nerve(\B_G) \to \Spectra_G$ making the diagram \eqref{eq:lifting K} commute. 
\end{proposition}

A proof of Proposition \ref{prop:K_G} is given in Appendix \ref{sec:simplicial stuff}. We note that as a consequence of the proof, if $X$ is an object of $\B_G$ then under the equivalence in \ref{item:K(A) is the right map}, the action of $g\in G$ on $K_G(X)$ permutes spheres, in the sense that it sends the sphere $\S_x$ corresponding to $x\in X$ to $\S_{gx}$ corresponding to $gx\in X$ via the identity. In other words, the  $G$-action on $K(X)\in\Spectra_G$ (and thus the $\Sq$-module structure) is determined by the $G$-action on the sets in $\B_G$.

\begin{remark}
Proposition \ref{prop:K_G} says that the outputs of $K$ can be taken to land in $\Spectra_G$ if the inputs are in $\B_G$. In particular, $K_G$ does not change the values of $K$.  Compare with \cite[Conjecture C]{DGML}, which would require a significant modification of $K$. 
\end{remark}

\subsection{The Khovanov spectra for links in $\D$ or $\A$}
\label{sec:XKh for links in D or A}
Let $X\in \{ \D, \A\}$, and let $L\subset X$ denote a link diagram with $N$ crossings.  In this section we summarize parts of the construction in \cite{LLS_Khovanov_Spectra} to define the Khovanov spectrum for $L$.

In Section \ref{sec:quantum annular homological background} we defined the Khovanov chain complex $\KC[X](L)$ by applying the relevant TQFT $\FKh[X]$ to the Bar-Natan complex $[[L]]$.  It is convenient to have the following alternative description, following \cite[Section 2.5]{LLS_Khovanov_Spectra}. 

\begin{definition}
    \label{def:IMC}
    Let $\cat{C}$ denote any pointed model category.  The \emph{iterated mapping cone} is a functor
    \[\Func(\cube^N,\cat{C}) \rar{\IMC_{\cube^N}} \cat{C}\]
    defined as follows.  Let $\cube^N_+$ denote the category obtained by adding a single object $*$ to $\cube^N$ and a unique morphism $u\to *$ whenever $u\neq \vec{1}:= (1,1,\ldots, 1)$. Note that $\cube^N$ is a full subcategory of $\cube_+^N$.  There is a pointed extension functor
    \[\Func(\cube^N,\cat{C})\rar{(-)_+} \Func(\cube_+^N,\cat{C})\]
    which extends $F\in\Func(\cube^N,\cat{C})$ to a functor $F^+ \in\Func(\cube^N_+,\cat{C})$ by $F^+ (x) = F(x)$ if $x\in \cube^N$ and $F^+ (x)=0_{\cat{C}}$ otherwise (and extends natural transformations correspondingly).  Then $\IMC_{\cube^N}$ is defined to be the composition
    \[\IMC_{\cube^N}(-) := \hocolim_{\cube_+^N} (-)_+.\]

    There is a related enlargement of the cube $\cube^N_\dagger := (\cube^1_+)^N$ which also contains $\cube^N$ as a full subcategory. Given a functor $F:\cube^N \to \cat{C}$, we can extend it to $F^\dagger : \cube^N \to \cat{C}$ by defining $F^\dagger$ to send any object not in $\cube^N$ to $0_{\cat{C}}$.
\end{definition}

 Both of these enlarged cubes appear in \cite{LLS-Burnside}, while \cite{LLS_Khovanov_Spectra} uses $\cube^N_\dagger$ and \cite{LLS_CK_Spectra} uses $\cube^N_+$. Homotopy colimits of functors out of these two enlargements (or analogous enlargements of related categories) are often equivalent; see \cite[Lemma 4.41, Lemma 5.27]{LLS-Burnside}. 

\begin{lemma}
    \label{lem:hocolim over plus vs dagger}
    Given any functor $F : \cube^N \to \Spectra$, there is a weak equivalence $\hocolim_{\cube^N_\dagger} F^\dagger \to \hocolim_{\cube^N_+} F^+$ which is functorial in $F$. The same holds if the target category of $F$ is $\Top_*$ or $\Top_*^G$, the categories of based topological spaces and based topological $G$-spaces, respectively, or $\Spectra_G$.
\end{lemma}

\begin{proof}
Consider the functor $J: \cube^N_\dagger \to \cube^N_+$ which is the identity on $\cube^N$ and sends any object not in $\cube^N$ to $* \in \cube^N_+$. We have $ F^+ \circ J = F^\dagger$.  Arguing as in \cite[Lemma 4.41]{LLS-Burnside}, one sees that $J$ is homotopy cofinal, and the statement follows from property (ho-4) in \cite[Section 4.2]{LLS-Burnside}.
\end{proof}

The above and \cite[Corollary 2.15]{LLS_Khovanov_Spectra}\footnote{Note, the category $P$ therein is our $\cube^1_+$.} imply that $\KC[X](L)$ is quasi-isomorphic to the iterated mapping cone of the  composition \[\cube^N \rar{} \BN(X) \rar{\FKh[X]} \Ab^\bullet,\]
where the first arrow indicates the \emph{cube of resolutions of $L$} (see Definition \ref{def:cube cat} and the ensuing discussion) and the bullet $\bullet$ indicates the proper grading as usual.

We now describe the relevant canonical thickenings.  The thickened cube $(\cube^N)'$ has the same objects as those of $\cube^N$ (vertices), and its morphism groupoids\footnote{Note that both $\cube^N$ and $(\cube^N)'$ have no non-empty multimorphism groupoids with more than one input; in this sense they are both ordinary categories being viewed trivially as a multicategories} $\Hom(v,w)$ have objects consisting of sequences of composable basic morphisms
\[v\rar{}u\rar{}\cdots\rar{} w.\]
Between any two such morphism objects within $\Hom(v,w)$ there is exactly one 2-morphism, indicating that they compose to give the same overall morphism in $\cube^N$.  Composition of morphism groupoids is based upon concatenations of sequences of basic morphisms.

Similarly we have the canonical thickenings $\divcob(X)'$ of the divided cobordism categories of Definition \ref{def:divcob} with the same objects (divided 1-manifolds), but morphisms consist of trees of composable morphisms (divided cobordisms) in $\divcob(X)$.  In this setting, there is a unique 2-morphism between two such trees if and only if the corresponding compositions (concatenations of divided cobordisms) are isotopic.

Now let $L\subset X$ denote a link diagram with $N$ crossings.  Vertices $v\in\cube^N$ then determine resolutions $L_v\subset X$ as in Section \ref{sec:The Khovanov and annular Khovanov functors}.  Following \cite{LLS_Khovanov_Spectra} we choose an arbitrary set of marked points called \emph{pox} on $L$ away from the crossings such that, for each $v\in\cube^N$, every component of $L_v$ contains at least one pox.   We now use the notation $L$ to refer to the data of the link diagram together with an ordering of crossings and choice of pox.  This abuse of notation is justified by the fact that differing auxiliary choices of crossing orderings and pox give rise to equivalent spectra at the end of the construction \cite[Proposition 3.26]{LLS_Khovanov_Spectra}.\footnote{The argument given there is for arbitrary tangles and generalizes immediately to links in either $\D$ or $\A$.} 

\begin{definition}\label{def:cobfunc for links}
The \emph{enriched cube of resolutions} for a link diagram $L\subset X$ is the strict 2-functor
\[(\cube^N)' \rar{\cobfunc_L} \divcob(X)'\]
defined as follows.
\begin{itemize}
    \item To an object $v\in (\cube^N)'$ we assign the divided 1-manifold
    \[ \cobfunc_L(v) := (L_v \,,\, A)\]
    where $A$ is the union of small (disjoint) neighborhoods of each pox and each 0-resolution in $L_v$.
    \item To a basic morphism $t=(v\rar{} w)$
    in $(\cube^N)'$ we assign a divided cobordism
    \[\cobfunc_L(t):=(\Sigma,\Gamma)\]
    consisting of single divided saddles from any 0-resolution neighborhood in $L_v$ to the corresponding 1-resolution neighborhood in $L_w$ (see \cite[Figure 3.2]{LLS_Khovanov_Spectra}).  Outside of these neighborhoods, we use the identity divided cobordism (including in the neighborhood of the pox).  More general morphisms (sequences of basic morphisms) are assigned sequences of divided cobordisms accordingly.
    \item Suppose $s, s'$ are $1$-morphisms in $(\cube^N)'$ from $v$ to $w$, and let $f$ be the unique $2$-morphism $s\to s'$. Then $\cobfunc_L(s)$ and $\cobfunc_L(s')$ are sequences of divided saddles arranged at different relative heights whose compositions are isotopic rel boundary and hence are related by a unique $2$-morphism in $\divcob(X)'$, so we set $\cobfunc_L(f)$ to be this $2$-morphism.\footnote{Note that in $\BN(X)$ isotopic cobordisms are equal, so that no enrichments and 2-morphisms were necessary in Section \ref{sec:The Khovanov and annular Khovanov functors}}
\end{itemize}
\end{definition}

\begin{definition}\label{def:VHKK}
For $X\in\{\D,\A\}$, the strict multifunctor $\divcob(X)'\rar{\VHKK[X]'}\Burn^\bullet$ is the functor induced by the lax 2-functor $\divcob(X)\rar{\VHKKlax[X]}\Burn^\bullet$ as guaranteed by \cite[Lemma 3.16]{LLS_Khovanov_Spectra}.  In particular, $\VHKK[X]'=\VHKKlax[X]$ on objects and basic multimorphisms in $\divcob(X)'$, ensuring that the following diagram commutes.
\[
\begin{tikzcd}
\divcob(X)' \ar[r,"{\VHKK[X]'}"] \ar[d,"\pi"] & \Burn^\bullet \ar[d,"\BtoAb"] \\
\divcob(X) \ar[r,"{\FKh[X]}"] & \Ab^\bullet
\end{tikzcd}
\]
\end{definition}

We remark here that the proof of \cite[Lemma 3.16]{LLS_Khovanov_Spectra} is for the case $X=\D$ only, but the argument is completely formal, using only the lax 2-functor structure and the respecting of the monoidal structure.  The case for $X=\A$ follows analogously, except that there is no monoidal structure to track. Thus $\VHKK[\A]'$ is a strict $2$-functor which we will sometimes refer to as a multifunctor by viewing $\divcob(\A)'$ as a multicategory with no multiple-input morphisms.

\begin{definition}\label{def:Khovanov spectra}
Let $X\in\{\D,\A\}$ and let $L\subset X$ be a poxed link diagram with $N$ (ordered) crossings.  Let $\cube^N\rar{\mathfrak{R}_L}\Spectra^\bullet$ 
denote the rectification of the composition
\[(\cube^N)'\rar{\cobfunc_L}\divcob(X)' \rar{\VHKK[X]'} \Burn^\bullet \rar{K} \Spectra^\bullet.\]
The \emph{Khovanov spectrum} $\XKh[X](L)$ is then defined to be
\[\XKh[X](L):= \sh^{-N_+} \IMC_{\cube^N}(\mathfrak{R}_L).\]
Here, $N_+$ is the number of positive crossings in $L$,  $\sh^{-1}$ is the shift functor from \cite[Definition 2.32]{LLS_Khovanov_Spectra}, and $\sh^{-n} : = (\sh^{-1})^n$ for $n \in \Z_{\geq 0}$.
\end{definition}

\begin{proposition}\label{prop:XKh recovers Kh for D and A}
The Khovanov spectrum of a link diagram $L\subset X$ recovers the Khovanov chain complex up to quasi-isomorphism:
\[\chainsfunc(\XKh[X](L))\simeq \KC[X](L).\]
\end{proposition}
\begin{proof}
For $X=\D$, this is precisely the $n=m=0$ case of \cite[Proposition 4.2]{LLS_Khovanov_Spectra}; see also Proposition \ref{prop:LLS lift platalg and CK bimod} later in this document. The proof when $X=\A$ is completely analogous.
\end{proof}

\subsection{The quantum annular Khovanov spectrum}
\label{sec:qaXKh}
Let $L\subset\Aq$ be a link diagram with $N$ crossings, let $\choice$ be a coherent choice of generators as in Definition \ref{def:coherent choice of generators}, and let $G=\langle q \rangle$ denote an arbitrary cyclic group.  Then the quantum annular complex $\qaKCchoice(L)$ can be viewed (up to quasi-isomorphism) as arising from the iterated mapping cone of the composition
\[\cube^N\rar{} \BNq \rar{\choiceFAq} \ZGmod.\]

The category $\divcobq'$ has the same objects (divided 1-manifolds) as  $\divcobq$, but $\Hom( Z_0, Z_m )$ consists of sequences of $G$-labelled basic morphisms (divided cobordisms)
\[ Z_0 \rar{q^{r_1}\Sigma_1} Z_1 \rar{q^{r_2}\Sigma_2} \cdots \rar{q^{r_m}\Sigma_m} Z_m,\]
where we have suppressed the arcs and divides from the notation to avoid clutter.  In this setting, there is a unique 2-morphism between two such sequences if and only if there is a 2-morphism between their corresponding compositions in $\divcobq$ - that is, if the concatenations are isotopic in a manner that treats the group labels properly.

 We then have the following analogs of Definitions \ref{def:cobfunc for links} and \ref{def:VHKK} in the quantum annular setting. For a link diagram $L\subset \Aq$, we will always assume that pox are disjoint from $\mq$.

\begin{definition}\label{def:cobfunc for links in Aq}
The \emph{cube of resolutions} for a link diagram $L\subset\Aq$ is the strict 2-functor
\[(\cube^N)' \rar{\cobfunc_L} \divcobq'\]
defined in the same manner as Definition \ref{def:cobfunc for links}.
\end{definition}

\begin{definition}\label{def:qVHKK}
The strict $2$-functor $\divcobq'\rar{\qVHKKchoice'}\ggBurn_G$ is the functor induced by the lax 2-functor $\divcobq\rar{\qVHKKlax} \ggBurn_G$ of Proposition \ref{prop:qVHKK lax}.
\end{definition}
As before, the existence of $\qVHKKchoice'$ is guaranteed by \cite[Lemma 3.16]{LLS_Khovanov_Spectra} whose formal proof immediately generalizes to the quantum annular setting.

\begin{definition}\label{def:qAKh spectrum}
Let $L\subset \Aq$ be a link diagram with $N$ crossings and let $G$ be an arbitrary cyclic group.  Rectifying the composition
\[(\cube^N)'\rar{\cobfunc_L}\divcobq' \rar{\qVHKKchoice'} \ggBurn_G \rar{K_G} \ggSpectra_G\] 
gives rise to a functor 
\[
\cube^N \rar{\rect_L} \ggSpectra_G.
\]
The \emph{quantum annular Khovanov spectrum} $\qaXchoice(L)\in \ggSpectra_G$ is defined to be 
\[
\qaXchoice(L) = \sh^{-N_+} \IMC_{\cube^N} (\mathfrak{R}_L).
\]
We write $\qaXchoice(L; j)$ for the summand of $\qaXchoice(L)$ in annular degree $j$.
\end{definition}

\begin{proposition}
\label{prop:invariance}
    Up to weak equivalence in $\ggSpectra_G$, the spectrum $\qaXchoice(L)$ is an invariant of the isotopy class of the annular link represented by $L$ that does not depend upon the choice of generators $\choice$.
\end{proposition}

\begin{proof}
    For a fixed annular link diagram $L$, that two choices of (coherent) generators $\choice$, $\choice'$ yield equivalent spectra follows from the fact that the two functors $\divcobq' \rar{{V}'_{\Aq,\Psi}} \ggBurn_G $ and  $\divcobq' \rar{{V}'_{\Aq,\Psi'}} \ggBurn_G $ are naturally isomorphic, as in the proof of \cite[Proposition 4.3]{AKW}. Alternatively, this follows from  Theorem \ref{thm:qTHH recovers qAX}  below, since the left-hand side of the equivalence therein does not depend on a choice of generators. 
    
    Next, we address isotopies. If $L, L'\subset \Aq$ are diagrams representing isotopic annular links, then $L$ and $L'$ may be related by a sequence of (1) Reidemeister moves disjoint from $\mq$ and (2) cyclic rotation of the diagrams. Invariance under (1) follows from standard arguments, for instance as in \cite[Section 6]{LS-Kh-spectra} or \cite[Section 2]{LLS_Khovanov_Spectra}. For (2), suppose $L$ and $L'$ are related by a rotation. We may assume that $L'$ is obtained from $L$ by moving either (i) an arc or (ii) a crossing across $\mq$, shown in \eqref{eq:rotating diagram}.
    \begin{equation}
    \label{eq:rotating diagram}
    \begin{aligned}
        \begin{tikzpicture}[scale=.5]
\draw (0,2) arc  (180:360:1 and 1); 

\draw[densely dashed] (-1,0) -- (3,0);
\node at (-1,0) {$\times$};
\draw[>=stealth,->] (2.8,0) -- (2.8,-.5);

\draw[->] (3.75,0) -- (5.25,0)  node[midway, above] {(i)};

\begin{scope}[shift={(7,0)}]
    \draw (0,0) arc  (180:360:1 and 1); 
\draw (0,0) -- (0,2);
\draw (2,0) -- (2,2); 
\node at (0,-2) {};
\draw[densely dashed] (-1,0) -- (3,0);
\node at (-1,0) {$\times$};
\draw[>=stealth,->] (2.8,0) -- (2.8,-.5);
\end{scope}

\begin{scope}[shift={(16,0)}]
\draw[ ] (0,2) .. controls (0,1.5) and (2,.5).. (2,0);

\draw[preaction={draw, line width=12pt, white}] (0,0) .. controls (0,.5) and (2,1.5) .. (2,2); 

\draw (0,0) -- (0,-2);
\draw (2,0) -- (2,-2);

\draw[densely dashed] (-1,0) -- (3,0);
\node at (-1,0) {$\times$};
\draw[>=stealth,->] (2.8,0) -- (2.8,-.5);

\draw[->] (3.75,0) -- (5.25,0)  node[midway, above] {(ii)};

\begin{scope}[shift={(7,0)}]
    \draw[ ]  (0,0) .. controls (0,-.5) and (2,-1.5) .. (2,-2);

\draw[preaction={draw, line width=12pt, white}] (0,-2) .. controls (0,-1.5) and (2,-.5).. (2,0); 

\draw (0,0) -- (0,2);
\draw (2,0) -- (2,2);

\draw[densely dashed] (-1,0) -- (3,0);
\node at (-1,0) {$\times$};
\draw[>=stealth,->] (2.8,0) -- (2.8,-.5);
\end{scope}
\end{scope}
\end{tikzpicture}
\end{aligned}
    \end{equation}
Note, for each  of these two cases, there are four possibilities, and we only depict one possibility in each case (for instance, in case (ii), the move may be viewed in reverse, or the crossing data can be switched). 

The argument is similar to the proof of \cite[Proposition 4.4]{AKW} but adapted to the present language. Let $N$ be the number of crossings in $L$ (and in $L'$). 
In each case, we construct a functor $(\cube^N \times \cube^1)' \rar{\cobfunc_{L,L'}} \divcobq'$ whose restriction to $(\cube^N\times \{0\})'$ and $(\cube^N\times \{1\})'$ is equal to 
$(\cube^N)'\rar{\cobfunc_L}\divcobq'$ and $(\cube^N)'\rar{\cobfunc_{L'}}\divcobq'$, respectively.

For a fixed choice of pox on $L$, we pick pox on $L'$ so that they are obtained by rotation, and we assume the ordering of crossings in $L'$ is obtained from the ordering of crossings of $L$ in the natural way. In case (ii), we also assume that the depicted crossing is first in the orderings. Let $\delta : 0\to 1$ be the unique non-identity morphism in $\cube^1$, and let $\varphi_{u,v} \in \cube^N$ be the unique morphism from $u$ to $v$. The functor $\cobfunc_{L,L'}$ is already determined on objects and on basic morphisms of the form $(\varphi_{u,v}, \id_i)$ for $i = 0,1$. 
For each $u\in \{0,1\}^N$, there is a divided cobordism $(\Sigma_u, \Gamma_u) : L_u \to L_u'$ traced out by the rotation. Moreover, these are invertible up to $2$-isomorphism, in the sense that there is some $k_u\in \Z$ such that\footnote{In fact, each $k_u \in \{-1,0,1\}$, since in each case the compositions are either equal to the identity or differ from the identity by moving a single cup or saddle through $\Mq$.} there are $2$-morphisms $\left( q^{k_u} (\b{\Sigma_u}, \b{\Gamma_u})  \right) \circ (\Sigma_u, \Gamma_u) \Rightarrow \id_{L_u}$ and $(\Sigma_u, \Gamma_u) \circ \left( q^{k_u} (\b{\Sigma_u}, \b{\Gamma_u})  \right) \Rightarrow \id_{L_u'}$ in $\divcobq$, where the bar denotes reflection through the thickening direction. 

In case (i), define 
\[
\cobfunc_{L,L'}(\varphi_{u,v}, \delta):= (\Sigma_v, \Gamma_v) \circ \cobfunc_{L}(\varphi_{u,v}).
\]
In case (ii), suppose that the situation is exactly as depicted in \eqref{eq:rotating diagram} (the other three cases are similar). Set $n_{u,v} = 0$ if $u_1=v_1$ and $n_{u,v} = -1$ if $u_1 = 0$ and $v_1 = 1$. Define 
\[
\cobfunc_{L,L'}(\varphi_{u,v}, \delta):= q^{n_{u,v}} (\Sigma_v, \Gamma_v) \circ \cobfunc_{L}(\varphi_{u,v}).
\]
Now  $\cobfunc_{L,L'}$ is defined on all basic morphisms, and we extend it to sequences of composable morphisms in the natural way. 

It remains to show that $\cobfunc_{L,L'}$ extends to $2$-morphisms.
Given a $2$-morphism $\alpha : (f_1, \ldots, f_n) \Rightarrow (g_1, \ldots, g_m)$ in $(\cube^N\times \cube^1)'$, if none of the $f_i$ or $g_j$ involve $\delta$, then $\alpha$ corresponds to rearranging heights of saddles, and thus there is a corresponding $2$-morphism in $\divcobq$. Otherwise, precisely one $f_i$ and one $g_j$ involves $\delta$ (that is, $f_i = (\varphi_{u,v},\delta)$). In case (i), saddles commute with the rotation cobordism up to isotopy, so there is a corresponding $2$-morphism in $\divcobq$. In case (ii), the existence of the required $2$-morphism is guaranteed by Figure \ref{fig:trace moves} (in particular the saddle in the top right) and the factor $q^{n_{u,v}}$. 

The composition $\qVHKKchoice' \circ \cobfunc_{L,L'}$ is a quasi-isomorphism, in the sense of \cite[Definition 3.25]{LLS_Khovanov_Spectra}, from $\qVHKKchoice' \circ \cobfunc_L$ to $\qVHKKchoice' \circ \cobfunc_{L'}$. Then the special case of \cite[Proposition 4.7]{LLS_Khovanov_Spectra} where $n=m=0$ finishes the argument. 
\end{proof}

\begin{proposition}
\label{prop:quantum annular V recovers the right thing}
The quantum annular Khovanov spectrum of a link diagram $L\subset \Aq$ recovers the quantum annular Khovanov chain complex up to quasi-isomorphism in $\Kom(\ggRmod[{{\Z[G]}}])$:
\[\chainsfunc(\qaXchoice(L)) \simeq \qaKCchoice(L).\]
\end{proposition}
\begin{proof}
Unwinding various definitions one proves that the diagram
\[
\begin{tikzcd}
(\cube^N)' \ar[r,"\cobfunc_L"] \ar[d] &
\divcobq' \ar[r,"\qVHKKchoice'"] \ar[d] &
\ggBurn_G \ar[d,"\BtoAb"] \\
\cube^N \ar[dr,bend right]&
\divcobq \ar[r,"\choiceFAq",swap]
&
\ggRmod[{{\Z[G]}}] \\
 & \BNq \ar[ru,"\choiceFAq",swap,bend right] & 
\end{tikzcd}
\]
commutes up to natural isomorphisms.  From that point the argument is analogous to the proof of \cite[Proposition 4.2]{LLS_Khovanov_Spectra}. First, since $C_*(\rect_L(v))$ is a wedge of copies of $\sphere$, \cite[Lemma 4.3]{LLS_Khovanov_Spectra}
provides a zig-zag of quasi-isomorphisms between $\chainsfunc \circ \rect_L$ and $H_0 \circ \rect_L$, where $H_0$ denotes the composition of $\chainsfunc$ followed by taking $0$-th homology. Note, a priori the maps involved in this quasi-isomorphism are in $\Kom(\Ab)$ rather than $\Kom(\Rmod[{{\Z[G]}}])$; however, examining the proof of \cite[Lemma 4.3]{LLS_Khovanov_Spectra} shows that they are induced by inclusion and projection maps, which are $\Z[G]$-linear. Next, we see that $H_0 ( \rect_L (v)) \cong \choiceFAq(L_v)$ as $\Z[G]$-modules. The map $H_0(\rect_L(v)) \to H_0(\rect_L(w))$ assigned to an edge $v\to w$ recovers the cobordism map $\choiceFAq(L_v) \to \choiceFAq(L_w)$ by construction. Finally, $\chainsfunc$ commutes with homotopy colimits by \cite[Proposition 2.31]{LLS_Khovanov_Spectra}, and applying $\IMC$ recovers the process of collapsing the cube of resolutions to a chain complex by \cite[Corollary 2.15]{LLS_Khovanov_Spectra}. 
\end{proof}

\section{Spectral tangle invariants, qTHH, and proof of Theorem \ref{thm: qthh recovers bpw spectrum}}\label{sec:spectral tangle invariants}
Throughout this section we fix $0\leq k \leq n$.  In Section \ref{sec:spectral platalg and tangle bimod} we review the constructions from \cite{LLS_CK_Spectra} of the spectral platform algebras $\platalgk$ and tangle bimodules over them using the language of Section \ref{sec:linear and spectral cats via multifunctors}.  In Section \ref{sec:aug quant HM functor} we then compute the quantum topological Hochschild homologies of their scalar extensions using the language of Section \ref{sec:qTHH via shape cats} together with Lemma \ref{lem:scalar ext properties}.

\subsection{The spectral platform algebra and tangle bimodules}
\label{sec:spectral platalg and tangle bimod}
The spectral platform algebra $\platalgk$ is defined as a chain lift of $\Platnk$, both viewed as multifunctors out of $\shape(\matchingsk)$.  Just as in Section \ref{sec:spectral link invts}, the construction utilizes groupoid enrichments and a composition of the form
\[\shape(\matchingsk)' \rar{\cobfunc_{\mathcal{A}}} \divcob(\D)' \rar{\VHKK'} \gBurn.\]
Similarly the spectral Chen-Khovanov bimodule $\X(T)$ for an $(n,n)$-tangle diagram\footnote{The discussion here can be generalized to the case of $(n,m)$-tangles, but we will only need to consider $(n,n)$-tangles moving forward.}  $T$
 with $N$ crossings utilizes a composition of the form
\[(\cube^N\tiltimes \shape(\matchingsk,\matchingsk))' \rar{\cobfunc_T} \divcob(\D)' \rar{\VHKK'} \gBurn,\]
where the symbol $\tiltimes$ indicates a special type of product to be described shortly.  Unlike in Section \ref{sec:spectral link invts} however, these compositions will be adjusted in a manner analogous to quotienting by the ideal $\Ideal{n}{k}$ in Section \ref{sec:Platform algebras and bimodules} before using rectification to complete the construction.  We provide only a summary here following \cite{LLS_CK_Spectra}; more details are provided in \cite{LLS_Khovanov_Spectra}.

\begin{figure}
\centering
\begin{tikzpicture}[xscale=.8, yscale=.7]

\draw[densely dashed] (0,3.25) -- (0,-.25); 

\begin{scope}
\clip (1.6,0) -- (.4,0) -- (.4,3) -- (1.6,3) -- cycle;

    \draw[very thick] (0,3) arc (90:-90:1.5 and 1.5); 
    
    \draw[very thick] (0,2) arc (90:-90:.75 and .5);     
\end{scope}

\begin{scope}
\clip(-1,0) -- (-.4,0) -- (-.4,3) -- (-1,3) -- cycle;
 \draw[very thick] (0,1) arc (90:270:.75 and .5);
  
    \draw[very thick] (0,3) arc (90:270:.75 and .5);
      
\end{scope}

\begin{scope}
\clip (-.4,-.25) -- (.4,-.25) -- (.4,3.25) -- (-.4,3.25) -- cycle;

  \draw[thick,densely dotted] (0,1) arc (90:270:.75 and .5);
  
    \draw[thick,densely dotted] (0,3) arc (90:270:.75 and .5);

    
    \draw[thick,densely dotted] (0,3) arc (90:-90:1.5 and 1.5); 
    
    \draw[thick,densely dotted] (0,2) arc (90:-90:.75 and .5);   
\end{scope}
\end{tikzpicture}
\caption{Active (\textbf{bold}) and inactive (dotted) arcs of $a_1\flip{a_2}$.}\label{fig:matching active arcs}
\end{figure}
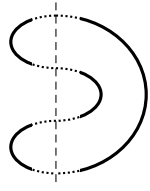

We begin by describing $\shape(\matchingsk)'$ and $\cobfunc_{\mathcal{A}}$.  Objects of $\shape(\matchingsk)'$ are pairs $(a_1,a_2)$ of platform matchings; $\cobfunc_{\mathcal{A}}(a_1,a_2)$ is defined to be the divided 1-manifold $(a_1\flip{a_2} \,,\, a_1\flip{a_2}\setminus B_\eps(\partial a_1))$.  In other words, we glue $a_1$ to the reflection of $a_2$, and declare a neighborhood of the gluing points to be \emph{inactive}.  See Figure \ref{fig:matching active arcs}.

Basic multimorphisms $t$ in $\shape(\matchingsk)'$ are `single vertex trees' with inputs 
\[
(a_1,a_2),(a_2,a_3),\dots, (a_{m-1},a_m)
\]
and output $(a_1,a_m)$; $\cobfunc_{\mathcal{A}}(t)$ is then the divided cobordism $(\Sigma,\Gamma)$  where
\[ a_1\flip{a_2} \sqcup a_2\flip{a_3} \sqcup \cdots \sqcup a_{m-1}\flip{a_m} \xrightarrow{\Sigma} a_1\flip a_m\]
consists of \emph{divided multisaddles} $\flip{a_i}a_i \xrightarrow{s(a_i)} \id_{2n}$ with divides $\Gamma$ as indicated in \cite[Figure 3.2]{LLS_Khovanov_Spectra}. The case $m=0$ is allowed, in which case there is a basic morphism $\varnothing \to (a,a)$ that is assigned the divided cup cobordism $\varnothing \to a \flip{a}$ as shown in \cite[Figure 3.2]{LLS_Khovanov_Spectra}. More general multimorphisms (trees) in $\shape(\matchingsk)'$ are assigned trees of divided cobordisms accordingly.  Then the unique 2-morphism between any two trees in $\shape(\matchingsk)'$ with the same domain and codomain is assigned the unique 2-morphism between the corresponding trees in $\divcob(\D)'$ indicating that various compositions of such multisaddles, arranged in any order of relative heights, are all isotopic rel boundary (\cite[Corollary 3.14, Proposition 3.15]{LLS_Khovanov_Spectra}).

The functor $\cobfunc_T$ is similar.  As in Section  \ref{sec:explicitly constructing the quasi-iso}, we let $\strands{T}$ denote the tangle obtained by adding $n-k$ and $k$ horizontal strands below and above $T$, respectively.  We then add pox to $\strands{T}$ so that each component of each resolution $\strands{T}_v$ contains at least one pox, as in Section \ref{sec:spectral link invts}.   The multicategory $(\cube^N\tiltimes \shape(\matchingsk,\matchingsk))'$ has three types of objects.  There are pairs $(a_1,a_2)$ and $(b_1,b_2)$ of crossingless matchings for tracking $\platalgk$ and its ability to act on either the left or the right; we set $\cobfunc_T(a_1,a_2):= \cobfunc_{\mathcal{A}}(a_1,a_2)$, and similarly for pairs $(b_1,b_2)$.  There are also objects of the form $(v,a \mid b)$ for $a,b\in\matchingsk$ and $v\in\cube^N$; we set $\cobfunc_T(v,a \mid b):=(a\strands{T}_v\flip{b},A)$ where $a\strands{T}_v\flip{b}$ is as in Figure \ref{fig:capping off platform tangle} with the set of active arcs $A$ consisting of the three following types of arcs:
    \begin{itemize}
        \item $a\setminus B_\eps(\partial a)$ and $\flip{b}\setminus B_\eps(\partial \flip{b})$, i.e. the matchings $a$ and $\flip{b}$ excluding a neighborhood of the gluing points, similar to the case of $\cobfunc_{\mathcal{A}}$;
        \item small (disjoint) neighborhoods of each pox in $\strands{T}_v$; and
        \item small (disjoint) neighborhoods of each 0-resolution in $\strands{T}_v$.
    \end{itemize}

Basic multimorphisms $t$ in $(\cube^N\tiltimes \shape(\matchingsk,\matchingsk))'$ are single vertex trees of one of two types. For the first type, there is a unique basic morphism with inputs 
\[
(a_1,a_2),(a_2,a_3),\ldots,(a_{m-1},a_m)
\]
and output $(a_1,a_m)$ (or similarly for $b$'s instead of $a$'s); for these types, $\cobfunc_T(t)=\cobfunc_{\mathcal{A}}(t)$ (including the case $m=0$).  For the second type, there is a unique basic morphism with inputs of the form
\[(a_1,a_2),(a_2,a_3),\dots, (a_{\ell-1},a_\ell), (v, a_\ell \mid b_1), (b_1,b_2),(b_2,b_3),\dots (b_{m-1},b_m)\]
and output $(w, a_1 \mid b_m)$ with $v\leq w$ in $\cube^N$.  For these types, $\cobfunc_T(t)$ consists of divided multisaddles between any two `pairs' (as in $\cobfunc_{\mathcal{A}}$) as well as single divided saddles from any 0-resolution neighborhood in $\strands{T}_v$ to the corresponding 1-resolution neighborhood in $\strands{T}_w$ (as in $\cobfunc_L$ from Definition \ref{def:cobfunc for links}).  More general multimorphisms (trees) and 2-morphisms are treated just as for $\cobfunc_{\mathcal{A}}$ and $\cobfunc_L$.

\begin{remark}\label{rmk:tilde product enrichment}
In \cite[Section 3.2.4]{LLS_Khovanov_Spectra}, a different groupoid enrichment of these types of products $\cube^N\tiltimes\shape$ (with fewer basic multimorphisms) is used. We instead follow \cite{LLS_CK_Spectra} which uses the canonical enrichment as described here; see also \cite[Remarks 3.4 and 3.21]{DGML}.
\end{remark}

Overall, the multifunctors $\cobfunc_{\mathcal{A}}$ and $\cobfunc_T$ can be summarized as follows.  Crossingless matchings $a_i,b_j$ and tangle resolutions $\strands{T}_v$ (if they are present) are glued together along their boundaries to form divided 1-manifolds such as $a_1\flip{a_2}$ or $a \strands{T}_v \flip{b}$.  Morphisms in the shape category are assigned (trees of concatenable) divided cobordisms based upon multisaddles $\flip{a_i} a_i \rightarrow \id_{2n}$ and single saddles from 0-resolutions to 1-resolutions.  Active arcs and divides have been positioned so that neighborhoods of the pox on $\strands{T}$ and any portion of a 1-manifold which has the potential to be involved in a saddle cobordism are considered active. Everything else is considered inactive.

Verifying that these are indeed  enriched multifunctors is treated in detail in \cite[Section 3.4]{LLS_CK_Spectra}.  Roughly speaking, the multimorphism structure is handled by the monoidal $\sqcup$ structure in $\divcob(\D)'$.  The 2-morphism structure utilizes the fact that any arrangement of a fixed set of divided (multi)saddles and divided cups at various heights is isotopic to any other.  The uniqueness of such an isotopy uses the divides on the various cobordisms (\cite[Lemma 3.1]{LLS_Khovanov_Spectra}); the inclusion of the pox makes sure that every connected component of any such cobordism is indeed divided. 

At this stage the compositions $\VHKK'\circ\cobfunc_{D}$ for $D\in\{\mathcal{A},T\}$ are modified into new functors $\func{G}_{D}$ by taking quotients with respect to so-called \emph{absorbing sub-functors} which parallel the process of quotienting by $\Ideal{n}{k}$ and $\mathcal{I}^k(T)$ in Section \ref{sec:Platform algebras and bimodules}; see \cite[Section 4]{LLS_CK_Spectra}.  The resulting functors $\func{G}_D$ can be described as follows.  If an object $O$ in the domain gives rise to any Type III circles, we define $\func{G}_D(O)$ to be the empty set in $\gBurn$.  If $O$ gives rise to Type II circles and $\VHKK'\circ\cobfunc_D(O) = \Gen(Z)$, we define $\func{G}_D(O)$ to be $(\Gen(Z)\setminus \Gen_{II}(Z))$, where $\Gen_{II}(Z)$ consists of all the Khovanov generators which label some type II circle by $X$.  Elements in the correspondences between such sets are discarded if and only if they map to or from such elements (via source and target maps).  The proof that this process gives rise to a well-defined multifunctor is discussed in \cite[Lemma 4.2]{LLS_CK_Spectra}, and follows from the original argument of \cite{Chen-Khovanov} showing that $\Ideal{n}{k}$ is an ideal (together with the fact that multiplication and comultiplication in the standard Khovanov basis are given by matrices with non-negative entries).  For the full details, see \cite[Section 4]{LLS_CK_Spectra}. 

We review the following constructions from \cite[Section 4.1]{LLS_Khovanov_Spectra}. Recall the rectification procedure discussed in Section \ref{sec:Groupoid enrichments}. 

\begin{definition}\label{def:spec plat alg and CK tangle spectrum}
Fix $n\geq 0$.  The \emph{spectral platform algebra} $\platalg=\bigvee_{0\leq k \leq n} \platalgk$ consists of summands $\platalgk$ defined as follows.  Rectifying the composition 
\[\shape(\matchingsk)' \rar{\func{G}_\mathcal{A}} \gBurn \rar{K} \gSpectra\]
gives a functor 
\begin{equation}
\label{eq:algebra multifunctor}
\shape(\matchingsk) \rar{\Xfunck[\mathcal{A}]} \gSpectra
\end{equation}
which defines the algebra $\platalgk$ as a spectral category (Definition \ref{def:spectral cat}) via the translation of Proposition \ref{prop:Alg as multifunctor}.

Similarly, the \emph{spectral Chen-Khovanov bimodule} $\X(T)=\bigvee_{0\leq k \leq n} \Xk(T)$ of an $(n,n)$-tangle diagram $T$ with $N$ crossings, $N_+$ of which are positive, consists of summands $\Xk(T)$  defined as follows.  Rectifying the composition
\[(\cube^N\tiltimes\shape(\matchingsk,\matchingsk))' \rar{\func{G}_T} \gBurn \rar{K} \gSpectra\]
gives rise to a functor $\cube^N\tiltimes\shape(\matchingsk,\matchingsk) \rar{\mathfrak{R}_T} \gSpectra$.

The multicategory $\cube^N\tiltimes\shape(\matchingsk,\matchingsk)$ has two copies of $\shape(\matchingsk)$ as full sub-multicategories (one for the left action and one for the right action). As discussed in \cite[Section 4.1]{LLS_Khovanov_Spectra}, $\mathfrak{R}_T$ restricted to each of these sub-multicategories agrees with $\Xfunck[\mathcal{A}]$ from \eqref{eq:algebra multifunctor}.  Define a multifunctor 
\[
\Xfunck(T) : \shape(\matchingsk,\matchingsk) \to \gSpectra
\]
as follows. As above, $\shape(\matchingsk,\matchingsk)$ contains two copies of $\shape(\matchingsk)$ as full sub-multicategories. Define $\Xfunck(T)$ to agree with $\mathfrak{R}_T$ (and hence with $\Xfunck[\mathcal{A}]$) on each of these copies. Next, for $a,b\in \matchingsk$, consider the functor 
\[
i_{a,b}: \cube^N \to \cube^N\tiltimes\shape(\matchingsk,\matchingsk)
\]
sending $v\in \cube^N$ to $(v,a \mid b)$ and sending a morphism $v\to w$ to the multimorphism (with one input) $(v, a \mid b) \to (w,a \mid b)$. Then  
\[
\Xfunck(T) (a \mid b) := \sh^{-N_+} \IMC_{\cube^N}( \mathfrak{R}_T \circ i_{a,b}).
\]
The value of $\Xfunck(T)$ on a multimorphism of the form 
\[
(a_1, a_2) ,\ldots (a_{m-1}, a_m), (a_m \mid b_1), (b_1, b_2), \ldots, (b_{\ell-1},b_\ell) \to (a_1 \mid b_\ell) 
\]
uses $\mathfrak{R}_T$ and naturality of homotopy colimits and of $\sh^{-1}$; see \cite[Section 4.1]{LLS_Khovanov_Spectra} for details. Finally, $\Xfunck(T)$ defines $\Xk(T)$ as a spectral bimodule (Definition \ref{def:bimodule over spectral cats}) via the translation of Proposition \ref{prop:Alg as multifunctor}.
\end{definition}

\begin{proposition}[{{\cite[Theorem 3]{LLS_Khovanov_Spectra}}}]
\label{prop:LLS lift platalg and CK bimod}
For any $0\leq k \leq n$, the spectral platform algebra recovers the platform algebra and the spectral Chen-Khovanov bimodule for a tangle $T$ recovers the Chen-Khovanov complex up to quasi-isomorphism:
\[ \chainsfunc(\platalgk)\simeq \Platnk \quad\text{ and }\quad \chainsfunc(\Xk(T))\simeq \CCKnk(T). \]
\end{proposition}

Recall from \eqref{eq:scalar ext for Sq} and \eqref{eq:quantum extensions of CK} the notion of extension of scalars on the spectral and algebraic sides, respectively.  
The following  says that extending scalars commutes with taking chains for the spectra of interest. 
\begin{proposition} There are quasi-isomorphisms
    \[
    \chainsfunc(\qplatalgk) \simeq \qPlatk{n}{k} \quad \text{ and }\quad \chainsfunc(\qXk(T))\simeq \qCCKnk(T).
    \]
\end{proposition}

\begin{proof}
    The argument is similar to the proof of \cite[Proposition 4.2]{LLS_Khovanov_Spectra}. We establish the second equivalence, and the first equivalence is similar.
    To simplify the notation, set 
    \[
\cat{A} = \platalgk, \ \ \cat{X} = \Xk(T),\ \  \func{X} = \Xfunck(T), \ \ \rect_{a,b} = \rect_T \circ i_{a,b},
    \]
   for $a,b\in \matchingsk$. We also will omit the shift $\sh^{-N_+}$ in what follows in light of \cite[Proposition 2.35]{LLS_Khovanov_Spectra}. We have
    \[
 {}_q\cat{X}
 =
 \left(
 \bigvee_{a,b\in \matchingsk} \func{X} (a\mid b)
 \right) \wedge \Sq
\cong
\bigvee_{a,b\in \matchingsk} \func{X} (a\mid b) \wedge \Sq
=
\bigvee_{a,b\in \matchingsk} \left(\IMC_{\cube^N}\rect_{a,b}\right) \wedge \Sq,
    \]
so
\[
\chainsfunc({}_q\cat{X})    \simeq \bigoplus_{a,b\in \matchingsk} \chainsfunc \left( (\IMC_{\cube^N} \rect_{a,b}) \wedge \Sq
\right).
\]
Since $\wedge$ respects homotopy colimits in each variable (\cite[Proposition 2.21]{LLS_Khovanov_Spectra}), we have a quasi-isomorphism
\[
(\IMC_{\cube^N}\rect_{a,b}) \wedge \Sq
\rar{\simeq}
\IMC_{\cube^N} 
\left(
\rect_{a,b} \wedge \Sq
\right)
\]
where $\rect_{a,b} \wedge \Sq$ is the composition $\cube^N \rar{\rect_{a,b}} \Spectra \rar{(-)\wedge \Sq} \Spectra_G$.
Since $\chainsfunc$ respects homotopy colimits (\cite[Proposition 2.31]{LLS_Khovanov_Spectra}), we  have a quasi-isomorphism 
\[
\chainsfunc
\left[
\IMC_{\cube^N} 
(\rect_{a,b} \wedge \Sq)
\right]
\rar{\simeq}
\IMC_{\cube^N} \chainsfunc
\left(
\rect_{a,b} \wedge \Sq
\right).
\]
By \cite[Proposition 4.3]{LLS_Khovanov_Spectra}, there is a zig-zag of quasi-isomorphisms between $\chainsfunc
\left(
\rect_{a,b} \wedge \Sq
\right)$ and $H_0
\left(
\rect_{a,b} \wedge \Sq
\right)$. 
Note that $H_0 \left(\rect_{a,b} \wedge \Sq \right)$ is isomorphic to the cube of resolutions used to build $\qCCKnk(T; a,b)$ since  $\rect_{a,b}(v)$ is a wedge of (the correct amount of) spheres for each $v\in \cube^N$, and edges in $\cube^N$ are, by construction, assigned maps that lift the Chen-Khovanov edge maps. 

Altogether we have a zig-zag of quasi-isomorphisms from $\chainsfunc({}_q \cat{X})$ to 
\[
\bigoplus_{a,b\in \matchingsk} \IMC_{\cube^N} H_0 \left(\rect_{a,b} \wedge \Sq \right)
\]
which is quasi-isomorphic to $\bigoplus_{a,b\in \matchingsk} \qCCKnk(T; a,b)$ as complexes of abelian groups by \cite[Corollary 2.15]{LLS_Khovanov_Spectra}. Let us verify that the left action is correct; the right action is analogous. The map
\[
\lambda : \cat{A}(c,a) \wedge \rect_{a,b}(v) \to \rect_{c,b}(v)
\]
furnished by $\rect_T$ for $a,b,c\in \matchingsk$ and $v\in \cube^N$, can be identified, upon taking chains, with the left action map at the algebraic level up to quasi-isomorphism. Then
\[
(\cat{A}(c,a) \wedge \Sq) \wedge_{\Sq} (\rect_{a,b}(v) \wedge \Sq) \cong 
(\cat{A}(c,a) \wedge \rect_{a,b}(v)) \wedge \Sq 
\rar{\lambda \wedge \id}
\rect_{c,b}(v) \wedge \Sq
\]
similarly recovers, upon taking chains, the left action map of the scalar-extended platform bimodule. 
 
\end{proof}

\subsection{The augmented quantum Hochschild-Mitchell functor to lift $\qXi^\choice$}
\label{sec:aug quant HM functor}
In this section we build an augmented quantum Hochschild-Mitchell functor
\[\HMshape^\omega(\shape(\matchingsk)) \rar{\qHMfunc{T}} \gSpectra_G,\]
in the sense of Definition \ref{def:qHM functor}, for the extended Chen-Khovanov bimodule $\qXk(T)$ over the extended platform algebra $\qplatalgk$. 

The construction of $\qHMfunc{T}$ will be based on a composition of functors through $\divcobq'$ and $\gBurn_G$. The spectrum $\qHMfunc{T}(\omega)$ will recover the quantum annular spectrum $\qaXchoice(\closure{T};n-2k)$, and the induced map $F$ of spectra with $G$-action from Corollary \ref{cor:aug qHM gives map on qTHH},
\[\qTHH(\qplatalgk;\qXk(T)) \rar{F}  \qaXchoice(\closure{T};n-2k),\]
will recover the quasi-isomorphism $\qXi^\choice$ of Proposition \ref{prop:qXichoice is a quasi-iso} upon applying the chains functor $\chainsfunc$, which together with Theorem \ref{them:Whitehead thm} proves Theorem \ref{thm: qthh recovers bpw spectrum}.  Verifying that $\qHMfunc{T}$ satisfies the required properties will also involve heavy use of Lemma \ref{lem:scalar ext properties}.  Since this is our main new result, we will go into more detail than in the previous section which largely summarized past work.  Note that, as indicated in \cite[Section 5]{LLS_CK_Spectra} via \cite[Lemma 4.5]{LLS_Khovanov_Spectra}, $\platalgk$ is pointwise cofibrant as a spectral category; by \cite[Ch. IV, Theorem 1.4]{schwedesymmetric} the scalar extension $\qplatalgk$ is also pointwise cofibrant.

\subsubsection{The quantum cobordism functor $\cobfunc_T^{H^\omega}$}\label{sec:func HMshape to qdivcob}
Having fixed a cyclic group $G=\langle q \rangle$ and $(n,n)$-tangle $T$ with $N$ crossings, we build a functor
\[(\cube^N \times \HMshape^\omega(\matchingsk))' \rar{\cobfunc_T^{H^\omega}} \divcobq'\]
similar to the functor $\cobfunc_T$ of Section \ref{sec:spectral platalg and tangle bimod}.

Recall that we had set $S^1  = [0,2\pi]/\partial$, $\Aq= S^1\times (-1,1)$ equipped with the seam $\mq = \{0\} \times (-1,1)$, and $\D = (0,1) \times (-1,1)$. We will use the embedding $\iota : \D \hookrightarrow \A$ given by $\iota(x,y) = (2\pi x, y)$; note that $\iota$ embeds $\D$ as the complement of $\mq$.

\begin{figure}
\centering 
\subcaptionbox{\label{fig:annular embedding for non-terminal objects}}[1 \linewidth]
{\begin{tikzpicture}
\draw[densely dashed] (.2,0) -- (9.8,0);
\draw[densely dashed] (.2,3) -- (9.8,3);
\node[below] at (.2,0){$0$};
\node[below] at (9.8,0){$2\pi$};

\draw[ultra thick,densely dashed,blue] (.2,0) to node[pos=.85,left]{$\mu_q$} (.2,3);
\draw[ultra thick,densely dashed,blue] (9.8,0) to node[pos=.85,right]{$\mu_q$} (9.8,3);

\begin{scope}[shift={(.5,0)}]
\draw (1,.75) rectangle (2,2.25);

\draw (1,2.25) arc (90:270:1 and .75);
\draw (2,.75) arc (-90:90:1 and .75);

\node at (1.5,1.5) {$\strands{T}_v$};
\node at (.5,1.5) {$a_m$};
\node at (2.5,1.5) {$\flip{a_0}$};
\end{scope}

 \draw (5,1.5) ellipse (1 and .75);
  \draw (5,2.25) -- (5,.75);
  \node at (4.5,1.5) {$a_{0}$};
  \node at (5.5,1.5) {$\flip{a_1}$};

  \begin{scope}[shift = {(4.5,0)}]
        \draw (4,1.5) ellipse (1 and .75);
  \draw (4,2.25) -- (4,.75);
  \node at (3.5,1.5) {$a_{m-1}$};
  \node at (4.5,1.5) {$\flip{a_m}$};
  \end{scope}

  \node at (6.75,1.5) {$\cdots$};
\end{tikzpicture}
}\\
\vskip2em
\subcaptionbox{\label{fig:annular embedding for terminal objects}}[.35\linewidth]
{\begin{tikzpicture}[xscale=.7]
\draw[densely dashed] (0,0) -- (5,0);
\draw[densely dashed] (0,3) -- (5,3);
\node[below] at (0,0){$0$};
\node[below] at (5,0){$2\pi$};

\draw[ultra thick,densely dashed,blue] (0,0) to node[pos=.85,left]{$\mu_q$} (0,3);
\draw[ultra thick,densely dashed,blue] (5,0) to node[pos=.85,right]{$\mu_q$} (5,3);

\draw (1.5,1) rectangle (3.5,2);

\node at (2.5,1.5) {${T}_v$};

\draw[thick] (0,2.5)--(4,2.5) to node[below,scale=.7,pos=.25]{$k$} (5,2.5);

\draw[thick] (3.5,1.5) to node[below,scale=.7,pos=.5]{$n$} (5,1.5);
\draw[thick] (0,.5)--(4,.5) to node[below,scale=.7,pos=.25]{$n-k$} (5,.5);
\draw[thick] (0,1.5) -- (1.5,1.5);
\end{tikzpicture}}
\subcaptionbox{\label{fig:HM functor last matching}}[.6\linewidth]
{\begin{tikzpicture}
\draw[densely dashed] (0,0) -- (3,0);
\draw[densely dashed] (0,3) -- (3,3);
\node[below] at (1,0){$0 = 2 \pi$};

\draw[ultra thick,densely dashed,blue] (1,0) to node[pos=.88,left]{$\mu_q$} (1,3);

\draw (1,2.25) arc (90:270:1 and .75);

\draw (1,2.25) --++ (2,0);
\draw (1,.75) --++ (2,0);

\node at (.5,1.5) {$a_j$};
\node at (2.5,1.5) {$\strands{T}_v$};

\draw (2,2.25) -- (2,.75);

\begin{scope}[shift= {(5,0)}]
    \draw[densely dashed] (0,0) -- (3,0);
\draw[densely dashed] (0,3) -- (3,3);
\node[below] at (1,0){$0 = 2 \pi$};

\draw[ultra thick,densely dashed,blue] (1,0) to node[pos=.88,left]{$\mu_q$} (1,3);

\draw (2,2.25) arc (90:270:.8 and .75);

\node at (1.7,1.5) {$a_j$};

\node at (2.5,1.5) {$\strands{T}_v$};

\draw (2,2.25) -- (2,.75);

\draw (2,2.25) -- (3,2.25);
\draw (2,.75) -- (3,.75);
\end{scope}

\draw[->] (3.5, 1.5) -- (4.5,1.5);
\end{tikzpicture}}
\caption{(A): evaluating $\cobfunc_T^{H^\omega}(v,(a_0,a_1,\cdots,a_m))\in\divcobq'$.  (B): evaluating $\cobfunc_T^{H^\omega}(v,\omega)\in\divcobq'$.  The divided 1-manifold $\closure{\strands{T}_v}$ is embedded into  $\A$ as shown.  The $k$ and $n-k$ indicate the extra strands above and below the resolution $T_v$. (C):  the cobordism pulling the last matching $a_j$ through the seam, used to define $\cobfunc_T^{H^\omega}$ on a basic morphism $(u,(a_0,\dots,a_m))\rar{}(v,(a_i,\dots,a_j))$ when $j\neq m, (a_i,\dots,a_j) \neq \omega$.  \label{fig:HM functor} }
\end{figure}
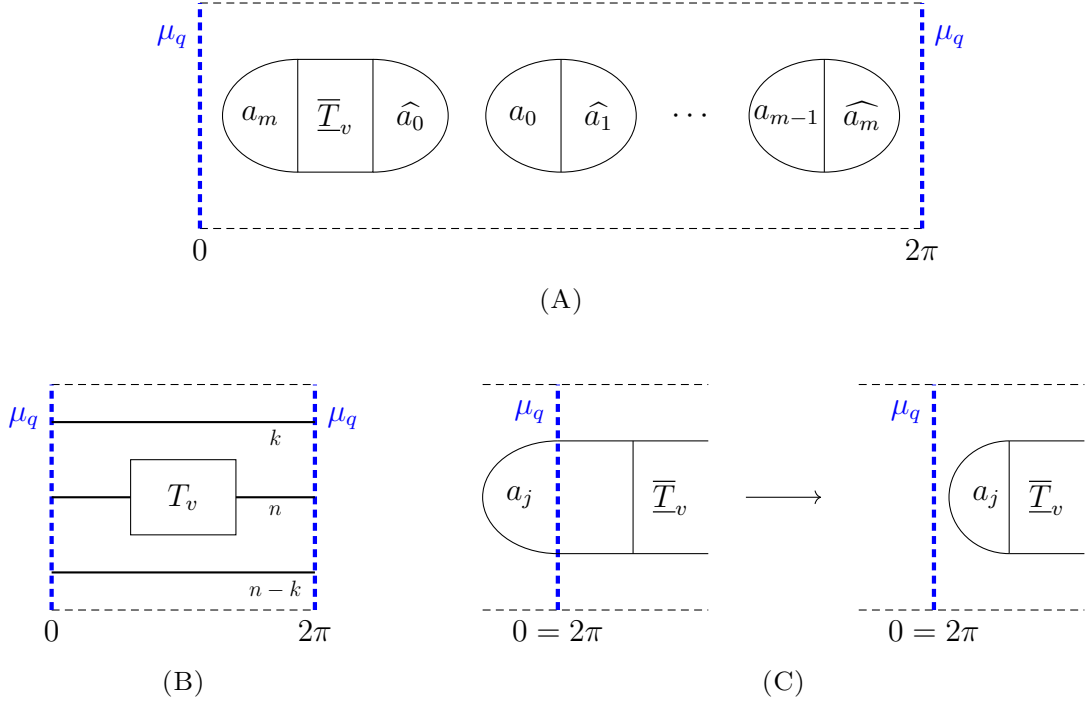

\begin{definition}\label{def:func HMshape to qdivcob}
Consider a fixed $(n,n)$-tangle diagram $T$ with $N$ (ordered) crossings and $0\leq k \leq n$. We let $\strands{T}$ denote the diagram obtained by adding $n-k$ and $k$ horizontal strands below and above $T$, respectively. Given a choice of pox on $\strands{T}$, the functor 
\[
(\cube^N \times \HMshape^\omega(\matchingsk))' \xrightarrow{\cobfunc_T^{H^\omega}} \divcobq'
\]
is defined as follows.

\begin{itemize}
    \item Objects of $(\cube^N \times \HMshape^\omega(\matchingsk))'$ are either pairs $(v,(a_0,\dots,a_m))$ of a vertex $v\in\cube^N$ and a tuple of platform matchings $a_i\in\matchingsk$, or they are so-called \emph{closure} objects $(v,\omega)$ with $\omega$ functioning as the `empty tuple'.     For a non-closure object $(v,(a_0,\cdots,a_m))$, first consider the divided 1-manifold 
    \[ (a_m \strands{T}_v \flip{a_0} \sqcup a_0\flip{a_1} \sqcup \cdots \sqcup a_{m-1}\flip{a_m}, A),\]
    in $\D$ with active arcs $A$ determined on each disjoint term as described in Section \ref{sec:spectral platalg and tangle bimod}, and then embed it into $\A$ via $\iota: \D\hookrightarrow \A$.  See Figure \ref{fig:annular embedding for non-terminal objects}. To a closure object $(v,\omega)$ we assign the divided 1-manifold $(\closure{\strands{T}_v},A)$.  The active arcs $A$ are small disjoint neighborhoods of each 0-resolution and each pox.
    See Figure \ref{fig:annular embedding for terminal objects}.
    \item A basic morphism in $(\cube^N\times \HMshape^\omega(\matchingsk))'$  is represented by
    \[(u,(a_0,\dots,a_m)) \to (v,(a_i,\dots,a_j)),\]
    corresponding to a morphism $u\rightarrow v$ in $\cube^N$ and a deletion of some number of platform matchings from $(a_0,\dots,a_m)$ to arrive at $(a_i,\dots,a_j)$ (again with $\omega$ acting as the empty tuple).  We assign it the triple $(\Sigma,\Gamma,q^r)$ as follows.  The divided cobordism $(\Sigma,\Gamma)$ consists of:
    \begin{itemize}
        \item Divided saddles from 0-resolutions in $\strands{T}_u$ to 1-resolutions in $\strands{T}_v$, as determined by $u\rightarrow v$.
        \item A divided multisaddle $\flip{a_\ell} a_\ell \rightarrow \id_{2n}$ for each deleted matching $a_\ell$ in the morphism $(a_1,\dots,a_m)\rightarrow(a_i,\dots,a_j)$. If $j=m$ then all of these multisaddles (and $\Sigma$) are disjoint from $\Mq$. If $j\neq m$ (this includes the case $(a_i,\dots,a_j)= \omega$) then the the multisaddle  $\flip{a_m} a_m \to \id_{2n}$ intersects $\Mq$ in the minimal way; if $j\neq m$ and further $(a_i,\dots,a_j)\neq \omega$, we also compose these multisaddles with  the cobordism formed by a planar isotopy which pulls the last remaining matching $a_j$ through $\mq$. See Figure \ref{fig:HM functor last matching}.
        \item Identity divided cobordisms away from these saddles.
    \end{itemize}
    If $j\neq m$ and $(a_i,\dots,a_j)\neq\omega$, we set $r:=-n$; this shift down by $n$ is  analogous to the shift in the Chen-Khovanov algebra, as in \eqref{eq:subalgebra of H^n}.   Otherwise, we set $r:=0$. 

    \item The morphism groupoid $\Hom( (u,(a_0,\dots,a_m)) , (v,(a_i,\dots,a_j)) )$ has objects consisting of sequences of composable basic morphisms
    \[(u,(a_0,\dots,a_m))\rar{}\cdots\rar{} (v,(a_i,\dots,a_j)).\]
    Such a sequence is sent to the corresponding sequence of divided cobordisms in $\divcobq'$.  Between any two such sequences in $\Hom( (u,(a_0,\dots,a_m)) , (v,(a_i,\dots,a_j)) )$, there is exactly one 2-morphism if and only if composing along each sequence yields the same morphism in $\cube^N \times \HMshape^\omega(\matchingsk)$.  A $2$-morphism in $(\cube^N \times \HMshape^\omega(\matchingsk) )'$  is then sent to the unique 2-morphism in $\divcobq'$ which indicates that the resulting sequences of divided cobordisms compose to give isotopic cobordisms.

\end{itemize}
\end{definition}

\begin{lemma}
Definition \ref{def:func HMshape to qdivcob} gives a well-defined strict 2-functor 
\[(\cube^N \times \HMshape^\omega(\matchingsk))' \rar{\cobfunc_T^{H^\omega}} \divcobq'.\]
\end{lemma}
\begin{proof}
In the canonical thickenings, compositions of 1-morphisms are formal concatenations, and thus $\cobfunc_T^{H^\omega}$ is functorial with respect to 1-morphisms by definition.  We must check that it is well-defined on 2-morphisms.  Let $t,t'$ be two 1-morphisms (sequences) in $(\cube^N\times \HMshape^\omega(\matchingsk))'$  related by a 2-morphism, so $t$ and $t'$ have the same starting and ending objects.  Suppose $t$ (respectively $t'$) contains $\ell$ (respectively $\ell'$) basic morphisms which delete the last element of a tuple without arriving at the empty tuple.  Then the sequence $\cobfunc_T^{H^\omega}(t)$ (respectively $\cobfunc_T^{H^\omega}(t')$) composes to a morphism of the form $(\Sigma,\Gamma,q^{-n\ell})$ (respectively $(\Sigma',\Gamma',q^{-n\ell'})$).  The divided cobordisms $(\Sigma,\Gamma)$ and $(\Sigma',\Gamma')$ are both isotopic to a third divided cobordism $(\Sigma'',\Gamma'')$ where all saddles are done at the same height and pulling through $\Mq$ is done afterwards.  These isotopies consist of adjusting heights of saddles away from $\Mq$, together with pulling $n$ saddles through $\Mq$ for every morphism that deletes the last element of a tuple except for the last such morphism.  As indicated in Figure \ref{fig:trace moves}, each saddle pull contributes a shift of $q$, indicating that 
\[\FAq(q^{-n\ell}\Sigma) = \FAq(q^{-n} \Sigma'') = \FAq(q^{-n\ell'}\Sigma')\]
as required for the existence of a (unique) 2-morphism between $\cobfunc_T^{H^\omega}(t)$ and $\cobfunc_T^{H^\omega}(t')$.
\end{proof}

Altogether then, the functor $(\cube^N \times \HMshape^\omega(\matchingsk))' \rar{\cobfunc_T^{H^\omega}} \divcobq'$ is defined on non-empty tuples similarly to the functor $\cobfunc_T$ of Section \ref{sec:spectral platalg and tangle bimod} via the embedding $\D\rar{\iota}\Aq$, with the empty tuple assigned the standard closure of our tangles in $\Aq$.  For morphisms, we embed our divided manifolds into $\Aq\times I$ in a manner that avoids $\Mq$, \emph{unless} we delete the final matching of a tuple (which corresponds to performing the last face map in the bar complex).  Note that morphisms of the form $m_a:(v,(a))\rar{} (v,\omega)$ are assigned (divided) multisaddles as in Figure \ref{fig:A map}, so that we have
\begin{equation}\label{eq:A map lifted}
    \BtoAb\circ \qVHKKchoice' \circ \cobfunc_T^{H^\omega} (m_a) = A_a^{\choice};
\end{equation}
see the proof of Proposition \ref{prop:quantum annular V recovers the right thing}.

\subsubsection{Modifying the composition $\qVHKKchoice'\circ \cobfunc_T^{H^\omega}$}
\label{sec:modifying qVHKK}

At this stage we modify the composition $\qVHKKchoice'\circ \cobfunc_T^{H^\omega}$ analogously to the adjustments in Section \ref{sec:spectral platalg and tangle bimod} and \cite[proof of Theorem 8]{LLS_CK_Spectra}.  Recall Definitions \ref{def:qVHKK Z lax} and \ref{def:qVHKK} for what follows.

\begin{definition}\label{def:modified qVHKK of cobfunc for HM construction}
The functor $(\cube^N \times \HMshape^\omega(\matchingsk))' \rar{\modifiedqHMcombo} \gBurn_G$ is defined as follows.
\begin{itemize}
    \item Let $O$ be a non-closure object in $(\cube^N \times \HMshape^\omega(\matchingsk))'$.  If $\cobfunc_T^{H^\omega}(O)$ contains any Type III circles, then $\modifiedqHMcombo(O):=\varnothing$.  Otherwise if $\qVHKKchoice'\circ \cobfunc_T^{H^\omega}(O)=\Gen(Z)\times G$, then set $\modifiedqHMcombo(O):=(\Gen(Z)\setminus\Gen_{II}(Z))\times G$ as in Section \ref{sec:spectral platalg and tangle bimod}.
    \item If $O=(v,\omega)$ is a closure object in $(\cube^N \times \HMshape^\omega(\matchingsk))'$, then set 
    \[
    \modifiedqHMcombo(v,\omega):=\Gen_-^+(\closure{\strands{T}_v})\times G
    \]
    where $\Gen_-^+$ indicates that we include only the Khovanov generators which are in annular degree zero, label the inner $n-k$ circles by $v_-$, and label the outer $k$ circles by $v_+$.  Compare to Equation \eqref{eq:C map on generators} for the $C$ map.  In this way, all objects are assigned sets of generators in annular degree zero, so that we may ignore this grading, viewing the functor as landing in $\gBurn_G$ via the quantum grading.
    
    \item If $m$ is any basic morphism in $(\cube^N \times \HMshape^\omega(\matchingsk))'$ with correspondence $\qVHKKchoice'\circ \cobfunc_T^{H^\omega}(m)= S \times G$, define $\modifiedqHMcombo(m):= (S\setminus \widetilde{S})\times G$ where $\widetilde{S}$ consists of all elements $e\in S$ such that either $s(e)$ or $t(e)$ were removed from the relevant objects above.  The source and target maps for $\modifiedqHMcombo(m)$ are then defined to be the restrictions of the source and target maps for $\qVHKKchoice'\circ \cobfunc_T^{H^\omega}(m)$.  More general sequences of basic morphisms in $(\cube^N \times \HMshape^\omega(\matchingsk))'$ are sent to the composition of the relevant basic morphisms in $\gBurn$ as implicit in \cite[proof of Theorem 8]{LLS_CK_Spectra}; see also Definitions \ref{def:VHKK} and \ref{def:qVHKK} as well as \cite[Lemma 3.16]{LLS_Khovanov_Spectra}.
    \item If $\tau$ is a 2-morphism in $(\cube^N \times \HMshape^\omega(\matchingsk) )'$, then $\modifiedqHMcombo(\tau)$ is the restriction of the 2-morphism $\qVHKKchoice'\circ \cobfunc_T^{H^\omega}(\tau)$.
\end{itemize}
\end{definition}

\begin{lemma}\label{lem:modified qHM combo is a 2functor}
Definition \ref{def:modified qVHKK of cobfunc for HM construction} gives a well-defined strict 2-functor
\[(\cube^N \times \HMshape^\omega(\matchingsk))' \rar{\modifiedqHMcombo} \gBurn_G.\]
\end{lemma}
\begin{proof}
Our modifications transforming $\qVHKKchoice'\circ \cobfunc_T^{H^\omega}$ into $\modifiedqHMcombo$ are entirely analogous to the modifications transforming $V_{HKK}\circ G^T$ into  $L^T$ described in \cite[proof of Theorem 8]{LLS_CK_Spectra}.  The key combinatorial checks that the modified 2-functor is well-defined are handled via \cite[Lemma 4.2, Lemma 5.4]{LLS_CK_Spectra}, and the same arguments apply here.
\end{proof}

The composition $K_G\circ \modifiedqHMcombo$ can then be rectified (Lemma \ref{lem:G rectification}) to yield  a functor 
\[
\cube^N \times \HMshape^\omega(\matchingsk)\rar{} \gSpectra_G,
\]
for which we can take the iterated mapping cone (Definition \ref{def:IMC}) to arrive at a functor
\begin{equation}\label{eq:rect qHM combo}
\rectifiedqHMcombo : \HMshape^\omega(\matchingsk)\rar{} \gSpectra_G.
\end{equation}

\subsubsection{Finishing the proof}
\label{sec:proof of main qTHH thm}

Let $L=\closure{T} \subset \Aq$ denote the annular link diagram obtained by closing $T$. The main point is that $\rectifiedqHMcombo$ is weakly equivalent to an
augmented quantum Hochschild-Mitchell functor
(Definition \ref{def:qHM functor}) for the bimodule $\qXk(T)$, and the induced map from Corollary \ref{cor:aug qHM gives map on qTHH} recovers the quasi-isomorphism $\qXi^\choice$ of Proposition \ref{prop:qXichoice is a quasi-iso} at the level of homology. Consider the composition 

\[(\cube^N)'\rar{\cobfunc_L}\divcobq' \rar{\qVHKKchoice'} \ggBurn_G \]
from  Definition \ref{def:qAKh spectrum}. There is a full sub-$2$-category  $(\cube^N \times \{\omega\} )'$ of $(\cube^N \times \HMshape^\omega(\matchingsk))'$ generated by objects of the form $(v,\omega)$ for $v\in \cube^N$. We can identify $(\cube^N \times \{\omega \} )'$ with $(\cube^N)'$ since $\omega$ is terminal in $\HMshape^\omega(\matchingsk)$. Letting $\modifiedqHMcomborestricted$ denote the restriction of $\modifiedqHMcombo$ to this sub-category, we see that forgetting the inner $n-k$ and outer $k$ circles yields a natural isomorphism from $\modifiedqHMcomborestricted$ to the part of $\qVHKKchoice'\circ \cobfunc_L$ supported in annular grading $n-2k$ (viewed as a functor to $\gBurn_G$).  Rectifying $K_G\circ \modifiedqHMcomborestricted$ results in a functor $ \cube^N \to \gSpectra_G$, and taking the iterated mapping cone yields
\[
\mathfrak{R}_T^\omega \in  \gSpectra_G.
\]
It follows that $\sh^{-N_+} \mathfrak{R}_T^{\omega}$ is naturally isomorphic to $\qaXchoice(L; n-2k)$, and thus the value of the already rectified functor
\[
\sh^{-N_+}\rectifiedqHMcombo(\omega) \in \gSpectra_G
\]
on $\omega$ is weakly equivalent to $\qaXchoice(L; n-2k)$ (via a zig-zag of weak equivalences as in \cite[Lemma 2.43]{LLS_Khovanov_Spectra}).

\begin{theorem}\label{thm:qTHH recovers qAX}
For any $(n,n)$-tangle diagram $T$ with annular closure $L=[T]$, any $0\leq k\leq n$, and any coherent choice of generators $\choice$, we have a zig-zag of weak equivalences of spectra with $G$-action
\[ \qTHH(\qplatalgk;\qXk(T) ) \simeq \qaXchoice(L;n-2k).\]
\end{theorem}
\begin{proof}
Comparing Definition \ref{def:modified qVHKK of cobfunc for HM construction} for $\modifiedqHMcombo$ with the constructions in Section \ref{sec:spectral platalg and tangle bimod} shows that, for any non-closure object $(v,(a_0,\dots,a_\ell))$, we have an isomorphism of objects (sets) in $\gBurn$
\[
\modifiedqHMcombo(v,(a_0,\dots,a_\ell)) \cong \func{G}_T(v,a_\ell\mid a_0) \times \func{G}_T(a_0,a_1) \times \cdots \times \func{G}_T(a_{\ell-1},a_\ell) \times G.
\]
(Recall that $\func{G}_T=\func{G}_\mathcal{A}$ on pairs separated by commas rather than separated by a bar.)  On basic morphisms $m$ in $(\cube^N\times \shape_{HM}^\omega(\matchingsk) )'$ which do not delete the last entry of the tuple, this isomorphism is natural in the following sense.  The morphism $m$ corresponds to a collection of multimorphisms $\til{m}_i$ in $\cube^N\tiltimes \shape(\matchingsk,\matchingsk)$ based upon a morphism in $\cube^N$ together with strings of consecutive deleted matchings from the tuple.  Then $\modifiedqHMcombo(m)$ acts on the $\func{G}_T$ factors according to the various $\func{G}_T(\til{m}_i)$. This follows mainly from the definition of $\cobfunc_T^{H^\omega}$ utilizing the embedding $\D\rar{\iota}\Aq$ disjoint from the $q$-seam $\mq$.  Morphisms which delete the final element of the tuple involve saddle cobordisms which pass through the seam, which are treated by $\qVHKKchoice'$ in accordance with $\choiceFAq$ which is described on such cobordisms in Lemma \ref{lem:FAq recovers deformed face map}.

After composing with $K_G$ and utilizing properties \ref{item:K(X) is wedge of sphere} and \ref{item:K(A) is the right map} of $K$, together with the natural isomorphisms provided by Lemma \ref{lem:scalar ext properties}, we see that the functor $K\circ\modifiedqHMcombo$, when restricted to objects and basic morphisms involving only a single vertex $v$ of the cube $\cube^N$, behaves like an augmented quantum Hochschild-Mitchell functor for the scalar extension $(K\circ \func{G}_T)\wedge \Sq$ restricted to the same vertex $v$.  Rectifying and taking iterated mapping cones along $\cube^N$ maintains this behavior up to zig-zags of natural weak equivalences so that $\sh^{-N_+}\rectifiedqHMcombo$ is naturally weakly equivalent to an augmented quantum Hochschild-Mitchell functor $\func{F}$ for $\qXk(T)$ (compare to \cite[Proof of Theorem 8]{LLS_CK_Spectra}).  The induced map $F$ of Corollary \ref{cor:aug qHM gives map on qTHH} fits in a  commuting diagram involving zig-zags (all entries interpreted in $\gSpectra_G$)
\begin{equation}\label{eq:zigzag with augHMmap}
\begin{tikzcd}
    \qTHH(\qplatalgk;\qXk(T) ) \ar[d,"F"] & \ar[l] \cdots \ar[r] & \hocolim_{\HMshape} \sh^{-N_+} \rectifiedqHMcombo|_{\HMshape} \ar[d,"R"] \\
    \func{F}(\omega) & \ar[l] \cdots \ar[r] & \sh^{-N_+}\rectifiedqHMcombo(\omega)
\end{tikzcd}
\end{equation}
where the map $R$ is constructed in the same manner as $F$ but for the functor $\modifiedqHMcombo$ instead of $\func{F}$.

Taking into account the further zig-zag of weak equivalences between $\sh^{-N_+}\rectifiedqHMcombo(\omega)$ with $\qaXchoice(L;n-2k)$ and comparing the combination of Equation \eqref{eq:A map lifted} and Definition \ref{def:modified qVHKK of cobfunc for HM construction} with the constructions in Section \ref{sec:explicitly constructing the quasi-iso} shows that $\chainsfunc(R)$ recovers the quasi-isomorphism $\qXi^\choice$ of Proposition \ref{prop:qXichoice is a quasi-iso} up to further zig-zags of weak equivalences.  Upon taking homology, all such zig-zags become isomorphisms, including those in Equation \eqref{eq:zigzag with augHMmap}, completing the proof.
\end{proof}

Combining Theorem \ref{thm:qTHH recovers qAX} with Proposition \ref{prop:quantum annular V recovers the right thing} moreover gives the following corollary, which could also be deduced as in \cite[Proposition 4.2]{LLS_Khovanov_Spectra}; compare with \cite[Proposition 7.5]{LLS_Khovanov_Spectra} and \cite[Proposition 5.1]{LLS_CK_Spectra}.
\begin{corollary}
    There is an isomorphism of $\Z[G]$-modules
\[
H_*\left(\qTHH(\qplatalgk;\qXk(T) ) \right) \cong \qhh(\qPlatk{n}{k};\qCCKnk(T)).
\]
\end{corollary}


\section{Equivalence with $\qakwX$}
\label{sec:comparison} 

We now fix $G$ to be a finite cyclic group with generator $q$. Let $L$ be an annular link with diagram $D$. A spectrum with $G$-action $\qakwX(D)\in \ggSpectra_G$ was introduced in \cite{AKW} whose equivariant stable homotopy type is an invariant of $L$. In this section we show that $\qaX(L)$ from Definition \ref{def:qAKh spectrum} and $\qakwX(L)$ are weakly equivalent in $\ggSpectra_G$. Note that both spectra are constructed using  a choice of generators $\choice$ which we typically omit from the notation; the result is independent of this choice by Proposition \ref{prop:invariance} and \cite[Theorem 5.10]{AKW}. The argument is similar to that of \cite[Section 8]{LLS-Burnside}; in particular, we define an intermediate $\hatX(L) \in \ggSpectra_G$ in Section \ref{sec:thickened diagram and comparison} and show that $\qaX(L) \simeq \hatX(L) \simeq \qakwX(L)$.

Recall the $2$-category  $\til{\B}_G$  from Definition \ref{def: finite Burnside}.
Given a choice of generators $\choice$ and a link diagram $D\subset \Aq$ with $N$ crossings, a strictly unitary lax $2$-functor (see Lemma \ref{def:lax 2-functor})
\[
F_{\mathfrak{q}}' : (\cube^N)^{\rm op}  \to \til{\B}_G
\]
was defined in \cite[Section 4.1]{AKW}. Note, in \cite{AKW} the functor was denoted just by $F_{\mathfrak{q}}$ (to be explained shortly), and the notation $\cube^N$ meant the opposite category to the one in the present paper; see also Remarks \ref{rmk:homology vs cohomology} and \ref{rmk:homology vs cohomology 2}. By taking opposites and using the isomorphism $\til\B_G \cong \til \B_G^{\rm op}$ that is the identity on objects and sends $X\lar{} A \rar{} Y$ to $Y\lar{} A \rar{} X$, we obtain a functor 
\begin{equation}
\label{eq:AKW Burnside functor}
F_{\mathfrak{q}} : \cube^N  \to \til\B_G.
\end{equation}
The dependence on $\choice$ and $D$ is not included in the notation. By \cite[Proposition 4.3]{AKW}, $F_\mathfrak{q}$ is independent of $\choice$ up to natural isomorphism.

Consider the  lax $2$-functor $F_\q^{\#} : \cube^N \to \til\B_G$ associated with $D$, defined as the composition of lax functors defined analogously to Definitions \ref{def:cobfunc for links in Aq} and \ref{def:qVHKK}, except without using canonical thickenings.

\begin{lemma}
\label{lem:AKW functor agrees with HKK}
   The functors $F_\q^{\#}$ and $F_\mathfrak{q}$ are naturally isomorphic.
\end{lemma}
\begin{proof}
The functors agree on objects and $1$-morphisms. 
By \cite[Lemma 8.5]{LLS-Burnside}, the Hu-Kriz-Kriz approach used to define $2$-morphisms associated to a genus $1$ surface agrees with the ladybug matching approach in \cite{LS-Kh-spectra, LLS-Burnside}.  Then \cite[Corollary 2.21]{AKW} implies that, when there is no choice to be made (in the case of a genus $1$ surface with contractible boundary and nontrivial winding number), it actually did agree with the ladybug matching made with the left pair. See also Remark \ref{rmk:intersection is -1}. 
\end{proof}

\subsection{From Burnside functors to spectra via little boxes}

We give a summary of the little box construction from \cite[Section 5]{LLS-Burnside} and its $G$-equivariant analogue from \cite[Section 5]{AKW}.  This builds a spectrum  with $G$-action from a strictly unitary lax $2$-functor $F: \cube^N \to \til{\B}_G$.  We shall refer to such functors as \emph{Burnside functors}.

A $k$-dimensional box is $\prod_{i=1}^k [a_i, b_i] \subset \R^k$. Given two $k$-dimensional boxes $B$ and $B'$, there is a canonical homeomorphism $B\rar{\sim} B'$ obtained by scaling and translating the ambient space $\R^k$. 
Fix an identification $S^k = [0,1]^k/\d([0,1]^k)$ so that $B/\d B$ is canonically identified with $S^k$ for any $k$-dimensional box $B$. 

Suppose we have a correspondence $X\xleftarrow{s} A \xrightarrow{t} Y$ in $\til{\B}$. Fix a collection of disjoint boxes $\{B_x\}_{x\in X}$ and let $F(\{B_x\},s,t)$ be the space of all collections of boxes $e = \{B_a\}_{a\in A}$ such that $B_a\subset B_{s(a)}$ for each $a\in A$ and $B_a \cap B_{a'} = \varnothing$ whenever $t(a) = t(a')$. Each $B_a \subset B_{s(a)}$ is called a \emph{little box}. 

Given $e\in F(\{B_x\},s,t)$, define the \emph{overlapping box map}
\begin{equation}
\label{eq:overlapping box map}
\Phi'(e,A) : \bigvee_{x\in X} S^k_x \to \prod_{y\in Y} S^k_y
\end{equation}
as follows. On the wedge summand corresponding to $x\in X$, $\Phi'(e,A)$ is given as a product of maps $\prod_{y\in Y} \Phi_y'(e,A)$ where $\Phi_y'(e,A)$ is the composition 
\[
  S^k_x  = B_x/\partial B_x \to B_x / ( B_x \setminus \bigcup_{a\in s^{-1}(x) \cap t^{-1}(y)} B_a^{\text{int}} ) = \bigvee_{a\in s^{-1}(x) \cap t^{-1}(y)} B_a/\partial B_a \rar{t} \bigvee_{y\in Y} S^k_y
\]
The first map is a quotient and the last map sends the sphere $B_a/\d B_a$  to the sphere $B_{t(a)}/ \d B_{t(a)} = S^k_{t(a)}$ via the canonical homeomorphism $B_a\cong B_{t(a)}$.

Let 
\[
E(\{B_x\},s)\subset  F(\{B_x\},s,t)
\]
be the subspace consisting of $e = \{B_a\}_{a\in A}$ where $B_a\subset B_{s(a)}$ and $B_{a} \cap B_{a'} = \varnothing$ whenever $a\neq a'$. A point $e=\{B_a\} \in E(\{B_x\},s)$ determines a \emph{disjoint box map} 
\begin{equation}
\label{eq:disjoint box map}
\Phi(e,A): \bigvee_{x\in X} S^k_x \to \bigvee_{y\in Y} S^k_y
\end{equation}
which is given on each wedge summand as the composition
\begin{equation}\label{eq:box map on each summand}
S^k_x = B_x/\d B_x \to B_x/ ( B_x\setminus \bigcup_{a\in s^{-1}(x)} B_a^{\text{int}} ) = \bigvee_{a\in s^{-1}(x)} B_a/\d B_a \xrightarrow{t} \bigvee_{y\in Y}S^k_y.
\end{equation}
 A map  $\bigvee_{x\in X} S^k_x \to \bigvee_{y\in Y} S^k_y$ of this form is said to \emph{refine the correspondence } $X \xleftarrow{s} A \xrightarrow{t} Y$.

The natural inclusion 
\[
\iota_X: \bigvee_{x\in X} S^k_x \to \prod_{x\in X} S^k_x
\]
(inclusion of the $k$-skeleton) is also the overlapping box map assigned to $X\lar{\id} X \rar{\id} X$ where each little box is the entire $B_x$. As discussed in \cite[Example 5.4]{LLS-Burnside},  if $e\in E(\{B_x\},s)$ then $\Phi'(e,A) = \iota_Y \circ \Phi(e,A)$. 

Now suppose $X\lar{s} A \rar{t} Y$ is a $1$-morphism in $\til{\B}_G$. The spaces $\bigvee_{x\in X} S^k_x$ and $\bigvee_{x\in Y} S^k_y$ have  natural $G$-actions given by permuting spheres according to the action of $G$ on $X$ and $Y$. At the level of boxes, $g\in G$ sends $S^k_x = B_x/\partial B_x$ to $S^k_{gx} = B_{gx}/\partial B_{gx}$ via the canonical homeomorphism $B_x\cong B_{gx}$. 
To ensure that overlapping and disjoint box maps are $G$-equivariant, define the subspaces 
\[
F_G(\{B_x\},s,t)\subset  F(\{B_x\},s,t), \ \ E_G(\{B_x\},s)\subset  E(\{B_x\},s)
\]
where $F_G(\{B_x\},s,t)$ consists of little boxes $\{B_a\}_{a\in A}$ such that the canonical homeomorphism $B_x \to B_{gx}$ restricts to the canonical homeomorphism $B_{a} \to B_{ga}$ for each $x\in X$, $a\in s^{-1}(x)$, and $g\in G$. The subspace $E_G(\{B_x\},s)$ is defined by the same condition.

\begin{lemma}
\label{lem:equivariant box maps facts}
Let $X\lar{s} A \rar{t} Y$ be a $G$-equivariant correspondence, and let $e\in F_G(\{B_x\}, s,t)$.
\begin{enumerate}
    \item The overlapping box map $\Phi'(e,A)$ is $G$-equivariant.
    \item If $e\in E_G(\{B_x\}, s)$, then the disjoint box map $\Phi(e,A)$ is $G$-equivariant.
    \item $F_G(\{B_x\}, s,t)$ and  $E_G(\{B_x\}, s)$ are $(k-2)$-connected.
\end{enumerate}
    
\end{lemma}

\begin{proof}
    Item (2) is \cite[Lemma 5.1]{AKW}. Item (1) is proven similarly. For item (3), the statement regarding $E_G(\{B_x\}, s)$ is \cite[Corollary 5.3]{AKW}. The proof of the statement regarding $F_G(\{B_x\}, s,t)$ is analogous: for a fixed $x\in X$, a collection of little boxes $\{B_a\}_{a \in s^{-1}(x)}$ in $B_x$ such that $B_a \cap B_{a'} = \varnothing$ whenever $t(a) = t(a')$ uniquely determines, via the canonical homeomorphism between boxes, a collection of little boxes $\{B_c\}_{c\in s^{-1}(gx)}$ in $B_{gx}$ satisfying the same disjointness property, for each $g\in G$. 
    Therefore we have a homeomorphism 
    \[
F_G(\{B_x\}, s,t) \cong F( \{ B_{[x]}\}_{[x]\in X/G}, s/G, t/G).
    \]
    The right-hand side is $(k-2)$-connected by \cite[Lemma 5.18]{LLS-Burnside}. 
\end{proof}

We write \emph{equivariant} overlapping and disjoint box maps to mean maps $\Phi'(e,A)$ and $\Phi(e,A)$ associated with little boxes $e\in F_G(\{B_x\},s,t)$ and $e\in E_G(\{B_x\},s)$, respectively. Note that the inclusion $\iota_X$ above is an equivariant overlapping box map.

Let $\Top_*^G$ denote the category of based $G$-spaces. We do not recall the definitions here but refer the reader to \cite[Section 4.2]{odd_khovanov_homotopy} for the notion of a \emph{homotopy coherent diagram} in $\Top_*^G$ as well as the homotopy colimit of such a diagram.

\begin{definition}\label{def:spatial refinement}
Let $\cat{C}$ be a small category and let $F:\cat{C} \to \til{\B}_G$ be a strictly unitary lax 2-functor. An (equivariant) \emph{$k$-dimensional spatial refinement} of $F$ is a homotopy coherent diagram $\til{F}_k:\cat{C} \to \Top^G_*$ satisfying
\begin{enumerate}[label= (\arabic*)]
\item\label{it:vertices are spheres} For any $u\in \cat{C}$, $\til{F}_k(u) = \bigvee_{x\in F(u)} S^k$.\\
\item\label{it:edges are box maps} For any sequence $u_0 \xrightarrow{f_1} \cdots \xrightarrow{f_m} u_m$ of composable morphisms in $\cat{C}$ and $t\in I^{m-1}$, the map
\[
\til{F}_k(f_m,\ldots, f_1)(t) : \bigvee_{x\in F(u_0)} S^k \to \bigvee_{y\in F(u_m)} S^k
\]
is an equivariant disjoint box map which refines the correspondence $F(f_m \circ \cdots \circ f_1)$. 
\end{enumerate}
\end{definition}

By \cite[Proposition 5.5]{AKW},  if $\cat{C}$ is a small category in which every sequence of composable non-identity morphisms has length at most $N$ and $F:\cat{C} \to \til{\B}_G$ is a strictly unitary lax 2-functor, then there is a $k$-dimensional spatial refinement of $F$ for every $k \geq N$.

Recall the enlarged cubes $\cube^N_+$ and $\cube^N_\dagger$ from Definition \ref{def:IMC}. 
Let $F:\cube^N \to \til{\B}_G$ be a Burnside functor, and let $\til{F}_k : \cube^N \to \Top_*^G$ be a $k$-dimensional spatial refinement of $F$.
Extend $\til{F}_k$ to homotopy coherent diagrams $\til{F}_k^+ : \cube^N_+ \to \Top_*^G$ and $\til{F}_k^\dagger : \cube^N_\dagger \to \Top_*^G$ by setting their value on any object not in $\cube^N$ to be a point.
Define the space
\begin{equation}\label{eq:realization}
\LR{F}_k := \hocolim  \til{F}_k^+.
\end{equation}
Since the homotopy coherent diagram $\til{F}_k$ takes values in $\Top_*^G$, the space $\LR{F}_k$ is again a based $G$-space and is well-defined up to equivariant homotopy equivalence by \cite[Lemma 5.6]{AKW}. By \cite[Lemma 5.27]{LLS-Burnside}, $\hocolim \til{F}_k^\dagger \simeq \LR{F}_k$, so one could work with the larger cube $\cube^N_\dagger$ instead.

Let $L\subset \Aq\times I$ be an annular link and let $D\subset \Aq$ be a diagram for $L$ with $N$ crossings, $N_+$ of which are positive. Consider the associated functor $F_\mathfrak{q} : \cube^N \to \til{\B}_G$ from \eqref{eq:AKW Burnside functor}.  Consider the space $\LR{F_\mathfrak{q}}_k$ assigned to 
a $k$-dimensional spacial refinement $\til{F}_{\mathfrak{q},k}$ of $F_\mathfrak{q}$ extended to $\cube^N_+$, as in \eqref{eq:realization}.

\begin{definition}
    \label{def:AKW spectrum}
    The quantum annular Khovanov spectrum of $D$ is the spectrum with $G$-action 
    \[
\qakwX(D) := \Sigma^{-k-N_+}\left( \Sigma^{\infty} \LR{F_\mathfrak{q}}_k\right),
    \]
   obtained by desuspending $k+N_+$ times the suspension spectrum of $\LR{F_\mathfrak{q}}_k$. 
\end{definition}

The following is a combination of results from \cite{AKW}. 
\begin{theorem}
Up to equivariant homotopy equivalence, $\qakwX(D)$ is an invariant of $L$.
Moreover, the homology of $\qakwX(D)$ is isomorphic to $\qaKh(L)$ as $\Z[G]$-modules. 
\end{theorem}

\begin{remark}
\label{rmk:homology vs cohomology 2}
    As discussed in Remark \ref{rmk:homology vs cohomology}, \cite{AKW} uses the functor $F_\mathfrak{q}'$ and builds a spectrum whose \emph{cohomology} recovers the quantum annular homology of $L$, viewed as a \emph{cohomology} theory by negating the homological gradings. That dualizing (passing from $F_\mathfrak{q}'$ to $F_\mathfrak{q}$ and from cohomology to homology) goes through as expected follows, for instance, from \cite[Proposition 5.7]{AKW}, which discusses the cellular \emph{chain} complex of $\LR{F}_k$. 

\end{remark}

\subsection{The thickened diagram and comparison with $\qakwX(L)$}
\label{sec:thickened diagram and comparison}

Now we recall the construction of the \emph{thickened diagram} from \cite[Section 4.4]{LLS-Burnside}. As an executive summary, given a strictly unitary lax $2$-functor $F: \cat{D} \to \til{\B}$, where $\cat{D}$ is a $1$-category viewed as a $2$-category with only identity $2$-morphisms, there is an associated category $\h{D}$ and \emph{strict} functors $\h{F}_k : \h{\cat{D}} \to \Top_*$, the target being the category of based topological spaces. Assembling the various $k$'s produces $\h{F} : \h{\cat{D}} \to \Spectra$. We discuss this construction and the relevant modifications below. 

\begin{definition}
    Let $\h{\cat{D}}$ be the category whose objects are pairs of composable morphisms $ u \rar{f} v \rar{g} w$ in $\cat{D}$ and whose  morphisms are commutative diagrams
\begin{equation}
\label{eq:morphism in hat category}
\begin{tikzcd}
    u  \ar[r,"f"] \ar[d,"\varphi"] & v \ar[r,"g"] & w \ar[d, "\theta"] \\
    u' \ar[r,"f'"] & v' \ar[u,"\psi"'] \ar[r,"g'"] & w'
\end{tikzcd}.
\end{equation}
Note the direction of the middle arrow. Composition is defined by stacking diagrams, $(\varphi', \psi', \theta') \circ (\varphi, \psi, \theta) := (\varphi'\circ \varphi, \psi\circ \psi', \theta' \circ \theta)$.  
\end{definition}

Given a strictly unitary lax $2$-functor $F: \cat{D} \to \til{\B}$ and $k\geq 1$,  $\h{F}_k: \h{\cat{D}} \to \Top_*$ is defined on objects as 
\begin{equation}
\label{eq:Fhat on objects}
\h{F}_k( u \rar{f} v \rar{g} w) = \bigvee_{a\in F(f)} \prod_{\substack{ b\in F(g) \\ s(b) = t(a)}} S^k.
\end{equation}
For a morphism  $(\varphi, \psi, \theta)$ as in \eqref{eq:morphism in hat category}, the map 
\[
\bigvee_{a\in F(f)} \prod_{\substack{ b\in F(g) \\ s(b) = t(a)}} S^k
\to 
\bigvee_{a'\in F(f')} \prod_{\substack{ b'\in F(g') \\ s(b') = t(a')}} S^k
\]
is defined as follows.
We have bijections 
\begin{align}
    F(f) & \cong F(\psi) \times_{F(v')} F(f') \times_{F(u')} F(\varphi) \label{eq:hat func on morphism 1} \\
      F(g') & \cong F(\theta) \times_{F(w)} F(g) \times_{F(v)} F(\psi)  \label{eq:hat func on morphism 2}.
\end{align}
Let $(y,a',x) \in F(\psi) \times_{F(v')} F(f') \times_{F(u')} F(\varphi) $ be the element corresponding to $a\in F(f)$. The map $\h{F}_k(\varphi, \psi, \theta)$ will send the summand corresponding to $a\in F(f)$ to the summand corresponding to $a' \in F(f')$. 
Define 
\begin{equation}
\label{eq:hat F on morphism step 1}
\prod_{\substack{ b\in F(g) \\ s(b) = t(a)}} S^k
\rar{\prod_b \Delta_b}
\prod_{\substack{ b\in F(g) \\ s(b) = t(a)}} 
\prod_{(z,b,y) \in F(\theta)\circ F(g) \circ F(\psi)}
S^k
\cong
\prod_{\substack{ b'\in F(g') \\
b' \leftrightarrow (z,b,y) \\ 
s(b) = t(a)} }
S^k,
\end{equation}
where $\Delta_b$ is the diagonal map, in the middle space we write $\circ$ rather than fiber products to avoid clutter, and the right homeomorphism uses \eqref{eq:hat func on morphism 2} to relabel the indexing sets. Note that 
\[
\{ b'\in F(g') \mid  b' \leftrightarrow (z,b,y), s(b) = t(a)  \}
\subset 
\{
b' \in F(g') \mid s(b') = t(a')
\}
\]
since $t(a') = s(y) = s(b')$. Define  $\h{F}_k(\varphi, \psi, \theta)$ on the summand corresponding to $a\in F(f)$ by extending \eqref{eq:hat F on morphism step 1} to 
\[
\prod_{\substack{ b\in F(g) \\ s(b) = t(a)}} S^k
\to 
\prod_{\substack{ b'\in F(g') \\ s(b') = t(a')}} S^k
\]
by mapping to the basepoint in the complementary factors.

\begin{definition}
\label{def:burnside functor to hat functor G}
    Given a strictly unitary lax $2$-functor $F:\cat{D} \to \til{\B}_G$, define  $\h{F}_k^G : \h{\cat{D}} \to \Top_*^G $ as follows. First, compose with the forgetful functor $\til{\B}_G \to \til{\B}$ to obtain $F' : \cat{D} \to \til{\B}$ and apply the above construction to form $\h{F}'_k : \h{\cat{D}} \to \Top_*$. The value of $\h{F}'_k$ on an object $u\rar{f} v \rar{g} w$ of $\h{D}$ (see \eqref{eq:Fhat on objects}) inherits a $G$-action by declaring that $q^k\in G$ sends the wedge summand $\prod_{\substack{ b\in F(g) \\ s(b) = t(a)}} S^k$ corresponding to $a\in F(f)$ to the wedge summand $\prod_{\substack{ b'\in F(g) \\ s(b') = t(q^k a)}} S^k$ corresponding to $q^k a \in F(f)$ by a map that permutes the factors, $S^k_{a,b} \rar{\id} S^k_{q^k a, q^k b}$. 
\end{definition}

\begin{lemma}
    The above definition gives a functor $\h{F}^G_k : \h{\cat{D}} \to \Top_*^G$. 
\end{lemma}
\begin{proof}
    One needs to check that the map assigned to a morphism in $\h{\cat{D}}$ by $\h{F}'_k$ is a map of $G$-spaces. We leave the straightforward check to the reader. 
\end{proof}

Assembling the functors $\h{F}_k^G$ over all $k\geq 0$ as in \cite[Section 4.4]{LLS-Burnside}, we obtain a functor $\h{F}^G: \h{\cat{D}} \to \Spectra_G$.
There are enlargements $\h{\cube^N}_+$ and $\h{\cube^N}_\dagger $  of $\h{\cube^N}$ analogous to those in Definition \ref{def:IMC}. Define $\h{\cube^N}_+$ by adding one object $*$ and a unique morphism $(u\to v\to w) \to *$ whenever $w\neq \vec{1} = (1,1,\ldots, 1)$. Also set $\h{\cube^N}_\dagger := ( \h{\cube^1}_+)^N$. Given a functor $F:\cube^{N}\to \til{\B}_G$, extend $\h{F}^G : \h{\cube^N} \to \Spectra_G$ to 
\[
(\h{F}^G)^+ : \h{\cube^N}_+ \to \Spectra_G  \ \ \text{ and } \ \ (\h{F}^G)^\dagger : \h{\cube^N}_\dagger \to \Spectra_G
\]
in the natural ways.

Now consider the Burnside functor $F_\mathfrak{q}$ \eqref{eq:AKW Burnside functor} associated to an $N$-crossing annular link diagram $D$ of a link $L$. Form $\h{F}_\mathfrak{q}^G : \h{\cube^N} \to \Spectra_G$ from Definition \ref{def:burnside functor to hat functor G} and its extensions  
$(\h{F}_\mathfrak{q}^{G})^+: \h{\cube^N}_+ \to \Spectra_G$ and $(\h{F}_\mathfrak{q}^{G})^\dagger : \h{\cube^N}_\dagger \to \Spectra_G$. 
By \cite[Lemma 4.41]{LLS-Burnside}, 
\[
\hocolim (\h{F}_\mathfrak{q}^{G})^+ \simeq \hocolim (\h{F}_\mathfrak{q}^{G})^\dagger.
\]
Note, the above equivalence can be taken to be $G$-equivariant by examining the proof of \cite[Lemma 4.41]{LLS-Burnside}. 
The \emph{realization} of $F_\mathfrak{q}$ is defined to be
\[
\lr{F_\mathfrak{q}} := \hocolim (\h{F}_\mathfrak{q}^{G})^+.
\]

Set 
\[
\hatX(D) := \Sigma^{-N_+} \lr{F_\mathfrak{q}}.
\]
Our goal in the remainder of this subsection is to show that $\qakwX(D)$ and $\hatX(D)$ are weakly equivalent in $\ggSpectra_G$. The argument follows \cite[Section 5.4]{LLS-Burnside} with modifications to account for $G$.

\begin{definition}
\label{def:arrow cat}
    Given a small category $\cat{C}$, its \emph{arrow category} $\Arr(\cat{C})$ has objects given by morphisms $u\rar{f} v$ in $\cat{C}$ and morphisms given by commutative squares 
    \[
\begin{tikzcd}
    u \ar[r, "f"] \ar[d,"\varphi"'] & v \ar[d,"\psi"] \\
    u' \ar[r, "f'"] & v'
\end{tikzcd}.
    \]
    Composition is defined by vertically stacking squares.  
\end{definition}

\begin{definition}[{\cite[Definition 5.33]{LLS-Burnside}}]
    Given a strictly unitary lax $2$-functor $F: \cat{C} \to \til{\B}$, define the functor (also strictly unitary and lax) $\vec{F}: \Arr(\cat{C}) \to \til{\B}$ as follows. For an object $u\rar{f} v$, the correspondence $F(f)$ is also a set, and we define $\vec{F}(u\rar{f} v) = F(f)$. Given a commutative square $(\varphi, \psi)$ as in Definition \ref{def:arrow cat},  define $\vec{F}(\varphi, \psi) = F(f' \circ \varphi) = F(\psi \circ f)$. Source and target maps are given by 
    \[
F(f) \lar{} F(\psi) \times_{F(v)} F(f) \lar{\alpha} F(f'\circ \varphi) \rar{\beta} F(f') \times_{F(u')} F(\varphi) \rar{} F(f')
    \]
    where $\alpha, \beta$ are the bijections supplied by $F$ and the unlabeled arrows are the natural projections. 
For a pair of composable morphisms 
    \[
\begin{tikzcd}
    u \ar[r, "f"] \ar[d,"\varphi"'] & v \ar[d,"\psi"] \\
    u' \ar[r, "f'"] \ar[d, "\theta"'] & v' \ar[d, "\rho"] \\
    u'' \ar[r, "f''"] & v''
\end{tikzcd}
    \]
    the $2$-morphism from $
    \vec{F}(\theta,\rho) \times_{F(f')} \vec{F}(\varphi, \psi) = F(\rho\circ f') \times_{F(f')} F(f'\circ \varphi)$ to $\vec{F}(\theta \circ \varphi, \rho\circ \psi) = F(\rho \circ f' \circ \varphi)$ is the composition 
    \begin{align*}
        F(\rho\circ f') \times_{F(f')} F(f'\circ \varphi) & \cong F(\rho)\times_{F(v')} F(f') \times_{F(f')} F(f')\times_{F(u')} F(\varphi) \\
        & \cong F(\rho) \times_{F(v')}  F(f') \times_{F(u')} F(\varphi) \\
        & \cong F(\rho \circ f' \circ \varphi).
    \end{align*}
\end{definition}

Given $F: \cat{C} \to \til{\B}_G$, we may view it as a functor landing in $\til{\B}$ by composing with the forgetful functor $\til{\B}_G\to \til{\B}$ and applying the above to form $\vec{F} : \Arr(\cat{C}) \to \til{\B}$. Since $F$ landed in $\til{\B}_G$, the values of $\vec{F}$ on objects and morphisms are also finite free $G$-sets.

\begin{lemma}
    Given a strictly unitary lax $2$-functor $F: \cat{C} \to \til{\B}_G$, the associated functor $\vec{F}$ lands in $\til{\B}_G$. 
\end{lemma}
\begin{proof}
    We need to verify that correspondences and $2$-morphisms assigned by $\vec{F}$ are $G$-equivariant. This is clear from the definitions and that $F$ lands in $\til{\B}_G$.
\end{proof}

The following is our analogue of \cite[Lemma 5.37]{LLS-Burnside}. Recall the functor $B: \h{\cat{C}} \to \Arr(\cat{C})$ from \cite[Definition 5.28]{LLS-Burnside}. 

\begin{lemma}
\label{lem:hat vs little cubes}
    Given a strictly unitary lax $2$-functor $F: \cube^N \to \til{\B}_G$ and $k\geq 2N$, let $\squig{F}$ be an equivariant $k$-dimensional spacial refinement of $\vec{F}$ (guaranteed to exist by \cite[Proposition 5.5]{AKW}). Consider the homotopy coherent diagrams $\squig{F} \circ B : \h{\cube^N} \to \Top_*^G$ and $\h{F}_k^G : \h{\cube^N} \to \Top_*^G$. There is a morphism of homotopy coherent diagrams $J_K : \squig{F} \circ B \to \h{F}_k^G$ which induces an isomorphism on $H_i$ for $i\leq 2k-1$ on each object. 
\end{lemma}

\begin{proof}
    The proof of \cite[Lemma 5.37]{LLS-Burnside} goes through with minor modifications. For objects $u,v\in \cube^N$ with $u\leq v$, let $\varphi_{u,v}$ denote the unique morphism from $u$ to $v$. We need to build a homotopy coherent diagram $J_k$ over $\h{\cube^N} \times \cube^1$ which restricts to $\squig{F} \circ B$ and $\h{F}_k^G$ on $\h{\cube^N} \times \{0\}$ and $\h{\cube^N} \times \{1\}$, respectively. This determines $J_k$ on all objects. Consider an object $u\to v \to w$ in $\h{\cube^N}$. We have 
    \[
(\squig{F}\circ B)(u \rar{\varphi_{u,v}} v \rar{\varphi_{v,w}} w)
=
\bigvee_{x\in F(\varphi_{v,w} \circ \varphi_{u,v}) } S^k 
\overset{(!)}{=}
\bigvee_{a\in F(\varphi_{u,v})} 
\bigvee_{ \substack{b\in F(\varphi_{v,w}) \\ s(b) =t(a)} }
S^k,
    \]
\[
\h{F}_k^G (u \rar{\varphi_{u,v}} v \rar{\varphi_{v,w}} w)
=
\bigvee_{a\in F(\varphi_{u,v})} 
\prod_{ \substack{b\in F(\varphi_{v,w}) \\ s(b) =t(a)} }
S^k.
\] 
Note, the identification marked $(!)$ is given by the $2$-morphism 
\[
\alpha : F(\varphi_{v,w} \circ \varphi_{u,v}) \Rightarrow F(\varphi_{v,w}) \times_{F(v)} F(\varphi_{u,v})
\]
provided by $F$. The action of $G$ on both sides of $(!)$ simply permutes spheres according to the action of $G$ on the indexing sets, so $(!)$ is $G$-equivariant since $\alpha$ is.  

    The value of $J_k$  on the morphism $\id_{u\to v \to w} \times \varphi_{0, 1}$ in $\h{\cube^N} \times \cube^1$ is defined to be the natural inclusion of the $k$-skeleton
    \[
\bigvee_{a\in F(\varphi_{u,v})} 
\bigvee_{ \substack{b\in F(\varphi_{v,w}) \\ s(b) =t(a)} }
S^k
\hookrightarrow
\bigvee_{a\in F(\varphi_{u,v})} 
\prod_{ \substack{b\in F(\varphi_{v,w}) \\ s(b) =t(a)} }
S^k
    \]
which is $G$-equivariant by the above observation.

The families of maps needed to define the homotopy coherent diagram $J_k$ are built inductively as in \cite[Lemma 5.37]{LLS-Burnside}. The only modification is that the overlapping box map in item (X-2) of their proof is chosen to be $G$-equivariant. The connectivity result in Lemma \ref{lem:equivariant box maps facts} ensures we can inductively build $J_k$. 
\end{proof}

\begin{theorem}
   For an annular link diagram $D\subset \Aq$, $\qakwX(D)$ and $\qaX(D)$ are equivalent as spectra with $G$-action. 
\end{theorem}

\begin{proof}
   With Lemma \ref{lem:hat vs little cubes} in hand, one can show that $\qakwX(D) \simeq \hatX(D)$ exactly as in the proof of \cite[Theorem 7]{LLS-Burnside}. It then suffices to demonstrate an equivalence $\qaX(D) \simeq \hatX(D)$. Recall that $\qaX(D)$ is defined as the iterated mapping cone of a functor
    \[
\rect_J : \cube^N \to \Spectra_G
    \]
    which is the rectification of a functor out of the nerve
    \[
J: \nerve((\cube^N)')  \to \Spectra_G.
    \]
    Extend $J$ to $J^\dagger : \nerve((\cube^N)')_\dagger \to \Spectra_G$ by setting $J^\dagger(*)$ to be a point. The rectified functor $\rect_{J^\dagger} : \cube^N_\dagger \to \Spectra_G$ is weakly equivalent to $(\rect_J)^\dagger$ since rectification commutes with restriction up to weak equivalence. We have
    \[
\qaX(D) \simeq 
\hocolim_{\cube^N_\dagger} (\rect_J)^\dagger 
\simeq 
\hocolim_{\cube^N_\dagger} \rect_{J^\dagger}
\simeq 
\hocolim_{\nerve((\cube^N)')_\dagger}  \rect_{J^\dagger} \circ \pi^\dagger
\simeq
\hocolim_{\nerve((\cube^N)')_\dagger} J^\dagger
    \]
    The first equivalence is due to Lemma \ref{lem:hocolim over plus vs dagger}; the second is due to the above equivalence $\rect_{J^\dagger} \simeq (\rect_J)^\dagger$; the third is due to the
    projection $\pi^\dagger : \nerve((\cube^N)')_\dagger \to \cube^N_\dagger$ being a weak equivalence and hence homotopy cofinal (see property (ho-4) in \cite[Section 4.2]{LLS-Burnside}); and the fourth is due to the weak equivalence $\rect_{J^\dagger} \circ \pi^\dagger \simeq J^\dagger$ of Lemma \ref{lem:G rectification}. 

    The proof of \cite[Proposition 8.15]{LLS-Burnside} gives a zig-zag of maps  between $\hocolim_{\nerve((\cube^N)')_\dagger} J^\dagger$ and  $\hatX(D)$ on the \emph{underlying} spectra, each of which is a weak equivalence in $\Spectra$. Each intermediate spectrum in the zig-zag naturally inherits a $G$-action, and one checks that each of the five maps in the zig-zag is $G$-equivariant. 
    
    For instance, in Items (2) and (3) of the proof of \cite[Proposition 8.15]{LLS-Burnside}, the functors $\func{G}$ and $K(\func{G}_p)$ therein\footnote{In \cite{LLS-Burnside} these are denoted $G$ and $K(G_p)$. We use a different font to distinguish from the group $G$.}  can be viewed as landing in $\Spectra_G$ as follows. Let $q^k\in G$ act on 
    \[
\func{G}( u \rar{f} v \rar{g} w) = \bigvee_{a\in F_\q(f)} \prod_{ \substack{b\in F_\q(g) \\ s(b) = t(a)}} K(\Sets)
 \ \text{ and } \ 
K(\func{G}_p)( u \rar{f} v \rar{g} w) = \prod_{a\in F_\q(f)} \prod_{ \substack{b\in F_\q(g) \\ s(b) = t(a)}} K(\Sets)
    \]
    by permuting factors. Specifically, $q^k$ sends the wedge summand corresponding to $a\in F_\q(f)$ to the wedge summand corresponding to $q^k a \in F_\q(g)$ in the way that permutes factors, sending the copy of $K(\Sets)$ corresponding to $b \in F_\q(g)$ to the copy corresponding to $q^k b$; $q^k$ acts on $K(\func{G}_p)( u \rar{f} v \rar{g} w)$ similarly. Each of the natural transformations $\h{F}_\q^G \to \func{G}$ and $\func{G}\to K(\func{G}_p)$ is then an object-wise equivalence in $\Spectra_G$. In Item (5), the map from 
    \[
\func{G}_p'( u \rar{f} v \rar{g} w) =  \coprod_{a\in F(f)} \prod_{ \substack{b\in F(g) \\ s(b) = t(a)}} K(\Sets)
\ \text{ to } \ 
(\vec{F}_\q \circ B)'( u \rar{f} v \rar{g} w) \prod_{c\in F_\q(g\circ f)} K(\Sets) 
    \]
    is induced by the bijection $F_\q(g\circ f) \to F_\q(g) \times_{F_\q(v)} F_\q(f)$, which is $G$-equivariant, and the identification in $\Permu$ of products and coproducts with the Cartesian product (see \cite[Property (EM-2)]{LLS-Burnside}).  
    
\end{proof}

\appendix

\section{Proofs of Lemma \ref{lem:G rectification} and
Proposition \ref{prop:K_G}}
\label{sec:simplicial stuff}

Here we prove Lemma \ref{lem:G rectification} and
Proposition \ref{prop:K_G} from Section \ref{sec:Groupoid enrichments}. These statements deal with simplicially enriched categories and functors between them, and we refer the reader to \cite{Kelly} for the necessary background. By a \emph{simplicial category} we mean a category enriched over the monoidal category $(\sset,\times)$ of simplicial sets. We write $\Map_\cat{C}(a,b)$ (rather than $\Hom$) for the simplicial set of morphisms from $a$ to $b$ in a simplicial category $\cat{C}$.

As defined in \cite[Section 1.3]{symmetric_spectra}, the category  $\Spectra$ of symmetric spectra is enriched over (pointed) simplicial sets by 
\[
\Map_\Spectra(X,Y) := \Hom_\Spectra(X\wedge \Delta[-]_+,Y), 
\]
where $\Delta[n] = \Hom_\Delta(-, [n])$ is the standard $n$-simplex and $(-)_+$ denotes adding a disjoint basepoint. Composition of $n$-simplices is the map
\begin{equation}
\label{eq:composition in simplicial spectra}
\Map_{\Spectra}(X,Y)_n \times \Map_{\Spectra}(Y,Z)_n \to 
\Map_{\Spectra}(X,Z)
\end{equation}
sending $(X\wedge \Delta[n]_+\rar{f_1} Y, Y\wedge \Delta[n]_+\rar{f_2} Z)$ to 
\[
X\wedge \Delta[n]_+ \rar{\id \wedge  \delta } X \wedge \Delta[n]_+ \wedge \Delta[n]_+ \rar{f_1\wedge  \id} Y\wedge \Delta[n]_+ \rar{f_2} Z
\]
where $\delta$ is the diagonal map. There is a map $X\wedge \Delta[n]_+ \to X$ induced by the unique map $[n]\to [0]$ in $\Delta$ together with the canonical isomorphism $X\wedge \Delta[0]_+ \cong X$.
Similarly, $\Spectra_G$ is a simplical category by declaring $\Map_{\Spectra_G}(X,Y) = \Hom_{\Spectra_G}(X\wedge \Delta[-]_+,Y)$, where $g\in G$ acts on $X\wedge \Delta[-]_+$ by $g\wedge \id$, and composition is defined in the same way as above. 

\begin{remark}
    Note that $\Spectra$ and $\Spectra_G$ are enriched over pointed simplicial sets $(\sset_*, \wedge)$. The simplicial category $\nerve((\cube^N)')$ does not have natural basepoints. When considering enriched functors from a simplicial category $\cat{C}$ into $\Spectra$ or into $\Spectra_G$, we view these latter two categories as enriched over simplicial sets $(\sset, \times)$. This is why in \eqref{eq:composition in simplicial spectra} we use $\times$ rather than $\wedge$. 
\end{remark}

Given a category $\cat{D}$, let $\u{\cat{D}}$ be the \emph{constant} simplicially enriched category where $\Map_{\u{\cat{D}}}(a,b)$ is the constant simplicial set at $\Hom_{\cat{D}}(x,y)$; that is, the set of $n$-simplices of $\Map_{\u{\cat{D}}}(a,b)$ is equal to $\Hom_{\cat{D}}(x,y)$ for all $n$, and each structure map is the identity.

For enriched categories $\cat{C}$ and $\cat{D}$, let $\EnrFunc(\cat{C},\cat{D})$ denote the category of enriched functors and enriched natural transformations; see \cite[Section 1.2]{Kelly}. 
If  $\cat{C}$ and $\cat{D}$ are simplicial categories, denote by $\cat{C}\times \cat{D}$ the simplicial category whose objects are $\ob(\cat{C}) \times \ob(\cat{D})$ and $\Map_{\cat{C}\times \cat{D}}((a_1,b_1), (a_2, b_2)) = \Map_\cat{C}(a_1, a_2)\times \Map_{\cat{D}}(b_1, b_2)$.

The following may be viewed as an enriched analogue of the usual adjunction isomorphism $\Func(\cat{C}, \Func(\cat{D}, \cat{E}) \cong \Func(\cat{C} \times \cat{D}, \cat{E})$ where $\cat{C}, \cat{D}, \cat{E}$ are categories. Recall that for a group $G$, $BG$ is the category with one object $*$ and $\Hom_{BG}(*,*) = G$.

\begin{lemma}
\label{lem:simplicial cat adjunction}
    For any simplicial category $\cat{C}$, there is an isomorphism of categories 
    \[
    \EnrFunc(\cat{C} \times \u{BG}, \Spectra) \cong \EnrFunc(\cat{C}, \Spectra_G).
    \]
\end{lemma}

\begin{proof}
    This follows from \cite[Section 2.3]{Kelly}. 
\end{proof}

We now prove Lemma \ref{lem:G rectification}. 

\begin{proof}[Proof of Lemma \ref{lem:G rectification}]
    Let $\til{J} : \cat{C}\times \u{BG} \to \Spectra$ be the enriched functor corresponding to $J:\cat{C} \to \Spectra_G$ under the isomorphism in Lemma \ref{lem:simplicial cat adjunction}. Let $\til{\pi} : \cat{C}\times \u{BG}\to (\cat{C}\times \u{BG})^0 = \cat{C}^0\times BG$ be the projection, which is a weak equivalence since $\pi:\cat{C} \to \cat{C}^0$ is by assumption.  
    Lemma \ref{lem:rectification} provides a functor $\rect_{\til{J}}: \cat{C}^0 \times BG \to \Spectra$ and a zig-zag of natural transformations from $\til{J}$ to $\rect_{\til{J}}\circ \til{\pi}$, each of which is an object-wise weak equivalence. Define $\rect_J : \cat{C}^0 \to \Spectra_G$ to be obtained from $\rect_{\til{J}}$ by adjunction, 
    \[
\Func(\cat{C}^0\times BG, \Spectra) \cong \Func(\cat{C}^0, \Func(BG, \Spectra)).
    \]
    Let $(\rect_{\til{J}} \circ \til{\pi})^* : \cat{C}\to \Spectra_G$ denote the enriched functor corresponding to $\rect_{\til{J}} \circ \til{\pi}$ under Lemma \ref{lem:simplicial cat adjunction}. 
    
    Passing the zig-zag of natural transformations through the isomorphism in Lemma \ref{lem:simplicial cat adjunction} yields a zig-zag of natural transformations from $J$ to $(\rect_{\til{J}} \circ \til{\pi})^*$.
 By examining the isomorphism, we see that each of the corresponding functors in the zig-zag $\cat{C} \to \Spectra_G$ is also an object-wise weak equivalence  (the isomorphism is, essentially, the identity on natural transformations). It remains to verify that $(\rect_{\til{J}} \circ \til{\pi})^* = \rect_J\circ \pi$, which is a straightforward unwinding of definitions. 
 
\end{proof}

Our final goal is to prove Proposition \ref{prop:K_G}. 
The following definition and lemma will be useful. 

\begin{definition}
\label{def:twisted correspondence} 
Suppose $X\lar{s} A \rar{t} Y$ is a $1$-morphism in $\B$ and $X$ and $Y$ are free $G$-sets (but $A$ does not necessarily satisfy property \ref{item:E}). For $g,h\in G$, let $A^{g,h}$ denote the \emph{twisted} correspondence given by $A^{g,h}_{y,x} = A_{g^{-1}y, h^{-1}x}$. We denote this schematically by $X \lar{hs} A \rar{gt} Y$. Of particular importance will be 
\[
X\lar{\id} X \rar{g} X,
\]
which corresponds to the twisted identity correspondence $1_X^{g,1}$. If $\alpha:A\Rightarrow B$ is a $2$-morphism in $\B$, let $\alpha^{g,h} : A^{g,h} \Rightarrow B^{g,h}$ be given by $
\alpha^{g,h}_{y,x} = \alpha_{g^{-1}y, h^{-1}x}$. We represent $\alpha^{g,h}$ schematically as 
\[
 \begin{tikzcd}[row sep = tiny]
     & A \ar[dl, "hs"'] \ar[dr, "gt"] \ar[dd, Rightarrow, "\alpha^{g,h}"] & \\
     X & & Y \\
        & B \ar[ul, "hs'"] \ar[ur, "gt'"'] & 
 \end{tikzcd}
\]
\end{definition}

Note that if $A$ is a $1$-morphism in $\B_G$, then $A^{g,h}$ is again a $1$-morphism in $\B_G$ from $X$ to $Y$. Moreover, there are equalities $A^{g,h} = A^{1,g^{-1}h} = A^{h^{-1}g,1}$, but we allow both $g$ and $h$ for additional flexibility. In particular, $A^{g,g} = A$ by \ref{item:E}, and also $\alpha^{g,g} = \alpha$.

\begin{lemma}
Suppose $X\lar{s} A \rar{t} Y$ is a $1$-morphism in $\B$ where $X$ and $Y$ are free $G$-sets and $g\in G$. The following equalities hold in $\B$. 
    \begin{equation}
        \begin{tikzcd}[row sep = small, column sep =small]
            & X \ar[dl, "g"'] \ar[dr, "\id"]\\
            X && X
        \end{tikzcd}
        =
            \begin{tikzcd}[row sep = small, column sep =small]
            & X \ar[dl, "\id"'] \ar[dr, "g^{-1}"]\\
            X && X
        \end{tikzcd}
    \end{equation}

\begin{equation}
\label{eq:twist on left}
\begin{tikzcd}[row sep = small, column sep =small]
    & X \ar[dl, "\id"'] \ar[dr, "g^{-1}"] & & A \ar[dl, "s"'] \ar[dr, "t"] & \\
    X & & X & & Y 
\end{tikzcd}
=
\begin{tikzcd}[row sep = small, column sep =small]
    & A \ar[dl, "g s"'] \ar[dr, " t"] & \\
    X & & Y 
\end{tikzcd}
\end{equation}

\begin{equation}
\label{eq:twist on right}
\begin{tikzcd}[row sep = small, column sep =small]
    & A \ar[dl, "s"'] \ar[dr, "t"] & & Y \ar[dl, "\id"'] \ar[dr, "g"] & \\
    X & & Y & & Y 
\end{tikzcd}
=
\begin{tikzcd}[row sep = small, column sep =small]
    & A \ar[dl, "s"'] \ar[dr, " gt"] & \\
    X & & Y 
\end{tikzcd}
\end{equation}

    \begin{equation}
\label{eq:strict G action}
    \begin{tikzcd}[row sep = small, column sep =small]
    & X \ar[dl, "\id"'] \ar[dr, "h"] & & X \ar[dl, "\id"'] \ar[dr, "g"] & \\
    X & & X & & X 
\end{tikzcd}
=
\begin{tikzcd}[row sep = small, column sep =small]
    & X \ar[dl, "\id"'] \ar[dr, "g h"] & \\
    X & & X 
\end{tikzcd}
\end{equation}
If moreover $A$ is a $1$-morphism in $\B_G$, then we have 
\begin{equation}
    \label{eq:twisting equivariant correspondence}
       \begin{tikzcd}[row sep = small, column sep =small]
            & A \ar[dl, "hs"'] \ar[dr, "gt"]\\
            X && Y
        \end{tikzcd}
        =
            \begin{tikzcd}[row sep = small, column sep =small]
            & A \ar[dl, "s"'] \ar[dr, "h^{-1}gt"]\\
            X && Y
        \end{tikzcd}\ 
            =
            \ 
            \begin{tikzcd}[row sep = small, column sep =small]
            & A \ar[dl, "g^{-1}hs"'] \ar[dr, "t"]\\
            X && Y
        \end{tikzcd}.
\end{equation}
\end{lemma}

\begin{proof}
    These are all straightforward. For instance, \eqref{eq:twist on right} follows from
    \[
(1_Y^{g,1} \circ A)_{y,x} = \bigcup_{y'\in Y} (1_Y)_{g^{-1}y,y'} \times A_{y',x} = A_{g^{-1}y,x} = A^{g,1}_{y,x}.
    \]
    Equation \eqref{eq:twisting equivariant correspondence} follows from property \ref{item:E}: $
A^{g,h}_{y,x} = A_{g^{-1}y,h^{-1}x} = A_{hg^{-1}y,x} = A_{y,x}^{h^{-1}g}$, and similarly $A^{g,h}_{y,x} = A_{y, g h^{-1} x} = A_{y,x}^{1, g^{-1} h}$. 
\end{proof}

\begin{proof}[Proof of Proposition \ref{prop:K_G}]
 Let $X$ be an object of $\B_G$. We are forced to define the underlying spectrum of $K_G(X)$ to be equal to $K(X)$ and the map $K(X) \to K(Y)$ assigned to a $1$-morphism $X\lar{} A \rar{} Y$ in $\B_G$ to be equal to $K(A)$. We need to (1) define the $G$-action on $K(X)$, (2) explain why $K(A)$ is a map in $\Spectra_G$, and (3) explain why $K_G$ is simplicially enriched. 

 Let us begin with (1). Given $g\in G$, recall the $1$-morphism $X\lar{\id} X \rar{g} X$ from Definition \ref{def:twisted correspondence}. Define the action of $g\in G$ on $K(X)$ to be given by the map
 \[
K(X\lar{\id} X \rar{g} X): K(X) \to K(X).
 \]
 This defines a $G$-action on $K(X)$ by \eqref{eq:strict G action}. For (2), note that \eqref{eq:twist on left}, \eqref{eq:twist on right}, and \eqref{eq:twisting equivariant correspondence} imply that if $X\lar{} A \rar{} Y$ is a $1$-morphism in $\B_G$, then $K(A) : K(X)\to K(Y)$ is $G$-equivariant when $K(X)$ and $K(Y)$ are equipped with the above $G$-actions. 

 For item (3), let $X,Y$ be objects in $\B_G$ and consider the composite 
 \begin{equation}
 \label{eq:KG is enriched}
\Map_{\nerve(\B_G)}(X,Y) \to \Map_{\nerve(\B)}(X,Y) \rar{K} \Map_{\Spectra}(K(X),K(Y)).
 \end{equation}
 The first map, induced by the forgetful functor, is evidently a map of simplicial sets, and so is the second map since $K$ is an enriched functor. There is a $G$-action on $\Map_{\nerve(\B)}(X,Y)$ given by sending a correspondence $X \lar{s} B \rar{t} Y$ (which may not satisfy property \ref{item:E}!) to the twisted correspondence $B^{g,g}= (X \lar{gs} B \rar{gt} Y)$  of Definition \ref{def:twisted correspondence}, and a sequence of composable $2$-morphisms $\u{\alpha} = (B_1 \overset{\alpha_1}{\Rightarrow} B_2 \overset{\alpha_2}{\Rightarrow} \cdots \overset{\alpha_n}{\Rightarrow} B_{n+1})$, which is a $n$-simplex in $\Map_{\nerve(\B)}(X,Y)$,
 to $g\cdot \u{\alpha} = (B_1 \overset{\alpha_1^{g,g}}{\Rightarrow} B_2 \overset{\alpha_2^{g,g}}{\Rightarrow} \cdots \overset{\alpha_n^{g,g}}{\Rightarrow} B_{n+1})$. We denote this schematically by 
 \[
 g\cdot \ \left(
\begin{tikzcd}[row sep=small]
    & B_1 \ar[dl,"s_1"'] \ar[dr, "t_1"] \ar[d, Rightarrow, "\alpha_1"] & \\
    X & \vdots  \ar[d, Rightarrow, "\alpha_{n}"] & Y \\
    & B_{n+1} \ar[ul, "s_{n+1}"] \ar[ur, "t_{n+1}"'] 
\end{tikzcd}
\right)
=
\begin{tikzcd}[row sep=small]
    & B_1 \ar[dl,bend right=10,"gs_1"'] \ar[dr,bend left=10, "gt_1"] \ar[d, Rightarrow, "\alpha_1^{g,g}"] & \\
    X & \vdots  \ar[d, Rightarrow, "\alpha_{n}^{g,g}"] & Y \\
    & B_{n+1} \ar[ul,bend left=10, "gs_{n+1}"] \ar[ur,bend right=10, "gt_{n+1}"'] 
\end{tikzcd}
 \]
Note that the forgetful map $\Map_{\nerve(\B_G)}(X,Y) \to \Map_{\nerve(\B)}(X,Y)$ lands in the $G$-fixed points $\Map_{\nerve(\B)}(X,Y)^G$, since if $A$ is a $1$-morphism in $\B_G$ then $A^{g,g} = A$, and similarly if $\alpha$ is a $2$-morphism in $\B_G$ then $\alpha^{g,g} = \alpha$. 
 
 Similarly, $G$ acts on $\Map_{\Spectra}(K(X),K(Y))$ by conjugation:
 \[
 g\cdot \left( K(X) \wedge \Delta[n]_+ \rar{f} K(Y) \right)
 :=
 K(X) \wedge \Delta[n]_+ \rar{g^{-1}\wedge \id} K(X) \wedge \Delta[n]_+ \rar{f} K(Y) \rar{g} K(Y).
 \]
We will argue that $\Map_{\nerve(\B)}(X,Y) \rar{K} \Map_{\Spectra}(K(X),K(Y))$ is $G$-equivariant. It follows that the composition \eqref{eq:KG is enriched} lands in the fixed points $\Map_{\Spectra}(K(X),K(Y))^G = \Map_{\Spectra_G}(K(X),K(Y))$, which will complete the proof. 

That $K$ is an enriched functor means, in particular, that the diagram
 \[
\begin{tikzcd}
    \Map_{\nerve(\B)}(X,X) \wedge \Map_{\nerve(\B)}(X,Y) \wedge \Map_{\nerve(\B)}(Y,Y) \ar[r] \ar[d, "K\wedge K \wedge K"] & \Map_{\nerve(\B)}(X,Y) \ar[d, "K"] \\
    \Map_{\Spectra}(K(X),K(X))\wedge \Map_{\Spectra}(K(X), K(Y)) \wedge \Map_{\Spectra}(K(Y),K(Y)) \ar[r] & \Map_{\Spectra}(K(X), K(Y))
\end{tikzcd}
 \]
commutes, where horizontal maps are composition. Equations \eqref{eq:twist on left} and \eqref{eq:twist on right} imply that $g\cdot \u{\alpha}$ is the image under the top horizontal map of 
 \[
 \left(
\begin{tikzcd}[row sep=small]
    & X \ar[dl,"\id"'] \ar[dr, "g^{-1}"] \ar[d, Rightarrow, "\id"] & \\
    X & \vdots  \ar[d, Rightarrow, "\id"] & X \\
    & X \ar[ul, "\id"] \ar[ur, "g^{-1}"'] 
\end{tikzcd},
\begin{tikzcd}[row sep=small]
    & B_1 \ar[dl,"s_1"'] \ar[dr, "t_1"] \ar[d, Rightarrow, "\alpha_1"] & \\
    X & \vdots  \ar[d, Rightarrow, "\alpha_{n}"] & Y \\
    & B_{n+1} \ar[ul, "s_{n+1}"] \ar[ur, "t_{n+1}"'] 
\end{tikzcd},
\begin{tikzcd}[row sep=small]
    & Y \ar[dl,"\id"'] \ar[dr, "g"] \ar[d, Rightarrow, "\id"] & \\
    Y & \vdots  \ar[d, Rightarrow, "\id"] & Y \\
    & Y \ar[ul, "\id"] \ar[ur, "g"'] 
\end{tikzcd}
\right).
 \]
Denote this triple by $(\u{\beta}, \u{\alpha}, \u{\gamma})$. Note that $\u{\beta}$ (resp. $\u{\gamma}$) is the image of the $0$-simplex $X \lar{\id} X \rar{g^{-1}} X$ (resp. $Y \lar{\id} Y \rar{g} Y$) under the unique map to the $n$-simplices. That $K$ is a map of simplicial sets implies that $K(\u{\beta})$ and $K(\u{\gamma})$ are equal to 
\[
K(X) \wedge \Delta[n]_+ \to K(X) \rar{g^{-1}} K(X)  \ \text{ and } \  K(Y) \wedge \Delta[n]_+ \to K(Y) \rar{g} K(Y),
\]
respectively. It follows that $K$ is $G$-equivariant, which completes the proof. 
\end{proof}

\section{Summary of notation}
\label{sec:summary of notation}

\renewcommand{\arraystretch}{1.5}
\begin{table}[H]
    \centering
     \begin{tabulary}{\textwidth}{l|J|l}
    \multicolumn{3}{c}{
\textbf{Cobordism categories}
}\\
    Notation & Meaning & Found in \\ \hline \hline
    $\D$ & The open disk $(0,1) \times (-1,1)$ &  Sec. \ref{sec:The Khovanov and annular Khovanov functors} \\ \hline
    $\A$ & The open annulus $S^1 \times (-1,1)$ & Sec. \ref{sec:The Khovanov and annular Khovanov functors} \\ \hline 
    $\Aq$ & The quantum annulus &  Def. \ref{def:quantum annulus} \\ \hline 
        $\BN(\D)$, $\BN(\A)$ &  Bar-Natan categories of the disk and annulus & Def. \ref{def:BN category} \\ \hline
        $\divcob(\D)$, $\divcob(\A)$ & Divided cobordism categories of the disk and annulus & Def \ref{def:divcob} \\ \hline
        $\BN(\Aq)$ & The quantum annular Bar-Natan category & Def. \ref{def:qBN cat} \\ \hline 
        $\divcobq$ & The quantum divided cobordism category of $\A_q$ &  Def. \ref{def:qdivcob of annulus} \\ \hline 
    \end{tabulary}
\end{table}


\renewcommand{\arraystretch}{1.5}
\begin{table}[H]
    \centering
     \begin{tabulary}{\textwidth}{l|J|l}
    \multicolumn{3}{c}{
\textbf{Domain shape categories}
}\\
    Notation & Meaning & Found in \\ \hline \hline
    $\cube^N$ & $\{0,1\}^N$, the $N$-dimensional cube category & Def. \ref{def:cube cat} \\ \hline
    $\shape(X)$, $\shape(X,Y)$ & Shape multicategories & Defs. \ref{def:shape(X) multicat}, \ref{def:shape gluing multicat} \\ \hline 
$\HMshape(S)$, $\HMshape^\omega(S)$ & Hochschild-Mitchell and augmented Hochschild-Mitchell shape categories &     Def. \ref{def:HM shape category} \\ \hline 

$\cube^N\tiltimes \shape(\matchingsk,\matchingsk)$ & A certain type of product multicategory & Sec. \ref{sec:spectral platalg and tangle bimod} \\ \hline
    \end{tabulary}
\end{table}


\begin{table}[H]
    \centering
    \begin{tabulary}{\textwidth}{l|J|l}
    \multicolumn{3}{c}{
\textbf{Algebraic constructions}
}\\
 Notation & Meaning & Found in \\ \hline \hline
 $G$ & A cyclic group with distinguished generator $q$ & \\ \hline 
$\KC(L)$, $\aKC(L)$ & The Khovanov and annular Khovanov chain complex  & Sec. \ref{sec:The Khovanov and annular Khovanov functors} \\ \hline 
$\qaKC(L)$, $\qaKh(L)$ & The quantum annular complex and its homology & Sec. \ref{sec:quantum annular homology defs} \\ \hline 
$\matchingsk$ & Crossingless matchings of $2n$ points with no matchings between the top $k$ and bottom $n-k$ points & Def. \ref{def:matchings and platforms}  \\ \hline
$\displaystyle \Plat{n} = \prod_{0\leq k \leq n} \Platk{n}{k}$ & The platform algebra & Def. \ref{def:platform algebra} \\ \hline 
$\FCK(T)$ & The platform bimodule for a planar tangle $T$ &  Def. \ref{def:platform bimodule} \\ \hline 
$\CCK(T)$ & The chain complex of platform bimodules for a tangle diagram $T$ & Sec. \ref{sec:Platform algebras and bimodules} \\ \hline 
$\qPlat{n}$, $\qPlatnk$, $\qCCK(T)$ & Scalar extensions of platform algebras and bimodules, obtained by applying $(-) \otimes \Z[G]$ & Eqn \eqref{eq:quantum extensions of CK} \\ \hline
    \end{tabulary}
\end{table}


\begin{table}[H]
    \centering
    \begin{tabulary}{\textwidth}{l|J|l}
    \multicolumn{3}{c}{
\textbf{Functors out of cobordism categories}
}\\
 Notation & Meaning & Found in \\ \hline \hline
 $\FKh, \FA$ & Khovanov and annular Khovanov TQFTs, out of $\BN(\D)$ and $\BN(\A)$. The same notations are also used for functors out of $\divcob(\D)$ and $\divcob(\A)$  & Defs. \ref{def:FKh on BN(D)}, \ref{def:FA on BN(A)}, Lem. \ref{lem:FKh and FA work on divcob}  \\ \hline 
 $\FAq$ & Quantum annular TQFT, out of $\BN(\A_q)$ & Sec. \ref{sec:quantum annular homology defs} \\ \hline 
 $\choice$ & A choice of generators & Defs. \ref{def:choice of generators}, \ref{def:coherent choice of generators} \\ \hline 
 $\choiceFAq$ & A version of $\FAq$ using $\choice$, out of $\BNq$ or $\divcobq$ & Def. \ref{def:FAq choice}  Lem. \ref{lem:FAq choice well-defined on divcob}\\ \hline 
 $\B_G$, $\B$ & The (free) $G$-Burnside $2$-category and its $G=\{1\}$ specialization & Def. \ref{def:G-Burn} \\ \hline 
 $\VHKKlax[X]$ & Lax functors $\divcob(X) \to \B$ lifting $\F_X$ for $X=\D, \A$ & Sec. \ref{sec:VHKKlax for D,A} \\ \hline 
 $\cat{C}'$ & Canonical groupoid enrichment of a multicategory $\cat{C}$ & Sec. \ref{sec:Groupoid enrichments} \\ \hline
 $\VHKK[X]'$ & The strict versions $\divcob(X)' \to \B$ of $\VHKKlax[X]$ & Def. \ref{def:VHKK} \\ \hline 
 $\qVHKKlax$ & Lax functor $\divcobq \to \B_G$, lifting $\choiceFAq$ & Def. \ref{def:qVHKK Z lax}, Prop. \ref{prop:qVHKK lax} \\ \hline 
 $\qVHKKchoice'$ & The strict version $\divcobq' \to \B_G$ & Def. \ref{def:qVHKK} \\ \hline 
    \end{tabulary}
\end{table}


\begin{table}[H]
    \centering
    \begin{tabulary}{\textwidth}{l|J|l}
    \multicolumn{3}{c}{
\textbf{Spectra}
}\\
 Notation & Meaning & Found in \\ \hline \hline
 $\Sq$ & The spherical group ring & Def. \ref{def:spherical group ring} \\ \hline
 $ \Spectra_G$ & The category $\Func(BG, \Spectra)$ of spectra with $G$-action, equivalent to $\Rmod[\Sq]$ & Def. \ref{def:S_G} \\ \hline
 $\displaystyle \platalg = \bigvee_{0\leq k \leq n} \platalgk$ & Spectral lifts of platform algebras & Def. \ref{def:spec plat alg and CK tangle spectrum} \\ \hline 
 $\displaystyle \X(T) = \bigvee_{0\leq k \leq n} \Xk(T)$ & Spectral lifts of platform bimodules for a tangle $T$ & Def. \ref{def:spec plat alg and CK tangle spectrum} \\ \hline 
 $\qplatalgk$, $\qXk(T)$ & Scalar extensions, obtained by applying $(-) \wedge \Sq$ & Eqn. \eqref{eq:scalar ext for Sq}  \\ \hline 
 $\qaX(L)$ & The quantum annular spectrum built using $\qVHKKchoice'$ & Def. \ref{def:qAKh spectrum} \\ \hline 
 $\qakwX(L)$ & The quantum annular spectrum from \cite{AKW} & Def. \ref{def:AKW spectrum} \\ \hline 
    \end{tabulary}
\end{table}

\bibliographystyle{alpha}
\bibliography{references}

\end{document}